\documentstyle[11pt,amsfonts,amssymb,fullpage]{amsart}
          \newcommand{\X}{\textstyle \scriptscriptstyle  X}
          \newcommand{\D}{\textstyle \scriptscriptstyle  D}
          \newcommand{\N}{\textstyle \scriptscriptstyle  N}
          \newcommand{\V}{\textstyle \scriptscriptstyle  V}
          \newcommand{\diam}{\operatorname{diam}}
          \newcommand{\supp}{\operatorname{supp}}
          \newcommand{\Spec}{\operatorname{Spec}}
          \newcommand{\osc}{\operatorname{osc}}
          \newcommand{\Vol}{\operatorname{\it Vol}}
          \newcommand{\Diff}{\operatorname{\rm Diff}}
          
           \newcommand{\smallplus}{\operatorname{\scriptscriptstyle +}}
          \newcommand{\smallminus}{\operatorname{\scriptscriptstyle -}}
     
           \newsymbol\lnleq 120A

          \def\Ker{\mathop{\rm Ker\, }\nolimits}
          \def\osc{\mathop{\rm osc\, }}
          \def\vol{\mathop{\rm vol}\nolimits}
          \def\dist{\mathop{\rm dist}\nolimits}
          \def\Im{\mathop{\rm Im\, }\nolimits}
          \def\tr{\mathop{\rm tr\, }\nolimits}
          \def\inj{\mathop{\rm inj\, }\nolimits}

   \newtheorem{specialtheorem}{Theorem}
    
 \newtheorem{theorem}{Theorem}[section]
 \newtheorem{corollary}[theorem]{Corollary}
 \newtheorem{lemma}[theorem]{Lemma}
 \newtheorem{proposition}[theorem]{Proposition}
 
 \theoremstyle{remark}
 \newtheorem{remark}[theorem]{\bf Remark}

 \newtheorem{definition}[theorem]{\bf Definition}

\renewcommand{\subsectionname}{$\mathsection$\kern0pt}

 \newcommand{\Roof}[2]{#2{#1}}
 \newcommand{\roof}[2]{#2{#1}}

\begin{document}

\makeatletter
\def\@sect#1#2#3#4#5#6[#7]#8{%
  \edef\@toclevel{\ifnum#2=\@m 0\else\number#2\fi}%
  \ifnum #2>\c@secnumdepth \let\@secnumber\@empty
  \else \@xp\let\@xp\@secnumber\csname the#1\endcsname\fi
 \ifnum #2>\c@secnumdepth
   \let\@svsec\@empty
 \else
    \refstepcounter{#1}%
    \protected@edef\@svsec{\ifnum#2<\@m
       \@ifundefined{#1name}{}{%
         \ignorespaces\csname #1name\endcsname\space}\fi
       \@nx\textup{%
      \@nx\@secnumfont
         \csname the#1\endcsname.}\enspace
    }%
  \fi
  \@tempskipa #5\relax
  \ifdim \@tempskipa>\z@ 
    \begingroup #6\relax
    \@hangfrom{\hskip #3\relax\@svsec}{\interlinepenalty\@M #8\par}%
    \endgroup
    \ifnum#2>\@m \else \@tocwrite{#1}{#8}\fi
  \else
  \def\@svsechd{#6\hskip #3\@svsec
    \@ifnotempty{#8}{\ignorespaces#8\unskip
       \@addpunct.}%
    \ifnum#2>\@m \else \@tocwrite{#1}{#8}\fi
  }%
  \fi
  \global\@nobreaktrue
  \@xsect{#5}}


\renewcommand{\subsectionname}{\textsection\ignorespaces}
\def\@secnumfont{\bfseries}




 \title[]{Scalar Curvature, Metric Degenerations and
 the Static Vacuum Einstein Equations on 3-Manifolds, I.}
 
 \author[]{Michael T. Anderson}
 
 \thanks{Partially supported by NSF Grants DMS 9505744 and 9802722} 

\maketitle
 
 \tableofcontents

 \setcounter{section}{-1}

 
 


 \section{\bf Introduction} \setcounter{equation}{0}


  In this paper, we prove that degenerations of
  sequences of Yamabe metrics on \hbox{3-manifolds} are
  modeled or described by solutions to the static
  vacuum Einstein equations. One underlying
  motivation to understand such degenerations is the
  question of existence of constant curvature metrics
  on \hbox{3-manifolds,} in other words with the
  geometrization conjecture of Thurston [Th2]. An
  approach towards resolving this conjecture via
  study of Yamabe metrics is outlined in [An1].
 
  Let ${\Bbb M}  $ denote the space of all smooth
  Riemannian metrics on a closed, oriented 3-manifold
  $M,$ and ${\Bbb M}_{1}$ the subset of metrics
  satisfying $\vol_{g}M = 1.$ Define the total scalar
  curvature or Einstein-Hilbert action ${\cal S} :
  {\Bbb M}  \rightarrow  {\Bbb R} $ by
    \begin{equation} 
     \label{e0.1} 
          {\cal S} (g) = v^{-1/3}\int_{M}s_{g}dV_{g}, 
    \end{equation}
  where $s_{g}$ is the scalar curvature of the metric
  $g,$ $dV_{g}$ is the volume form associated with $g$
  and $v$ is the volume of $(M, g).$ The critical
  points of ${\cal S} $ are {\it  Einstein metrics},
  i.e. metrics satisfying the equation
 \begin{equation} 
     \label{e0.2} 
   -\nabla{\cal S}|_{g} =
              v^{-1/3}z = v^{-1/3} \left(r -  \frac{s}{3}\cdot  g\right) = 0,
     \end{equation}
  where $r$ is the Ricci curvature of $g$ and $z$ is
  the traceless Ricci curvature. In dimension 3, (and
  only in this dimension), Einstein metrics are
  exactly the metrics of constant curvature.

  There is a well-known minimax procedure to obtain
  critical values of ${\cal S} .$ First, the solution
  to the Yamabe problem [Y], [Tr], [Au1], [Sc1]
  implies that in the conformal class $[g]$ of any
  metric $g\in{\Bbb M}_{1},$ there is a metric
  $\bar{g} \in  {\Bbb M}_{1}$ which realizes
  the infimum $\mu [g]$ of ${\cal S}|_{[g]\cap{\Bbb
  M}_{1}},$ i.e.
        \begin{equation} 
          \label{e0.3} s_{\roof{g}{\bar}} =
            \mu [g] \equiv  \inf \ {\cal S}|_{[g]\cap{\Bbb M}_{1}}.  
           \end{equation}
  Such metrics are called {\it Yamabe metrics}. Let
  ${\cal C} $ denote the space of Yamabe metrics in
  ${\Bbb M} ,$ (of arbitrary volume), and ${\cal
  C}_{1}$ the subset of unit volume Yamabe metrics;
  any metric in ${\cal C} $ can be scaled uniquely to
  a metric in ${\cal C}_{1}.$ In dimension 2, the
  space ${\cal C}_{1},$ when divided by the action of
  the diffeomorphism group, corresponds to the
  Teichm\"uller or moduli space of a surface,
  and the functional $s_{\roof{g}{\bar}}: {\cal
  C}_{1} \rightarrow  {\Bbb R} $ is a constant
  function, (by the Gauss-Bonnet theorem). In
  dimension 3, (and above), this is an infinite
  dimensional space, and the functional is not
  constant.  A simple but important comparison
  argument of Aubin [Au1] implies that, for any
  conformal class [g] on any $M,$ 
 	\begin{equation} \label{e0.4} 
 		\mu [g] \leq  \mu (S^{3}, g_{can}),
 	\end{equation} 
  where $g_{can}$ is the canonical
  constant curvature metric on $S^{3}$ of volume 1.
  Define the {\it Sigma constant} $\sigma (M)$ of $M$ by 
 	\begin{equation} \label{e0.5} 
 	    \sigma (M) = \sup_{{\cal C}_{1}}\mu [g].  
 	\end{equation}
 
  This is a smooth invariant of the 3-manifold $M,$
  which should act very much like an Euler
  characteristic for 3-manifolds, when appropriately
  normalized. In case $\sigma (M) \leq  0,$ it is
  easy to prove, c.f.~[Bes, Prop.4.47], that any
  Yamabe metric $g_{o}\in{\cal C}_{1}$ such that
  $s_{g_{o}}=\sigma (M)$ is Einstein. In case
  $\sigma (M) >  0,$ this has been conjectured to
  be true, c.f. [Bes, Remark4.48], but remains still
  unknown.
 
  Of course, an arbitrary closed, oriented 3-manifold
  does not admit an Einstein metric; this is the case
  if for instance $M$ has an essential 2-sphere.
  Thus, if $M$ is an arbitrary closed 3-manifold and
  if $\{g_{i}\}$ is a maximizing sequence of unit
  volume Yamabe metrics on $M,$ that is
   \begin{equation} \label{e0.6} 
        s_{g_{i}} \rightarrow \sigma (M), 
   \end{equation} 
  then in general this sequence must degenerate in some manner.
 
  More specifically, if $\{g_{i}\}$ satisfies the
  uniform curvature bound
 	  \begin{equation} \label{e0.7}
 	  \int_{M}|z_{g_{i}}|^{2}dV_{g_{i}} \leq  \Lambda ,
 	  \end{equation}
  for some $\Lambda  <  \infty ,$ then as outlined
  in [An1], the structure of a suitable subsequence
  of $\{g_{i}\}$ and its limit can be used to
  geometrize the 3-manifold $M,$ (at least when $\sigma
  (M) \leq  0).$ This is accomplished essentially via
  the Cheeger-Gromov theory of convergence and
  collapse of Riemannian manifolds [C], [G, Ch.8],
  [CG1,2]. Since an arbitrary 3-manifold cannot be
  given a geometric structure, a maximizing sequence
  $\{g_{i}\}$ of unit volume Yamabe metrics must in
  general satisfy
 	\begin{equation} \label{e0.8}
 		\int_{M}|z_{g_{i}}|^{2}dV_{g_{i}} \rightarrow
 		\infty , \ \ \ \ \ 
                {\rm as}\  i \rightarrow  \infty .  
         \end{equation}
  For the purposes of this paper, we will say that a
  sequence $\{g_{i}\}$ satisfying (0.8) {\it
  degenerates.}  One then would like to relate the
  geometry of a degenerating sequence of Yamabe
  metrics with some topological structure on the
  underlying manifold $M$ that is the cause for the
  degeneration, c.f. [An1, \S 4-6] for further
  discussion.
 
  The main purpose of this paper is to show that the
  degenerations of such a sequence $\{g_{i}\}$ are
  described by solutions to the static vacuum
  Einstein equations. Further, under natural
  conditions, the degenerations correspond to
  non-trivial solutions of these equations.
 
  Some further background is needed to explain this;
  we refer to \S 1-\S 2 for further details. Let $L$
  be the linearization of the scalar curvature
  function at $g$ given by
          \begin{equation} \label{e0.9} 
           L(\alpha ) = \left. \frac{d}{dt}s(g+t\alpha )\right|_{t=0} 
                      =  -\Delta \tr\alpha  + \delta\delta\alpha  \,
                         - \langle r, \alpha\rangle , 
          \end{equation} 
 c.f. [Bes, Ch.1K]. The $L^{2}$ adjoint $L^{*}$ of $L$ is given by
 	\begin{equation} \label{e0.10} 
               L^{*}(f) = D^{2}f -\Delta f\cdot  g - f\cdot  r.
         \end{equation} 
 For any $g\in{\cal C} ,$ one thus has
 a natural splitting
 	\begin{equation} \label{e0.11} 
 	   T_{g}{\Bbb M}  = T_{g}{\cal C}\oplus N_{g}{\cal C} , 
 	\end{equation}
 where 
 	$$T_{g}{\cal C}  = \{\alpha\in T_{g}{\Bbb M} :
 	L(\alpha ) = const\},  
         \ \ \ \ \ 
          N_g{\cal C}  = \left\{\beta\in
 	T_{g}{\Bbb M} : \beta  = L^{*}(q), \int q = 0 \right\}, $$
 are the (formal) tangent and normal spaces to
  ${\cal C} $ in ${\Bbb M} .$
 
  In particular, the $L^{2}$ gradient $\nabla{\cal
  S}  = -  v^{- 1/3}z$ of the scalar curvature
  functional ${\cal S} $ splits at $g$ as
        \begin{equation} \label{e0.12} 
          z = z^{T} + z^{N},
        \end{equation} 
 where $-  v^{- 1/3}z^{T}$ is the
 gradient of the scale invariant functional
 $v^{2/3}\cdot  s = {\cal S}|_{{\cal C}}$ on ${\cal
 C} .$
 
  Another natural splitting of $T{\Bbb M} ,$
  differing from (0.10) by just a 1-dimensional
  factor, is \begin{equation} \label{e0.13}
 T_{g}{\Bbb M}  = \Ker L\oplus  \Im L^{*}.
 \end{equation} As in (0.12), we have then an $L^{2}$
 orthogonal sum \begin{equation} \label{e0.14} z =
 \xi  + L^{*}f, \end{equation} for some $\xi\in \Ker
 L$ and $f\in C^{\infty}(M).$ We set \begin{equation}
 \label{e0.15} u = 1+f.  \end{equation}

  Next, we describe briefly the {\it  static vacuum
  Einstein equations}. These are equations
 	\begin{gather*} \label{e0.16}  
                         \     \ \ \ hr  = D^{2}h, \tag{0.16}\\
 	         \Delta  h = 0, 
 	\end{gather*}
     \addtocounter{equation}{1}
\noindent  
for a pair $(g, h)$ consisting of metric $g$ and
  positive harmonic potential function $h,$ defined on
  an open 3-manifold $N.$ These equations have been
  extensively studied in classical general
  relativity. They imply that the 4-manifold $X = N\times_{h} S^{1},$ 
   with warped product metric $g_{\X} =
  g_{\N} + h^{2}d\theta^{2}$ is Ricci-flat, i.e.~a~vacuum 
 solution to the Einstein equations in
  dimension 4, c.f. \S 1.3.

\bigskip
 
  Of course, a flat metric, with $h$ an affine
  function, is a solution of (0.16). It is proved in the Appendix, c.f. also Theorems 3.2 and 3.3, that the only complete solution to the static vacuum equations (0.16) with potential $h > 0$ everywhere is a flat metric, with $h = const$. The most important or ``canonical''
  (non-trivial) solution to the static vacuum
  equations is the Schwarzschild metric
  $g_{s},$
 \begin{equation} 
       \label{e0.17} g_{s} = (1 - 2mt^{-1})^{-1}dt^{2} + t^{2}ds^{2}_{S^{2}},
 \end{equation}
  defined on $[2m,\infty )\times S^{2},$
   with $h = (1 - 2mt^{-1})^{1/2};$ 
   here the mass $m$ is a positive constant. 
   This metric is asymptotically flat, (for large t), 
   and the locus $\Sigma  = \{h = 0\} = \{t = 1\} $ 
   is a round totally geodesic 2-sphere, of
  radius $2m.$ Physically, this represents the surface
  of an (idealized) static black hole.

\bigskip
 
  Finally, given a sequence $\{g_{i}\}\in{\Bbb
  M}_{1},$ a sequence $r_{i}\in{\Bbb R} $ converging
  to 0, and points $x_{i}\in M,$ the blow-up sequence
  is the pointed sequence of Riemannian manifolds 
  $(M, \Roof{g}{\bar}_{i}, x_{i})$ with
  $\Roof{g}{\bar}_{i} = r_{i}^{-2}\cdot  g_{i}.$ For
  instance, in case the curvature tensor $R_{i}$ of $(M, g_{i})$ 
  is everywhere bounded by $r_{i}^{-2},$ the
  Cheeger-Gromov theory implies that a subsequence of
  the pointed sequence $(M, \Roof{g}{\bar}_{i},
  x_{i})$ either converges, modulo diffeomorphisms,
  uniformly on compact subsets, to a limit complete
  Riemannian manifold $(N, \bar{g}, x),$ or
  collapses along a sequence of $F$-structures to a
  lower dimensional space.

\bigskip

  A simplified version of the main result of this
  paper is the following:
 \begin{specialtheorem}
   Let $\{g_{i}\}$ be a sequence of unit volume
   Yamabe metrics on a closed oriented 3-manifold M,
   with $s_{g_{i}} \geq  -  s_{o} \geq  -\infty ,$
   satisfying the following condition: \\
  (i).	 (Non-Collapse). There is a constant $\nu_{o} \geq
 	   0$ such that 
        \begin{equation} \label{e0.18} 
        \vol B_{x}(r) \geq  \nu_{o}r^{3}, 
        \end{equation} 
      for any geodesic ball 
      $B_{x}(r) \subset  (M, g_{i}), r \leq  \diam(M, g_{i}).$\\
  \indent { \bf  (I).}  
 		 Then given any sequence of points $x_{i}\in (M, g_{i})$ 
 		 for which
 		 $\int_{B_{x_{i}}(r_{o})}|z_{g_{i}}|^{2}dV_{g_{i}}
 		 \rightarrow  \infty ,$ for some $r_{o} >  0,$
 		 there is a blow-up sequence $(M, g_{i}' , x_{i}),$
 		$$g_{i}'  = \rho_{i}^{-2}\cdot  g_{i}, \rho_{i}
 		\rightarrow  0, $$ such that a subsequence converges
 		to a locally defined solution of the static vacuum
 		Einstein equations.\\
\indent {\em\bf (II).}  
 		 In addition to (i), suppose the following
 		 conditions hold:   \\ 
  (ii). 
 	There is a constant $K< \infty ,$ 
         such that \begin{equation} \label{e0.19}
 	\int_{M}|z^{T}_{g_{i}}|^{2}dV_{g_{i}} \leq  K.
 	\end{equation} 
  (iii). 
       For $u$ as in (0.15), there is
 	a constant $\delta_{o} >  0$ such that
 	the sequence $\{g_{i}\}$ degenerates on the domain
 	$\displaystyle   U_{\delta_{o}} = \left\{x\in (M, g_{i}):
 	\frac{|u_{i}(x)|}{\sup|u_{i}|} \geq  \delta_{o}\right\},$
 	i.e.  
         \begin{equation} \label{e0.20}
  	       \int_{U_{\delta_{o}}}|z_{g_{i}}|^{2}dV_{g_{i}}
 	       \rightarrow  \infty .  \end{equation}
 
  Then there are points $y_{i}\in U_{\delta_{o}}$
  with $|z_{g_{i}}|(y_{i}) \rightarrow  \infty ,$ and
  a blow-up sequence $(M, g_{i}' , y_{i})$ defined as
  above, such that a subsequence converges to a
  non-trivial, (in particular non-flat), locally
  defined solution to the static vacuum Einstein
  equations.  
     \end{specialtheorem}
 
  We first make some comments on the hypotheses and
  conclusions. Condition (i) is used to rule out the
  possibility of collapse of the blow-up sequence.
  This condition will be weakened in \S 3, but some
  version of it is essential for Theorem A. Without
  any lower bound on the volumes of (small) geodesic
  balls, one cannot expect blow-ups to converge at
  all, see \S 4.1 for further discussion. Theorem
  A(I) gives then a general relationship between
  degenerations of Yamabe metrics and solutions,
  possibly trivial, of the static vacuum equations.
 
  Conditions (ii) and (iii) are required to obtain a
  {\it  non-trivial}  limit solution. For instance,
  in case $\sigma (M) < 0,$ condition (ii)
  prevents the function $u,$ which basically serves as
  the potential function $h$ in the static vacuum
  solution (0.16), from going to 0 in $L^{2}.$ The
  condition (ii) is not as strong an assumption as it
  might at first appear. For example, it is proved in
  Theorem 2.10 that if $\{g_{i}\}$ is an arbitrary
  sequence of unit volume Yamabe metrics on $M$ for
  which there exist points $x_{i}\in M$ and arbitrary
  numbers $r_{o} \geq  0,$ $R_{o}  <   \infty ,$
  such that \begin{equation} \label{e0.21}
 	\int_{B_{x_{i}}(r_{o})}
 	  |r_{g_{i}}|^{2} \leq R_{o}, \ \ \ \ \ 
             \vol B_{x_{i}}(r_{o}) \geq  R_{o}^{-1},
 	\end{equation} 
 then (0.19) holds, with $K$ depending
 only on $r_{o},\ R_{o}.$
 
  Note that (0.21) requires only that $\{g_{i}\}$
  have locally bounded $(L^{2,2})$ geometry in some
  ball of small but fixed size in $M,$ and yet it
  implies the global conclusion (0.19). This novel
  feature of some of the global aspects of Yamabe
  metrics plays a central role in this paper.
 
  Finally, the condition (iii) that $\{g_{i}\}$
  degenerates on $U_{\delta_{o}},$ for some
  $\delta_{o} >  0,$ and not just somewhere on 
  $(M, g_{i})$ as in (0.8), is also essential to obtain a
  non-trivial blow-up limit solution. Examples
  discussed in \S 6 will illustrate the necessity of
  this condition.
 
  A much more precise formulation of Theorem A will
  be proved in Theorem 3.10 below, where the base
  points $y_{i}$ and scale factors $\rho_{i}$ are
  constructed from the geometry of $(M, g_{i}).$
  Further, it is shown in \S 5 that the limit
  solutions are complete in a natural sense and some
  initial results on their asymptotic geometry are
  obtained.

\bigskip

  Theorem A implies that blow-up limits of
  degenerations in $U_{\delta_{o}}$ of sequences of
  Yamabe metrics satisfy a strong, in particular a
  determined, set of equations. This is quite
  surprising, since Yamabe metrics themselves satisfy
  no such rigid equations; the equation defining a
  Yamabe metric is highly underdetermined. Note also
  that Theorem A holds for quite general sequences of
  Yamabe metrics; for example no assumption is made
  that the sequence $\{g_{i}\}$ is a maximizing
  sequence for ${\cal S}|_{{\cal C}}.$
 
  It has long been known that static vacuum solutions
  on an open 3-manifold $N$ are closely tied with
  2-spheres in $N,$ occurring at the event horizon or
  boundary $\Sigma $ of a black hole in general
  relativity. This is already seen in the
  Schwarzschild metric (0.17), and is apparent in the
  work of Hawking; c.f. [HE, Ch.9.3] and \S 1.3 for
  further discussion. More generally, many solutions
  of the static vacuum equations are asymptotically
  flat, and thus have natural 2-spheres near
  infinity. These remarks illustrate a basic tie
  relating the degeneration of Yamabe metrics with
  the underlying topology of the manifold mentioned
  above.
 
  On the other hand, the first two examples discussed
  in \S 6 show that the restriction to
  $U_{\delta_{o}}$ in Theorem A is essential. If
  $\{g_{i}\}$ degenerates only on the complement of
  $U_{\delta_{o}},$ for any given $\delta_{o} >  0,$ then the blow-ups describing the degeneration
  may satisfy the static vacuum equations only
  trivially; see~\S~3.2. We also present examples in
  \S 6 of degenerating sequences which do satisfy the
  hypotheses of Theorem A(II), showing that this
  result is indeed applicable.
 
  It is certain global features of Yamabe
  metrics which lead to the necessity of condition
  (iii) in Theorem A. In \S 7, we consider the relation of this condition with the condition that a sequence of Yamabe metrics is Palais-Smale for the
  functional ${\cal S}|_{{\cal C}}.$
  (Recall that a sequence is Palais-Smale for a
  functional if the gradient tends to 0, in a
  suitable (weak) topology; this is not to be
  confused with Condition $C$ of Palais-Smale, which
  is much too strong to be of use here). In \S 4.2, it is shown that
  condition (iii) follows for sequences of Yamabe
  metrics satisfying a natural strengthened version
  of the definition of Palais-Smale sequence for
  ${\cal S}|_{{\cal C}},$ at least when 
  $\sigma (M)  <   0.$
 
  In any case, in order to effectively describe the
  structure of degenerations of sequences of Yamabe
  metrics on all of $M,$ (not only on
  $U_{\delta_{o}}),$ one needs to restrict to special
  classes of sequences. In addition to the reason
  above, this is also needed in case the sequence
  collapses, as discussed in \S 4.1. Some such
  preferred sequences are considered already in
  [An2, \S 8], and we will discuss these and related
  special sequences in more detail in [AnII].
 
  This paper is partly intended as the third in a
  series of works on the geometrization conjecture of
  Thurston, the previous works being [An1, 2].

\medskip

\noindent
 {\bf Acknowledgments:}
  I would like to thank Dennis Sullivan and Mike
  Freedman, who organized and led seminars in Stony
  Brook and La Jolla on some of the topics treated
  here and the sequel papers. Thanks also to Richard
  Hamilton, Steve Kerckhoff, Yair Minsky, Richard
  Schoen, Jon Wolfson, and other participants in
  these seminars for their comments, and to Jeff
  Cheeger for his abiding support.
 
  Many thanks also to Matt Gursky for pointing out to
  me the well-known, (but not to me at the time),
  splitting (0.13) during some conversations we had
  in Summer '94.
 
  I'm most grateful to Cliff Taubes for his
  careful reading of the manuscript and suggestions
  for its improvement, and for his continued support
  and interest in this program.  

  Finally, many thanks to Ramon del Castillo and especially
  my wife, Myong-hi Kim, for carrying out the huge task of 
  converting the original paper into a \TeX manuscript.
 \section{\bf Background Material.} 
\subsection{}
 \setcounter{equation}{0} 
  Throughout the paper, $M$ will denote a compact, connected, oriented 3-manifold, without boundary. Let ${\Bbb M} $ denote the space of smooth
  $(C^{\infty})$ Riemannian metrics on $M.$ This is an open cone, (the
  positive definite cone), in the linear space
  $S^{2}(M)$ of smooth symmetric bilinear forms. For
  any $g\in{\Bbb M} ,$ one has a natural
  identification \begin{equation} \label{e1.1}
 T_{g}{\Bbb M}  = S^{2}(M) , \end{equation} of the
 tangent space of ${\Bbb M} $ at $g.$ Let ${\cal M} $
 denote the space of isometry classes of all
 $C^{\infty}$ Riemannian metrics; thus ${\cal M} $ is
 the quotient of ${\Bbb M} $ by the action of the
 diffeomorphism group $\Diff(M)$ on ${\Bbb M} .$ It is
 well known, [Bes, Ch.4B], that the tangent space
 $T_{g}{\Bbb M} $ splits as 
 	\begin{equation}
 	\label{e1.2} T_{g}{\Bbb M}  = \Im \delta^{*} \oplus \Ker \delta , 
 	\end{equation} 
 where 
 $\delta: S^{2}(M) \rightarrow  \Omega^{1}(M)$ 
 is the divergence operator (w.r.t. $g$) on $S^{2}(M),$ given by 
 $\delta~(\alpha )=-(D_{e_{i}}\alpha) (e_{i},\cdot);$ here
 $\{e_{i}\}$ is an orthonormal basis of $TM.$ The
 operator 
 $\delta^{*}: \Omega^{1}(M) \rightarrow S^{2}(M)$ 
  is the formal adjoint of $\delta ,$ given by 
 $\delta^{*}(\omega ) = \frac{1}{2}{\cal L}_{\omega}g; {\cal L} $ 
 denotes the Lie derivative,
 and we are identifying vector fields and $1$-forms via
 the metric. In the splitting (1.2), the factor 
 $\Im \delta^{*}$ is the tangent space to the orbit of
 $\Diff(M)$ at $g$, while the factor $\Ker \delta $ is
 naturally identified with the tangent space of the
 quotient ${\cal M} ,$ at least when $g$ has no
 isometries.
 
  There is a natural (weak) Riemannian metric on
  ${\Bbb M} ,$ the $L^{2}$ metric, given for $\alpha
  , \beta  \in  T_{g}{\Bbb M} $ by
 	\begin{equation} \label{e1.3} 
         \langle\alpha , \beta \rangle~=\int_{M}\tr_{g}(\alpha\circ\beta ) dV_{g},
 	\end{equation}
  where $\alpha , \beta $ are considered as linear
  maps of $TM,$ via $g,$ and $dV_{g}$ denotes the volume
  form of $g.$ It is easily verified that this $L^{2}$
  metric is invariant under the action of $\Diff(M),$
  and thus passes to the quotient to give the $L^{2}$
  metric on ${\cal M} .$ Further, the splitting (1.2)
  is orthogonal w.r.t. the $L^{2}$ metric.
 
  The subsets of ${\Bbb M} $ and ${\cal M} $
  consisting of metrics of volume 1 on $M$ will be
  denoted by ${\Bbb M}_{1}$ and ${\cal M}_{1}.$ The
  discussion above is also valid for ${\Bbb M}_{1}$
  and ${\cal M}_{1}.$
 
  The (weak) $L^{2}$ Riemannian metric does not give
  a good (smooth) topology to ${\Bbb M} .$ For our
  purposes, it is natural to put the $L^{2,2}$
  topology on ${\Bbb M} ,$ i.e. the topology given by
  the $L^{2,2}$ Riemannian metric, defined on each
  $T_{g}{\Bbb M} $ by 
     \begin{equation} \label{e1.4}
       || h ||_{L^{2,2}} 
       = \int_{M}|h|^{2}+|Dh|^{2}+|D^{2}h|^{2}\ dV_{g}\,\strut^{1/2} , 
     \end{equation}
  for $h\in T_{g}{\Bbb M} .$ Here the norms $|\cdot |,$ 
  and covariant derivative $D$ are taken with
  respect to the metric $g\in{\Bbb M} .$ Thus, the
  $L^{2,2}$ distance between two metrics $g, g' $ is
  the infimum of the lengths of smooth curves joining
  $g$ and $g' ;$ the length of a curve $\gamma $ is
  defined by 
    \begin{equation} \label{e1.5} 
       L(\gamma ) = \int  \left|\left|\frac{d\gamma}{dt}\right|\right|_{L^{2,2}
                   (\gamma (t))}dt.  \end{equation} 
    This norm corresponds
 formally to the Sobolev space $L^{2,2}$ of functions
 with two weak derivatives in $L^{2}.$ In fact, it is
 not difficult to verify, c.f. [E, p.21], that the
 $L^{2,2}$ metric topology above induces the same
 topology as the $L^{2,2}$ topology defined by local
 coordinates, i.e.
 \begin{equation} \label{e1.6} 
        \dist_{L^{2,2}}(g, g') <   \varepsilon  
        \ \Leftrightarrow \  |g_{ij} - g'_{ij}|_{L^{2,2}} 
           <   \varepsilon' , 
 \end{equation}
 where $| \cdot  |_{L^{2,2}}$ is the Sobolev topology on
 functions on bounded domains in ${\Bbb R}^{3},$ and
 the components $g_{ij}, g'_{ij}$ are taken with
 respect to an arbitrary fixed coordinate atlas of $M;$
 here $\varepsilon'  = \varepsilon' (\varepsilon )$ is small
 if $\varepsilon $ is small, and vice versa. Although
 the two topologies defined by (1.5) and (1.6) are
 the same, the two metrics are very different in the
 large on ${\Bbb M} .$
 
  The Sobolev embedding theorem for $L^{2,2}$
  functions in dimension 3 reads
 \begin{equation} \label{e1.7} L^{2,2} \subset
 L^{1,6} \subset  C^{1/2} ; \end{equation}
  this is understood to apply to functions of compact
  support say in the unit ball in ${\Bbb R}^{3}.$
  Here $C^{\alpha}$ denotes the space of
  H\"older continuous functions with
  H\"older exponent $\alpha,\ L^{p,q}$ the
  Sobolev space of functions with $p$ weak
  derivatives in $L^{q}.$ Thus, at any $g\in{\Bbb M} ,$ 
  one has an estimate of the form
 	\begin{equation} \label{e1.8}
 		 c_{s}(g)||\alpha||_{C^{1/2}(g)} \leq
 	||\alpha||_{L^{2,2}(g)}, 
 	\end{equation}
  for any 2-tensor $\alpha .$ The Sobolev constant
  $c_{s}$ depends strongly on the metric $g,$ so that
  the estimate (1.8) is not uniform on ${\Bbb M} .$
 
  It is useful to compare the $L^{2,2}$ metric with a
  somewhat weaker metric. Thus, define the $T^{2,2}$
  metric to be the Riemannian metric given on each
  $T_{g}{\Bbb M} $ by the norm
   \begin{equation} \label{e1.9} 
     ||h||_{T^{2,2}} = \int_{M}|h|^{2}+|Dh|^{2}+|D^{*}Dh|^{2}\ dV_{g}\,\strut^{1/2}, 
   \end{equation}
  where $D^{*}$ is the $L^{2}$ adjoint of $D;$
  $(D^{*}D$ is the so-called rough Laplacian of $g$).
  Since $D^{*}D$ is an elliptic operator, by elliptic
  regularity, there is a constant $C = C(g)$ such
  that
 	\begin{equation} \label{e1.10} 
          ||h||_{T^{2,2}} \leq ||h||_{L^{2,2}} \leq  C(g)\cdot ||h||_{T^{2,2}} .
 	\end{equation}
  Using (1.10), it can be shown that the topology
  defined by the $T^{2,2}$ metric is the same as the
  topology defined by the $L^{2,2}$ metric, c.f.
  again [E,p.21]. However the $L^{2,2}$ and $T^{2,2}$
  metrics are far from being uniformly equivalent;
  the constant $C,$ which by definition is the constant
  on which one has $L^{2}$ elliptic estimates, c.f.
  [GT, Thm.8.8], cannot be chosen independent of $g$.
 
  We return to the study of ${\Bbb M} $ in more
  detail in \S 7.

\subsection{}
  The solution to the Yamabe problem by Yamabe [Y],
  Trudinger [Tr] Aubin [Au1] and Schoen [Sc1],
  implies that in each conformal class
  $[g]\subset{\Bbb M}_{1}$ of metrics there is a
  preferred metric, namely a {\it  Yamabe metric} ,
  i.e. a metric of constant scalar curvature which
  minimizes the total scalar curvature functional
  ${\cal S} $ (0.1) restricted to the conformal class
  $[g]\cap{\Bbb M}_{1}.$
 
  We let ${\cal C} $ denote the subset of ${\Bbb M} $
  consisting of all Yamabe metrics, and ${\cal
  C}_{1}$ the subset of unit volume Yamabe metrics;
  ${\cal C} $ is obtained from ${\cal C}_{1}$ by
  rescaling. If $g$ is a fixed background metric in
  $[g]\cap{\Bbb M}_{1}$ and $\Roof{g}{\bar}$ is a
  Yamabe metric in $[g]\cap{\Bbb M}_{1},$ then one
  may write 
	 \begin{equation} \label{e1.11}
	 \Roof{g}{\bar} = \psi^{4}\cdot  g, 
	 \end{equation}
 where $\psi $ is a smooth positive function on $M,$
 with $L^{6}$ norm 1, w.r.t. $dV_{g},$ corresponding
 to $\vol_{\Roof{g}{\bar}}M = \vol_{g}M.$ The equation
 that $\Roof{g}{\bar}$ have constant scalar curvature
 $\bar{s} = \inf {\cal S}|_{[g]}$ is
 	\begin{equation} \label{e1.12} \psi^{5}\cdot
 	\bar{s} = - 8\Delta\psi  + \psi\cdot s.  
 	\end{equation}
 
  In general, although there always exists at least
  one Yamabe metric, it is unknown to what extent
  they are unique in their conformal class. It is
  well-known, and follows easily from the maximum
  principle applied to (1.12), that Yamabe metrics
  $\Roof{g}{\bar}\in{\cal C}_{1}$ are unique in their
  conformal classes, when the scalar curvature
  $\bar{s}$ of $\Roof{g}{\bar},$ necessarily
  constant, is non-positive. In fact, in this case,
  $\Roof{g}{\bar}$ is the unique metric of constant
  scalar curvature in $[g]\cap{\Bbb M}_{1}.$ c.f.
  [Sc2].
 
  Thus, in case $\sigma (M) \leq  0,$ any unit volume
  Yamabe metric has non-positive scalar curvature and
  so is unique. Further, it is well-known, c.f.
  [Bes, Cor.4.49], [Sc2] that if $\sigma (M) \leq  0,$
  any metric of zero scalar curvature is flat, so
  that, with the exception of flat 3-manifolds, any
  Yamabe metric on $M$ has negative scalar curvature.
  Conversely, if $M$ is a flat 3-manifold, then $M$
  admits no metrics of positive scalar curvature,
  c.f. [Sc2], so that $\sigma (M) = 0.$ Thus, on a
  flat 3-manifold, any flat metric realizes the Sigma
  constant (0.5).
 
  On the other hand, if $\sigma (M) >  0,$ i.e.
  $M$ admits a metric $g$ of positive scalar
  curvature, then there may be many (necessarily)
  positive constant scalar curvature metrics in
  $[g]\cap{\Bbb M}_{1}.$
 
  We note that in case $\sigma (M) \leq  0,$ the
  space ${\cal C} $ of Yamabe metrics is an infinite
  dimensional submanifold of ${\Bbb M} .$ More
  generally, if $g\in{\cal C} ,$ and $- s_{g}/2$ is
  not in the spectrum of the Laplacian, then an
  $L^{2,2}$ neighborhood of $g$ in ${\cal C} $ is an
  infinite dimensional submanifold of ${\Bbb M} ;$
  c.f. [Bes, Ch.4F], [Ks].
 
\subsection{}
  As noted in the Introduction, the static vacuum
  Einstein equations
 	\begin{equation} \label{e1.13} 
         hr = D^{2}h, \ \ \ \ \  \Delta h = 0, \end{equation}
  will play a fundamental role in the study of the
  degeneration of Yamabe metrics. Here $D^{2}h$ is
  the Hessian of h, i.e. $D\nabla h,$ where $D$ is
  the covariant derivative and $\nabla $ the
  gradient. The Laplacian will always be defined as
  $\Delta h = \tr D^{2}h,$ (sum of second derivatives), 
  so that $\Delta $ has non-positive spectrum.
 
  These are equations for a pair $(g, h),$ defined on
  an open 3-manifold $N$. From the
  maximum principle, it is obvious that there are no
  (non-trivial) solutions of (1.13) on closed
  manifolds.
 
  If one considers the 4-manifold $X = N\times S^{1},$ with
  warped product metric $g_{\X} = g + h^{2}d\theta^{2},$ 
  then the equations for the Ricci
  curvature $r_{\X}$ on $(X, g_{\X})$ are
 	\begin{equation} \label{e1.14} 
 	r_{\X}(H,H) = r(H,H) -  \frac{1}{h}D^{2}h(H,H),  
 	\ \ \ \ \  r_{\X}(V, V) = -
 	\frac{1}{h}\cdot \Delta h\cdot |V|^{2},
 	\end{equation}
  where $H$ is tangent to $N,$ $V$ is tangent to
  $S^{1}.$ c.f. [Bes, 9.106]. Thus the static vacuum
  Einstein equations (1.13) are exactly the equations
  that the Ricci curvature vanish on $(X, g_{\X}).$
 
  Using the regularity theory of elliptic equations,
  c.f. [GT, Ch.8], it is quite standard to prove that
  $L^{2,2}$ weak solutions of (1.13) are
  $C^{\infty},$ in fact real-analytic, in harmonic
  coordinates, in any domain on which $h$ is bounded
  away from 0, (where the equations degenerate). Here
  one uses the fact that, to leading order, the Ricci
  curvature is given by the Laplacian of the metric
  in harmonic coordinates. Alternately, one may use
  regularity estimates for the Einstein metric (1.14)
  to prove smoothness of weak solutions, c.f.~[DK].
 
  Of course, a flat metric, with $h$ an affine
  function, gives a trivial solution to the static
  vacuum equations. The canonical solution is given
  by the Schwarzschild metric (0.17). One has the
  following remarkable `black-hole uniqueness
  theorem', proved by physicists.

 \begin{theorem} \label{t 1.1.}
    {\rm [I], [Ro], [BM].}
     Let $(N, g, h)$ be a solution to
    the static vacuum equations (1.13), complete up
    to the locus $\Sigma  = \{h=0\},$ in the sense
    that the metric completion $\Roof{N}{\bar}$ of
    $(N, g)$ is given by $N\cup\Sigma .$ Suppose
    further that the metric $g$ extends smoothly to
    $N\cup\Sigma ,$ and that $( N, g, h)$ is {\sf
    asymptotically flat}, so that, outside a large
    compact set in $\Roof{N}{\bar},$ the metric $g$
    is given in a chart by 
 	\begin{equation} \label{e1.15} 
 	       g_{ij} = \left( 
                         1+\frac{2m}{t}
                         \right)
                          \delta_{ij} + \gamma_{ij},  
 	       \ \ \ \ \ \ 
 	       h = 1 -  \frac{m}{t} + O(t^{-2});
 	\end{equation} 
 here, $t = |x|,\ |D^{k}\gamma_{ij}| = O(t^{- 2- k})$ and $m >  0.$
 
  If $\Sigma$ is a compact, possibly disconnected
  surface, then $(N, g, h)$ is the Schwarzschild metric of mass $m$
  for some $m>0.$  
 \end{theorem} 
  It is shown in [An4] that the assumption that $(N, g, h)$ 
   is asymptotically flat is usually
  superfluous. Except in some rather rare situations, (which may occur however), it is a consequence of the
  assumptions that $N\cup\Sigma $ is complete and
  $\Sigma $ is compact, (possibly singular).
 
  The locus $\Sigma  = \{h=0\}$ is called the {\it
  event horizon}  in general relativity, and plays a
  special role. From the equations (1.13), note that
  if $g$ is smooth up to $\Sigma ,$ then it follows
  immediately that $D^{2}h = 0 $ on $\Sigma$. In particular, $|\nabla h| = const. $ on each component of $\Sigma $, and it its straightforward to prove that $|\nabla h| > 0$ everywhere on $\Sigma$, c.f. [An4, Rmk.1.5]. Hence, each component of $\Sigma$ is a totally geodesic surface, with $|\nabla h| = const. > 0.$
 
  There are however many other solutions to the
  static vacuum equations (1.13). Usually, these are
  either singular at the event horizon $\Sigma ,$ or
  not complete away from $\Sigma ,$ or both. As
  indicated above the equations (1.13) may be viewed
  as an elliptic system, (in harmonic coordinates for
  example), away from $\Sigma ,$ but the equations
  formally degenerate at $\Sigma .$
 
  We point out though that any weak $L^{2,2}$ static vacuum solution on a domain $D$, with $L^{2,2}$ potential function $h$, is smooth, in fact real-analytic, in the interior of $D$, without any assumptions that $h$ is bounded away from $0$ in $D$. Hence, in this case, $\Sigma \cap D$ is a smooth, totally geodesic surface, and $g$ is smooth up to $\Sigma$, as above. Since no use will be made of this sharpening of the local regularity mentioned prior to Theorem 1.1, we do not include a proof, except to  note that it is a straightforward exercise in elliptic regularity techniques, on the associated Ricci-flat 4-manifold.

  A large and interesting class of explicit solutions
  of the static vacuum equations are given by the
  {\it Weyl solutions}, where $(N, g)$  itself is a
  warped product of the form 
  $$ N = V\times_{f}S^{1}, \ \ \ \ \ g = g_{\V} + f^{2}d\theta^{2}, $$ 
  and $(V, g_{\V})$ is a
 Riemannian surface, c.f. [EK, \S 2.3-9], [Kr, Ch.16-18]; 
 these metrics are discussed in much more
 detail in [An4]. For metrics of this form, the
 function $$r = f\cdot  h $$ is harmonic on $(V, g_{\V}).$ 
 Let $z$ be the (locally defined) harmonic
 conjugate of $r$ on $V.$ Then the function
	 \begin{equation} \label{e1.16} 
                  \nu  = \log h
	 \end{equation} 
 is an axially symmetric harmonic
 function on a domain in ${\Bbb R}^{3},$ where ${\Bbb
 R}^{3}$ is given cylindrical coordinates
 $(r,z,\theta ), \theta\in [0,2\pi ].$ The set $I =
 \{\nu =-\infty\},$ (considered as
 $\cap_{n}\nu^{-1}(-\infty ,$ n)), is usually assumed
 to be non-empty, (c.f. Theorem 3.2(I)), and in most
 circumstances, (but not always), is a subset of the
 axis $A = \{r=0\}.$ Note that the event horizon
 $\Sigma $ corresponds to the locus $I$. The metric $g$
 is given in these cylindrical coordinates by
 \begin{equation} \label{e1.17} 
 g = h^{-2}(e^{2\lambda}(dr^{2}+dz^{2}) + r^{2}d\theta^{2}), 
 \end{equation} 
 where $\lambda $ is a solution to the integrability equations
 $$\lambda_{r} = r(\nu_{r}^{2}-\nu_{z}^{2}),\ \ \ \ \ 
 \lambda_{z} = 2r\nu_{r}\nu_{z}. $$

  Conversely, given any axially symmetric harmonic
  function $\nu $ on a domain in ${\Bbb R}^{3},$ the
  equations above determine $\lambda $ (locally) up
  to a constant and the metric (1.17) gives a
  solution to the static vacuum equations with
  $S^{1}$ symmetry. 

   Thus, Weyl solutions are
  completely determined locally by an axially
  symmetric harmonic function $\nu $ on a maximal
  domain $D$ in ${\Bbb R}^{3}.$
 
  There is a large variety of behaviors in the global
  geometry of Weyl solutions, c.f. [An4]. Probably
  the most interesting class of metrics are those for
  which the potential $\nu $ in (1.16) is a globally
  subharmonic function on ${\Bbb R}^{3}.$ For these,
  one may use the value distribution theory of
  subharmonic functions on domains in ${\Bbb R}^{3}$
  to analyse the geometry of the associated Weyl
  metric.
 
  If $\nu $ is subharmonic on ${\Bbb R}^{3}$ and
  bounded above, say $\sup \nu  = 0,$ then the Riesz
  representation theorem c.f. [Ha, Thms. 3.9, 3.20],
  implies that $\nu $ may be globally represented as
 	\begin{equation} \label{e1.18} 
 	       \nu (x) = -\int_{{\Bbb R}^{3}}
                          \frac{1}{|x-\xi|}d\mu_{\xi},
 	\end{equation} 
 where $d\mu_{\xi}$ is a positive
 Radon measure on ${\Bbb R}^{3},$ the {\it  Riesz
 measure}  of $\nu .$ For such axi-symmetric
 functions, one has the characterization $$\supp \, d\mu
 = \Roof{I}{\bar} \Leftrightarrow  \supp \, d\mu
 \subset  A. $$

  The Riesz measure $d\mu $ in (1.18) encodes all the
  geometric information on the structure of such Weyl
  solutions. The fact that the Riesz measure is
  positive implies that such solutions have positive
  mass, in the sense of general relativity. Note that
  $\nu  \rightarrow  \sup \nu  = 0$ as the distance
  to the event horizon goes to $\infty .$
 
  On the other hand, there are solutions with
  negative mass. One may just take the function $-\nu
  $ in (1.18) for example, so that the potential is
  superharmonic. Note that in this case, $u$ is
  unbounded {\it  above}  on $\supp d\mu ,$ while $u$
  tends to its infimum at large distances to $\supp\, d\mu .$ 
  Of course in this case, there is no event
  horizon.

\medskip

  It is worthwhile to list a few standard examples of
  Weyl solutions and their corresponding measures. In
  these examples, the metric is given globally by
  (1.17), and the mass is assumed positive.
 
  From the point of view of the Riesz measure,
  perhaps the simplest example is the measure $d\mu $
  given by a multiple of the Dirac measure at some
  point on $A,$ so that $\nu  = - m/t,\ t(x) = |x|,$
  is a multiple of the Green's function on ${\Bbb
  R}^{3}.$ This gives rise to the Curzon solution
  [Kr, (18.4)],
      $$g_{c} = 
          e^{2m/t} (e^{- m^{2}r^{2}/t^{4}}(dr^{2} 
                + dz^{2}) + r^{2}d\theta^{2} ). $$
  The event horizon $\Sigma $ corresponds formally to
  \{0\} and $u = e^{- m/t}.$ This metric has a
  complicated singularity at the origin.
 
  The Schwarzschild metric (0.17) is a Weyl metric,
  with measure $d\mu  = \frac{1}{2}dA,$ where $dA$ is
  the standard Lebesgue measure on a finite interval,
  say $[- m, m]$ in $A$. This gives $\nu $ of the form
     $$ 
 	    \nu  = \frac{1}{2} \ln 
    \frac{r_{\smallplus} + r_{\smallminus} - 2m}
         {r_{\smallplus} + r_{\smallminus} +2m},
               \ \ \ {\rm where} \ \ \ 
           r_{\scriptscriptstyle\pm}^{2} = r^{2} + (z \pm  m)^{2}. 
     $$
 
 \noindent
  As mentioned before, the event horizon $\Sigma $
  here is a smooth totally geodesic 2-sphere of
  radius $2m.$
 
  It is easy to see that a Weyl solution $(N, g)$
  generated by a potential $\nu $ as in (1.18) for
  which $\supp\, d\mu  = \Roof{I}{\bar}$ is a compact
  subset of the axis $A,$ is asymptotically flat, in
  the sense of (1.15). Further, the simplest or most
  natural surfaces enclosing any finite number of
  compact components of $\Roof{I}{\bar},$ and
  intersecting $A$ outside $\Roof{I}{\bar},$ are
  2-spheres in $N.$ Of course if $\supp\, d\mu\subset A$ is
  non-compact, then the solution cannot be
  asymptotically flat.
 
  Note that (sub)-harmonic functions of the form
  (1.18) form a convex cone. In particular, one thus
  has a natural linear superposition principle for
  Weyl solutions. For example, one may choose the
  measure $d\mu  = \frac{1}{2}dA$ on two, or any
  number of disjoint intervals $\{I_{j}\}$ on the
  axis $A,$ provided the integral in (1.18) is finite.
  These correspond to solutions with `multiple black
  holes'. Although such solutions are essentially
  smooth up to the axis $A,$ they do not define smooth
  Weyl solutions $(N, g).$ There are cone
  singularities, (called struts in the physics
  literature), along geodesics (corresponding to
    $A\setminus \{I_{j} \} )$ joining the 2-spheres of $\Sigma ,$
  so that the metric $g$ is not locally Euclidean
  along such curves. Nevertheless, the curvature of
  such metrics is uniformly bounded everywhere. Of
  course, the black hole uniqueness theorem, Theorem
  1.1, also implies that such solutions cannot be
  smooth everywhere, when the number of intervals is
  finite.
 
  It seems to be unknown whether there are any smooth
  Schwarzschild type metrics with infinitely many
  black holes, although it is natural to conjecture
  that such solutions do not exist, c.f. [An4].
 
  On the other hand, there are Weyl solutions, c.f.
  [Sz], which are smooth and complete everywhere away
  from the event horizon $\Sigma ,$ for which $\Sigma
  ,$ (corresponding to $I \subset $ A), consists of
  any number, including infinity, of components; here
  the Weyl metric is highly singular at the event
  horizon.

\smallskip

  It is often useful in analysing the behavior of
  static vacuum solutions to consider the conformally
  equivalent metric 
 \begin{equation} \label{e1.19}
 \Roof{g}{\tilde} = u^{2}\cdot  g \end{equation} on
 $N,$ c.f. [EK, \S 2-3.5] and compare with (1.17). An
 easy calculation, c.f. [Bes, p.59] shows that the
 Ricci curvature $\Roof{r}{\tilde}$ of
 $\Roof{g}{\tilde}$ is given by 
	 \begin{equation}
	 \label{e1.20} \Roof{r}{\tilde} = 2 (d\log u)^{2} \geq  0, 
	 \end{equation} \hspace*{-2mm}
 and 
	 \begin{equation}
\hspace*{-2mm}
	 \label{e1.21} {\widetilde \Delta}\log u = 0.
	 \end{equation}
 
  This allows one to use methods and results on
  spaces with non-negative Ricci curvature and the
  behavior of harmonic functions on such spaces. Note
  that the metric $\Roof{g}{\tilde}$ is necessarily
  singular at the event horizon, even if the event
  horizon $\Sigma $ is smooth in $(N, g).$

\subsection{}
  We briefly summarize some of the main aspects of
  the theory of {\it convergence/degeneration} of metrics
  under uniform curvature bounds, but refer to the
  primary sources for further details.
 
  Let $V_{i}$ be a sequence of (possibly open)
  manifolds and let $\{\gamma_{i}\}$ be a sequence of
  smooth Riemannian metrics on $V_{i}.$ The sequence
  $(V_{i}, \gamma_{i})$ is said to converge to $(V, \gamma )$ 
  in the $C^{1,\alpha}$ topology, if first
  $V_{i}$ is diffeomorphic to $V,$ for all $i$
  sufficiently large, and second there exist
  diffeomorphisms $\phi_{i}: V \rightarrow  V_{i},$
  such that the pull-back metrics
  $\phi_{i}^{*}\gamma_{i}$ converge to the metric
  $\gamma $ in the $C^{1,\alpha}$ topology on $V.$ This
  means that there is a smooth coordinate atlas on
  $V$ for which the component functions of
  $\{\phi_{i}^{*}\gamma_{i}\}$ converge to the
  component functions of $\gamma ;$ here the
  convergence is with respect to the usual
  $C^{1,\alpha}$ topology for functions defined on
  domains in ${\Bbb R}^{3}.$ Note that this notion is
  well-defined for metrics $\gamma $ which are not
  necessarily $C^{\infty};$ for instance, it may well
  be that the limit metric $\gamma $ is only a
  $C^{1,\alpha}$ metric on $V.$
 
  In exactly the same way one defines convergence
  with respect to other topologies, for instance the
  weak or strong $L^{2,2}$ topology. It follows from
  the Sobolev embedding theorem (1.7) that
  convergence in the strong $L^{2,2}$ topology
  implies convergence in the $L^{1,6}$ and $C^{1/2}$
  topologies. Further, convergence in the weak
  $L^{2,2}$ topology implies convergence in the
  $L^{1,p}$ and $C^{\alpha}$ topologies, for $p < $ 6 
  and $\alpha  <  \frac{1}{2}.$

\bigskip

  We require the following definitions for a Riemannian
  3-manifold $(M,g)$, c.f. [An2, \S 3].
 
\begin{definition}
  {\bf (I)} The $\mu ${\it -volume radius}  
   at $x$ is given by
 	\begin{equation} \label{e1.22} 
 	\nu (x) = \sup  \left\{
 		   r: \frac{ \vol(B_{y}(s)\cap B_{x}(r))}{s^{3}}
 	     \geq \mu\omega_{3}, 
           \ {\rm for \ all\ } y\in B_{x}(r), \ s\leq r,
 		       \right\} 
 	\end{equation}
  where $\omega_{3}$ is the volume of the unit ball
  in ${\Bbb R}^{3}.$ The parameter $\mu $ is chosen
  to be an arbitrary but fixed small number, e.g.
  $10^{-2},$ so that we will suppress the parameter
  $\mu .$
\end{definition}
 
  We note that trivially one has the bound
 \begin{equation} \label{e1.23} 
   \nu (x) \leq  \left(\mu\omega_{3})^{-1}\cdot  \vol M\right)^{1/3}.
 \end{equation}
 {\bf (II).}
   The $L^{2,2}$ {\it harmonic radius}  $r_{h}(x),$ at $x$ 
   is the radius of the largest geodesic ball about $x$
   in which there exist harmonic coordinates $u_{i}:
   B_{x}(r_{h}(x)) \rightarrow  {\Bbb R} ,$ with
   respect to which the metric components $g_{ij}$
   satisfy
 	\begin{equation} \label{e1.24} 
 	    e^{-C}\delta_{ij} \leq  g_{ij} \leq  e^{C}\delta_{ij},   
 	   {\ \ \ \ \ \ \rm as \ bilinear \ forms,}
 	\end{equation} 
  and 
 	\begin{equation}
 		\label{e1.25}
 		r_{h}^{1/2}||\partial^{2}g_{ij}||_{L^{2}(B_{x}(r_{h}(x)))}
 		\leq  C.  
        \end{equation}
 {\bf (III).}
  The $L^{2}$ {\it  curvature radius}  $\rho (x)$ is
  the radius of the largest geodesic ball at $x$ such
  that for $y\in B_{x}(\rho (x)),$ and $D_{y}(s) =
  B_{x}(\rho (x))\cap B_{y}(s),$ one has the bound
 	\begin{equation} \label{e1.26}
 	\frac{s^{4}}{\vol D_{y}(s)}   
                \int_{\D_{y}(s)}|r|^{2} \leq  c_{o}, 
         \end{equation}
 for any $s \leq  \rho (x).$ 
 Note that the left-hand side of (1.26) is not necessarily a
 monotonic function of $s$, for a fixed $y$, thus the
 need to vary the center point, as in (1.22). Define
 the harmonic radius $r_{h},$ or $r_{h}(N)$ of $(N, g)$
 by $r_{h}(N) = \inf_{x\in N}r_{h}(x),$ and similarly
 for the $L^{2}$ curvature radius $\rho (N)$ of $N.$ If
 $N$ is a complete flat manifold, note that $\rho (N) = \infty .$
 
  The constant $C$ in (1.24)--(1.25) is an arbitrary
  but fixed parameter that may be taken to be 1.
  Similarly, $c_{o}$ is a free parameter, but will be
  chosen to be a fixed sufficiently small number, say $10^{-3},$
  throughout the paper. Note that the bounds
  (1.24)-(1.26) are invariant under rescaling of the
  metric (one also rescales the coordinate functions
  $u_{i}).$ In particular, these radii and the volume
  radius scale, i.e. behave under rescalings of the
  metric, as distance functions do. Observe from the
  definition that if $y\in B_{x}(\rho (x)),$ then
	 \begin{equation} \label{e1.27} 
		\rho (y) \geq \dist(y, \partial B_{x}(\rho (x))).  
	 \end{equation}
 The same estimate holds for $r_{h}.$ Thus, $\rho $
 and $r_{h}$ are Lipschitz functions, with Lipschitz
 constant 1. By the Sobolev embedding theorem,
 $L^{2,2} \subset  C^{1/2},$ so that one also has
 $C^{1/2}$ control of the metric components on
 $B_{x}(r_{h}(x)).$
 
  It is important to note that the $L^{2,2}$ harmonic
  radius and the $L^{2}$ curvature radius are
  continuous with respect to convergence in the
  strong $L^{2,2}$ topology. Thus, if 
  $(M_{i}, g_{i}) \rightarrow  (M, g)$ in the 
  $L^{2,2}$ topology,
  with $x_{i} \rightarrow  x,$ then 
	 \begin{equation} \label{e1.28} 
		 \lim_{i\to\infty} r_{h}(x_{i},
		g_{i}) = r_{h}(x, g), \ \ \  {\rm and}\ \ \ 
		\lim_{i\to\infty} \rho_{g_{i}}(x_{i}) =
		\rho_{g}(x).  
	 \end{equation}
  This is more or less obvious for $\rho $ and is
  proved for $r_{h}$ in [An3, Lemma 2.2], c.f. also
  [AC,Prop.1.1]. However, (1.28) is not true if the
  convergence is only in the weak $L^{2,2}$
  topology.
 
  There is an obvious relation between $r_{h}$ and
  $\rho ,$ namely
 	\begin{equation} \label{e1.29} \rho (x) \geq c\cdot  r_{h}(x), 
 	\end{equation}
  where $c = c(C, c_{o}).$ To see this, we may
  assume (by rescaling) that $r_{h}(x) = 1.$ The
  bound (1.24) then implies upper and lower bounds on
  $\vol B_{x}(1),$ while the bound (1.25) provides a
  bound on the $L^{2}$ norm of curvature on
  $B_{x}(1).$ This shows $\rho (x) \geq  c,$ as
  required.
 
  On the other hand, the opposite inequality, $\rho
  (x) \leq  c\cdot  r_{h}(x),$ is not true, as seen
  for example on compact flat manifolds, or more
  generally on manifolds which are highly collapsed
  on the scale of their curvature, that is for which
  $\nu< < \rho .$ On the other hand, under the
  presence of a lower bound on $\nu (x),$ one does
  obtain a bound of the form $\rho (x) \leq  c\cdot
  r_{h}(x),$ provided $\rho (x)$ is not too large,
  i.e. $\rho (x) \leq  K,$ where $c = c(\nu (x), K)$, c.f. [An2, \S 3],[An3].
 
  The natural bounds on sequences of Yamabe metrics
  that we obtain are bounds on the $L^{2}$ norm of
  (components of) the curvature. The following
  results, proved in [An2, \S 3], summarize the
  behavior of sequences of metrics on 3-manifolds,
  having a uniform $L^{2}$ bound on curvature. These
  are extensions of the fundamental $(L^{\infty})$
  Cheeger-Gromov theory of convergence and collapse
  of Riemannian manifolds, see [C], [G], [CG1,2].
 
 \begin{theorem} \label{t 1.3.}
   Let $\{g_{i}\}$ be a sequence of metrics in ${\Bbb
   M}_{1}(M)$, where $M$ is a closed 3-manifold. Suppose there is a uniform bound
 	\begin{equation} \label{e1.30}
 	   \int_{M}|r_{g_{i}}|^{2}dV_{i} \leq  \Lambda  < \infty .  
 	\end{equation}
  Then there is a subsequence, also called
  $\{g_{i}\},$ and diffeomorphisms $\psi_{i}$ of $M$
  such that exactly one of the following occurs:
 \begin{enumerate} \renewcommand{\labelenumi}{\upshape(\bf\Roman{enumi})}
   \item
  (Convergence) The metrics $\psi_{i}^{*}g_{i}$
  converge in the weak $L^{2,2}$ topology to a
  $L^{2,2}$ metric $g_{o}$ on $M.$
   \item  
  (Collapse) The metrics $\psi_{i}^{*}g_{i}$
  collapse $M$ along a sequence of orbit structures
  of a graph manifold structure. In particular, $M$
  is necessarily a graph manifold. The metrics
  $\psi_{i}^{*}g_{i}$ collapse each orbit ${\cal
  O}_{x},$ (namely a circle or torus), of a sequence
  of orbit structures to a point, as $i \rightarrow
  \infty ,$ i.e. $\diam_{\psi_{i}^{*}g_{i}}{\cal
  O}_{x} \rightarrow  0,$ $\forall x\in M.$
 \item  
   (Cusps) There is a maximal open domain $\Omega ,$ 
   (not necessarily connected), such that the
  metrics $\psi_{i}^{*}g_{i}$ converge, uniformly on
  compact subsets in the weak $L^{2,2}$ topology, to
  an $L^{2,2}$ metric $g_{o}$ defined on $\Omega ,$
  of volume $\leq  1.$ Any smooth compact domain $K
  \subset  \Omega $ embeds in $M$ and for
  sufficiently large $K \subset  \Omega ,$ the
  complement $ M\setminus  K $  has the structure of a graph
  manifold, part of which is collapsed along a
  sequence of orbit structures as in II.
  \end{enumerate}
 \end{theorem}

  Graph manifolds are 3-manifolds which are unions of
  Seifert fibered spaces glued along toral boundary
  components; thus, they have naturally defined
  embedded circles or tori, see [CG1], [An1,2]. From
  the Sobolev embedding theorem, the convergence in
  cases I and III above is also in the $L^{1,p}$
  and $C^{\alpha}$ topologies, $p < $ 6, $\alpha
  <  \frac{1}{2}.$

  The distinction between these three cases is
  determined by the behavior of the volume radius
  $\nu_{i}.$ Case I occurs if $\nu_{i}$ is
  uniformly bounded below, $\nu_{i} \geq  \nu_{o}
  \geq  0,$ Case II occurs if $\nu_{i}\rightarrow 0$
  everywhere on $M,$ while in Case III, there are
  regions in $(M, g_{i})$ where $\nu_{i}$ is bounded
  below and regions where it goes to 0.
 
  The following elementary result shows that $\{(M, g_{i})\}$ 
  degenerates in the sense of (0.8) only if 
  $\rho_{i}(x_{i}) \rightarrow  0,$ 
  for some    $x_{i}\in M.$
 
 \begin{lemma} \label{l 1.4.}
   Let $g$ be a unit volume metric on $M$ for which
   there is a uniform lower bound on the $L^{2}$
   curvature radius $\rho (M),$ say 
 	\begin{equation}
 		\label{e1.31} \rho (M) \geq  \rho_{o}.
 	\end{equation} 
 Then there is an explicit constant
 $\Lambda  = \Lambda (\rho_{o})$ 
   such that
 	\begin{equation} \label{e1.32} 
 	       \int_{M}|z|^{2}dV \leq  \Lambda .  
 	\end{equation}
 \end{lemma} 
\begin{pf}
  From the definition of $\rho ,$ we have
 $$
 \int_{B_{x}(\rho_{o})}|z|^{2}dV 
 \leq
 c_{o}\frac{\vol B_{x}(\rho_{o})}{\rho_{o}^{4}}, 
 $$
 for any geodesic ball 
 $B_{x}(\rho_{o})\subset (M, g).$ 
 Choose a maximal family of disjoint balls
 $B_{x_{k}}, k = 1,...,m$ of radius $\rho_{o}/4$ in
 $(M, g);$ thus the corresponding balls of radius
 $\rho_{o}/2$ cover $M.$ Now observe that there is a
 uniform upper bound $N$ on the multiplicity of this
 covering, independent of $\rho_{o}, M.$ In fact,
 since $c_{o}$ in (1.26) is sufficiently small,
 $(c_{o} = 10^{-3}$ suffices), the geometry of
 $B(\rho_{o}/2)$ is very close to that of a Euclidean
 ball, (or its quotient by Euclidean isometries), so
 that the Besicovitch covering theorem holds, c.f.
 [M, Theorem 2.7]. Thus
 $$\int_{M}|z|^{2}dV 
    \leq  
      N\sum \int_{B_{x_{k}}({\rho_{o}\over2})}|z|^{2}dV 
    \leq 
      \frac{Nc_{o}}{\rho_{o}^{4}}\sum \vol B_{x_{k}}
                   \left({\rho_{o}\over 2} \right)
    \leq
       C\frac{Nc_{o}}{\rho_{o}^{4}}\sum \vol B_{x_{k}}
                   \left({\rho_{o}\over 4} \right) 
    \leq
 	C\frac{Nc_{o}}{\rho_{o}^{4}} \vol M.$$\\
\end{pf}

  There are also local versions of Theorem 1.3 which
  will be frequently used below. While there are
  related results which hold in the collapse case,
  c.f. [An2, \S 3], we consider only the
  non-collapsing situation here; 
  c.f. [An2, Remark 3.6] for the following result.
 
 \begin{theorem} \label{t 1.5. }
   Let $(U_{i}, g_{i}, x_{i})$ be a pointed sequence of smooth
   Riemannian 3-manifolds such that
 $$
 	\rho_{i}(x_{i}) \geq  
         \rho_{o}  >  
         0,
         \ \ \ \ \ 
 	\nu_{i}(x_{i}) \geq  \nu_{o}  >  0, 
         \ \ \ \ \ 
         \diam U_{i}
 	\leq  D <  \infty , $$ 
 and 
 $$\dist(x_{i}, \partial U_{i}) \geq  \delta , $$
  for some arbitrary positive constants $\rho_{o},
  \nu_{o},$ D, $\delta .$ Then for any given
  $\varepsilon  >  0$ 
  sufficiently small, there are smooth domains 
  $V_{i}\subset U_{i}$ with 
  $\varepsilon /2 \leq  \dist(\partial V_{i}, \partial U_{i}) \leq \varepsilon ,$ 
  such that a subsequence of $(V_{i}, g_{i})$ 
  converges, modulo diffeomorphisms, to a limit 
  $L^{2,2}$ Riemannian manifold $(V, g_{o}).$
  The convergence is in the weak $L^{2,2}$ topology.
  In particular, the limit domain $V$ embeds in the
  domains $U_{i}.$ 
   \end{theorem} 
 
  In studying the degeneration of a sequence
  $\{g_{i}\}$ of Yamabe metrics on $M,$ we will {\it blow-up}  
  the metrics in neighborhoods of points
  where the curvature of $g_{i}$ goes to infinity,
  i.e. consider the behavior of the rescaled sequence
	 \begin{equation} \label{e1.33} 
	     g_{i}'  = \rho (x_{i})^{-2}\cdot  g_{i}, 
	 \end{equation}
  when $\rho (x_{i}) \rightarrow  0;$ of course, 
  $\rho (x_{i})$ is the $L^{2}$ curvature radius of $g_{i}$
  at $x_{i}.$ Theorem 1.5 will be used to examine the
  behavior of this sequence. In considering limits of
  $\{g_{i}'\},$ one must always consider based
  limits, i.e. limits w.r.t. a sequence of base
  points; the points $\{x_{i}\}$ will always be
  chosen to be the base points. Thus, assuming for
  instance that $\nu' (x) \geq  \nu_{o},$ for all
  $x\in B_{x_{i}}' (1),$ for some $\nu_{o} \geq  0,$
  the pointed sequence $\{(B_{x_{i}}' (1), g_{i}' ,
  x_{i})\}$ has a subsequence converging weakly in
  $L^{2,2}$ to a limit $(B' , g' ,$ x),\ x = lim
  $x_{i},$ with limit $L^{2,2}$ metric $g' .$ The
  convergence $g_{i}'  \rightarrow  g' $ is also
  understood to be in the pointed Gromov-Hausdorff
  topology [G, Ch.5A].
 
  Any smooth compact subdomain of $B' $ is naturally,
  but not canonically, embedded as a (very small)
  domain in $(M, g_{i})$ and the structure of the
  limit $(B' , g' )$ mirrors the very small scale
  behavior of $(M, g_{i})$ near $x_{i}$ as
  $i\rightarrow  \infty .$

\bigskip

  Finally, some general remarks. Since we are
  constantly dealing with the behavior of sequences
  $\{g_{i}\}, \{g_{i}'\}$ etc., we will often drop
  the subscript $i$, or prime, in order to simplify
  notation, when there is no danger of confusion. The
  main point is to establish uniform estimates,
  independent of $i$. Similarly, we will often pass to
  subsequences to obtain convergence, without always
  indicating the specific subsequence. A sequence
  $\{\alpha_{i}\}$ is said to sub-converge if a
  subsequence converges.  

 \section{\bf Initial Global Estimates for Yamabe Metrics}
 \setcounter{equation}{0}
  In this section, we derive a number of simple
  global relations on the behavior of the gradient of
  the (normalized) scalar curvature functional, i.e.
  the traceless Ricci curvature $- z,$ as well
  various components of z, on the space of Yamabe
  metrics ${\cal C} .$

  For a given smooth metric $g\in{\cal C} ,$ consider
  the operator $L: S^{2}(M) \rightarrow C^{\infty}(M)$ 
  giving the derivative or linearization of the scalar 
   curvature function at $g,$ i.e.
 	\begin{equation} \label{e2.1} L(\alpha ) = s'
 	(\alpha ) = \frac{d}{dt}s_{(g+t\alpha )}.
 	\end{equation}
  This operator has been classically studied, for
  instance in general relativity, and is given, c.f.
  [Bes, Ch.1K], by 
	   \begin{equation} \label{e2.2}
		L(\alpha ) = -\Delta \tr\alpha  +
		\delta\delta\alpha \  -  \langle r, \alpha \rangle, 
	   \end{equation}
  The $L^{2}$ adjoint of L, (w.r.t. the metric $g$),
  $L^{*}: C^{\infty}(M) \rightarrow  S^{2}(M),$ is
  the expression 
      \begin{equation} \label{e2.3}
 	L^{*}(h) = D^{2}h -  \Delta h\cdot  g - h\cdot  r.  
      \end{equation}
  This is an (overdetermined) elliptic operator, and
  thus by general elliptic theory there is a
  splitting of $T_{g}{\Bbb M} ,$ orthogonal w.r.t.
  the $L^{2}$ metric, of the form 
 	\begin{equation}
 	\label{e2.4} T_{g}{\Bbb M}  = \Im L^{*} \oplus  \Ker L.  
 	\end{equation} 
 To our knowledge, the splitting
 (2.4) first appeared in work of Berger-Ebin [BE],
 where it is attributed to Fadeev and Nirenberg, c.f.
 also [Bes, Ch.4F].
 
  Of course one sees immediately that
 \begin{equation} \label{e2.5} - r = L^{*}(1),
 \end{equation} corresponding to the fact that the
 $L^{2}$ gradient of the functional 
 $v^{-1}\int s\, dV,$ 
 (at a Yamabe metric), is given by $-r.$
 
  On the other hand, one may decompose the traceless
  Ricci tensor $z\ (=z_{g})$ with respect to this
  splitting and write 
 \begin{equation} \label{e2.6} 
           z = L^{*}f + \xi , \end{equation}
  where $\xi\in \Ker L$ and $f$ is a smooth function on
  $M.$ Clearly $\xi $ and $L^{*}f$ are uniquely
  determined by z. If 
 \begin{equation} \label{e2.7}
        \Ker L^{*} = \{0\}, 
 \end{equation} then $f$ is also
 uniquely determined by $z.$
 
  To examine $\Ker L^{*},$ suppose $h\in \Ker L^{*},$
  so that $D^{2}h -  \Delta h\cdot  g -
  h\cdot  r = 0.$ Taking the trace, one obtains
 	\begin{equation} \label{e2.8} 
   	       - 2 \Delta h - sh = 0, 
         \end{equation} 
 i.e. $h$ is an eigenfunction
 of the Laplacian, with eigenvalue $-\frac{s}{2}.$
 Suppose first that $\sigma (M) \leq  0.$ Then it
 follows, see \S 1.2, that $s  <   0,$ (unless $M$
 is a flat 3-manifold), and since the Laplacian has
 non-positive spectrum, the only solution of (2.8) is
 $h = 0.$ In other words, $\sigma (M) \leq  0$
 implies $\Ker L^{*} = 0,$ except in the special case
 that $(M, g)$ is a flat 3-manifold, where $\Ker L^{*}$
 consists of the constant functions. In case 
 $\sigma (M) >   0,$ \ $ \Ker L^{*}$ may well be non-zero; 
 see \S 3.6 and \S 6.4 for further discussion.
 
\noindent
   \begin{remark} It is worth emphasizing at this point
  that although $z$ is locally determined by $g,$
  neither $L^{*}f$ nor $\xi $ in (2.6) are locally
  determined by $g,$ (in contrast to the terms in
  (2.5)). They both depend on the global geometry of
  $(M, g),$ in particular on the global volume or scale
  of $(M, g).$ Thus they cannot be expressed locally in
  terms of the metric and its derivatives; see
  Theorem 2.10 for a simple concrete illustration of
  this.
 
  On the other hand, in analogy to (2.5), note that
  $z = H^{*}(- 1),$ where $H^{*}(f) = D^{2}f -
  \Delta f\cdot  g -  f\cdot  z$ is the
  $L^{2}$ adjoint of the operator $H,\ H(\alpha ) =
  -\Delta \tr \alpha  + \delta\delta\alpha\  -
  \langle  z, \alpha\rangle ,$ giving the derivative or
  linearization of the functional $v^{2/3}\cdot  s$
  on ${\cal C} .$ 
 \end{remark} 
 
  Returning to the general discussion, we may add the
  equations (2.5) and (2.6); this gives the $L^{2}$
  splitting for the metric $g,$ (considered as an
  element of $T_{g}{\Bbb M} ),$ i.e.
    \begin{equation} \label{e2.9} 
       -\frac{s}{3}\cdot  g = L^{*}(1+f) + \xi , 
    \end{equation} 
  so that, setting
   $u = 1+f$ gives 
  \begin{equation} \label{e2.10}
     L^{*}(u) + \xi  =-  \frac{s}{3}\cdot  g.
   \end{equation}
  Note that at least in the case $\sigma (M)\leq 0, $
   $g$ is never in $\Im L^{*}$ unless $(M, g)$ is
  Einstein. To see this, suppose there is an $h$ such
  that 
  $$ L^{*}h = g. 
  $$ 
 Taking the trace gives
 $$
    2\Delta h + sh = - 3. 
 $$
  It follows that $h+\frac{3}{s}$ is an eigenfunction
  of $\Delta,$ with eigenvalue
  $-\frac{s}{2},$ so that the argument above implies
  $(M, g)$ must be Einstein. The same reasoning shows
  that if $\sigma (M) \leq  0,$ then $z\in \Im L^{*}$
  only if $(M, g)$ is Einstein; this is also
  conjectured to be the case when $\sigma (M) > 0$, 
  c.f. [Bes, Remark 4.48].
 
  Thus, the pair $(f, \xi )$ measure in a certain
  global way the deviation from $g$ being an Einstein
  metric on $M.$
 
  Taking the $L^{2}$ norm of both sides of (2.10)
  gives the following simple but important estimate.

 \begin{theorem} \label{t 2.2.}
    For $f$ and $\xi $ defined as above, and 
   $u = 1+f,\ v = \vol M,$ one has
 \begin{equation} \label{e2.11} 
       \int|L^{*}u|^{2} + \int|\xi|^{2} = \frac{s^{2}}{3}\cdot   v.
 \end{equation}
 \end{theorem}
 \hfill $\blacksquare$\hspace*{4.5em}

 Thus, one has apriori bounds on the $L^{2}$ norms of
 $L^{*}(u)$ and $\xi .$ Referring back to (2.6), it
 follows that the component of $z$ in $\Ker L$ is
 apriori bounded in $L^{2}.$ Thus, $z$ is
 uncontrolled in $L^{2},$ say on a sequence
 $g_{i}\in  {\cal C}_{1}$ with scalar curvature
 bounded below, only in the direction $\Im L^{*}.$
 
  The splittings (2.6), (2.10) immediately give
  corresponding trace equations 
 	\begin{equation}
 	\label{e2.12} 2\Delta f + sf = \tr \xi  ,
 	\end{equation}
 	\begin{equation} \label{e2.13} 
           2\Delta u + su = \tr \xi  + s.  
         \end{equation}
  These equations easily lead to the following
  identities.
 \begin{lemma} \label{l 2.3.}
    For an arbitrary Yamabe metric $g \in  {\cal C} ,$ 
    one has the relations
 	\begin{equation} \label{e2.14} 
          \int|\xi|^{2} = - \frac{s}{3}\int \tr \xi  , 
 	\end{equation} and
 	\begin{equation} \label{e2.15} 
           \int \tr \xi  =  s\int f.  
 	\end{equation}
 \end{lemma} 
  \begin{pf}
   For $\xi\in \Ker L,$ note by (2.2) that $\int \langle r, \xi \rangle = 0.$ 
   Hence
 	$$
          \int \langle z, \xi \rangle  =\int \langle r, \xi \rangle -  \frac{s}{3}\int
 	  \langle g, \xi \rangle= -  \frac{s}{3}\int \tr \xi . 
        $$
 But also 
        $$
           \int \langle z, \xi \rangle   =\int \langle L^{*}f, \xi \rangle
           +\int|\xi|^{2} = \int|\xi|^{2}. 
        $$ 
  This gives (2.14).
 
  The equation (2.15) follows immediately from the
  trace equation (2.12) by integrating over $M$.\\
 \end{pf}

  Of course (2.14) and (2.15) combine to give
 \begin{equation} \label{e2.16} \int|\xi|^{2} = -
 \frac{s^{2}}{3}\int f .  \end{equation} Note that if
 $g\in{\cal C}_{1},$ then (2.11) implies that the
 $L^{2}$ norm of $\xi $ is bounded by
 $\frac{s^{2}}{3}.$

\bigskip

\begin{definition} 
  For $g\in{\cal C}\cap{\Bbb M}_{1},$ define the
  quantities $\delta  \in $ [0, 1] and $\lambda  \in
  $ [0, 1] by
 \begin{equation} \label{e2.17} \delta  = -
 \int_{M}f, \ \ \ \ \   \lambda  = 1 -  \delta  = \int_{M}u.
 \end{equation}
\end{definition} 
  The behavior of the quantities $\delta $ and
  $\lambda $ plays a fundamental role in the
  discussion. For $g\in{\cal C} $ not of unit volume,
  $\delta , \lambda $ are defined by the averages in
  (2.17); thus, they are invariant under scaling.

\bigskip

  One immediate consequence of this discussion is
  obtained by taking the $L^{2}$ inner product of (2.10) with $z,$
  to give the interesting bound
 	\begin{equation} \label{e2.18} 
 	\int u|z|^{2} = \int \langle\xi , z\rangle 
 		      = \int|\xi|^{2} 
                       \leq 
                       \frac{s^{2}}{3}\cdot  v.  
 	\end{equation}
 In particular, if one had an estimate of the form $u
 \geq  u_{o} >  0$ on $(M, g),$  then (2.18) gives an
 $L^{2}$ bound on $z;$ compare with \S 0.

\bigskip

  Another important $L^{2}$ orthogonal splitting of
  $T_{g}{\Bbb M} $ relevant to the study of the
  functional ${\cal S}|_{{\cal C}}$ is given by
 	\begin{equation} \label{e2.19} 
 	T_{g}{\Bbb M}  = T_{g}{\cal C}  \oplus  N_{g}{\cal C} ,
 	\end{equation} 
 where 
 	\begin{equation} \label{e2.20}
 		T_{g}{\cal C}  = \{\chi\in T_{g}{\Bbb M} : L(\chi )
 			       = const.\},    
 		       \ \ 
 		N_{g}{\cal C}  = \left\{\tau\in T_{g}{\Bbb M} : \tau  
 			       = L^{*}(h), \int_{M}h = 0\right\},
 	\end{equation}
  as in (0.10). At least in the case where $\sigma (M) \leq  0,$
  i.e. where Yamabe metrics are the unique metrics 
  of constant scalar curvature in their conformal class, 
  the spaces in (2.19) are the tangent and normal spaces 
  to the space of Yamabe metrics ${\cal C} .$ 
  In fact, using the splitting
  (2.19), one may show that ${\cal C} $ is an
  infinite dimensional submanifold of ${\Bbb M} ,$
  c.f. [Bes, Ch.4F], [Ks]. In case however 
  $ - {s/2} \in Spec 
             \left(   \Delta       \right),$ 
  the spaces $T_{g}{\cal C} $ and $N_{g}{\cal C} $ should only
  be considered as formal tangent spaces to ${\cal C}
  .$ (Recall that $\Delta $ has non-positive spectrum).
 
  Note that the splittings (2.4) and (2.19) differ by
  only 1-dimensional factors. We may then decompose
  $z$ also with respect to this splitting, and write
	 \begin{equation} \label{e2.21} 
		z = z^{T} + z^{N}.
	 \end{equation} 
 The component $-  z^{T},$ more
 precisely $- v^{-1/3}\cdot  z^{T},$ is the $L^{2}$
 gradient of the functional 
	 \begin{equation} \label{e2.22} 
	      v^{2/3}\cdot  s: {\cal C} \rightarrow  {\Bbb R} , 
	 \end{equation}
  (again assuming 
   $- s/2 \notin Spec \left(\Delta \right)$). 
  One sees this by
  observing that $- v^{-1/3}\cdot  z$ is the $L^{2}$
  gradient of the (volume normalized) total scalar
  curvature functional ${\cal S} : {\Bbb M}
  \rightarrow  {\Bbb R} $ in (0.1) and $z^{T}$ is the
  $L^{2}$ projection of $z$ onto ${\cal C} .$
 
  Now the space $T{\cal C} $ may be further
  decomposed into 
 \begin{equation} \label{e2.23}
 T{\cal C}  = T{\cal L}\oplus N{\cal L}  = \Ker L
 \oplus (\Im L^{*}\cap T{\cal C} ), \end{equation}
 where 
 \begin{equation} \label{e2.24} 
 T{\cal L}  = \{\xi\in T{\cal C} : L(\xi ) = 0\},  
   \ \ \ \  
 N{\cal L}  = \{\phi\in T{\cal C} : \phi =L^{*}(h),\  {\rm some}\ h\}.
 \end{equation} Clearly, the factor $N{\cal L}  =
 (\Im L^{*}\cap T{\cal C} )$ is 1-dimensional.
 
  Observe that $T{\cal L} $ corresponds to the
  tangent space to the level sets ${\cal L} $ of the
  functional $s = {\cal S}|_{{\cal C}}$ while $N{\cal
  L} $ corresponds to the normal space of ${\cal L} $
  in $T{\cal C} .$ The level sets ${\cal L} $ of $s$
  and the level sets ${\cal H} $ of $v^{2/3}\cdot  s$
  are obviously hypersurfaces in ${\cal C} $ at
  non-critical points of these functionals. These of
  course coincide on ${\cal C}\cap{\Bbb M}_{1},$ but
  do {\it  not}  coincide off ${\Bbb M}_{1}.$ Note
  also that by definition, $z^{T}$ is $L^{2}$
  orthogonal to the tangent spaces of ${\cal H} .$
 
  The splittings (2.4), (2.19) and (2.23) are all
  compatible with the splitting (1.2), since constant
  scalar curvature is a diffeomorphism invariant, see
  also [BE].
 
  In particular, the $L^{2}$ projection of z,
  $z^{T}\in  T{\cal C} $ can be split further as
 \begin{equation} \label{e2.25} z^{T} = L^{*}k + \xi
 , \ \ \ \ {\rm where}\  L^{*}k\in N{\cal L} , \  \xi\in T{\cal L} .
 \end{equation}
 
 \begin{lemma} \label{l 2.5.}
   The function $k$ is characterized by the basic
   property that
 \begin{equation} \label{e2.26} 
    LL^{*}k = const. = -\frac{1}{v}\int|z^{T}|^{2}.  
 \end{equation}
 \end{lemma} 
 \begin{pf}
  Apply the operator $L$ to both sides of (2.25) to
  obtain $$L(z^{T}) = LL^{*}k, $$ since $L(\xi ) = 0. $
  Since $z^{T}\in T{\cal C} , L(z^{T})$ is a
 constant function and
 $$\int L(z^{T})dV = \int \langle L^{*}1, z^{T}\rangle dV =
 -\int \langle r, z^{T}\rangle dV =  $$
 $$= -\int|z^{T}|^{2}dV -  \int \langle z^{T},z^{N}\rangle dV 
     -\frac{s}{3}\int \tr z^{T}dV. $$
  The components $z^{T}$ and $z^{N}$ are $L^{2}$
  orthogonal. Also $\int \tr z^{T}$ is the derivative
  of the functional ${\cal S}|_{{\cal C}}$ in the
  direction of the metric $g,$ i.e. its derivative
  under homothetic changes of the metric. Since the
  functional is scale-invariant, this derivative is
  0, and the result follows.\\
 \end{pf}
 
  A straightforward calculation shows that the
  operator $LL^{*}$ is given by 
 \begin{equation} \label{e2.27} 
         LL^{*}v = 2\Delta\Delta v + 2s\Delta v
                  -  \langle D^2v, r \rangle  + v|r|^{2}.  
 \end{equation}
  Clearly, when $\Ker L^{*} = \{0\},$ the space of
  functions $\phi $ such that $LL^{*}\phi  = c,$ for
  some constant c, is 1-dimensional. This gives the
  following relation between the functions $u$ and $k.$
 \begin{proposition} \label{p 2.6.}
    For a Yamabe metric $g\in{\cal C} ,$ with $s_{g}  <   0,$ 
   one has the relation 
          \begin{equation} \label{e2.28} 
 		 \frac{u}{\lambda} = - \frac{k}{\delta}.  
 	 \end{equation} 
  In general, the relation (2.28) holds $\mod \Ker L^{*}.$
 \end{proposition} 
\begin{pf}
  We have 
 	\begin{equation} \label{e2.29} 
 	  LL^{*}u = L\left(-\frac{s}{3}\cdot  g -  \xi \right) 
 		  = -\frac{s}{3}L(g) 
 		  = \frac{s^{2}}{3}; 
 	\end{equation} 
 the last equality just corresponds to the fact 
 that varying the scalar curvature $s$ in the direction 
 of $g,$ i.e. by a homothety, just changes $s$ by a constant. 
 Thus, by (2.26) and (2.29), both $LL^{*}k = \alpha $ and
 $LL^{*}u = \beta $ are constant. It follows that 
    $$ LL^{*}\left(\frac{k}{\alpha} -  \frac{u}{\beta}\right) = 0. $$
 If $ \Ker L^{*} = 0,$ we see that 
 $$ \frac{u}{\beta} =
      \frac{k}{\alpha}, 
 $$
  so that $\alpha $ and $\beta $ are the mean values
  of $k$ and $u$ respectively. From (2.6) and (2.25),
  we have \begin{equation} \label{e2.30} z^{N} =
 L^{*}(f- k), \end{equation} so that, in particular,
 since $z^{N}\in N{\cal C} ,$ \begin{equation}
 \label{e2.31} \int f = \int k = -  \delta .
 \end{equation} This implies (2.28).\\
\end{pf}

\bigskip

  Using (2.29), it is easy to see that the function
  $u$ minimizes the $L^{2}$ norm of $\Im L^{*}$ among
  functions with the same mean value, i.e.
 \begin{proposition} \label{p 2.7.}
   The function $u/\lambda ,$ (or $k/\delta ),$ is
   characterized uniquely by the fact that
 \begin{equation} 
       \label{e2.32} 
          \int|L^{*}(u/\lambda)|^{2}dV \leq  \int|L^{*}\phi|^{2}dV, 
       \end{equation}
  for all functions $\phi $ on $(M, g)$ with mean value 1.  
  \end{proposition} 
 \begin{pf}
 Write 
 $$
 \int|L^{*}\phi|^{2}dV = 
      \int|L^{*}(\phi  - u/\lambda )+L^{*}(u/\lambda )|^{2}dV  = $$
 $$=   \int|L^{*}(\phi  -  u/\lambda )|^{2}  dV
          +\int|L^{*}(u/\lambda )|^2 dV - 
  2\int \langle L^{*}(\phi  -  u/\lambda ), L^{*}(u/\lambda )\rangle dV. $$ 
 But 
    $$\int \langle L^{*}(\phi  -  u/\lambda ), L^{*}(u/\lambda )\rangle dV 
         =\int\langle (\phi  - u/\lambda ), LL^{*}(u/\lambda )\rangle dV = 0, $$
  where the last equality follows from the fact that
  $LL^{*}u$ is constant by (2.29), and $\phi $ and
  $u/\lambda $ have the same mean value.\\
\end{pf}
  The equation (2.28) contains the basic relation
  between the functions $f$ and $k$, relating $z$ and
  $z^{T}.$ It is useful to derive several further
  relations.
 \begin{lemma} \label{l 2.8.}
   The following identities hold for any Yamabe
   metric $(M, g)$:
	 \begin{equation} \label{e2.33} \xi  +
		 \frac{s}{3}\cdot  g = \lambda (z^{T}+
		 \frac{s}{3}\cdot  g), 
          \end{equation}
 \begin{equation} \label{e2.34} 
      \int|L^{*}k|^{2} =
      \frac{s^{2}}{3}\frac{\delta^{2}}{\lambda}\cdot  v,
 \end{equation}
	 \begin{equation} \label{e2.35} 
	 \int|z^{T}|^{2} = \frac{s^{2}}{3}\frac{\delta}{\lambda}\cdot  v.
	 \end{equation}
 \end{lemma} 
\begin{pf}
  First, recall that 
   $$
   L^{*}u = -\xi  -  \frac{s}{3}\cdot  g  \ \ \ \ \ {\rm and} \ \ \ \ \ L^{*}k
                = z^{T}-  \xi  . 
   $$
       From (2.28), one then obtains
 $$ \frac{\delta- 1}{\delta}z^{T} = \frac{\delta- 1}{\delta} \xi  
   -  \xi  -  \frac{s}{3}\cdot  g = -
    \frac{1}{\delta}\xi  -  \frac{s}{3}\cdot  g, 
 $$
  so that from Definition 2.4, 
	 \begin{equation}
	 \label{e2.36} \xi  + \frac{\delta\cdot  s}{3}\cdot
	 g = \lambda z^{T}.  
	 \end{equation} 
 which implies (2.33). 
 Thus the tensors 
 $\xi  + \frac{s}{3}g$ and
 $z^{T}+ \frac{s}{3}g$ are proportional. If $\xi_{o}$
 and $z_{o}^{T}$ denote the trace-free parts of $\xi
 $ and $z^{T},$ i.e. $\xi_{o} = \xi  -
 \frac{tr\xi}{3}\cdot  g,$ then one also has from
 (2.33) that $$\xi_{o} = \lambda\cdot  z_{o}^{T}. $$

  Next, from (2.28) and (2.29), we find
	 \begin{equation} 
	    \label{e2.37} LL^{*}k = -\frac{s^{2}}{3}\frac{\delta}{\lambda},
	 \end{equation}
  Multiply (2.37) by $k$ and integrate over $M,$ using
  (2.31) to obtain (2.34). Finally from (2.16), we
  have $\int|\xi|^{2} = \frac{s^{2}}{3}\delta\cdot
  v,$ so that (2.25) and (2.34) give (2.35).\\
\end{pf} 

  We now turn to the trace equation (2.12) in order
  to obtain estimates on $f,$ or $u.$
 \begin{proposition} \label{p 2.9.}
    For any unit volume Yamabe metric $g\in{\cal C}_{1},$ 
    with $s_{g}  <   0,$ one has the estimate 
 	   \begin{equation} \label{e2.38}
 		||f||_{T^{2,2}}  \leq  C, 
 		\end{equation} 
 where $C$ is a constant depending only on $|s_{g}|.$
 \end{proposition} 
\begin{pf}
  Recall from (1.9), that the $T^{2,2}$ norm is given by 
 $$                 
    ||f||_{T^{2,2}}^{2} =   ||f||_{L^{2}}^{2} 
                         + ||df||_{L^{2}}^{2} 
                         + ||\Delta f||_{L^{2}}^{2} . 
 $$
 Multiplying (2.12) by $f$ and integrating gives
   \begin{equation} \begin{array}{rl}
   \displaystyle
    2\int|df|^{2} -  s\int f^{2} 
                     & \displaystyle = -\int f\cdot \tr\xi   \\
                     & \displaystyle \leq  \tfrac{1}{2} |s|
                       \displaystyle \int f^{2} 
                     + \tfrac{1}{2} |s|^{-1}
                       \displaystyle \int (\tr\xi )^{2}, 
  \end{array} \end{equation} 
 so that, since $ s  <  0,$  
   \begin{equation} 
          2\int|df|^{2} 
            + \tfrac{1}{2}|s|\displaystyle \int f^{2}
      \leq   \displaystyle \tfrac{1}{2} |s|^{-1}
             \displaystyle \int (\tr\xi )^{2} 
              \leq \tfrac{1}{2}|s|\delta , 
   \end{equation} 
  where the last inequality follows from (2.16).
  Thus, one has an apriori bound 
 $$ ||f||_{L^{1,2}}^{2} \leq  C = C(|s|). $$
  Returning to the trace equation (2.12), it also
  then follows that 
 	 $$||\Delta f||_{L^{2}}^{2} \leq  C = C(|s|). $$
\end{pf} 

  We note that the estimate (2.38) is false in case
  $s_{g} >  0,$  see \S 3.6 for further discussion,
  while for $s_{g} = 0,$ it is borderline. Namely if
  $g$ is a Yamabe metric with $s_{g} = 0,$ then
  either $\sigma (M) >  0$ or $\sigma (M) = 0.$ In
  the former case, it is obvious that $f \equiv  - 1$
  and $\xi  = 0$ is one solution of the $z$-splitting
  equation (2.6) and hence it is the only solution
  since $\Ker L^{*} = 0$ when $s_{g} = 0.$ Thus, in
  this case we have $u \equiv  0.$ The splittings
  (2.6), (2.10) contain no essential information in
  this situation; they are equivalent to the identity
  (2.5). If however $\sigma (M) = 0,$ then $g$
  realizes the Sigma constant on $M,$ and thus $g$ is
  Einstein, with $r = 0.$ Thus, the solution to (2.6)
  is given by $f = const.$ and $\xi  = 0,$ where the
  constant for $f$ is however undetermined.
 
\bigskip

  The bound (2.38) implies that there is an apriori
  bound for $f$ in the $T^{2,2}$ norm, in case
  $\sigma (M)~\leq~0.$ On the other hand, we will
  see later in \S 3 that it is far from true that one
  has uniform bounds for $||D^{2}f||_{L^{2}}$ , i.e.
  uniform bounds for the $L^{2,2}$ norm of $f.$ This
  will be the case when one has non-flat blow-up
  limits. Thus, one has a breakdown of (uniform)
  elliptic regularity on $\{g_{i}\},$ c.f. (1.10).
 
  We refer to \S 4.2 for some further discussion on
  the $L^{2}$ behavior of $u$ or $f.$

\bigskip

  Next, consider the $L^{2}$ bound on $z^{T} = L^{*}k
  + \xi .$ From Theorem 2.2, $\xi $ is uniformly
  controlled in $L^{2},$ while $k$ is a function in
  the 1-dimensional space $\Im L^{*}\cap  T{\cal C} .$
  Thus, one would expect to be able to control
  $z^{T}$ in a natural way. This is given by the
  following general result, which illustrates in a
  simple manner the highly global nature of the
  $L^{2}$ projection operator onto $T{\cal C} .$

 \begin{theorem} \label{t 2.10.}
   Let $g$ be a unit volume Yamabe metric on $M.$
   Suppose there is a point $x\in M,$ and arbitrary
   but fixed constants $\rho_{o},\ \nu_{o} \geq  0$
   such that 
 	   \begin{equation} \label{e2.39} 
 	       \rho (x) \geq  \rho_{o},  \ \ \ \ \  \nu (x) \geq  \nu_{o}.
 	   \end{equation} 
   Then there is a constant $K = K(\rho_{o}, \nu_{o})$ such that 
 	\begin{equation} \label{e2.40} 
 	     \int_{M}|z^{T}|^{2}dV_{g} \leq  K.
 	\end{equation}
 \end{theorem} 
\begin{pf}
  Since $z^{T}\in T_{g}{\cal C} $ is the $L^{2}$
  orthogonal projection of $z$ onto $T{\cal C} ,$ we
  have \begin{equation} \label{e2.41}
 \int_{M}|z^{T}|^{2}dV_{g} = \inf_{\int\phi =0}
 \int_{M}|z -  L^{*}\phi|^{2}dV_{g}, 
\end{equation}
  where the infimum is over all smooth functions
  $\phi $ on $M$ with 0 mean value on $(M, g);$ recall
  $N{\cal C}  = \Im L^{*}\phi ,$ over mean value 0
  functions on $(M, g).$
 
  A straightforward computation gives
 \begin{equation} \label{e2.42}
    \int_{M}|z -
    L^{*}\phi|^{2}dV_{g}=\int_{M}\lbrace (1+\phi
    )^{2}|z|^{2}+|D^{2}\phi|^{2}+(\Delta\phi
    )^{2}-  \frac{2}{3}s|d\phi|^{2}+2z(d\phi ,d\phi
    )+\frac{s^{2}}{3}\phi^{2}\rbrace dV_{g}
 \end{equation}
 	$$= \int_{M}\lbrace (1+\phi )^{2}|z|^{2} -
 	|D^{2}\phi|^{2} + 3(\Delta\phi )^{2} -
 	\frac{4}{3}s|d\phi|^{2} +
 	\frac{s^{2}}{3}\phi^{2}\rbrace dV_{g}, 
 	$$
  where the last inequality follows from the use of the
  Bochner-Lichnerowicz formula to eliminate the term $z(d\phi, d\phi)$.
 
  Given $x\in M$ satisfying (2.41), choose a function
  $\phi ,$ with $\phi  \geq  - 1$ everywhere on $M,$
  $\phi  \equiv  - 1$ on $ M\setminus B,$  where 
  $B = B_{x}(\rho_{o}/2)$ 
   and with 0 mean value on $(M, g)$.
  Note that since $\vol B \geq  c\cdot
  \min((\rho_{o},\nu_{o}))^{3},$ so that $B$ has a
  definite proportion of the volume of $(M, g),$ one
  may choose such a $\phi $ so that $\sup\ \phi  \leq  H,$
  where $H$ depends only on $\rho_{o},\ \nu_{o}.$
  Now in the ball $B$ the geometry of $g$ is
  controlled in the $L^{2,2}$ topology, see the
  discussion concerning (1.29). It follows that one
  may choose $\phi $ so that the $L^{2,2}$ norm of
  $\phi $ is also uniformly controlled in $B$.
 
  For such a choice of $\phi ,$ it is clear that the
  expression on the right in (2.44) is uniformly
  controlled by $\rho_{o}$ and $\nu_{o}.$ The
  estimate (2.42) then follows immediately from
  (2.44).\\
\end{pf} 
 \begin{remark} \label{r 2.11.}
   In a related vein, suppose $\{g_{i}\}$ is a
   sequence of unit volume Yamabe metrics on $M,$ and
   $\{x_{i}\}$ is a sequence of points in $M$
   satisfying (2.39). Then the local estimate
		 \begin{equation} \label{e2.43}
			 \int_{B_{i}}|z^{T}|^{2}dV_{g_{i}} \rightarrow  0 
	         	 \ \ \ {\rm as}\ i \rightarrow  \infty , 
		 \end{equation} 
 implies the global estimate
 		\begin{equation} \label{e2.44}
 			\int_{M}     |z^{T}|^{2}dV_{g_{i}} \rightarrow  0 
 			\ \ \ {\rm as}  \  i \rightarrow  \infty , 
 		\end{equation}
  where $B_{i} = B_{x_{i}}(\rho_{o}).$ To see this,
  from (2.25) and (2.26), we have 
 	 $ L(z^{T}) = -||z^{T}||_{L^{2}}. $
 Multiplying this by a suitable cutoff function  
 supported in $B_{i}$  implies
 	$$  \vol_{g_{i}}B_{i}\cdot ||z^{T}||_{L^{2}}  
 	   \leq
 	       c\int \langle L^{*},  z^{T}\rangle
 	   \leq
 	       c\int_{B_{i}}|z^{T}|^{2}, 
 	$$
  which gives the result.  
 \end{remark} 
 \begin{remark} \label{r 2.12. }
  From (2.35), we see that in case $\sigma (M) \neq  0,$ 
  a bound on $||z^{T}||_{L^{2}}$ is equivalent to
  a bound on $\lambda $ away from 0. More precisely,
  if $\{g_{i}\}$ is a sequence of unit volume Yamabe
  metrics with scalar curvature bounded away from 0
  (and $-\infty ),$ then $||z^{T}||_{L^{2}}$ remains
  uniformly bounded exactly when $\lambda ,$ the mean
  value of $u$,  remains bounded away from 0.
 
  On the other hand, this is certainly not the case
  when $\sigma (M) = 0$ or more generally when
  $s_{g_{i}} \rightarrow  0.$ Suppose for instance
  $(M, g_{o})$ is a flat 3-manifold, so that $g_{o}$
  realizes $\sigma (M).$ Let $g_{t}$ be a smooth
  curve of unit volume Yamabe metrics on $M$ through
  $g_{o},$ and let $\alpha  = \frac{dg_{t}}{dt}$
  satisfy $||\alpha||_{L^{2}(g_{t})} = 1.$ We have
 	$$\frac{d}{dt}s_{g_{t}}  
           = \int_{M} \langle z^{T}, \alpha \rangle     \ 
            \leq \left( \int_{M}|z^{T}|^{2}\right)^{1/2}, 
 	$$ 
  so that 
 	$$\left(
           \frac{d}{dt}s_{g_{t}} \right)^{2} \ 
         \leq 
         \int_{M}|z_{g_{t}}^{T}|^{2}. 
         $$ 
 Since $s_{g_{t}}$ is a smooth function of $t$, 
 with $s_{g_{o}} = 0,$ the
 maximal value, it is clear that 
 $$- s_{g_{t}} <  <   -\frac{d}{dt}s_{g_{t}}, $$ 
 as $t \rightarrow  0.$ Thus, it follows that
 	\begin{equation} \label{e2.45}
 		\frac{1}{s^{2}}\int_{M}|z^{T}|^{2}dV \rightarrow
 		\infty  \ \ \ \ \ {\rm as}\  t \rightarrow  0.  
 	\end{equation}
 By (2.35), this implies that 
 \begin{equation} \label{e2.46} 
 	\lambda  \rightarrow  0,
 \end{equation}
  as $t \rightarrow  0.$ Further, by the proof of
  Proposition 2.8, $|\nabla u| \rightarrow  0$ in $L^{2}$
  and since $u$ converges smoothly to its limit here,
  we have 
     $$u \rightarrow  0 \ \ \ \ \ {\rm in}\  C^{o}. $$ 
  Note that in this example, one still has of course
 $||z^{T}||_{L^{2}} \leq  C, $  in fact
 $||z^{T}||_{L^{2}} \rightarrow  0$ as 
 $t \rightarrow  0. $ 
 \end{remark} 
 
 \section{\bf Existence of Non-Flat Blow-Ups.  }
 \setcounter{equation}{0}
  In this section, we will analyse the degeneration
  of sequences of Yamabe metrics using the structural
  results in \S 2. We recall (Lemma 1.4), that a
  sequence of unit volume Yamabe metrics $\{g_{i}\}$
  degenerates, in the sense that the curvature
  becomes unbounded in $L^{2},$ only if the $L^{2}$
  curvature radius $\rho_{i}(M)$ of $(M, g_{i})$ goes to $0.$
 
\subsection{} 
   As indicated in \S 0, the static vacuum Einstein equations
 \begin{equation} \label{e3.1} 
   ur = D^{2}u,\ \ \ \ \  \Delta u = 0, 
  \end{equation}
  play the fundamental role in the understanding of
  degenerations of sequences of Yamabe metrics. It is
  immediate from the definition that these equations
  are equivalent to the equation
 \begin{equation} 
       \label{e3.2} L^{*}u = 0,
 \end{equation}
  on scalar-flat manifolds. Obviously, there are no
  non-trivial solutions of the equations (3.1) on
  closed manifolds, although (3.1) or (3.2) may have
  locally defined solutions.
 
  We discuss briefly the relation between the
  equation (3.2) and Einstein metrics on 4-manifolds,
  since it does not seem to appear in the
  literature.
\setcounter{theorem}{-1} 
 \begin{proposition} \label{p 3.0.}
    Let $g$ be a metric of constant scalar curvature
    on a 3-dimensional domain $\Omega ,$ with smooth
    boundary $\partial\Omega .$ Then any solution of
    (3.2) on $\Omega $ with $u >  0 $ on $\Omega $
    gives an Einstein 4-manifold of the form
 \begin{equation} \label{e3.3} 
          X^{4} = \Omega \times_{u}S^{1}, 
 \end{equation} 
 with scalar curvature
 $s_{\X} = 2s_{g}.$
 \end{proposition} 
\begin{pf}
 Let  $u > 0$ be a solution to (3.2) in $\Omega.$  
  For any form $\eta \in  T_g {\Bbb M}$  of compact support in $\Omega ,$
  we then have
 $$
   \int u \cdot  L(\eta )\, d V_{g}  = \int \langle L^{*} u, \eta \rangle  dV_{g} = 0. 
 $$ 
 Now consider the warped
 product $X^{4} = \Omega  \times_{u}S^{1},$ with metric
 given by 
 \begin{equation} \label{e3.4} g_{\X} = g +
 u^{2}d\theta^{2}.  
 \end{equation}
  The volume form of $g_{\X}$ is given by
  $udV_{g}\wedge d\theta ,$ and one computes, c.f.
  [Bes, Ch.9J], that the scalar curvature $s_{\X}$ of
  $g_{\X}$ is given by \begin{equation} \label{e3.5}
 s_{\X} = s -  2\frac{\Delta u}{u}.
 \end{equation} 
 Thus, applying the divergence theorem
      $$\int_{X}s_{\X}dV_{X} 
         = \int_{\Omega}\left(s - 2\frac{\Delta u}u\right)u\ dV_{g} 
         = \int_{\Omega}u\cdot  s\, dV_{g} -
          2\int_{\partial\Omega}\langle du, \nu\rangle dA_{g}, $$
  where $\nu $ is the outward unit normal.
 
  Consider a compactly supported 1-parameter
  variation of $g_{\X},$ i.e. a curve of metrics
  $g_{\X}(t)$ on $X,$ with $g_{\X}(0) = g_{\X},$ of the
  form $g_{\X}(t) = g_{t} + u_{t}^{2}d\theta^{2}.$ We
  suppose also that $\vol_{g_{\X}(t)}(X) =
  \vol_{g_{\X}}(X).$ Then, for $\alpha  =
  \frac{d}{dt}g_{\X}(t)|_{t=0},$ we have at $t = 0,$
                    $$\frac{d}{dt}  \int_{X}s_{\X}(t)\  dV_{t} =
                    \frac{d}{dt}    \int_{\Omega}s_{t}u_{t}\  dV_{t} -
 2\frac{d}{dt}\int_{\partial\Omega} 
           \langle du_{t}, \nu \rangle_{g_{t}}\ dA_{g_{t}}. 
                   $$
  Now the boundary term vanishes, since the variation
  is of compact support. For the first integral, we have 
 $$ \frac{d}{dt}\int_{\Omega}s_{t}u_{t}dV_{t}
       =\int_{\Omega}s' udV +s\int_{\Omega}(udV)' , $$
  where we have used the fact that $s$ is constant.
  The second integral here vanishes, since the
  variation is volume preserving, while for the first
  integral 
          \begin{equation} \label{e3.6}
 	    \int_{\Omega}s' udV  = \int_{\Omega}L(\alpha )udV =
 	    \int_{\Omega}\langle L^{*}u, \alpha\rangle dV = 0.
 	\end{equation}
  Thus the gradient $\nabla{\cal S} $ of the
  unnormalized scalar curvature functional on metrics
  on $X$ vanishes when paired with all compactly
  supported, volume preserving, variations for which
  $S^{1}$ acts by isometries. By the so-called
  symmetric criticality principle, c.f.
  [Bes, Thm.4.23], it follows that $\nabla{\cal S}  = 0 $
  on $(X, g_{\X}),$ i.e. $g_{\X}$ is an Einstein
  metric. The scalar curvature $s_{\X}$ is determined
  by (3.5).
 
  Alternately, the equations (1.14) can be used to
  show that the equation (3.2) is equivalent to the
  condition that $g_{\X}$ is Einstein.\\
\end{pf}

  Of course a similar relation is valid in all
  dimensions. Conversely, the proof shows that if $(X, g_{\X})$ 
  is an Einstein 4-manifold of the form
  (3.4), then (3.2) holds. We note that the Einstein
  metric on $X$ may (or may not) have singularities
  if $\Omega $ includes points where $u = 0.$
 
  The equations (1.14) for the Ricci curvature of the
  4-manifold $(X, g_{\X}),$ together with the basic
  identity (2.10), show that the Ricci curvature
  $r_{\X},$ or the Einstein tensor $r_{\X}
  -\frac{1}{2}s_{\X}g_{\X},$ of $X$ may be expressed
  solely in terms of $\xi , s$ and $u.$ In the language
  of general relativity, $(X, g_{\X})$ may be viewed
  as a (Riemannian) space-time with a matter or field
  term involving only $\xi,\ s$ and $u.$ In regions
  where $u$ is bounded away from 0 and $\infty ,$
  this term is apriori bounded in $L^{2}.$

\subsection{} 
  Although (3.1) has no non-trivial global solutions
  on compact manifolds, it is closely related to
  local degenerations of Yamabe metrics. Thus,
  suppose $\{x_{i}\}$ is a sequence in $(M, g_{i})$
  such that 
 	\begin{equation} \label{e3.7}
 		\rho_{i}(x_{i}) \rightarrow  0.  
 	\end{equation}
 Consider the rescaled metrics 
 	\begin{equation} \label{e3.8} 
 		 g_{i}'  = \rho_{i}(x_{i})^{-2}\cdot g_{i}.  
 	\end{equation} 
 Throughout \S 3, we make the
 following (weak) {\it non-collapse assumption}  on
 sequences $\{g_{i}\}$ of Yamabe metrics: for some
 arbitrary, but fixed $\nu_{o} >  0,$
 	\begin{equation} \label{e3.9} 
 		\nu_{i}(x) \geq \nu_{o}\cdot \rho_{i}(x).  
 	\end{equation}
  The inequality (3.9) is scale-invariant, and it
  implies that in the scale where $\rho (x) \sim  1,$
  as in (3.8), that one has a uniform lower bound on
  the volume of geodesic $r$-balls, for $r \leq  1,$ in
  terms of Euclidean volumes; see \S 4.1 for some
  further remarks on the collapse case.
 
  In particular, for the pointed sequence $(M, g_{i}'
  , x_{i})$ above, it follows from Theorem 1.5 that a
  subsequence of $\{g_{i}'\}$ converges, modulo
  diffeomorphisms, in the weak $L^{2,2}$ topology to
  a limit metric $g' ,$ defined on a unit ball
  $B_{x}' (1),$ centered at $x = \lim x_{i}.$
 
  The following Proposition proves Theorem A(I).
 \begin{proposition} \label{p 3.1.}
    Let $(M, g_{i}, x_{i})$ be a pointed sequence of
    Yamabe metrics on a closed 3-manifold $M$ satisfying (3.7) and (3.9)
    with $s_{g_{i}} \geq  - s_{o},$ for some 
     $s_{o}  <   \infty .$ 
     Then the blow-up limit $(B_{x}' (1), g' , x) $
     is an $L^{2,2}$ weak solution of
    the static vacuum Einstein equations (3.1).
 \end{proposition} 
\begin{pf}
   Each metric $g_{i}$ has an associated splitting
   (2.10), i.e.  
        \begin{equation} \label{e3.10}
            L^{*}u + \xi  = -\frac{s}{3}\cdot  g, 
        \end{equation}
 where we have dropped the subscript $i.$ Such a
 splitting holds also for the metrics $g_{i}' $ in
 (3.8), with the same function $u = u_{i},$ since $u$
 is scale invariant. The scaling properties of
 curvature imply 
 \begin{equation} \label{e3.11}
	 s_{i}'  = \rho_{i}(x_{i})^{2}\cdot   s_{i}
	 \rightarrow  0, 
 \end{equation} since $s_{i}$ is
 uniformly bounded, (c.f. (0.4)), and $\rho_{i}(x_{i}) \rightarrow  0,$ 
 by (3.7). Further, since $\xi $ scales as
 curvature, 
	 \begin{equation} \label{e3.12}
	 \int_{M}|\xi_{i}'|^{2}dV_{i}'  =
	 \rho_{i}(x_{i})\cdot \int_{M}|\xi_{i}|^{2}dV_{i} \rightarrow  0, 
	 \end{equation} since
 $\int_{M}|\xi_{i}|^{2}dV_{i}$ is uniformly bounded,
 see (2.11). We emphasize that (3.12) uses in an
 essential way that $M$ is 3-dimensional.

  Thus, it follows from (3.10)-(3.12) that
 \begin{equation} \label{e3.13}
  (L' )^{*}u_{i} \rightarrow  0, \ \ \ \ \
  {\rm  strongly} \  {\rm in}\  L^{2}(M, g_{i}' ), 
  \end{equation}
  and hence the limit $L^{2,2}$ metric $g' $
  satisfies \begin{equation} \label{e3.14}
  (L' )^{*}u = 0.  \end{equation} in $B_{x}' (1).$
 Here, to be precise, we need to examine the limiting
 behavior of $\{u_{i}\}.$ Suppose first that
 $\{u_{i}\}$ is bounded in $L^{2}(B_{i}), B_{i} =
 B_{x_{i}}' (1).$ Since then $\Delta u_{i}
 \rightarrow  0$ in $L^{2}(M, g_{i}' ),$ it follows
 from standard elliptic estimates, c.f. [GT,
 Thm.8.8], that $\{u_{i}\}$ is uniformly bounded in
 $L^{2,2}$ in $B_{i}(1-\delta ) = B_{x_{i}}'
 (1-\delta ),$ for any given $\delta  >  0.$ Hence
 $\{u_{i}\}$ sub-converges to a limit function $u
 \in  (L^{2,2})_{loc}$ on $B_{x}' (1)$ and (3.14)
 holds weakly, i.e. when paired with smooth 2-tensors
 of compact support in $B_{x}' (1).$
 
  If instead $||u_{i}||_{L^{2}(B_{i})} \rightarrow
  \infty $ as $i \rightarrow  \infty ,$ we just
  renormalize $u_{i}$ by setting $\Roof{u}{\bar}_{i}
  = u_{i}/||u_{i}||_{L^{2}(B_{i})}.$ When
  renormalizing (3.10) by the same factor, all the
  terms become even smaller and the argument proceeds
  as before. This renormalization process will recur
  several times throughout \S 3.\\
\end{pf}

  From the discussion in \S 1.3, weak $L^{2,2}$
  solutions of the static vacuum equations are smooth,
  (at least in regions where the potential function $u$ does
  not vanish).

\medskip
  The following examples of static vacuum solutions,
  although trivial, are important in understanding
  the structure of the arguments to follow.

\medskip
\noindent
 {\bf Examples of Static Vacuum Solutions:}
 
  Super-trivial solutions $(N, g, u)$: $u \equiv  0,$
  $(N, g)$ arbitrary.
 
  Trivial solutions: $(N, g, u) = ({\Bbb R}^{3}, g_{o}, u_{o}),$ 
  where $g_{o}$ is a flat metric on ${\Bbb R}^{3}$ and 
  $u_{o}$ is a constant or affine function. 
  Similarly, one may have such solutions on
  flat quotients of ${\Bbb R}^{3}.$
 
\medskip
  Note that super-trivial solutions give no
  information whatsoever about the Riemannian
  manifold $(N, g).$ Thus, in order for Proposition 3.1
  to be of any use, one must study the sequence $(M, g_{i})$ 
  of Yamabe metrics away from the locus
  where $u_{i}$ approaches 0, see also Remark 3.15(i).

\medskip
   We have the following characterization of the
  trivial or flat solution, generalizing a classical
  result of Lichnerowicz [Li], (which assumes that $u$ is 
  asymptotically constant).
 
 \begin{theorem} \label{t 3.2. (I). } {\bf (I)}
   Let $(N, g, u)$ be a complete solution to the static
   vacuum equations (3.1), i.e. $(N, g)$ is a complete
   Riemannian manifold. If $u >  0$ on $N,$ then $N$
   is flat, and $u$ is constant.
 
 {\bf (II)}
   Let $(N, g, u)$ be a solution of (3.1) and let 
   $U \subset  N$ be any domain with smooth boundary on
   which $u >  0.$ If $t(x) = \dist_{N}(x, \partial
   U),$ for $x \in $ U, then there is an absolute
   constant $K  <   \infty $ such that
     \addtocounter{equation}{1}
 	\begin{equation} \label{e3.16} 
 		  |z|(x) 
 		   \leq 
 		   \frac{K}{t^{2}(x)},   \ \ \ \ \
 		   and \ \ \ \ \
 		   (u^{-1}|\nabla u|)(x) 
 		   \leq
 		\frac{K}{t(x)}.  
 	\end{equation}
  The constant $K$ does not depend on the domain $U,$
  (provided $u >  0$ on $U),$ or on the static
  vacuum solution $(N, g).$
 \end{theorem} 
 \begin{pf}
  The proof is deferred to the Appendix, since the
  methods do not bear directly on the main discussion
  to follow. (The proof of (3.16) is similar though
  to the proof of Theorem 3.3 below).
 
  The relation between the non-existence (of
  non-trivial solutions) in (I) is closely related to
  the existence of the pointwise curvature estimate
  (3.16) in (II). This situation occurs frequently in
  geometric P.D.E.'s and statements of the form (I)
  and (II) are often equivalent. Thus, given (I), one
  obtains (II) by a basically standard scaling
  argument. Conversely, (II) immediately implies (I)
  since the function $t \equiv  \infty $ in this
  case. \\ 
 \end{pf}
 
  As noted before, the canonical solution of the
  static vacuum equations (3.1) is the Schwarzschild
  metric (0.17). This metric is characterized by the
  conditions in Theorem 1.1.

\subsection{} 
  In this subsection, we show that the curvature is
  controlled in $L^{2}$ in regions of $(M, g)$ where
  the level sets of $u = 1+f$ do not come too close
  together. Let 
   $$L^{c} = \{x\in M: u(x) = c\} 
    \ \ \ 
   {\rm and}
    \ \ \ U^{c} = \{x\in M: u(x) >  c\}, $$ 
 denote the $c$-level and super-level sets of $u.$

 \begin{theorem} \label{t 3.3. }
   Let $g$ be a Yamabe metric on a closed 3-manifold $M$, of volume 1,
   satisfying $s_{g} \geq  -  s_{o} >  -\infty $
   and (3.9), i.e. $\nu (x) \geq  \nu_{o}\cdot \rho
   (x).$ Given any constant $c >  0,$ there is a
   constant $\rho_{o} = \rho_{o}(c, s_{o}, \nu_{o}) >  0,$ 
   such that, for any $x\in U^{c},$
	 \begin{equation} \label{e3.17} \rho (x) \geq
	 \rho_{o}\cdot  \min\{1, \dist(x, L^{c})\}.
	 \end{equation}
 \end{theorem} 
\begin{pf}
  Note that the estimate (3.17) is exactly the
  $L^{2}$ analogue of the estimate (3.16), but
  pertains to quite general Yamabe metrics while
  (3.16) holds only for the much more rigid class of
  static vacuum solutions. The proof of (3.17)
  reduces to that of (3.16) or Theorem 3.1(I) by
  taking blow-up limits.
 
  Thus, we assume (3.17) is false, and will derive a
  contradiction. 
  Given $c > 0, $  
  $s_{o}  <  \infty ,$ 
  if (3.17) does not hold, then there is a
  sequence of Yamabe metrics 
   $\gamma_{k}$ on $M_{k}$
  such that $\vol_{\gamma_{k}}M_{k} = 1,$
  $s_{\gamma_{k}} \geq  - s_{o},$ and
 	\begin{equation} \label{e3.18}
 		\frac{\rho_{k}(x_{k})}{\min\{v^{1/3},
 		\dist_{\gamma_{k}}(x_{k}, L^{c})\}} \rightarrow  0,
 		\ \ \ \ \ {\rm as}\  k \rightarrow  \infty , 
 	\end{equation} 
 for some sequence of points $x_{k}\in U^{c}=U^{c}(k);$ 
 here, we replace the constant 1 in (3.17) by 
 $v^{1/3} = \vol_{\gamma_{k}}M_{k}^{1/3} = 1,$ so that the
 equation (3.18) is scale invariant. Choose points
 $y_{k}$ realizing the minimum of the ratio in
 (3.18); it follows that $\rho_{k}(y_{k}) \rightarrow
  0,$ as $k \rightarrow  \infty .$ We rescale the
 metrics to make $\rho_{k}(y_{k})$ of size 1. 
 Thus, define metrics $\gamma'_{k} = \rho_{k}(y_{k})^{-2}\gamma_{k}$ 
 and consider the sequence of pointed Riemannian manifolds 
 $(U^{c}, \gamma'_{k}, y_{k}).$ By construction
 	\begin{equation} \label{e3.19}
 	 \rho'_{k}(y_{k}) = 1, 
 	\end{equation} 
 and
 	\begin{equation} \label{e3.20}
 	\vol_{\gamma'_{k}}M_{k} \rightarrow  \infty ,
        \ \ \ \ \  
 	\dist_{\gamma'_{k}}(y_{k}, L^{c}) \rightarrow  \infty , 
       \ \ \ \ \ 
       {\rm as}\  k \rightarrow  \infty .  
 	\end{equation} 
 Also, for any sequence $z_{k}\in U^{c},$ the estimate
 	\begin{equation} \label{e3.21}
 	 \rho'_{k}(z_{k}) \geq
 	 \rho'_{k}(y_{k})\frac{\dist_{\gamma_{k}'}(z_{k},
 	 L^{c})}{\dist_{\gamma_{k}'}(y_{k}, L^{c})} ,
 	\end{equation}
  follows from the fact that the ratio (3.18) is
  scale invariant and is minimized at $y_{k}.$ Thus,
  the sequence $(U^{c}, \gamma'_{k}, y_{k})$ has
  $L^{2}$ curvature radius uniformly bounded below,
  at points within a bounded but arbitrary distance
  from $y_{k}.$
 
  Next, recall that the volume radius (1.15) scales
  as a distance. By (3.9), (3.21) and scale
  invariance, we have a uniform lower bound $\nu_{k}'
  (z_{k}) \geq  \nu_{o}\cdot \rho' (z_{k}) \geq
  \nu_{o}\cdot  C(\dist(z_{k},y_{k})).$ Thus the
  sequence $(U^{c}, \gamma'_{k}, y_{k})$ cannot
  collapse anywhere. From Theorem 1.5, it follows
  that this pointed sequence has a subsequence which
  converges uniformly on compact subsets, in the
  {\it  weak}  $L^{2,2}$ topology, to a limit
  $L^{2,2}$ Riemannian manifold $(N, \gamma' ,$ y).
  The estimate (3.20) further implies that the metric
  $\gamma' $ is complete; the distance to the
  boundary of $U^{c},$ namely $L^{c},$ goes to
  infinity as $k \rightarrow  \infty .$
 
  From the arguments in Proposition 3.1, c.f.
  (3.10)-(3.13), we see that
  \begin{equation} \label{e3.22}
        (L' )^{*}u_{k} \rightarrow  0, \ \ \ \ \ {\rm strongly \ in }
         L^{2}(\gamma_{k}' ).  
  \end{equation}
 
 \noindent
  Now by definition, $u_{k} \geq  c$ in $U^{c}.$ If
  $u_{k}(y_{k}) \rightarrow  \infty ,$ consider the
  function 
        \begin{equation} \label{e3.23}
		 \Roof{u}{\bar}_{k}(z_{k}) =
		 \frac{u_{k}(z_{k})}{u_{k}(y_{k})}, 
     \end{equation} so
 that $\Roof{u}{\bar}_{k} \geq  0,$
 $\Roof{u}{\bar}_{k}(y_{k}) = 1.$ Otherwise, let
 $\Roof{u}{\bar}_{k} = u_{k}.$ The trace equation
 associated to (3.22) reads 
 \begin{equation}
 	\label{e3.24}
 	 \Delta' \Roof{u}{\bar}_{k} +
 	 s'_{k}\Roof{u}{\bar}_{k} =
 	 (u_{k}(y_{k}))^{-1}(\tr\xi'_{k} + s'_{k}).
 \end{equation}
  Note that the right-hand side of (3.24) goes to 0
  in $L^{2}.$ It follows from the lower bound
  $\bar{u}_{k} \geq  0$ and the
  DeGiorgi-Nash-Moser estimates, c.f. [GT, Thms
  8.17, 8.18, 8.22], that the oscillation of
  $\Roof{u}{\bar}_{k}$ is uniformly bounded in
  $B_{y_{k}}' (\tfrac{1}{2}).$ Applying the same
  considerations to neighboring balls of radius
  $\tfrac{1}{2}\rho' ,$ it follows that the
  oscillation of $\Roof{u}{\bar}_{k}$ is uniformly
  bounded in $B_{z_{k}}' (\frac{1}{2}\rho' (z_{k}))$
  for all $z_{k}$ within uniformly bounded distance
  to $y_{k}.$ This, together with standard $L^{2}$
  estimates for elliptic equations, c.f. [GT,
  Thm.8.8] imply that a subsequence of
  $\{\Roof{u}{\bar}_{k}\}$ converges, in the weak
  $L^{2,2}$ topology, with respect to the metrics
  $\gamma'_{i},$ to a limit $L^{2,2}$ function
  $\Roof{u}{\bar}$ on $N.$ Further, by (3.22) and the
  fact that (3.22) is preserved under the
  renormalization (3.23), the limit pair $(\gamma' , \bar{u})$ 
 satisfies
	  \begin{equation} \label{e3.25}
	   {L'}^{*}\bar{u} = 0, \ \ \ \ \   \Delta' \bar{u}  = 0, 
	  \end{equation}
  so that $(N, \gamma', \bar{u} )$ is a weak
  solution of the static vacuum Einstein equations.
 
  Note that the (weakly) harmonic function
  $\Roof{u}{\bar}$ satisfies $\Roof{u}{\bar} \geq  0$
  everywhere and $\Roof{u}{\bar}(y) = 1.$ By the
  (weak) maximum principle, c.f. [GT, Thm.8.1], in
  fact $\bar{u} >  0$ everywhere. As noted
  in \S 1.3, when the potential $\bar{u} > 0,$ 
  elliptic regularity implies that $L^{2,2}$
  weak solutions of (3.25) are $C^{\infty}.$ From
  Theorem 3.2(I), it follows that $\gamma' $ is flat,
  and $\Roof{u}{\bar}$ is constant.
 
  On the other hand, we have the estimate (3.19) for
  the sequence $(\gamma'_{k}, y_{k}).$ The fact that
  the convergence of $\gamma'_{k}$ to $\gamma' $ is
  only in the weak $L^{2,2}$ topology does not imply
  that (3.19) passes continuously to the limit.
  However, we will show in Theorem 3.4 below that in
  fact the convergence $\gamma'_{k} \rightarrow
  \gamma' $ is in the strong $L^{2,2}$ topology. By
  (1.28), the radius $\rho (x)$ is continuous in the
  strong $L^{2,2}$ topology and one thus obtains the
  estimate $$\rho' (y) = 1, $$ on the limit. This
 however contradicts the fact that $\gamma' $ is
 flat; a complete flat manifold $(N, g)$ clearly has $\rho (x)
 =\infty ,$ at any $x.$
 
  Hence the proof of Theorem 3.3 is completed with
  the proof of the following result, which will also
  be used frequently in the work to follow. \\
\end{pf}
 
 \begin{theorem}[\bf Strong Convergence] \label{t 3.4. (Strong Convergence).}
   Let $\{g_{i}\}$ be a sequence of Yamabe metrics,
   (not necessarily of volume 1), with associated
   splittings 
 	  \begin{equation} \label{e3.26}
 	 L^{*}u_{i} + \xi_{i} = -  \frac{s_{i}}{3}\cdot g_{i} .  
 	 \end{equation}
 Suppose that for all i, there is a constant $d > 0$  such that
 	\begin{equation} \label{e3.27} 
 		r_{h}(x_{i}, g_{i})
 		\geq  1, 
                 \ \ {\rm and }\ \
                r_{h}(y_{i}, g_{i}) \geq  d,  \ \  \forall
 		y_{i}\in\partial B_{x_{i}}(1),  
 	\end{equation}
  where $r_h$ denotes the $L^{2,2}$ harmonic radius. Suppose also that $\xi_{i} \rightarrow  0$ strongly
  in $L^{2}(B_{x_{i}}(1)), s_{i} \geq  - s_{o},$ for
  some $s_{o},$ and that there is a constant $u_{o}
  >  0$ such that 
	  \begin{equation} \label{e3.28}
		 u_{i} \geq  u_{o}\ \ 
		 {\rm  on\ } B_{x_{i}}((1+\tfrac{1}{2}d)) .
	 \end{equation}

  Then on $B_{x_{i}}(1),$ a subsequence of
  $\{g_{i}\}$ converges strongly in the $L^{2,2}$
  topology to a limit metric $g_{o}$ on $B_{x}(1),$
  where $x = \lim x_{i}.$
 \end{theorem} 
\begin{pf}
  Let $B = B_{i} = B_{x_{i}}(1)$ and set $B'  =
  B_{x_{i}}(1+\frac{d}{2}).$ By (3.27) and Theorem
  1.5, we may assume that a subsequence of
  $\{g_{i}\}$ converges weakly in the $L^{2,2}$
  topology to a limit metric $g_{o}$ on $B' .$ Thus,
  there are coordinates $\{y_{k}\}$ on suitable balls
  $D\subset B_{x_{i}}(1+d),$ covering $B' ,$ such
  that the functions 
           \begin{equation} \label{e3.29}
             (g_{i})_{kl} \rightarrow  (g_{o})_{kl},
           \end{equation} 
  weakly in $L^{2,2}(D).$ As in the
 proof of Theorem 3.3, (by dividing by $u(x_{i})$ if
 necessary, see the discussion following (3.23)), we
 may assume that the functions $u_{i}$ are uniformly
 bounded, in $L^{\infty}(B' ),$ and that (3.28)
 holds, possibly with a different constant $u_{o}.$
 In particular, a subsequence converges strongly in
 $L^{2}$ to a limit function $u.$ Thus, from the trace
 equation
 	\begin{equation} \label{e3.30}
 		\Delta_{i}u_{i} = \tr\xi_{i} -  s_{i}u_{i} + s_{i}, 
 	\end{equation}
  and the arguments above, together with the
  assumption on $\xi , \Delta_{i}u_{i}$
  converges strongly in $L^{2}$ to a limit $L^{2}$
  function $\Delta_{o}u$ on $B' $ with $u\in T^{2,2}(B' ).$
 
  We now use the $L^{2}$ estimates for elliptic
  equations of the form (3.30), see [GT, Thm 8.8].
  This gives 
     \begin{equation} \label{e3.31}
 		||D^{2}u_{i}||_{L^{2}(B)} 
 		\leq
 		C(||\Delta_{i}u_{i}||_{L^{2}(B' )} +
 		||u_{i}||_{L^{2}(B' )}), 
 	\end{equation}
  where $C$ depends on the $L^{2,2}$ norm of the
  metrics $g_{i}$ and the constant $d$ in (3.27).
  Here the norms and derivatives $D^{2}$ are taken
  with respect to the fixed coordinates $\{y_{k}\}$
  and limit metric $g_{o}.$ From this, it then
  follows first that 
 	 \begin{equation} \label{e3.32}
		||D^{2}u_{i}||_{L^{2}(B)} \leq  C, 
 	\end{equation}
  so that a subsequence of $\{u_{i}\}$ converges
  weakly in the $L^{2,2}$ topology to $u.$ By the
  Sobolev embedding theorem, $\{u_{i}\}$
  (sub)-converges strongly in the $L^{1,2}$ topology.
  Repeating the estimate (3.31) on $u -  u_{i}$ gives
 	\begin{equation} \label{e3.33} 
		||D^{2}(u- u_{i})||_{L^{2}(B)} 
		\leq  
		C(||\Delta_{i}(u-
		u_{i})||_{L^{2}(B' )} + ||(u- u_{i})||_{L^{2}(B')}).  
 	\end{equation}
  Clearly $||(u- u_{i})||_{L^{2}(B' )}\rightarrow  0.$ 
  Write $\Delta_{i}(u- u_{i}) =
  (\Delta_{i}u-\Delta_{o}u)
    +(\Delta_{o}u-\Delta_{i}u_{i}).$
  From preceding arguments,
  $\Delta_{i}u_{i}\rightarrow\Delta_{o}u$
  strongly in $L^{2}(B' ).$ Consider then the
  sequence $\Delta_{i}u-\Delta_{o}u.$
  We have in the $y_{k}$-coordinates
	 \begin{equation} \label{e3.34} 
          \Delta_{i}u =
	  g_{i}^{kl}\partial_{kl}u + \partial
	  g_{i}^{kl}\partial_{k}u.  
	 \end{equation}
  Since $\{g_{i}\}$ is uniformly bounded in
  $L^{2,2},$ a subsequence converges strongly in
  $L^{1,2}\cap L^{\infty}$ to $g_{o}.$ This shows
  that $\Delta_{i}u$ converges strongly in
  $L^{2}$ to $\Delta_{o}u.$ Thus, (3.33)
  implies that
	 \begin{equation} \label{e3.35} 
		||D^{2}(u- u_{i})||_{L^{2}(B)}\rightarrow  0.  
	\end{equation}
  If $D_{i}^{2}$ denotes the Hessian with respect to
  the metrics $g_{i},$ the same reasoning as above on
  the Laplacian in (3.34) then also gives
	 \begin{equation} \label{e3.36}
	 D_{i}^{2}u_{i}\rightarrow  D^{2}u,   \ \ \ \ \
	 {\rm  strongly\ in }\   L^{2}(B).  
	 \end{equation} 
 Now return to the splitting equation (3.26), which we write as
	 \begin{equation} \label{e3.37} 
		   u_{i}r_{i} = D_{i}^{2}u_{i} -  \Delta_{i}u_{i} 
		   + \xi_{i} + \frac{s_{i}}{3}\cdot  g_{i}.  
	 \end{equation}
  By the preceding arguments, the right hand side of
  (3.37) converges strongly in $L^{2}$ to its limit
  on $B$ and thus 
	 \begin{equation} \label{e3.38}
	 \int_{B}|u\cdot r_{o} -  u_{i}\cdot r_{i}|^{2}dV_{o} \rightarrow  0, 
	  \end{equation}
  where the norm is taken in the $g_{o}$ metric and
  $r_{o}$ is the weak $L^{2}$ limit of the Ricci
  curvature $r_{i}$ of $g_{i}.$ Since $u_{i}
  \rightarrow  u$ in $L^{2,2}\cap C^{o}$ on $B$ and
  since $u_{i} \geq  u_{o}$ by (3.28), it follows
  that
	 \begin{equation} \label{e3.39} 
	   \int_{B}|r_{o} - r_{i}|^{2}dV_{o} \rightarrow  0, 
	 \end{equation}
  so that the Ricci curvature of $\{g_{i}\}$
  converges strongly in $L^{2}$ to the Ricci
  curvature of $g_{o}$ on $B.$ The same reasoning shows
  that $r_{i}$ converges strongly to $r$ on a
  thickening, say $B(1+\frac{d}{4})$ of $B.$ This
  implies the strong $L^{2,2}$ convergence of $g_{i}$
  to $g_{o},$ via the equation for the Ricci
  curvature in harmonic coordinates on $B,$ see [An3,
  p.434] for the details here.  \\
 \end{pf}
 \begin{remark}
 \label{r 3.5(i):} {\bf (i):}
   If one drops the assumption (3.27) on the behavior
   of $r_{h}$ at the boundary, so that the sequence
   $(B_{x_{i}}(1), g_{i})$ satisfies
   $r_{h}(x_{i},g_{i}) = 1$ with $u_{i} \geq  u_{o}$
   on $B_{x_{i}}(1)$ and $\xi_{i} \rightarrow  0,$
   $s_{i} \geq  - s_{o},$ then Theorem 3.4 implies
   that a subsequence converges strongly in the
   $L^{2,2}$ topology on $B_{x_{i}}(1-\delta ),$ for
   any fixed $\delta >  0,$ to a limit metric
   $g_{o}$ on $B_{x}(1-\delta ).$
 
 {\bf (ii):}
   Returning to Theorem 3.3, note that the estimate
   (3.17), together with Lemma 1.4, implies that if
   $u_{i} \geq  u_{o} >  0$ on $M,$  then
 	\begin{equation} \label{e3.40}
 		\int_{M}|z_{i}|^{2}dV_{i} \leq  C/u_{o}.
 	\end{equation}
  This has already been proved by elementary means in (2.18).
 
 {\bf (iii):}
  We note that exactly the same proof as Theorem 3.3
  shows the estimate (3.17) holds in the region
  $U_{-c} = \{x\in M: u < - c  <   0\},$ 
  i.e.
  $\rho (x) \geq  a\cdot  \dist(x, L^{-c}),$ for $x\in
  U_{-c}.$ 
  \end{remark}

\subsection{} 
   The discussion in \S 3.3 implies that $\rho
   (x_{i})$ can become arbitrarily small, under the
   sequence $\{g_{i}\},$ only if the 0-level set
   $L^{o},$ or an $\varepsilon$-level set
   $L^{\varepsilon}$ for $\varepsilon $ very small, comes
   arbitrarily close to $x_{i}.$
 
  In this subsection, we will begin the study of the
  geometry of the degeneration of $(M, g_{i})$ in a
  neighborhood of the 0-level $L^{o}.$ On the other
  hand, from the remarks above, the 0-level itself
  must be avoided; see also Remark 3.15(i).
 
  The treatment of the three possible cases $\sigma
  (M)  <   0,$ $\sigma (M) = 0 $ and $\sigma (M) > 0 $
  is somewhat different. Some indication for the
  need to treat these cases separately is already
  apparent from \S 2. Thus, throughout \S 3.4, we
  assume 
 	$$ \sigma (M)  <   0, 
 	$$
  or more precisely that $\{g_{i}\}$ is a sequence of
  unit volume Yamabe metrics on $M$ with scalar
  curvature bounded away from $0$  and $-\infty ,$ i.e.
 	\begin{equation} \label{e3.41} 
 	     -\infty   <   - s_{o} \leq  s_{g_{i}} \leq  - s_{1}  <   0.
 	\end{equation} 
 The cases $\sigma (M) = 0$ and
 $\sigma (M) >  0$ are treated in \S 3.5 and \S
 3.6 respectively. As in Theorem A(II), c.f. (0.19),
 we also assume 
	 \begin{equation} \label{e3.42}
	 \int_{M}|z^{T}|^{2}dV \leq  K, 
	 \end{equation} and
 \begin{equation} \label{e3.43} 
    \nu (x) \geq
 \nu_{o}\cdot \rho (x); 
 \end{equation}
  note that (3.43) is more general than (0.18).
 
  To set the stage for the considerations to follow,
  let $T_{i} = \max u_{i}$ and choose $x_{i}$ such
  that
 \begin{equation} \label{e3.44} 
	   |u_{i}(x_{i})/T_{i} -  1| \rightarrow  0, 
	 \ \ \ \ \ 
	 {\rm as}\ i \rightarrow  \infty .
 \end{equation}
  Thus, we are considering points as far away as
  possible, in terms of the function $u,$ from the
  0-levels of $u.$ By (2.35), as discussed in Remark
  2.12, (3.41) and (3.42) imply one has a uniform
  lower bound
 \begin{equation} \label{e3.45} 
          T_{i} \geq \lambda_{o} \geq  0.  
\end{equation}

\noindent
 {\bf Note:}
  The bound (3.45) is the only place in \S 3.4, (and
  also in \S 3.6), where the assumption (3.42) is
  used. In \S 3.5, corresponding to the case $\sigma
  (M) = 0,$ (3.42) is used in a stronger way.

\bigskip

  Consider $\rho (x_{i}),$ the $L^{2}$ curvature
  radius at $x_{i},$ of course with respect to the
  metric $g_{i}.$ If there exists a uniform lower
  bound 
         \begin{equation} \label{e3.46} 
             \rho (x_{i}) \geq  \rho_{o}, \ \   \forall i, 
          \end{equation}
  then Theorem 1.5 provides a complete description of
  the possible behavior of $\{g_{i}\}$ in
  $B_{x_{i}}(\rho_{o}).$ Suppose instead, for a
  possibly different choice $\{x_{i}^{1}\},$ that
  both 
\begin{equation} \label{3.47}
  \begin{array}{c}
        |u_{i}(x_{i}^{1})/T_{i} -  1| \rightarrow  0, \ {\rm and} \\[4mm]
                  \rho (x_{i}^{1}) \rightarrow  0, 
             \ \  {\rm as}\  i
 \rightarrow  \infty . 
  \end{array}
\end{equation}
  Thus, the metrics $g_{i}$ are degenerating in a
  (progressively smaller) neighborhood of $x_{i}^{1}$
  and this degeneration cannot be described by
  Theorem 1.5.
 
  To understand the degeneration of $(M, g_{i})$ at
  or near $\{x_{i}^{1}\},$ rescale the metrics
  $g_{i}$ by $\rho (x_{i}^{1})^{-2}.$ Thus, set
	 \begin{equation} \label{e3.48} 
		g_{i}^{1} = \rho (x_{i}^{1})^{-2}\cdot  g_{i}, 
	  \end{equation}
  and let $\rho^{1}$ denote the $L^{2}$ curvature
  radius with respect to the rescaled metrics
  $g_{i}^{1}.$ One then has
	 \begin{equation} \label{e3.49} 
		\rho^{1}(x_{i}^{1}) = 1, \ \ \ \ \ {\rm for \ all}\  i, 
	  \end{equation}
  so that Theorem 1.5 and (3.43) imply the
  sub-convergence of $\{g_{i}^{1}\}$ in
  $B_{x_{i}^{1}}^{1}(1).$ Thus we may assume that
  $\{g_{i}^{1}\}$ converges, weakly in the $L^{2,2}$
  topology, to a limit $L^{2,2}$ metric $g^{1},$
  defined on $B^{1} = B_{x^{1}}^{1}(1),$ with base
  point $x^{1}  = \lim x_{i}^{1}.$

 \begin{proposition} \label{p 3.6.}
   The limit $(B^{1}, g^{1}, x^{1})$ is a flat
   solution of the static vacuum equations (3.1),
   with potential function 
	 \begin{equation} \label{e3.50} 
	 \Roof{u}{\bar} = \lim (u_{i}/T_{i})
	 \equiv  1 \ \  {\rm on}\  B^{1}.  
	 \end{equation}
 \end{proposition} 
\begin{pf}
  The splitting (2.10) of $g_{i},$ when rescaled,
  gives the splitting for $g_{i}^{1}.$ Thus
 \begin{equation}
\begin{array}{c}
	 (L^{*})^{1}(u_{i}) + \xi_{i}^{1} = -
	 \frac{s_{i}^{1}}{3}\cdot  g_{i}^{1}, \\[4mm]
	 \Delta^{1}u_{i} + s_{i}^{1}u_{i} =
	 \tr\xi_{i}^{1} + s_{i}^{1}. 
 \end{array}
 \end{equation}
  As in the proof of Proposition 3.1, we have
  $\xi_{i}^{1} \rightarrow  0$ in $L^{2}(M, g^1_{i} ),\ 
  s_{i}^{1} \rightarrow  0$ in $L^{\infty}(M, g^1_{i}).$ 
  While it may be possible apriori that
  $u_{i}(x_{i}^{1}) \rightarrow  +\infty ,$
  (corresponding to the possibility that $T_{i}
  \rightarrow  \infty ),$ as in the proof of Theorem
  3.3, (or Proposition 3.1) we renormalize equation
  (3.51) by dividing $T_{i}$ to obtain the
  (renormalized) limit function $\Roof{u}{\bar}$ on
  $B^{1}.$ As noted above, $T_{i}$ is bounded away
  from 0 by (3.45), (a consequence of the standing
  assumption (3.42)), so that the renormalizations of
  the terms $\xi_{i}^{1},\ s_{i}^{1}$ go to~$0$ at least
  as fast as before.
 
  Thus, on $B^{1},$ the limit metric $g^{1}$
  satisfies the static vacuum Einstein equations 
 \begin{equation}\label{3.52}
\begin{array}{c}
  (L^{*})^{1}(\bar{u})  = 0, \\[4mm]
 \Delta^{1}\Roof{u}{\bar} = 0, 
\end{array}
 \end{equation}
  where $\Roof{u}{\bar} = \lim u_{i}/T_{i}.$  As in
  the proof of Theorem 3.3, the convergence of
  $\Roof{u}{\bar}_{i} = u_{i}/T_{i}$ to the limit is
  in the weak $L^{2,2}$ topology, and since
  $\Roof{u}{\bar}_{i}(x_{i}^{1}) = 1,$ the limit
  $\Roof{u}{\bar}$ is a weakly harmonic $L^{2,2}$
  function, which is not identically~$0.$\\ 
 
 As discussed in \S 1.3, the equations (3.52) imply
  in particular that the metric $g^{1}$ is smooth
  away from the 0-locus of $\Roof{u}{\bar}.$ Now by
  (3.44), $x^{1},$ the center of the ball $B^{1},$
  satisfies
	 \begin{equation} \label{e3.53} 1 =
		 \Roof{u}{\bar}(x^{1}) \geq  \Roof{u}{\bar}(y),
	 \end{equation}
  for all $y\in B^{1}.$ Since $\Roof{u}{\bar}$ is
  harmonic, it follows that $\Roof{u}{\bar} \equiv $
  1 and thus $(B^{1}, g^{1})$ is flat.\\
 \end{pf}
 
\bigskip

  The convergence $g_{i}^{1} \rightarrow  g^{1}$ is
  in the strong $L^{2,2}$ topology, on
  $B_{x_{i}^{1}}^{1}(s),$ for any fixed $s <  1,$
  by Theorem 3.4. It thus follows that all of the
  curvature of $\{g_{i}^{1}\}$ in $L^{2}$ is
  concentrating on the boundary $\partial
  B_{x_{i}^{1}}^{1}(1).$ This implies that the limit
  $(B^{1}, g^{1})$ itself does not yet effectively
  model the degeneration of $g_{i}$ near $x_{i}^{1};$
  the base points $x_{i}^{1}$ must be altered
  slightly.
 
  The fact that some or all of the curvature in
  $L^{2}$ is concentrating on $\partial B_{x_i^1}^1$ 
  implies, (and is equivalent~to),
			 $$\rho^{1}(y_{i}) \rightarrow  0, $$
  for some sequence $y_{i}\in\partial
  B_{x_{i}^{1}}^{1}(1),$ or, what is same,
  $\rho_{i}(y_{i})  <  <   \rho_{i}(x_{i}^{1}).$
  Thus the curvature is blowing up at $y_{i}$ much
  faster than at $x_{i}^{1}.$ This is the first
  indication that the curvature of $(M, g_{i})$ blows
  up at many different scales.
 
  The preceding remarks lead naturally to the
  following definitions.
 \begin{definition}
   A point $y$ in a complete Riemannian manifold $(N, \gamma)$ is
    {\it $(\rho ,c)$ buffered}  if  $c >  0$ and
  	 \begin{equation} \label{e3.54}
 		   \frac{\rho (y)^{4}}{\vol B_{y}(\rho
 		   (y))}   \int_{B_{y}((1-c)\cdot \rho
 		   (y))}|r|^{2} \geq  c\cdot  c_{o}, 
 	 \end{equation}
  where $c_{o}$ is the constant in the definition of
  $\rho ,$ c.f. (1.26).
 
  Similarly, $y$ is
  {\it  strongly}   {\it $(\rho ,d)$~buffered} if
   $d >  0$ and 
	 \begin{equation} \label{e3.55} 
	   \rho (z) \geq  d\cdot \rho (y), 
	 \end{equation}
  for all 
    $z \in \partial B_{y} (\rho ( y ) ).$

  A sequence of points ${y_{i}}$ in the manifolds $(N, \gamma_{i})$ is $(\rho,c)$ buffered, or strongly $(\rho,c)$ buffered, if each $y_i$ is.

 \end{definition}
 
  The buffer constants $c$ and $d$ are arbitrary small parameters, as is the parameter $c_o$; their precise values, beyond being small, are not important. From \S 1, $c_o$ is a fixed small number throughout the paper, e.g. $c_o = 10^{-3}$. The buffer constant $c$, (more important than $d$), might be allowed to vary over small numbers, say $0 < c < 10^{-3}$, but in any given discussion, the value of c will be fixed.

  Note that the strong $(\rho ,d)$
  buffered condition appears in Theorem 3.4 (strong
  convergence), and is invariant under scaling. Thus,
  if $\{g_{i}\}$ is a sequence of unit volume Yamabe metrics on a closed 3-manifold $M$ with
  $\rho (y_{i}) \rightarrow  0,$ then the blow up
  metrics $g_{i}'  = \rho (y_{i})^{-2}\cdot  g_{i}$
  converge strongly in $L^{2,2}$ on $B_{y_{i}}' (1),$ to a non-flat limit,
  provided the sequence $\{y_{i}\}$ is strongly
  $(\rho ,d)$ buffered, and $u_{i}$ is bounded away
  from 0 on $B_{y_{i}}'(1).$ Thus, this condition
  prevents {\it any}  of the curvature from
  concentrating in $L^{2}$ (with normalized measure)
  on the boundary.
 
  Similarly, the $(\rho ,c)$ buffered condition
  (3.54) is invariant under scaling and
  prevents {\it all}  of the curvature from
  concentrating in $L^{2}$ on the boundary. Thus, if
  $c$ is close to 0, then a definite (small)
  percentage of the curvature in $L^{2}$ is in the
  ball of radius $(1-c).$  (The situation where $c$ is
  close to $1$ will never arise here). In particular,
  from Theorem 3.4, any limit $B_{y}' (1)$ of a
  $(\rho ,c)$ buffered sequence $(B_{y_{i}}'(1), g_{i}', y_{i})$ of Yamabe metrics as above, with $u_i$ bounded away
  from 0 in $B_{y_{i}}' (1)$, cannot be flat, since one has strong convergence everywhere in the interior.
  Conversely, again by Theorem 3.4, if $u_i$ is bounded away from $0$ on $B_{y_{i}}'(1)$ and the limit is not flat, then the sequence is $(\rho ,c)$ 
  buffered, for some $c > 0$, c.f. also Remark 3.5(i).
 
  Note that from the definition, it is obvious that if a given sequence is $(\rho, c)$ buffered, then it is $(\rho, c')$ buffered, for any $c' < c$, and similarly for strongly buffered sequences.

  The next Lemma formalizes the relation between these notions, for Yamabe metrics.

\begin{lemma} \label{l 3.8.}
  Let $g$ be a unit volume Yamabe metric on a closed 3-manifold $M$, and $y$ a base point in $M$ which is strongly $(\rho, d)$ buffered. Suppose $u(y) = 1$, the non-collapse assumption (3.43) holds at y, and the blow-up $(B_{y}'(1), g', y)$, $g' = \rho(y)^{-2} \cdot g$ is $\epsilon$-close to a static vacuum solution, in the sense that $||L^*u||_{L^2} \leq \epsilon$ in $B_{y}'(1)$, as in (3.13).

\indent
   Under these assumptions, there exists $\varepsilon_{o} > 0$ such that if $\varepsilon \leq \varepsilon_{o}$, and if either

\indent
 (i) $u_{i}$ is bounded away from $0$ in $B_{y{i}}'(1)$, or

\indent
  (ii) $|u_{i} - 1| \leq \tfrac{1}{4}$ in $B_{y{i}}'(\tfrac{1}{2})$,

\noindent
then there is a constant $c = c(d, c_o, \nu_o) > 0$ such that $y$ is $(\rho, c)$ buffered.
\end{lemma}

\begin{pf}

  The proof is by contradiction. Thus, if the conclusion is not true, there exists a pointed sequence $(M, g_i, y_i)$ satisfying the assumptions with $\varepsilon = \varepsilon_i \rightarrow 0$, such that the sequence ${y_i}$ is not $(\rho, c)$ buffered. Let $(B'(1), g', y)$ be a weak $L^{2,2}$ limit of (a subsequence of) $(B_{y_{i}}'(1), g_{i}', y_i)$, so that $(B'(1), g')$ is an $L^{2,2}$ static vacuum solution. The existence of such limits follows, as previously in \S 3.4, from the non-collapse assumption (3.43).

   By the discussion above regarding strong convergence, it suffices to prove that $u_i$ is bounded away from $0$ on $B_{y_{i}}'(1)$ and the limit is not flat. In the case of assumption (i), this first condition is obvious, while the second follows since $\rho$ is continuous under strong $L^{2,2}$ convergence. Hence one has a contradiction in this case.

  With regard to (ii), suppose first the limit $(B'(1), g')$ were flat. Then from the static vacuum equations, the limit function $u$ is an affine function on $B'(1)$. The assumption (ii) then implies that $u$ is uniformly bounded away from $0$ in $B'(1)$, in fact $u \geq \frac{1}{2}$. Further, since ${y_{i}}$ is strongly $(\rho, d)$ buffered, the functions $u_i$ are uniformly controlled in $L^{2,2}$, and hence controlled in $C^{\alpha}$, $\alpha < \tfrac{1}{2}$, in $B_{y_{i}}'(1+\tfrac{d}{2})$. It follows that ${u_{i}}$ is uniformly bounded away from $0$ in $B_{y_{i}}'(1)$. As in Case (i) above, this implies that the full limit $B'(1)$ could not have been flat.

  Since the full limit $(B'(1), g')$ is not flat, it follows from the real-analyticity of static vacuum solutions discussed in \S 1.3, and (ii) again, that the limit is not flat in any domain in $B'(1)$ on which $u$ is bounded away from $0$; in particular, this is the case on $B'(\tfrac{1}{2})$. Hence, as above, $u_i$ is bounded away from $0$ in $B_{y_{i}}'(\tfrac{1}{2})$, and the strong convergence to this limit $B'(\tfrac{1}{2})$ implies that the sequence is $(\rho, c)$ buffered, for some $c > 0$.

\end{pf}

  The buffer constant $c = c(d, c_o, \nu_o)$ in Lemma 3.8 could be explicitly estimated in terms of $d$, $c_o$ and $\nu_o$; however, we will have no need to do this. Of course, $c$ may well be much smaller than $d$. In the following, we will usually use the contrapositive of Lemma 3.8, namely that if a given sequence ${y_i}$ is not $(\rho, c)$ buffered, for a given $c > 0$, then it is not strongly $(\rho, c)$ buffered if either (i) or (ii) in Lemma 3.8 hold; here $d = d(c, c_o, \nu_o)$ may be large compared with $c$, but may made arbitrarily small by choosing $c$ sufficiently small.

\medskip
 
  Now we return to the situation preceding Proposition 3.6. Of course the points $x_{i}^{1}$ in (3.47) above
  are not $(\rho , c)$ buffered, for any choice of $c >  0.$ 
  Observe that all the positive level sets
  $L^{s},\ s >  0,$ of $\Roof{u}{\bar}_{i}$ have
  points converging to $\partial
  B_{x_{i}^{1}}^{1}(1).$ For otherwise, by Theorem
  3.3, $\{x_{i}^{1}\}$ would be a (strongly) buffered
  sequence,  a contradiction.
 
  Let $U = U(i)$ be the component of
  $(u_{i}/T_{i})^{-1}((0, 2]) \subset  M$
  containing $x_{i}^{1}.$ Let $L^{k} = L^{k}(i),\  k \geq~1,$ 
  denote the level set
  $(u_{i}/T_{i})^{-1}(2^{-k+1})$ in $U(i),$ (possibly
  having many components). By a slight perturbation
  of the values, we may assume that $L^{k}$ are
  smooth hypersurfaces. Let $U^{k} = U^{k}(i)$ be the
  set $u_{i}^{-1}[2^{-k+1},2]\subset U.$ The choice
  of the factor 2 here is for convenience, and may be
  replaced by any fixed number $> 1.$
 
  In Theorem 3.10 below, we will construct a 
  $(\rho ,c)$-buffered sequence from an initial sequence
  $\{x_{i}^{1}\}$ satisfying conditions similar to
  (3.47) above, in order to produce a non-flat blow
  up limit. This is done essentially by `descending
  down' the level sets $L^{k}$ of $u_{i}$, for $i$ fixed.
 
  To understand the underlying idea and motivation of
  the proof of Theorem 3.10, it is useful to point
  out that the descent down the levels of $u$ must
  require an infinite process. This is shown by the
  following remark and example. Remark 3.9 and Example 3.9 below are not logically necessary in this respect, and may be skipped if preferred.

 \begin{remark} \label{r 3.9.}
   Recall again that the metrics $g_{i}^{1}$ converge
   uniformly on compact subsets of
   $B_{x_{i}^{1}}^{1}(1)$ to the flat metric, and
   $u_{i}/T_{i} \rightarrow  1$ uniformly on such
   compact subsets. It is clear that there exists a
   sequence $\{x_{i}^{2}\}\in U^{2}\cap
   B_{x_{i}^{1}}^{1}(1)$ such that
 	\begin{equation} \label{e3.56} 
 		\rho^{1}(x_{i}^{2}) \rightarrow  0,  \ \
 		{\rm i.e.}\  \rho (x_{i}^{2})  <  <   \rho
 		(x_{i}^{1}), 
 	\end{equation}
  and consequently 
    $$
       \dist_{g_{i}^{1}}(x_{i}^{2},
 	\partial B_{x_{i}}^{1}(1)) \rightarrow  0. 
    $$
  There are many possible choices of $\{x_{i}^{2}\}$
  satisfying (3.56).
 
  For reasons that will be apparent below, we require
  further that $\{x_{i}^{2}\}$ satisfy the following:
  $r_{i}^{2} \equiv  \dist_{g_{i}^{1}}(x_{i}^{2},
  x_{i}^{1}) \rightarrow  1 $ as $i \rightarrow
  \infty ,$ 
  and
 	\begin{equation} \label{e3.57} 
 	|u_{i}(x) - u_{i}(x_{i}^{1})| 
 	\leq  
 	\mu_{o}\cdot u_{i}(x_{i}^{1}), \  \  \forall x\in
 	B_{x_{i}^{1}}^{1}(r_{i}^{2}).  
 	\end{equation}
  Here $\mu_{o}$ is a fixed constant, $0  < \mu_{o}  <   \frac{1}{2},$ 
  for example $\mu_{o} = \frac{1}{4}.$ 
  Again it is clear that there are many
  possible choices of such points $x_{i}^{2}.$ If
  $T_{i} = \max u_{i} \rightarrow  \infty ,$ (3.57)
  is assumed to apply to the renormalized functions
  $u_{i}/T_{i},$ as before. The factor $T_{i}$ will
  be ignored below.
 
  Now rescale the metrics $g_{i}^{1}$ further to make
  $\rho (x_{i}^{2}) = 1,$ i.e. set
 	$$g_{i}^{2} = \rho^{1}(x_{i}^{2})^{-2}\cdot
 	g_{i}^{1} = \rho (x_{i}^{2})^{-2}\cdot  g_{i}. $$
  Consider then the sequence of pointed Riemannian
  manifolds $(B_{x_{i}^{2}}^{2}(1), g_{i}^{2},
  x_{i}^{2}).$ This has $L^{2}$ curvature radius
  equal to 1 at $x_{i}^{2},$ and as above, one may
  apply Theorem 1.5. The collapse case is ruled out
  by assumption (3.43).
 
  Thus a subsequence of $\{g_{i}^{2}\}$ converges, in
  the weak $L^{2,2}$ topology on
  $B_{x_{i}^{2}}^{2}(1)$, and in the strong $L^{2,2}$
  topology on $B_{x_{i}^{2}}^{2}(s),$ any $s  < 1 $
  away from the locus $\{|u_{i}| \leq\varepsilon\}$
  for any given $\varepsilon  > 0,$  to a limit
  $(B^{2}, g^{2}, x^{2}).$ Necessarily, this limit is
  a solution of the static vacuum Einstein equations
  (3.1), with limit function $u^{2}.$ Note that since
  $u_{i}(x_{i}^{2})\geq  1- 2\mu_{o}$ for $i$ large,
  the limit function $u^{2}$ is not identically~$0.$
 
  However, we claim that the limit $(B^{2}, g^{2})$ is also
  flat and $u^{2}$ is constant. To see this, consider first the
  `half-ball' $B_{x_{i}^{2}}^{2}(1)\cap
  B_{x_{i}^{1}}^{2}(R_{i}^{2})$, where $R_{i}^{2}=
  \dist_{g_{i}^{2}}(x_{i}^{2}, x_{i}^{1})$. Note that in this scale, $R_{i}^{2} \rightarrow \infty $ as $i \rightarrow \infty $ by (3.56)-(3.57). Since the curvature of $B_{x_{i}^{1}}^{1}(1)$ is bounded, and $g_{i}^{2}$ is a 'blow-up' of $g_{i}^{1}$, the curvature of $(D_{i}^{2}, g_{i}^{2})$ goes to $0$ uniformly in $L^2$. Hence the limit half-ball $(D^2, g^2)$ is flat.

  To extend this to the full ball $B^2$, the static vacuum equations imply that $D^{2}u^{2} = 0,$ 
  so that $u^{2}$ must be an affine function on
  its domain. We claim that the limit harmonic
  function $u^{2}$ on $D^{2}$ is constant, i.e.
 	\begin{equation} \label{e3.58} 
 	u^{2} = const. \geq 1 -  \mu_{o}.  
 	\end{equation}
 
  To see this, return to the balls
  $B_{x_{i}^{2}}^{2}(1)$ in the scale $g_{i}^{2}.$ In
  this scale, the previous ball
  $(B_{x_{i}^{1}}^{1}(1), g_{i}^{1})$ is expanded
  into a ball $B_{x_{i}^{1}}^{2}(R_{i})$ of radius
  $R_{i} \rightarrow  \infty $ as $i \rightarrow
  \infty .$ Let $\gamma_{i}^{2}$ be a minimizing
  geodesic from $x_{i}^{2}$ to $x_{i}^{1}$ in
  $B_{x_{i}^{1}}^{2}(R_{i}),$ parametrized by
  arclength w.r.t. $g_{i}^{2}.$ Note that the length
  of $\gamma_{i}^{2}$ becomes arbitrarily large. Now
  suppose that $u_{i}^{2}$ converges to a
  non-constant affine function $u^{2}$ on the limit
  $B^{2}.$ Note that all the (horo)-balls
  $B_{\gamma_{i}(t)}^{2}(t),$ centered at
  $\gamma_{i}(t)$ and containing $x_{i}^{2}$ at the
  boundary, also converge to flat balls in the limit,
  for any given $t$. (Here we recall from the
  definition of $\rho ,$ c.f.(1.27), that
  $\rho^{2}(x) \geq  \dist_{g_{i}^{2}}(x, \partial
  B_{x_{i}^{1}}^{2}(R_{i})).$ It follows that the
  functions $u_{i}$ approximate the unique
  non-constant affine function $u^2$ in
   $B^2_{\gamma(t)}(t)$ extending the limit
   function $u^2$ in $B^2.$ 
  From the assumption (3.57) on the choice
  of $x_i^2,$ it follows that $|u_i -1| \leq 2\mu_o$
   everywhere in $B^2_{\gamma_i(t)} (t)$
   and hence the same holds on the limit.
  But clearly a non-constant affine function 
  cannot have this property. This proves second claim above.
 
  Now $(B^2, g^2, u^2)$ is an $L^{2,2}$ static vacuum solution, and from the above, $(D^2, g^2)$ is flat, with $u^2$ a positive constant. Since, from the discussion in \S 1.3, $L^{2,2}$ static vacuum solutions are real-analytic away from the 0-locus of the potential, it follows that $(B^2, g^2)$ itself is flat, and $u^2$ is constant in $B^2$, as claimed.
 
  Thus the geometry of the sequence
  $(B_{x_{i}^{2}}^{2}(1), g_{i}^{2}, x_{i}^{2})$ and
  that of the potential function $u_{i}$ on
  $B_{x_{i}^{2}}^{2}(1)$ is exactly the same as on
  the previous sequence $(B_{x_{i}^{1}}^{1}(1),
  g_{i}^{1}, x_{i}^{1}).$ Clearly one may repeat this
  process an arbitrary finite number of times; just
  return to (3.56)-(3.58) and raise the indices
  successively by~1. However, each limit $g^{k}$
  constructed in this way is flat with constant
  potential function and hence the sequence
  $\{x_{i}^{k}\}$ is not $(\rho , c)$ buffered, for
  any fixed $c > 0.$ Thus, no finite iteration of this process leads to a $(\rho, c)$ buffered sequence.
 \end{remark} 
    It is worthwhile to examine the construction in
  Remark 3.9 on a specific example.

 \smallskip
 {\bf Example 3.9}
    Let $(N, g)$ be the Schwarzschild metric (0.17), the
  canonical solution of the static vacuum equations.
  Let $x_{j}^{1}$ be a divergent sequence of points
  in $N$, so that $t(x_{j}^{1})  \rightarrow  \infty ,$
  where $t$ is the distance to the event horizon.
  Note that the curvature $z$ of $(N, g)$ decays on the
  order of $t^{-3}$ as $t \rightarrow  \infty $ and
  $u(x_{j}^{1}) \rightarrow  1  =  \max u.$
 
  Now blow-down the metric $g$ based at
  $\{x_{j}^{1}\},$ i.e. form the `tangent cone at
  infinity' based at $\{x_{j}^{1}\},$ by considering
  the metrics $g_{j}^{1} = \rho
  (x_{j}^{1})^{-2}\cdot  g.$ The metrics $g_{j}^{1}$
  converge to the flat metric on $B_{x^{1}}^{1}(1),
  x^{1}= \lim x_{j}^{1}$ and the functions
  $u_{j}^{1}$ given by the restriction of $u$ to
  $B_{x_{j}^{1}}^{1}(1)$ converge to the constant
  function $1.$
 
  Of course the sequence $(B_{x_{j}^{1}}^{1}(1),
  g_{j}^{1}, x_{j}^{1})$ is not $(\rho ,$ c)
  buffered, for any $c > 0.$ 
  A short computation, using the cubic curvature 
  decay and the definition of $\rho ,$ shows that 
  $$\rho (x_{j}^{1}) \sim
 t(x_{j}^{1}) -  (t(x_{j}^{1}))^{1/3}. $$
 
  Thus, as in the construction in Remark 3.9,
  consider a second sequence of points $x_{j}^{2},$
  with $t(x_{j}^{2}) \sim  (t(x_{j}^{1}))^{1/3}.$
  Again if one blows-down the metric $g$ based at
  $\{x_{j}^{2}\},$ then the metrics $g_{j}^{2} = \rho
  (x_{j}^{2})^{-2}\cdot  g$ converge to the flat
  metric and the functions $u_{j}^{2},$ now given by
  the restriction of $u$ to $B_{x_{j}^{2}}^{2}(1)$
  again converge to the constant function 1.
 
  This process may be repeated any finite number of
  times, so that by induction, points $x_{j}^{k}$ are
  chosen satisfying $$t(x_{j}^{k}) \sim
 (t(x_{j}^{1}))^{1/3^{k}}.  $$ However one always
 obtains flat limits in this manner. Observe that for
 any given choice of $k$, the functions $u_{j}^{k}$
 defined as above, always converge to the constant
 function 1. In fact $u(x_{j}^{k}) \sim \,  1 -
 t^{-1}(x_{j}^{k}) \rightarrow  1$ as $j \rightarrow
 \infty ,$ for any fixed $k.$ It is clear that to
 obtain a non-flat limit, (namely the original
 Schwarzschild metric), one must increase the descent
 level $k$ in $\{x_{j}^k\}$ as a function of $j$,
 with $k = k(j) \rightarrow  \infty $ as $j
 \rightarrow  \infty $. One obtains a non-flat limit when the base points $x_{i}^{k}$ satisfy $u(x_{i}^{k}) = 1 - \delta$, for any fixed $\delta > 0$. In this particular example, even though $k(j) \rightarrow \infty $ as $j \rightarrow \infty $, the potential function $u_{j}^{k(j)}$ does not converge to the 0-function as $j \rightarrow \infty$.

 \medskip 
 
  We now return to the main issue of constructing a
  $(\rho ,c)$ buffered sequence from some initial
  sequence whose blow-up limits are flat, with
  constant potential function.
 \begin{theorem}[\bf Existence of Non-Flat Blow-ups] \label{t 3.10.} 
   Let $\{g_{i}\}$ be a sequence of unit volume
   Yamabe metrics on a closed 3-manifold $M,$ satisfying the non-collapse
   assumption (3.43) and scalar curvature bound
   (3.41), i.e.  
 	  \begin{equation} \label{e3.59}
 	   -\infty   <   - s_{o} \leq  s_{g_{i}} \leq  - s_{1}  <   0.  
 	   \end{equation} Suppose there are points
 $x_{i}\in (M, g_{i})$ and a number $\delta_{o} > 0$ 
 such that 
 	\begin{equation} \label{e3.60} \rho
 	(x_{i}) \rightarrow  0, 
 	 \ \ {\rm as }\ i\rightarrow\infty ,
 	 \ \ {\rm and}\ \ 
         u_{i}(x_{i})/T_{i} \geq  \delta_{o} >  0,
 	\end{equation}
  where $\displaystyle T_{i}\equiv  \max_{M} u_{i}$ satisfies
  (3.45), i.e. for some $\lambda_{o} >  0,$
 \begin{equation} \label{e3.61} T_{i} \geq
 \lambda_{o}.  \end{equation}

  Then for any fixed $c >  0$ sufficiently small,
  there exists a sequence of $(\rho ,c)$-buffered
  base points $\{y_{i}\},$ satisfying the following
  conclusions: $u_{i}(y_{i}) \in  (0,T_{i}),
  \dist_{g_{i}}(x_{i}, y_{i}) \rightarrow  0,$ and the
  rescalings
 	\begin{equation} \label{e3.62} 
         g_{i}'  = \rho_{i}^{-2}(y_{i})\cdot  g_{i}, \ \ \ 
         \rho_{i}(y_{i}) \rightarrow  0, 
       \end{equation}
  of $\{g_{i}\}$ sub-converge, in the weak $L^{2,2}$
  topology based at $\{y_{i}\}$ on a domain $D,$ to a
  limit smooth metric $g' $ on $D,$ based at $y'  =
  \lim y_{i}.$ The pair (D, $g' )$ is a non-flat,
  (and {\it non-super-trivial),} solution to the static
  vacuum Einstein equations, with limit potential
  function $\Roof{u}{\bar}$ constructed from
  $\{u_{i}\}.$ The limit domain $D$ is naturally
  embedded in $M$ via the metric $g' $ and the
  convergence is in the strong $L^{2,2}$ topology
  away from the locus $\{\bar{u} = 0\}.$
 \end{theorem} 
\begin{pf}
  The proof proceeds in several steps, organized
  around several lemmas. Recall again that the
  estimate (3.61) follows from a bound on the $L^{2}$
  norm of $z^{T}$ under the bounds (3.59), c.f.
  (3.45).
 
  First, as previously in Proposition 3.1 and Remark
  3.9 for example, it is necessary to renormalize
  $u_{i}$ to $\Roof{u}{\bar}_{i} = u_{i}/T_{i},$ at
  least if $T_{i} \rightarrow  \infty .$ We will
  assume this is done, so that $\sup \bar{u}_{i} = 1,$ and neglect the bar
  notation on $\bar{u}_{i}.$
 
\bigskip
\noindent

{\bf Step I.}
  To begin, the sequence $\{x_{i}\}$ satisfying
  (3.60) must be altered to a new sequence of points
  $\{x_{i}^{1}\},$ still satisfying (3.60), but with
  (locally) near-maximal u-values, (compare with
  (3.47)).
 
 \begin{lemma} \label{l 3.11.}
   {\bf (Initial Degeneration Points).} 
   Under the assumption (3.60), there exists a subsequence
  $\{i'\}$ of \{i\}, relabeled to \{i\}, and base
  points $x_{i}^{1},$ with $\dist_{g_{i}}(x_{i},
  x_{i}^{1}) \rightarrow  0$ which still satisfy
  (3.60), such that for any point $q_{i}\in
  B_{x_{i}^{1}}(\rho (x_{i}^{1})),$ 
  \begin{equation} \label{e3.63}
  u_{i}(q_{i}) \leq  u_{i}(x_{i}^{1}) + \varepsilon_{i}, 
  \end{equation} 
  for some sequence
 $\varepsilon_{i} = \varepsilon (x_{i}^{1}) \rightarrow  0$
 as $i \rightarrow  \infty .$
 \end{lemma} 
\begin{pf}
  The proof is based on the fact that $\rho $ is a
  Lipschitz function with Lipschitz constant 1, c.f.
  (1.27). Thus given a sequence $x_{i}$ satisfying
  (3.60), suppose (3.63) does not hold in
  $B_{x_{i}}(\rho (x_{i})).$ Then there exist
  $q_{i}\in B_{x_{i}}(\rho (x_{i}))$ such that
  $u(q_{i}) \geq  u(x_{i}) + \varepsilon_{o},$ for some
  $\varepsilon_{o} >  0.$ 
  But $\rho (q_{i}) \leq 2\rho (x_{i}) \rightarrow  0$ 
  as $i \rightarrow \infty ,$ so that $\{q_{i}\}$ also satisfies
  (3.60). Clearly this procedure may be repeated any
  given number of times. Since $\max u_{i} = 1,$ it
  follows that for any given $\varepsilon_{o} >  0,$
  there exist $x_{i}^{\varepsilon_{o}}$ such that, for
  all $i$ sufficiently large, (3.63) holds with
  $\varepsilon_{o}$ in place of $\varepsilon_{i}$ and $\rho
  (x_{i}^{\varepsilon_{o}}) \rightarrow  0$ as $i
  \rightarrow  \infty .$ Now replace $\varepsilon_{o}$
  by a sequence $\varepsilon_{j} \rightarrow  0$ and
  choose a suitable diagonal subsequence $j_{i}$ of
  $\{i,j\}$ with $\varepsilon_{j_{i}} \rightarrow  0$
  sufficiently slowly as $i \rightarrow \infty.$
  Then the base points $x_i^1\equiv x_i^{\varepsilon_{j_i}}$
  satisfy (3.60) and (3.63).\\ 
 \end{pf}

  From now on, we work with the subsequence
  $\{x_{i}^{1}\}$ from Lemma 3.11. Of course the
  choice of $\{x_{i}^{1}\}$ may be far from unique.
  Even more, sequences satisfying (3.60) may also be
  far from unique. The initial points $x_{i}^{1}$
  should not be confused with the points in (3.47),
  whose existence was assumed and not proved. Lemma
  3.11 is the only place where the hypothesis (3.60)
  is used in the proof of Theorem 3.10.
 
  Now, as in (3.48), consider the rescaled metrics
 \begin{equation} 
       \label{e3.64} g_{i}^{1} = \rho (x_{i}^{1})^{-2}\cdot  g_{i}.  
 \end{equation}
  Given the hypotheses (3.59) and (3.61), the
  discussion preceding Proposition 3.6 shows that a
  subsequence of $\{g_{i}^{1}\}$ converges on
  $B_{x_{i}^{1}}^{1}(1)$ to a limit, which is a non
  super-trivial solution $(B^{1}(1), g^{1}, u^{1}, x^{1})$ 
  of the static vacuum equations. (With the exception 
  of the non-collapse hypothesis (3.43),
  this is the only place where the remaining
  hypotheses, (3.59) and (3.61), are used in the proof
  of Theorem 3.10). For simplicity, we restrict ourselves in the following to this convergent subsequence. From (3.63), one sees as in
  Proposition 3.6, c.f.(3.53), that the limit is
  flat, with $u^{1} = const.$ If necessary, we
  renormalize $\{u_{i}\}$ again from now on so that
  $u^{1} \equiv  1$ in the limit.
 
  In the following steps, we show that for each $i$
  sufficiently large, there is a first $k_{o} =
  k_{o}(i)$ and points $x_{i}^{k_{o}}\in U^{k_{o}}$
  such that the sequence of points $x_{i}^{k_{o}}$
  are $(\rho ,$ c) buffered in the sequence 
  $(M, g_{i}),$ for a specified small $c > 0.$ 
  Here $U^{k}$ and $L^{k}$ are defined as preceding Remark
  3.9. This is a quantitative version of the argument
  in Remark 3.9, and rests basically on the fact that
  $(M, g_{i}, u_{i})$ is smooth, for each fixed $i.$
 
  The construction of such points will take place with a choice for the parameter $c$ for a $(\rho, c)$ buffered sequence, and a choice for a parameter $\varepsilon$ measuring the distance of the blow-up metrics to a static vacuum solution. As will be seen in the proof, the choice of $c$ and $\varepsilon$ is arbitrary, provided each is sufficiently small, and further $\varepsilon$ is sufficiently small, depending on $c$. While it is possible to make specific numerical choices for $c$ and $\varepsilon$, the work involved is cumbersome and of no specific value, and so we will not do so. The final specification of $c$ and $\varepsilon$ will be made in Step IV below, after a preliminary specification in Step II.

  From the discussion following (3.64), given an $\varepsilon > 0$, there is an $i_{o} = i_{o}(\varepsilon )$ such that for all $i \geq
  i_{o},$ the blow up $(B_{x_{i}^{1}}^{1}(1),
  g_{i}^{1}, x_{i}^{1})$ is $\varepsilon$-close to a
  flat static vacuum limit. By $\varepsilon$-close, we mean that the $L^2$ norm of $L^{*}u_{i}$ on $B_{x_{i}^{1}}^{1}(1)$  is less than $\varepsilon$. Since $u_{i}$ is $C^{o}$ (in fact $L^{2,2})$ close to 1
  at points in $B_{x_{i}^{1}}^{1}(1)$ an arbitrary fixed
  distance away from the boundary, this implies that the $L^2$ norm of the curvature $r$ is less than $\varepsilon'$ = $\varepsilon'(\varepsilon)$ in balls an arbitrary but fixed distance away from the boundary. In particular, by Theorem 3.4, $\rho^{1}$ is very small, (depending on $\varepsilon$), somewhere near the boundary of $B_{x_{i}^1}^{1}(1)$. (Further, by
  Theorem 3.3, some portion of the level set $L^{2}$
  is very close to the boundary of
  $B_{x_{i}}^{1}(1).$
 
\bigskip
\noindent
{\bf Step II.}
  For any given $i \geq i_o = i_o(\varepsilon)$ fixed, and base points $x_{i}^{1}$ as above, we now inductively construct new base points $x_{i}^{2}$, $x_{i}^{3}$, etc, chosen to satisfy quantitative versions of (3.56)-(3.57). For a given choice of $c$, $c_o = 10^{-3}$ and $\nu_o = 10^{-2}$, let $d = d(c, c_o, \nu_o)$ be the strong buffer constant determined by $c$, $c_o$ and $\nu_o$, as discussed following Lemma 3.8. For the moment, (until Step IV), we require only that $c$ be chosen sufficiently small so that $d \leq \frac{1}{5}$ and $d + 2c \leq \frac{1}{4}$.

\smallskip
\noindent
  {\bf $k^{\rm th}$-Level Inductive Hypothesis:}
  For any given $i \geq i_o$, suppose base points $x_{i}^{j}\in U^{j}$, have been constructed for each $j$, $1 \leq j \leq k$, up to a given $k \geq 2$. We require that each $x_{i}^{j}$ is not $(\rho ,c)$ buffered for $j \leq  k- 1$, and to be chosen to satisfy the following conditions:

 \begin{equation} \label{e3.65} 
        x_{i}^{j}\in B_{x_{i}^{j-1}}(\rho (x_{i}^{j-1})), 
  \end{equation}
 	\begin{equation} \label{e3.66} 
 	c\cdot \rho (x_{i}^{j-1}) 
 	\leq  
 	\dist_{g_{i}}(x_{i}^{j}, 
 	 \partial B_{x_{i}^{j-1}}(\rho (x_{i}^{j-1}))) \leq  2c\cdot
 	\rho (x_{i}^{j-1}), 
 	\end{equation}
 \begin{equation} \label{e3.67} 
 \rho (x_{i}^{j}) \leq  d_{1}\cdot \rho (x_{i}^{j-1}), d_{1} = d + 2c, \ \ {\rm and}
 \end{equation} 
  \begin{equation} \label{e3.68} 
 |u_{i}(x) - u_{i}(x_{i}^{j-1})| 
 \leq  
 	\tfrac{1}{4}\cdot u_{i}(x_{i}^{j-1}), 
 	 \ \ \ \ \ 
 	 \forall x\in B_{x_{i}^{j-1}}((1- c)\rho (x_{i}^{j-1})).  
 \end{equation}

  The conditions (3.65)-(3.66) on the placement of $x_{i}^{j}$, given $x_{i}^{j-1}$, can obviously always be satisfied. With regard to the next two conditions, condition (3.68) on the behavior of the potential, and the fact that $x_{i}^{j-1}$ is not $(\rho, c)$ buffered, would imply that $x_{i}^{j-1}$ is not strongly $(\rho, d)$ buffered provided Lemma 3.8 is applicable. If this is the case, it follows that there is a point $z_{i}^{j} \in B_{x_{i}^{j-1}}(\rho(x_{i}^{j-1}))$, with $\rho(z_{i}^{j}) \leq d\cdot \rho(x_{i}^{j-1})$, and hence by (1.27), there is a point $x_{i}^{j}$ satisfying (3.65)-(3.66) together with (3.67). In other words, (3.67) is a consequence of (3.65), (3.66) and (3.68) to the extent that Lemma 3.8 is applicable. We will show below that indeed it is, even when $u_{i}(x_{i}^{j-1})$ is very small.

  It is important to note that the terminal base point $x_{i}^{k}$ may or may not be $(\rho ,c)$ buffered. Further, by the paragraph preceding Step II, there exist base points $x_{i}^{2}$ satisfying (3.65)-(3.68) with $k = 2$
  whenever $i \geq  i_{o}$, so that one may start the induction process at $k = 2$.
 
  To normalize the discussion, we scale the metrics
  $g_{i}$ at each $x_{i}^{j}$ so that $\rho_{i}^{j}
  = 1,$  i.e. set $g_{i}^{j} = \rho
  (x_{i}^{j})^{-2}\cdot  g_{i},$ and renormalize
  $u_{i}$ by setting 
  \begin{equation} \label{e3.69}
 u_{i}^{j} = u_{i}/u(x_{i}^{j}).  
 \end{equation} 
 With this normalization, (3.65)-(3.68) are equivalent to
 the statements 
 \begin{equation} \label{e3.70}
 x_{i}^{j}\in B_{x_{i}^{j-1}}^{j-1}(1),
 \end{equation}
	 \begin{equation} \label{e3.71} 
	 c \leq \dist_{g_{i}^{j-1}}(x_{i}^{j}, 
	 \partial B_{x_{i}^{j-1}}^{j-1}(1)) 
	 \leq  2c, 
	 \end{equation}
	 \begin{equation} \label{e3.72} 
	 \rho^{j-1}(x_{i}^{j})
	 \leq  d_{1}, \ \ {\rm and}, 
	 \end{equation} 
	 \begin{equation}
		 \label{e3.73} 
		 |u_{i}^{j-1}(x) -  1| \leq  \tfrac{1}{4},\ \ \ \ \ 
		 \forall x\in B_{x_{i}^{j-1}}^{j-1}(1- c).
         \end{equation}
 
  The inductive step is the following:

\smallskip  
\noindent 
  {\bf Inductive Step.} Suppose $x_{i}^{k}$ is not $(\rho ,c)$
  buffered. Then there exist base points
  $x_{i}^{k+1}$ satisfying the conditions
  (3.65)-(3.68), with $k+1$ in place of $k$.
 
  The next two steps are concerned with the proof of the inductive step.

\smallskip
\noindent
{\bf Step III.}
  In order to carry out the inductive step from $k$
  to $k+1,$ we need to analyse the geometry of
  $g_{i}^{k}$ in the unit ball
  $B_{x_{i}^{k}}^{k}(1).$ Note that the following
  Lemma does not assume any buffer condition on
  ${x_{i}^{k}}$, as a sequence in $i$, for any $k$ fixed.
 
 \begin{lemma} \label{l 3.12.}
   For any $i \geq i_{o}$ and some $k \geq 2$, suppose base points $x_{i}^{k}$ have been constructed satisfying the $k^{\rm th}$-level
   inductive hypothesis. Then the triple
   $(B_{x_{i}^{k}}^{k}(1), g_{i}^{k}, u_{i}^{k})$
   satisfies the estimates
	 \begin{equation} \label{e3.74} 
	   u_{i}^{k}r_{i} = D^{2}u_{i}^{k} + o(i,k), 
	  \ \ \ \ \ \ 
	  \Delta u_{i}^{k} = o(i,k), 
	 \end{equation}
  where $o(i,k)$ denotes terms which become arbitrarily
  small in $L^{2}(B_{x_{i}^{k}}^{k}(1)),$ (in fact in
  $L^{2}(M, g_{i}^{k})),$ if either $i$ or $k$ is
  sufficiently large.  
  \end{lemma} 
\begin{pf}
   We return to the equation (2.10) for $u = u_{i},$
   i.e.  $$ur -  D^{2}u = -\Delta u\cdot  g 
            + \xi  + \frac{s}{3}\cdot  g, $$
  applied to $g_{i}^{k}$ in $B_{i}^{k}.$ Divide, i.e.
  renormalize, this equation by $\upsilon_{i,k}
  \equiv  u_{i}(x_{i}^{k}).$ By (3.66) and (3.68),
  $\upsilon_{i,k} \geq  2^{-k}.$ Let $u_{i}^{k} =
  u_{i}/\upsilon_{i,k}$ as in (3.69), so that
  $u_{i}^{k}(x_{i}^{k}) = 1.$ This gives
 	\begin{equation} \label{e3.75} u_{i}^{k}r -
 	D^{2}u_{i}^{k} = -\Delta u_{i}^{k}\cdot  g +
 	\upsilon_{i,k}^{-1}\cdot \xi  +
 	\upsilon_{i,k}^{-1}\cdot \frac{s}{3}\cdot  g.
 	\end{equation}
  Now we claim that the term
  $\upsilon_{i,k}^{-1}\cdot \xi $ is arbitrarily
  small in $L^{2},$ if either $i$ or $k$ is
  sufficiently large. To see this, recall from
  Theorem 2.2 that 
     $||\xi||_{L^{2}(M, g_{i})} \leq C.$ 
  As noted before in (3.12), in any rescaling
  $g_{i}'  = \rho_{i}^{-2}\cdot  g_{i},$ one has then
 	\begin{equation} \label{e3.76} 
 		  ||\xi||_{L^{2}(M, g_{i}' )} 
 		  \leq  
 		  C\cdot \rho_{i}^{1/2}.
 	\end{equation} 
 Now by (3.72), in the sequence
 $g_{i}^{j},$ for $i$ fixed and any 
 $j  <   k,$ we have 
       $$\rho_{i}^{j} \leq  d_{1}\cdot \rho_{i}^{j-1}, $$
 so that 
       $$\rho_{i}^{k} \leq  d_{1}^{k}\cdot \rho_{i}^{1}. $$ 
 It follows that 
      $$||\upsilon_{i,k}^{-1}\xi||_{L^{2}(M, g_{i}^{k})}
          \leq  
 	  C\cdot  2^{k}\cdot  d_{1}^{k/2}\rho
 	(x_{i}^{1})^{1/2}  \rightarrow  0, 
     $$
  as either $i$ or $k \rightarrow  \infty ,$ since 
  $d_{1} = d + 2c  \leq   1/4$. The same
  reasoning applies to the last term in (3.75)
  involving the scalar curvature, since the scalar
  curvature of $(M, g_{i})$ is uniformly bounded in
  $L^{\infty}.$ Thus, the last two terms in (3.75)
  are arbitrarily small, if either $i$ or $k$ is
  sufficiently large.
 
  Similarly, from Proposition 2.9, the $L^{2}$ norm
  of $\Delta u$ is uniformly bounded on $(M, g_{i})$ 
  and hence again the same reasoning implies
  that the $L^{2}$ norm of $\Delta u_{i}^{k}$
  goes to 0 uniformly on $(M, g_{i}^{k}),$ as either 
  $i \rightarrow  \infty $ or $k \rightarrow  \infty .$ 
  Thus the right hand side (3.75) is arbitrarily
  small whenever either $i$ or $k$ is sufficiently
  large.\\
\end{pf}
 
  Lemma 3.12, together with the non-collapse
  assumption (3.43) implies that for any $i \geq i_o = i_{o}(\varepsilon)$ and $k \geq 2$, the metrics
  $g_{i}^{k}$ are
  $\varepsilon$-close to a static vacuum solution in the sense that the metrics satisfy the static vacuum equations with a (volume normalized) $L^2$ error of at most $\varepsilon$. Observe that this $\varepsilon$ is the same as the $\varepsilon$ chosen initially at the end of Step I. By passing to limits, it follows that in fact $g_{i}^{k}$ on $B_{i}^{k}$ is $\varepsilon' = \varepsilon`(\varepsilon)$ close in the weak $L^{2,2}$ topology to an actual static vacuum solution; however, we will not use this fact. By the
  inductive hypothesis, the predecessors 
  $g_{i}^{j},\ j \leq  k- 1,$ 
  are not $(\rho ,c)$ buffered at their base points, and so
  are all $\delta_{1}$-close to a flat metric in
  $B_{i}^{j}(1- c)$ in the sense that

\begin{equation} \label{e3.77}
\bigl( \int_{B_{i}^{j}(1-c)} |r|^2 \bigr)^{1/2} \leq \delta_{1};
\end{equation}
here $\delta_{1} = \delta_{1}(c, c_o) = 2\pi c\cdot c_{o}$. (Again it can be shown that this implies that $(B_{i}^{j}(1-c), g_{i}^{j})$ is $\delta_{1}'$-close, $\delta_{1}' = \delta_{1}'(\delta_{1}, \nu_o)$, in the weak $L^{2,2}$ topology, to a flat metric, but this will not be used). Clearly, $\delta_{1}$ may be made arbitrarily small by choosing the initial buffer parameter $c$ sufficiently small. Further, by (3.73), observe that the potential function $u_{i}^{j}$ is $\frac{1}{4}$-close to a constant function in the $C^{0}$ topology in $B_{i}^{j}(1-c)$.

\bigskip
\noindent

{\bf Step IV.}
  We are now in position to prove the inductive step. The next result is the main point in this respect. For this, we restrict somewhat further the initial choice of $c$, and $\varepsilon$. Thus, we require for the remainder of the proof that $\varepsilon < \varepsilon_{o}$, where $\varepsilon_{o}$ is the constant from Lemma 3.8, and also $\varepsilon < \frac{1}{10}\delta_{1}$. In addition, $c$ will be chosen sufficiently small so that $d_{1} = d + 2c \leq \frac{1}{100}$, and $\delta_{2} = 2c_{s}\cdot \delta_{1} \leq \frac{1}{100}$, where $c_{s}$ is the Sobolev constant of the embedding $L^{2,2} \subset C^{0}$ in ${\Bbb R}^3$, (c.f. (3.79) below).

 \begin{lemma} \label{l 3.13.}
    Under the conditions above on $c$ and $\varepsilon$, for any given $i \geq i_{o} = i_{o}(\varepsilon)$, suppose the base point $x_{i}^{k}$ is not $(\rho ,c)$ buffered. Assume also that the non-collapse assumption (3.43) holds. Then the potential function $u_{i}^{k}$ is $\frac{1}{4}$-close
    to the constant function 1 on
    $B_{x_{i}^{k}}^{k}(1- c)$ in the $C^{0}$
    topology, i.e. (3.73) holds with $k+1$ in place
    of k.
 \end{lemma} 
\begin{pf}
  The proof is essentially the same as the proof of
  the similar issue in Remark 3.9, c.f. (3.58). Since
  $x_{i}^{k}$ is not $(\rho ,c)$ buffered and (3.43)
  holds, the metric $g_{i}^{k}$ is both $\varepsilon$-close 
  to a static vacuum solution, and
  $\delta_{1} = \delta_{1}(c)$ close to a flat
  metric, in the sense defined above, on
  $B^{k}(1- c) \equiv  B_{x_{i}^{k}}^{k}(1- c)$.

Let $v_{i}^{k} = sup_{B^{k}(1-c)}|u_{i}^{k}|$, and recall that $u_{i}^{k}(x_{i}^{k}) = 1$. It follows from Lemma 3.12 that

\begin{equation} \label{e3.78}
\bigl(\int_{B^{k}(1-c)}|D^{2}u_{i}^{k}|^{2} \bigl)^{1/2} \leq 
\bigl(\int_{B^{k}(1-c)}(u_{i}^{k})^{2}|r|^{2} \bigl)^{1/2} + \varepsilon
\leq 2\delta_{1} \cdot v_{i}^{k}.
\end{equation}
\noindent
  This, together with Sobolev embedding, c.f. [GT, Ch.7], shows that the potential $u_{i}^{k}$ is close in the $C^0$ topology to an affine function $\alpha = \alpha_{i,k}$ on $B^{k}(1-c)$, in that

\begin{equation} \label{e3.79}
 	  |u_{i}^{k} -  \alpha|
 	\leq  
 	\delta_{2}v_{i}^{k}.  
 	\end{equation} 
\noindent
  The constant $\delta_{2}$ is given by $\delta_{2} = 2c_{s}\cdot \delta_{1} \leq \frac{1}{100}$, where $c_{s}$ is the Sobolev embedding constant; $c_{s}$ is a fixed numerical constant in dimension 3, independent of all parameters. This estimate also, of course, applies on domains inside $B^{k}(1-c)$, with the same $\alpha$ and with $v_{i}^{k}$ replaced by the supremum of $|u_{i}^{k}|$ on the domain. In particular, it follows that $sup|u_{i}^{k}|$ is close to $sup|\alpha|$ on any such domain.

  Now to understand more precisely the form of $\alpha$ and the size of $v_{i}^{k}$, we first consider the validity of (3.79) on the predecessor domain, in the $g_{i}^{k}$ scale. Thus, consider the previous point $x_{i}^{k-1},$ of
  distance $R_{i}^{k}$ from $x_{i}^{k}$ in the
  $g_{i}^{k}$ scale. By (3.71)-(3.72), one computes
  $d_{1}^{-1}(1- 2c) \leq  R_{i}^{k} \leq  c^{-1}(1- c),$
  the important point being that $R_{i}^{k}$ is large, 
  since $d_{1}$ is small. Consider also the
  ball $B^{k-1}(R_{i}^{k}) \equiv
  B_{x_{i}^{k-1}}^{k}(R_{i}^{k})$ containing a
  portion (roughly half) of the ball $B^{k}(1- c).$
  Using the inductive assumption, the $L^2$ norm of the curvature of $g_{i}^{k}$ on $B^{k-1}(R_{i}^{k})$ is much smaller than $\delta_{1}$, since $g_{i}^{k}$ is much larger than $g_{i}^{k-1}$. This implies that (3.79) holds on $B^{k-1}(R_{i}^{k})$, i.e.

\begin{equation} \label{e3.80}
  |u_{i}^{k} - \alpha'| \leq \delta_{2}w_{i}^{k},
\end{equation}
where $\alpha'$ is an affine function on $B^{k-1}(R_{i}^{k})$ and $w_{i}^{k}$ is the supremum of $|u_{i}^{k}|$ on $B^{k-1}(R_{i}^{k})$. We may write $\alpha' = a' + b't$, where $t$ is an affine Euclidean coordinate with $t(x_{i}^{k}) = 0$.

    By the induction hypothesis (3.73), the function
 	$$u_{i}^{k} =
 	\frac{u_{i}(x_{i}^{k-1})}{u_{i}(x_{i}^{k})}u_{i}^{k-1}
 	$$
  is $2\cdot \frac{1}{4}$ close to the constant function 1 in
  $B^{k-1}(R_{i}^{k})$, and hence $w_{i}^{k} \leq \frac{3}{2}$. Together with (3.80), this implies that $\alpha' = a' + b't$ satisfies $|b'| \leq  (R_{i}^{k})^{-1} < 2d_{1}$. Thus the oscillation of $\alpha'$ is small. Further $|1 - \alpha'| \leq 2\delta_{2}$.

   It follows that on the 'half-ball' $D_{i}^{k} \equiv B^{k}(1-c) \cap B^{k-1}(R_{i}^{k})$, the function $u_{i}^{k}$ satisfies

\begin{equation} \label{e3.81}
 |u_{i}^{k} - 1| \leq 2(\delta_2 + d_1).
\end{equation}

   Now return to the estimate (3.79). It follows from (3.81) that $sup_{D^{k}}|u_{i}^{k}|$ is close to $1$. Hence, by (3.79) and the discussion in the paragraph following it, one may calculate that $v_{i}^{k} = sup_{B^{k}(1-c)}|u_{i}^{k}| \leq 5$. It then follows from the estimates (3.78)-(3.79) applied to the union $B^{k}(1-c) \cup B^{k-1}(R_{i}^{k})$ that $\alpha$ and $\alpha'$ are close. Working out concretely the degree of closeness then gives the estimate

\begin{equation} \label{e3.82}
  |u_{i}^{k} - 1| \leq 10(\delta_2 + d_1),
\end{equation}
on $B^{k}(1-c)$. The result then follows from the specifications of $\delta_2$ and $d_1$ above.\\
\end{pf}
 
  We are now in position to argue the inductive
  step.
 
 \begin{lemma} \label{l 3.14.}
   Suppose $x_{i}^{k}$ satisfies the $k^{\rm th}$-level inductive
   hypothesis (3.65)-(3.68). If the base point
   $x_{i}^{k}$ is not $(\rho ,c)$ buffered, then
   there exist points $x_{i}^{k+1}$ satisfying the $(k+1)^{\rm st}$-level inductive hypothesis, i.e.
   (3.65)-(3.68) with $k+1$ in place of $k$
   everywhere.
 \end{lemma} 
\begin{pf}
  Given the previous work, this is now essentially
  obvious. There are clearly (many) base points $x_{i}^{k+1}$ satisfying (3.65)-(3.66) with $k+1$ in place of $k$. Since $x_{i}^{k}$ is not $(\rho ,c)$
  buffered, by Lemma 3.13 $u_{i}^{k}$ is $\frac{1}{4}$-close to the constant function 1 in $B_{x_{i}^{k}}^{k}(1- c)$, so that (3.68), (or (3.73)) holds, with $k+1$ in place of $k$. By Lemma 3.12, or more precisely its contrapositive, $x_{i}^{k}$ is not strongly $(\rho, d)$ buffered, and so there exist $x_{i}^{k+1}$ such that (3.67) also holds, with $k+1$ in place of $k$.\\
\end{pf}
 
{\bf Step V.}
  Now, since $(M, g_{i})$ is a smooth compact Riemannian
  manifold, $\rho (M, g_{i})$ cannot be arbitrarily
  small for any given $i \geq i_o$ i.e. there exists $\rho_{o} =
  \rho_{o}(i)$ such that 
   $\rho (M, g_{i}) \geq \rho_{o} >  0.$ 
  Starting with the initial sequence 
  $x_{i}^{1}$ from Lemma 3.11, and any $i \geq i_o$, the inductive
  process above constructs points $x_{i}^{k}, k \geq 2,$ 
  satisfying $\rho (x_{i}^{k}) \leq  d_{1}\rho (x_{i}^{k- 1}),$ 
  and continues indefinitely as long
  as $x_{i}^{k}$ is not $(\rho ,c)$ buffered. It
  follows that there is a first level $k_{o} =
  k_{o}(i)$ such that $y_{i} \equiv  x_{i}^{k_{o}}$
  is $(\rho ,c)$ buffered. Observe that the discussion in Remark 3.9 proves that necessarily $k_{o}(i) \rightarrow \infty$ as $i \rightarrow \infty$. The sequence $(B_{i}^{k_{o}}(1), g_{i}^{k_{o}}, y_{i},
  u_{i}^{k_{o}})$ satisfies the equation (3.75), with
  the right side going to 0 uniformly in $L^{2}(M, g_{i}^{k_{o}})$ 
  as $i \rightarrow  \infty .$ 
  Thus, (c.f. the discussion following Definition 3.7), a
  subsequence converges to a non-flat limit solution
  $(D, g' )$ to the static vacuum equations, with
  (possibly renormalized) potential function
  $\bar{u} = \lim u_i^{k_o}.$ (Here we are
  using again the non-collapse assumption (3.43)).
  The limit domain $D$ may be chosen for instance
  to be an open domain in the limit
  $B_{y'}(1)$ of $B_i^{k_o}(1)$ which has smooth and compact closure in $B_{y'}(1)$. This completes
  the proof of Theorem 3.10.\\
\end{pf}
 
  We conclude this subsection with a number of
  remarks.
 
 \begin{remark} {\bf (i).}\label{r 3.15.(i).}
   At this point, if not well before, one may wonder
   why not just blow-up the sequence $\{g_{i}\}$ at
   base points $\{q_{i}\}$ realizing the minimum
   $\rho (M)$ of $\rho (x), x\in\{(M, g_{i})\},$ in
   which case $\{q_{i}\}$ is strongly $(\rho , 1)$
   buffered. In this case, a limit of the metrics 
 $$
  g_{i}'  = \rho (M)^{-2}\cdot  g_{i} 
 $$
  will automatically be complete and at least
  $L^{2,2}$ smooth.
 
  The main problem is that, from the work above,
  there is every indication that this minimum occurs
  at, or arbitrarily near the 0-levels of $u.$ The
  $z$-splitting of $\{g_{i}'\}$ has the form
 $$u_{i}z_{i}'  = (D^{2})' u_{i} -
 \Delta_{i}' u_{i} + \xi_{i}'  -
 \frac{s_{i}'}{3}\cdot  g_{i}' . $$ 
 Although, as before, the last three terms 
 on the right go to 0, it may well happen that 
 $u_{i}$ converges to 0 uniformly when blowing up 
 at $\{q_{i}\},$ so that in
 the limit this equation gives just $0 = 0,$ i.e. a
 super-trivial solution, and thus no information
 whatsoever regarding the limit metric. Examples
 exhibiting exactly this behavior will be discussed
 in \S 6. Further, without the descent construction
 as in the proof of Theorem 3.10, there may be no
 means to renormalize $u$ to obtain a non-trivial
 solution.
 
 {\bf (ii).}
  It is clear that all the preceding arguments are
  unchanged if $u$ is replaced by $- u,$ so that the
  construction of Theorem 3.10 may be carried out
  either ``down'' the $u$-levels from a local maximum, as
  in (3.63), or ``up'' the $u$-levels from a local minimum,
  provided such a local minimum is bounded away from
  zero.
 
 {\bf (iii). }
  We point out that it is not possible in the
  construction above to descend down to other levels
  besides the 0-level, e.g. trying to descend only to
  the level $u = \frac{1}{2}$ by renormalizing the
  differences $\frac{1}{2}+2^{-k-1}$ and
  $\frac{1}{2}+2^{-k}.$ This is due essentially to
  the linear, and not affine, nature of the equation
  (3.71) in $u.$
 
 {\bf (iv).}
  The construction in Theorem 3.10 has been referred
  to as a `descent' construction down the levels of
  $u = u_{i},$ since this is what occurs in all known
  examples, (as in Example 3.9), and since one
  certainly descends to $u$-levels less than $\max
  u_{i}/T_{i} = 1.$ However, the construction itself
  does not guarantee such a descent, i.e. that
 	\begin{equation} \label{e3.83} 
 	u_{i}(x_{i}^{j})
 	 <   u_{i}(x_{i}^{j-1}) \ \ \ \ \ {\rm for} \ j \leq  k_{o}.
 	\end{equation} 
 Only the control (3.68) or (3.73) is
 obtained on the behavior of $u.$
 
  On the other hand, although it has not been used in
  the proof, it is worth noting that at any stage in
  the inductive construction, the level sets $L^{s}$
  of $u_{i}^{j},\ j  <   k_{o},$ for any given $s$
  with $0  <   s  <   10^{-1},$ are somewhere
  $\delta  = \delta (d)$ close to $\partial
  B_{x_{i}^{j}}^{j}(1)$ in the $g_{i}^{j}$ metric;
  here $\delta $ may made arbitrarily small by
  choosing $d = d(c)$ sufficiently small. This
  follows from the fact that $x_{i}^{j}$ is not
  $(\rho ,c)$ buffered, and hence not strongly $(\rho
  ,d)$ buffered, by Lemma 3.8.
 
  Using this remark, it is possible to define
  variants of the inductive construction in Theorem
  3.10 for which (3.83) does hold, and for which the
  same conclusions are valid. Such variations will
  not be used here however.  \end{remark} 
 
 \begin{remark} {\bf (i).}\label{r 3.16(i).}
   The construction in Theorem 3.10 is very local.
   Thus, given any initial sequence $\{x_{i}^{1}\}$
   satisfying the conclusions of Lemma 3.11, there is
   a buffered sequence $\{y_{i}\}$ satisfying the
   conclusions of Theorem 3.10 very close to
   $x_{i}^{1};$ in fact
   $\dist_{g_{i}}(x_{i}^{1},y_{i}) \leq  2\rho
   (x_{i}^{1}).$
 
  However, there may be many possible (essentially)
  distinct choices of the initial sequence
  $\{x_{i}^{1}\}$ in $(M, g_{i})$ satisfying (3.60)
  and (3.63). Note also that the choice of $y_{i}$ is
  by no means unique, since the choice of the
  intermediate points $x_{i}^{j}$ may allow
  considerable freedom.
 
 {\bf (ii).}
  The limit solution $(D, g' )$ in Theorem 3.10 is
  only locally defined, and further may live in the
  region of $M$ where $u$ is approximately 0, since
  the potential function $\Roof{u}{\bar}$ is obtained
  by renormalizing, possibly infinitely many times,
  the original sequence $\{u_{i}\}.$ The geometry of
  the limit solutions $(D, g' )$ will be studied in
  more detail in \S 5.  \end{remark} 
 
 \begin{remark} \label{r 3.17.}
   The proof of Theorem 3.10 requires Lemma 3.11 to
   obtain the initial sequence $\{x_{i}^{1}\}$ whose
   blow-up limit is a flat static vacuum solution,
   with limit potential function $u^{1} = const.  >  0.$ 
   On the other hand, if one is given some
   initial sequence $\{p_{i}^{1}\}$ satisfying this
   latter property, (not necessarily satisfying the
   local maximum condition (3.63)), then the
   construction in Theorem 3.10 may be applied
   starting at $\{p_{i}^{1}\},$ without any further
   changes, to reach the same conclusions.
 \end{remark}

  There are instances where the hypothesis (3.60) in
  Theorem 3.10 can be deduced or replaced by suitable
  assumptions on the local behavior of the functions
  $u_{i}.$ We give two examples of this below.
 
 \begin{proposition} \label{p 3.18.}
   Let $(M, g_{i})$ be a sequence of unit volume
   Yamabe metrics satisfying the non-collapse
   assumption (3.43). Suppose there exists a $\delta >  0,$ 
   a sequence $r_{i} \rightarrow  0,$ and a
   sequence of base points $\{x_{i}\}\in (M, g_{i}),$
   such that on $B_{i} = B_{x_{i}}(r_{i}),$
 \begin{equation} \label{e3.84} 
     \osc_{B_{i}} (u_{i}/T_{i}) \geq  \delta. 
       \end{equation} 
       Then (3.60) holds, with $\delta $ in place of
 $\delta_{o}.$
 \end{proposition} 
\begin{pf}
  It suffices to show that $\rho (x_{i}) \rightarrow 0$ 
  as $i \rightarrow  \infty .$ If there exists
  $\rho_{o} > 0$ such that $\rho (x_{i}) \geq \rho_{o},$ 
  (in a subsequence), then $L^{2}$
  estimates on the trace equation (2.13), c.f. [GT,
  Thm.8.8], imply that $\{u_{i}/T_{i}\}$ is uniformly
  bounded in $L^{2,2}(B_{x_{i}}(\rho_{o}/2)).$
  Sobolev embedding then implies that $u_{i}/T_{i}$
  is uniformly bounded in $C^{1/2},$ so that in
  particular the oscillation of $u_{i}/T_{i}$ on
  $B_{x_{i}}(r)$ is bounded by $r^{1/2}.$ This
  contradicts (3.84).\\
\end{pf}
 
  Somewhat more generally, the hypotheses (3.60) and
  (3.61) in Theorem 3.10 can be replaced by the
  following local condition on the behavior of $u =
  u_{i}$ near its 0-level set. Let $t(x) =
  \dist_{g_{i}}(x, L^{o}).$ If $u_{i} > 0$
  everywhere, $L^{o}$ must be replaced by
  $L^{\varepsilon_{i}},$ for a suitable sequence
  $\varepsilon_{i} \rightarrow  0.$
 
 \begin{proposition} \label{p 3.19.}
   Let $(M, g_{i})$ be a sequence of unit volume
   Yamabe metrics satisfying (3.59) and the
   non-collapse assumption (3.43). Suppose there
   exists a sequence of base points $\{x_{i}\}\in (M, g_{i}),$ 
   with $t(x_{i}) \rightarrow  0$ as $i \rightarrow  \infty ,$ 
   such that 
      \begin{equation}
              \label{e3.85} u_{i}(x_{i}) \gg  t^{1/2} 
              \ \ \ \ \ {\rm as}\  i
	 \rightarrow  \infty , 
      \end{equation} and
 \begin{equation} \label{e3.86} 
     t|\nabla log u_{i}|(q_{i})
       <   
     \tfrac{1}{2}, 
 \end{equation}
  for all $q_{i}\in B_{x_{i}}(K_{i}t(x_{i}))$
  satisfying $\frac{1}{2}t(x_{i}) \leq  t(q_{i})
  \leq  K_{i}t(x_{i});$ here $K_{i}$ is any  given
  sequence with $K_{i} \rightarrow  \infty $ as $i
  \rightarrow  \infty $ and it is assumed that $t$
  achieves the value $K_{i}/2$ in
  $B_{x_{i}}(K_{i}t(x_{i})).$ Then the conclusions of
  Theorem 3.10 hold for some $(\rho ,c)$ buffered
  sequence $y_{i}$ with $dist_{g_{i}}(x_{i},y_{i})
  \rightarrow  0$ as $i \rightarrow  \infty .$
 \end{proposition} 
\begin{pf}
  First, note that as in Proposition 3.18, (3.85)
  forces $\rho (x_{i}) \rightarrow  0,$ since if
  $\rho (x_{i})$ were bounded away from 0, $L^{2}$
  elliptic estimates applied to the trace equation
  (2.13) as above imply that $u_{i}$ is bounded in
  $L^{2,2} \subset  C^{1/2},$ which violates (3.85).
 
  Consider the behavior of the pointed sequence 
  $(M, g_{i}' , x_{i}, \Roof{u}{\bar}_{i}), g_{i}'  =
  t(x_{i})^{-2}\cdot  g_{i}.$ We claim that Lemma
  3.12 holds on $B_{x_{i}}' (1).$ Thus, for
  $\Roof{u}{\bar}_{i} = u_{i}/u(x_{i}),$ (compare
  with (3.75)), the renormalized equation for $u$
  takes the form 
 \begin{equation} \label{e3.87}
 	 \bar{u}_{i}r -  D^{2}\bar{u}_{i} =
 	-\Delta \Roof{u}{\bar}_{i}\cdot  g +
 	u(x_{i})^{-1}\cdot \xi  + u(x_{i})^{-1}\cdot
 	\frac{s}{3}\cdot  g.  
 \end{equation} 
  As in (3.76),
 in this scale, we have
 $$u(x_{i})^{-1}||\xi||_{L^{2}} 
 \leq  C\cdot u(x_{i})^{-1}t(x_{i})^{1/2} \rightarrow  0,  \ \ \ 
 {\rm as} \ i \rightarrow  \infty , 
 $$ 
 where the last estimate follows from (3.85). 
 The same reasoning as in Lemma 3.12 shows that 
 the full right side of (3.87) also
 goes to $0$, as $i \rightarrow  \infty ,$ which proves
 the claim.
 
  Next we claim that there exists a constant
  $\rho_{o} >  0$ such that 
  \begin{equation} \label{e3.88} 
  \rho (x_{i}) \geq  \rho_{o}\cdot t(x_{i}).  
 \end{equation} 
 For in the scale $g_{i}' $
 where $t' (x_{i}) = 1,$ note that
 $\Roof{u}{\bar}_{i}(x_{i}) = 1,$ and by (3.86),
 $\Roof{u}{\bar}_{i}$ is bounded away from 0 in the
 region where $\frac{1}{2} \leq  t'  \leq  2;$ of
 course here $t' (p) = \dist_{g_{i}'}(p, L^{o}).$
 Hence, (3.88) follows from Theorem 3.3 applied to
 the triple $(M, g_{i}' , \Roof{u}{\bar}_{i}).$ Since
 $\rho' (x_{i})$ is bounded away from 0, blow-up
 limits based at $x_{i}$ exist. Further, the
 preceding discussion shows that all blow-up limits
 based at $x_{i}$ converge to a static vacuum
 solution with potential function $\Roof{u}{\bar}.$
 
  Now if a blow-up limit $(M, g_{i}' , x_{i})$
  happens to be non-flat, we are done (the assertion
  is proved in this case). Suppose instead that the
  blow-up limit is flat. Then the limit function
  $\Roof{u}{\bar}$ is either a non-constant affine
  function or $\Roof{u}{\bar}\ =\ const.\, = 1.$ Hence
  it suffices to prove that the latter case holds,
  since then all the arguments following Lemma 3.11
  can be carried over without change, c.f. Remark
  3.17.
 
  Suppose the limit $\Roof{u}{\bar}$ is a
  non-constant affine function. The product (3.86) is
  scale-invariant, and hence on the limit, one
  obtains 
     $$t'|\nabla log\bar{u}|(q) \leq \tfrac{1}{2}, $$ 
  for all $q$ in $B_{x_{i}}' (K_{i})$
  satisfying $t' (q) \geq  \frac{1}{2}.$ 
  The assumption on $t$ and $K_i$ implies that on the limit, $t' $
  is an unbounded function. However, it is obvious
  that no non-constant affine function satisfies this
  bound.\\ \end{pf}

  To conclude this subsection, it is clear that
  Theorem 3.10 proves Theorem A(II) in case 
  $\sigma (M) <   0,$ or more precisely in case the scalar
  curvature of $\{g_{i}\}$ satisfies (3.59).
 
  In the following two subsections, we complete the
  proof of Theorem A(II) and Theorem 3.10 in the
  cases $\sigma (M) = 0 $ and $\sigma (M) >  0.$ 
   As will be seen, this requires some, but as it turns out, no truly essential changes.
 
\subsection{} 
   Suppose $$\sigma (M) = 0, $$ or more precisely
 suppose $\{g_{i}\}$ is a sequence of unit volume
 Yamabe metrics on $M$ satisfying \begin{equation}
 \label{e3.89} s_{g_{i}}\rightarrow  0.
 \end{equation} As before, we require that the
 assumptions (3.42) and (3.43) hold uniformly on
 $\{g_{i}\}.$

  Of course, if $\lambda $ does not converge to 0, so
  that (3.45) holds, then all of the results of \S
  3.4 still hold, with identical proofs, at least
  when $s_{g_{i}} \leq  0.$ However, as discussed in
  Remark 2.12, one often expects that $\lambda
  \rightarrow  0$ when (3.89) holds. In fact, if
  $\sigma (M) >  0 $ and $s_{g_{i}} \rightarrow  0,$ 
  then as discussed following Proposition 2.9,
  there are many unit volume Yamabe metrics $g_{i}$
  which do not degenerate at all, and for which $f
  \rightarrow  - 1,\  u \rightarrow  0,$ and $\xi
  \rightarrow  0$ smoothly.
 
  Thus, suppose that $\lambda  =
  \lambda_{i}\rightarrow  0$ on $(M, g_{i}).$
  Assuming $\{g_{i}\}$ degenerates in the sense of
  (0.8), since $u$ may go to $0$ in $L^{2},$ one might
  not expect to produce any non-trivial blow-up
  solutions of the static vacuum equations by descent
  down the $u$-levels. However, when $\lambda
  \rightarrow  0,$ we also have $\delta  = -\int k
  \rightarrow  1,$ by (2.17) and (2.28). In
  particular, $\sup(- k)$ is uniformly bounded away
  from 0 on $\{g_{i}\},$ compare with (3.61).
 
  Given these preliminaries, we prove Theorem 3.10
  with the assumption (3.89) replacing (3.59).
 
 \begin{pf*}{\bf Proof of Theorem 3.10: $\sigma (M) = 0$.}
  \
 
  We assume that $\{g_{i}\}$ is a degenerating
  sequence of unit volume Yamabe metrics on
  $U_{\delta_{o}}\subset M $ , satisfying (3.89),
  together with the $L^{2}$ bound on $z^{T}$ (3.42),
  (or (0.19)), and the non-collapse assumption
  (3.43). For simplicity, assume first that
 	\begin{equation} \label{e3.90} 
 		   s_{g_{i}} \leq  0.
 	\end{equation}
 
  The splitting equation (2.25) for $z^{T}$ may be
  written in the form 
  \begin{equation} \label{e3.91}
        kr = D^{2}k -\Delta k\cdot  g -  z^{T} + \xi .  
   \end{equation}
  By (3.42), (and also (2.11)), the term $z^{T}-\xi $
  is uniformly bounded in $L^{2},$ as is the term
  $\Delta k,$ (by taking the trace of
  (2.25)).
 
  This equation, and its associated trace equation,
  now have exactly the same form and bounds as the
  $u$-equation (2.10), with $- k$ in place of $u$ and
  $z^{T}-\xi $ in place of $\xi +\frac{s}{3}g.$ As
  noted above, the mean value of $- k,$ namely
  $\delta ,$ converges to 1 as $i \rightarrow  \infty
  .$
 
  Exactly the same analysis as in \S 3.3 and \S 3.4
  may now be carried out on the function $- k,$ or
  more precisely $- k/\sup(- k)$ in place of $u/\sup u.$
   Note that from (2.28),
 \begin{equation} 
       \label{e3.92} \frac{- k}{\sup(- k)} = \frac{u}{\sup u}, 
 \end{equation}
  so that the main hypothesis (3.60) holds with $- k$
  in place of $u.$ Similarly, analogous to (3.61),
  $\sup(- k)$ is uniformly bounded away from 0, since
  its mean value converges to 1. Thus, one produces a
  $(\rho ,c)$ buffered sequence with associated limit
  giving a non-flat static vacuum solution, with
  limit potential function $\Roof{k}{\bar}$ coming
  from the geometry of $\{k_{i}\}.$ In fact, by
  (3.92), the potential function can also be obtained
  from renormalizations of the $\{u_{i}\}.$
 
  The proof in case 
 	\begin{equation} \label{e3.93}
 	s_{g_{i}} \geq  0,  \ \
       {\rm and}\ \ s_{g_{i}} \rightarrow  0,
 	\end{equation} 
 is exactly the same as in the case
 (3.90) with one exception. Namely, if $s_{g_{i}} >  0,$ 
  it is possible that $\Ker L^{*} \neq  0$ on
 $(M, g_{i})$ so that the functions $u$ or $k$ may
 not be uniquely defined. In this situation, one just
 takes any choice for $- k$ in $\{- k+\Ker L^{*}\}.$
 Since the choice is always renormalized by its
 supremum, this has no effect on any of the
 arguments. This situation will be discussed in more
 detail in \S 3.6, where it occurs more naturally.
 
  This completes the proof of Theorem 3.10 and
  Theorem A in case $\sigma (M) = 0,$ i.e. (3.89)
  holds.\\
 \end{pf*}
 
\subsection{} 
  Finally, we turn to the case where 
 $$\sigma (M) >  0, $$ 
 or more generally where $\{g_{i}\}$ is a
 sequence of unit volume Yamabe metrics with
 	\begin{equation} \label{e3.94} 
 	s_{g_{i}} \geq  s_{o} > 0.  
 	\end{equation}

  As noted in \S 2, such metrics may satisfy
 \begin{equation} \label{e3.95} 
 \Ker L^{*} \neq  0 .
 \end{equation} 
 For instance, on the standard sphere
 $S^{3},$ the $1^{\rm st}$ eigenfunctions of the
 Laplacian form exactly $\Ker L^{*},$ see also \S 6.4.
 Note that the condition (3.95) implies, by taking
 the trace, that 
 \begin{equation} \label{e3.96} 
  - s/2 \in  \Spec\Delta .  
  \end{equation} 
 The condition (3.95) is equivalent to the statement that
 the linearization $L$ of $s$ at a Yamabe metric $g$
 is not surjective onto $C^{\infty}(M, {\Bbb R} ),$
 and has been studied by several authors, c.f. [Bo],
 [Kb2], [La]. In particular, the last two authors
 show that there are many (even conformally flat)
 Yamabe metrics satisfying (3.95) which are not
 Einstein.
 
  Nevertheless, even for metrics satisfying (3.95),
  all the splittings (2.6), (2.10) (2.21) and (2.25)
  are valid and defined as before. The associated
  functions $f,$ $u$ and $k$ however are obviously not
  uniquely defined; they are unique only modulo $\Ker L^{*}.$ 
  Related to this, as noted following
  Proposition 2.9, there is no apriori bound on the
  $T^{2,2}$ norm or even the $L^{2}$ norm of $f.$
 
  On the other hand, if $g$ is a Yamabe metric with
  $s_{g} >  0,$ then by the definition of $g$ - as
  realizing $\inf {\cal S} $ in the conformal class
  $[g]$ - one has 
  \begin{equation} \label{e3.97}
     s_{g}\cdot  \left(
 		\int\psi^{6}dV_{g}
 		\right)^{1/3} 
 		 \leq  
 		 \int (8|d\psi|^{2} + s\psi^{2})\, dV_{g}, 
 		 \end{equation}
  for any positive smooth function $\psi $ on $M,$ cf.
  [Bes,4.28], [LP]. It follows that (3.97) holds for
  all smooth functions $\psi $ on $M.$ Thus on $(M, g)$
  one has a Sobolev inequality, uniform on any
  sequence of Yamabe metrics for which $s_{g} > 0$
  is uniformly bounded away from $0.$
 
  This has two immediate consequences. First, we
  have:
 \begin{lemma} \label{l 3.20.}
    Let $(M, g_{i})$ be a sequence of unit volume
    Yamabe metrics on $M$ satisfying (3.94). Then
    there is a constant $\nu_{o} = \nu_{o}(s_{o})$
    such that the volume radius satisfies
 	\begin{equation} \label{e3.98} 
 	\nu_{i}(x) \geq \nu_{o} , 
 	\end{equation} 
 for all $x \in (M, g_{i}).$ Thus, the sequence $\{g_{i}\}$ cannot
 collapse anywhere, c.f. \S 1.4.  
 \end{lemma} 
\begin{pf}
  It is a standard fact that the Sobolev inequality
  (3.97) gives rise to a uniform lower bound on the
  volumes of geodesic balls in comparison to the
  volumes of Euclidean balls, i.e. to the ratio
  $\vol(B_{x}(r))/r^{3},$ provided 
  $\vol B_{x}(r)  <  \frac{1}{2}\, \vol (M, g_{i});$ 
  we refer to [Ak] for example for a proof. 
  By definition of the volume radius, this gives (3.98). 
  Similarly, it follows from the definition of collapse in \S 1.4 that
  (3.98) prevents collapse of the sequence $(M, g_{i})$ 
  at any sequence of base points.\\
\end{pf}
 
  Hence, in the case of uniformly positive scalar curvature, 
  condition (i) of Theorem A is automatically satisfied.
 
  Second, from the trace equation (2.12),
 \begin{equation} 
  \label{e3.99} 2\Delta 
  f + sf = \tr\xi , 
 \end{equation} 
 and the fact that the
 right side is uniformly bounded in $L^{2},$ see
 Theorem 2.2, it follows from standard elliptic
 estimates using only the Sobolev inequality (3.97),
 c.f. [GT, Thm.8.15], that
 \begin{equation} \label{e3.100} 
 \sup \, |f| 
 \leq c(s_{o})
 \left\{ \left(
  \int f^{2}\, dV_{g})^{1/2}
  \right) +1\right\} .
 \end{equation}
  Of course, the $T^{2,2}$ norm of $f$ is also
  bounded by the $L^{2}$ norm of $f;$ see the proof of
  Proposition 2.9.
 
\medskip
  We now complete the proof of Theorem 3.10, with
  (3.94) in place of (3.59).
 
\smallskip

\noindent
\begin{pf*}{\bf Proof of  Theorem 3.10: $\sigma (M) >  0$.}
  \  

  Assume that $\{g_{i}\}$ is a degenerating sequence
  of unit volume Yamabe metrics on $U_{\delta_{o}}
  \subset  M,$ satisfying (3.94), together with the
  $L^{2}$ bound on $z^{T}$ (3.42), (or (0.19)). By
  Lemma 3.20, no non-collapse assumption is needed.
 
  As noted in Remark 2.12, the bound (3.94) implies
  that $\lambda  = \lambda_{i}$ is bounded away from
  0, and hence $T_{i} = \sup u_{i} \geq \lambda_{o},$ 
  for some $\lambda_{o} > 0;$ thus
  one has the bound (3.61) in this context.
 
  If $\Ker L^{*} = 0$ for all $g_{i},$ then all of the
  previous arguments in \S 3, in particular the proof
  of Theorem 3.10, hold without any changes.
 
  Suppose instead that $\Ker L^{*} \neq  0$ on some
  subsequence of $\{g_{i}\}.$ In this case, simply
  make some choice for the function $u = u_{i} \mod
  \Ker L^{*}.$ Note that $\lambda  = \int u$ is
  independent of any choice. As indicated above in \S
  3.5, recall that throughout the previous work in \S
  3, the potential $u$ is initially renormalized by
  its maximum, i.e. only the behavior of $u/{\sup u}$ is
  considered. By (3.100) this renormalization is
  equivalent to the renormalization of $u$ by its
  $L^{2}$ norm. Of course if $||u_{i}||_{L^{2}}
  \rightarrow  \infty ,$ then the initial structural
  equation (2.10) becomes
 $$L^{*} \left(
         \frac{u_{i}}{||u_{i}||_{L^{2}}}
         \right)
 \rightarrow  0, 
 \ \ {\rm in} \  L^{2}(M, g_{i}), 
  \ \ {\rm as}\  i
 \rightarrow  \infty , $$ 
 which is even stronger than
 the initial $L^{2}$ bound on $L^{*}(u),$ compare
 with the renormalization in the proof of Proposition
 3.1. The proofs of all the previous results in \S 3
 follow exactly as before, without any changes.
 
  This completes the proof of Theorem A and Theorem
  3.10 in all cases.\\
\end{pf*}
 \begin{remark} \label{r 3.21.}
   We point out that the assumption that the metrics
   $\{g_{i}\}$ in Theorem A or Theorem 3.10 are
   defined on a fixed 3-manifold $M$ has not been
   used. In fact, these and all previous results hold
   under the same assumptions on arbitrary sequences
   of closed oriented Riemannian 3-manifolds $(M_{i},
   g_{i}).$ \end{remark} 
 
 \section{\bf Remarks on the Hypotheses.}
 \setcounter{equation}{0}
  In this section, we make some further comments on
  the hypotheses of Theorem A. First, we consider the
  non-collapse hypothesis (i) in Theorem A, i.e.
  (0.18) or more generally (3.43), and afterwards the
  degeneration hypothesis (iii) in Theorem A, i.e.
  (0.20). The hypothesis (ii), i.e. the $L^{2}$ bound
  on $z^{T}$ is already discussed in Theorem
  2.10-Remark 2.12. 

\subsection{} 
   By Lemma 3.20, the non-collapse assumption (3.43)
   is only needed in the case $\sigma (M) \leq  0.$
   From some perspectives however, this assumption is
   perhaps the most crucial in Theorem~3.10 or Theorem~A.
 
  The basic difficulty in handling the case of
  collapse as opposed to convergence is already seen
  in the discussion in \S 3.2, for instance in the
  proof of Proposition 3.1. Thus, consider a blow-up
  sequence $g_{i}'  = \rho (x_{i})^{-2}\cdot  g_{i},$
  with $\rho (x_{i}) \rightarrow  0.$ If this
  sequence collapses, that is 
  \begin{equation}
        \label{e4.1} 
    \nu (x)  <  <   \rho (x),
 \end{equation} 
 then it can be proved that the collapse is along a sequence of injective
 F-structures on $B_{x_{i}}' (1-\varepsilon_{i}),$ for
 some sequence $\varepsilon_{i} \rightarrow  0$ as $i
 \rightarrow  \infty .$ We may then unwrap the
 collapse, (i.e. resolve the degeneration), by
 considering the universal covers
 $\widetilde{B}_{x_{i}}' (1-\varepsilon_{i}).$ This
 sequence does not collapse anywhere. Thus, a
 subsequence of $(\widetilde{B}_{x_{i}}'
 (1-\varepsilon_{i}), g_{i}, x_{i})$ converges to a
 limit $L^{2,2}$ Riemannian manifold
 $(\widetilde{B}_{x}' (1), g' , x);$ see [An2, \S
 3] for further details here.
 
  The $L^{2}$ curvature radius is essentially
  invariant under coverings. However, in the process
  of passing to the universal cover, we lose the
  property that
 	\begin{equation} \label{e4.2}
 	\int|\xi_{i}'|^{2}dV_{i}'  = \rho (x_{i})\cdot
 	\int|\xi_{i}|^{2}dV_{i} \rightarrow  0,
 	\end{equation}
  see (3.12). In passing to the universal cover, the
  apriori $L^{2}$ bound (2.11) on $\xi $ is lost. The
  property (4.2) was crucial in the case of
  convergence. The fact that $\xi'  \rightarrow  0$
  in $L^{2}$ as in Proposition 3.1 implies that the
  limit is a solution of the static vacuum Einstein
  equations. Without this property, the limit does
  not satisfy any particular equation.
 
  Recall that, under the assumption (3.43),
 	\begin{equation} \label{e4.3} 
 	\nu_{i} \geq \nu_{o}\cdot \rho_{i}, 
 	\end{equation}
  the previous results of \S 3 hold in
  $U_{\delta_{o}}$ for an {\it  arbitrary}  sequence
  of Yamabe metrics $\{g_{i}\}$ with a bound on the
  $L^{2}$ norm of $z^{T}.$ Even though the sequence
  $\{g_{i}\}$ satisfies no particular equation (PDE),
  (besides being a Yamabe metric), the blow-up limits
  modeling degenerations of $\{g_{i}\}$ do satisfy a
  strong PDE, namely the static vacuum equations. It
  seems quite unreasonable to believe this without
  some assumption like (4.3). Note that the integral
  bound $\int|\xi|^{2} \leq  \frac{s^{2}}{3}$ from
  (2.11), which gives rise to (4.2) in blow-ups,
  becomes less and less meaningful in regions where
  $(M, g_{i})$ is more and more collapsed, since the
  bound may then come primarily from the volume
  collapse and not reflect any particular behavior of
  $\xi $ itself.
 
  This difficulty can be overcome if one had apriori
  $L^{\infty}$ control on $\xi $ in place of an
  $L^{2}$ bound, (since such bounds are invariant
  under coverings), or if the sequence of Yamabe
  metrics satisfies some other (stronger) P.D.E. This
  latter will be the path taken in [AnII].

\subsection{} 
  Next we make some remarks on the hypothesis (iii)
  of Theorem A. Suppose $\{g_{i}\}$ is a degenerating
  sequence of unit volume Yamabe metrics on $M,$ so
  that $||z_{g_{i}}||_{L^{2}(M)} \rightarrow  \infty .$ 
  It is natural to ask if then (iii) is satisfied.
  In general, that is for arbitrary sequences, the
  answer is no, see \S 6 for further discussion.
 
  Suppose for instance $\{g_{i}\}$ is a unit volume
  maximizing sequence for the functional
  $v^{2/3}\cdot  s,$ so that \begin{equation}
 \label{e4.4} s_{g_{i}} \rightarrow  \sigma (M).
 \end{equation} We will see in \S 7, c.f. Theorem
 7.2 and Lemma 7.4, that any such sequence can be perturbed
 slightly if necessary, for instance in the $T^{2,2}$
 topology, so that the sequence is still maximizing,
 and so that the gradient $- z^{T} = \nabla
 (v^{2/3}\cdot  s|_{{\cal C}})$ goes to 0 in the
 natural dual topology, i.e. at $g_{i},$
 	\begin{equation} \label{e4.5}
 	||z^{T}||_{T^{-2,2}(T{\cal C} )}  
 	=  \sup\left\{ \left| 
 	  \int_{M} \langle z^{T}, \alpha\rangle dV
 	  \right| : \alpha\in T{\cal C}  
 	  {\rm\  and\ } ||\alpha||_{T^{2,2}}\leq 1
 	  \right\}
 	\rightarrow  0, 
 \end{equation}
  as $i \rightarrow  \infty .$ In other words, the
  sequence $\{g_{i}\}$ is Palais-Smale for the
  functional ${\cal S}|_{{\cal C}}$ w.r.t. the
  $T^{2,2}$ metric on ${\Bbb M}$.     In particular,
  there are many such Palais-Smale sequences.
 
  Consider now sequences $\{g_{i}\}\in{\cal C}_{1}$
  which satisfy the stronger condition that
 	\begin{equation} \label{e4.6}
 		  ||z^{T}||_{T^{-2,2}(T{\Bbb M} )}  \rightarrow   0;
 	\end{equation}
  in other words, the condition (4.5), but where
  $\alpha $ is no longer constrained to lie in
  $T{\cal C} .$ Such sequences will be called {\it
  strongly Palais-Smale}. It follows easily from the
  definition that then also 
 \begin{equation} \label{e4.7} 
 \tr z^{T} \rightarrow  0 \ \ \ \ \ { \rm in\ } T^{-2,2}.
 \end{equation} Under these circumstances, one has
 the following result.
 \begin{proposition} \label{p 4.1.}
    Let $\{g_{i}\}$ be a sequence of unit volume
    Yamabe metrics satisfying the estimate (4.7), and
    suppose $\sigma (M)  <   0,$ or just 
    $s_{g_{i}} \leq  s_{o}  <   0.$ Then for the associated
    functions $u = u_{i},$ one has
 \begin{equation} \label{e4.8}
 \left|\left|\tfrac{u-\lambda}{\lambda}\right|\right|_{L^{2}} 
  \rightarrow 0,  \ \ {\rm as \ \, }i \rightarrow  \infty , 
 \end{equation}
  where the $L^{2}$ norm is w.r.t. the metric $g =
  g_{i}.$
 \end{proposition} 
\begin{pf}
   Consider the trace equation (2.13) for u, that is,
   using (2.33), 
  \begin{equation} \label{e4.9}
 2\Delta \left(u - \lambda \right) + s\left(u-\lambda \right) =
 \lambda \tr z^{T}.  
 \end{equation} Set $v =
 \frac{u-\lambda}{\lambda}$ and let $\phi  =
 \phi_{i}$ be the solution on $(M, g) = (M, g_{i})$
 to the equation \begin{equation} \label{e4.10}
 \Delta\phi  -  \phi  =
 \frac{v}{||v||_{L^{2}}}.  \end{equation}
  Since 1 is not in the spectrum of $\Delta ,
  (\Delta $ is a non-positive operator), this
  equation has a unique solution. To estimate the
  $T^{2,2}$ norm of $\phi ,$ pair equation (4.10)
  with $\phi $ and integrate by parts to obtain
 \begin{equation} \label{e4.11} 
       \int|d\phi|^{2} + \int\phi^{2} = -  \int\phi  \frac{v}{||v||_{L^{2}}}
              \leq  \left(\int\phi^{2}\right)^{1/2}.  
 \end{equation}
 Further, squaring both sides of (4.10) gives
 \begin{equation} \label{e4.12} \int
 (\Delta\phi )^{2} + 2|d\phi|^{2}+\phi^{2} =
 1.  \end{equation} Thus, the $T^{2,2}$ norm of $\phi
 $ is bounded by 1.
 
  It follows that $\phi $ is an admissible test
  function in (4.7), so that from (4.9), one obtains
      $$\int 2\phi\Delta v + \phi sv \rightarrow 0. $$ 
   Integrating by parts gives 
   $$\int 2v  \left[  \frac{v}{||v||_{L^{2}}}+\phi 
 	       \right]
 	       + \phi sv \rightarrow  0, 
   $$ or, 
  \begin{equation} \label{e4.13} 
  2||v||_{L^{2}} + \int (2+s)\phi v
 \rightarrow  0.  
 \end{equation} 
 But from (4.10)
     $$||v||_{L^{2}}
     \left(\int|d\phi|^{2}+\int\phi^{2}\right) 
     = - \int\phi v, 
     $$ 
   so that, 
  $$2||v||_{L^{2}} -
 (2+s)||v||_{L^{2}}||\phi||_{L^{1,2}} \rightarrow  0.
 $$ 
 This implies 
 \begin{equation} \label{e4.14}
 ||v||_{L^{2}}(2 -  (2+s)||\phi||_{L^{1,2}})
 \rightarrow  0.  
 \end{equation} 
 Thus either
 $||v||_{L^{2}}\rightarrow  0,$ as required, or
 	\begin{equation} 
                \label{e4.15} ||\phi||_{L^{1,2}} \rightarrow  \frac{2}{2+s}.  
         \end{equation}
  This of course implies $s > - 2,$ and since
  $\sigma (M)  <   0,$  $\frac{2}{2+s} > \frac{2}{2+\sigma (M)} > 1,$ 
  (for all $i$),
  which with (4.15) contradicts (4.11). Hence the
  result follows.\\
\end{pf}
 
  Note that this argument makes strong use of the
  $T^{2,2}$ norm; it is not (likely to be) valid
  w.r.t. the $L^{2,2}$ norm. It is quite easy to see
  that there are counterexamples to Proposition 4.1
  when $\sigma (M) >  0,$ c.f. \S 6.4. If $\sigma
  (M) = 0,$ or more precisely $\{g_{i}\}$ is a
  sequence of Yamabe metrics with $s_{g_{i}}
  \rightarrow  0,$ then this result is borderline.
  Namely, the same argument as above gives either the
  conclusion (4.8), or (from (4.15),
  $||d\phi||_{L^{2}} \rightarrow  0$ and
  $||\phi||_{L^{2}} \rightarrow  1,$  with $\int\phi = 0.$ 
  The latter case is of course impossible if
  $\lambda_{1},$ the lowest non-zero eigenvalue of
  $\Delta ,$ is bounded away from 0, but
  without some assumption of this kind, it is not
  clear if (4.8) follows from (4.7).
 
  Now in case (4.8) holds, it is easy to see that
  condition (iii) of Theorem A follows from the
  non-collapse assumption (0.18) or (3.43) and the
  general degeneration assumption (0.8), (in place of
  (0.20)). In fact, (compare with Lemma 3.11 and
  (3.47)), we have
 \begin{proposition} \label{p 4.2.}
    Let $(M, g_{i})$ be a sequence of unit volume
    Yamabe metrics satisfying the non-collapse
    assumption (3.43) together with (4.8).

  If $\{g_{i}\}$ degenerates on $M,$ i.e. if there are
  points $p_{i}\in (M, g_{i})$ such that $\rho
  (p_{i}) \rightarrow  0,$ then there are points
  $x_{i}\in (M, g_{i})$ such that 
 	\begin{equation}
 	\label{e4.16} |u_{i}(x_{i})/T_{i} -  1| \rightarrow 0, 
 	 \ \ {\rm as}\ i \rightarrow  \infty , 
        \end{equation} 
 and
 	\begin{equation} \label{e4.17} \rho (x_{i})
 	\rightarrow  0, 
 	 \ \ {\rm as \ } i \rightarrow  \infty ,
 	\end{equation}
  where $T_{i} = \sup u_{i}.$
 \end{proposition} 
\begin{pf}
  Note that given (4.8), we do not assume $\sigma (M)  <   0. $
  Note also that $u / T = ( \frac{u}{\lambda}) / (\frac{T}{\lambda}).$
 
  First, suppose $T/\lambda  = T_{i}/\lambda_{i}
  \rightarrow  \infty .$ We claim that any sequence
  $\{x_{i}\}$ satisfying (4.16) must then satisfy
  (4.17). For if not, then for some $\rho_{o}  >  0 $
   and for some sequence $x_{i}$ as above,
 \begin{equation} \label{e4.18} 
 \rho (x_{i}) \geq \rho_{o}.  
 \end{equation}
 From (4.18) and the non-collapse assumption, it
 follows that $B = B_{x_{i}}(\rho_{o})$ has $L^{2,2}$
 bounded geometry. Thus, by Theorem 2.10, $z^{T}$ is
 uniformly bounded in $L^{2}.$ Consider the trace
 equation (4.9),
 \begin{equation} \label{e4.19} 
 2\Delta \bigl(\frac{u-\lambda}{\lambda}\bigr) 
               + s\bigl(\frac{u-\lambda}{\lambda}\bigr) 
 = \tr z^{T}.
 \end{equation} 
 The right side of (4.19) is uniformly
 bounded in $L^{2},$ so that $L^{2}$ elliptic
 estimates [GT, Thm.8.8] imply a bound
 \begin{equation} \label{e4.20} 
 ||v||_{L^{2,2}(D)}
 \leq  C[||v||_{L^{2}(B)}+1], 
 \end{equation}
  for $D\subset\subset B$ and $v =
  \frac{u-\lambda}{\lambda}.$ The right side of
  (4.20) is bounded by (4.8) so that Sobolev
  embedding on $B$ implies that $\sup u/\lambda  =
  T/\lambda $ is bounded. Thus, under the assumptions
  above, (4.18) cannot hold.
 
  The same argument as above holds if only $\limsup
  T_{i}/\lambda_{i} >  1.$ Namely, this
  assumption together with (4.20) and Sobolev
  embedding then implies that, in some subsequence,
  $u/\lambda  = u_{i}/\lambda_{i} \geq  1 + \varepsilon ,$ 
  for some $\varepsilon > 0,$ on a ball
  $D\subset B$ whose volume is uniformly bounded
  below, if (4.18) held. This is impossible since
  $u_{i}/\lambda_{i} \rightarrow  1$ in $L^{2}$ by
  (4.8). Thus, in the following, we may assume
 	\begin{equation} \label{e4.21} 
 	\lim T_{i}/\lambda_{i}
 	= 1 .  
 	\end{equation}

  For a given $t  <   1,$ let $U_{i} = U_{i}(t) =
  \{x\in (M, g_{i}): u_{i}(x)/\lambda_{i} \geq  t\}.$
  Suppose that there were $\rho_{o} > 0 $ such
  that 
  \begin{equation} \label{e4.22} 
  \rho (x) >  \rho_{o}, \ \ \ \ \  \forall x\in U_{i}.  
  \end{equation} 
  Since $u_{i}/\lambda_{i} \rightarrow  1$ in $L^{2},$ it
 follows that $\vol(M\setminus U_{i}) \rightarrow 0.$  Then,
 arguing as above on (4.19), one also has
 $u_{i}/\lambda_{i} \rightarrow  1$ pointwise on
 compact subsets of each $B = B_{x_{i}}(\rho_{o}),\ 
 x_{i}\in U_{i}.$ 
 Observe that these balls must cover $M.$ 
 For if not, then there exist balls
 $B_{q_{i}}(\rho_{o})$ disjoint from $U_{i}.$ On the
 other hand, by the non-collapse assumption (3.43),
 such balls have a definite volume, and hence
 $B_{q_{i}}(\rho_{o})$ must have intersected $U_{i}.$
 It follows (from (1.27)), that 
 $$\rho (x) \geq \rho_{o}/2,  \ \  \forall x\in M, $$ 
 which contradicts the
 assumption $\rho (p_{i}) \rightarrow  0,$ for some
 $p_{i}\in (M, g_{i}).$ This means that (4.22) does
 not hold, so that for all $t < 1,$ there exist
 $x_{i}(t) \in  U_{i}(t)$ such that $\rho (x_{i}(t))
 \rightarrow 0$  as $i \rightarrow  \infty .$ Since
 $t$ is arbitrary, one may choose a sequence $t_{i}
 \rightarrow  1$ suitably slowly as $i \rightarrow
 \infty $ to give (4.16) and (4.17).
  \\ \end{pf}

  We note that Proposition 4.2 does not require
  $\lambda $ to be bounded away from 0, or any bound
  on the $L^{2}$ norm of $z^{T}.$
 
  Also, note that if $T_{i}/\lambda_{i} \rightarrow
  \infty ,$ then the proof above does not require the
  hypothesis (4.8), at least when $\sigma (M)  < 0 $
  Hence, in this case, the main hypothesis (3.60)
  of Theorem 3.10 is automatically satisfied if
  $\{g_{i}\}$ degenerates somewhere on $M.$
 
  Although the estimate (4.8) follows from (4.7) when
  $\sigma (M) < 0,$ it is not clear how to
  construct, on general manifolds $M,$ sequences of
  Yamabe metrics satisfying (4.6) or (4.7), even
  though the Palais-Smale condition (4.6) is easily
  realized. In \S 6.5, we construct examples of
  sequences satisfying (4.6); in fact these sequences
  satisfy much stronger conditions, c.f. also the
  discussion in \S 7.
 
 
 \begin{remark} \label{r 4.3.}
   We will not focus here on any applications of
   Theorem 3.10, but mention one kind of potential
   application. Suppose $\{g_{i}\}$ is a sequence of
   unit volume Yamabe metrics, for instance a
   maximizing sequence, satisfying the conditions of
   Theorem A(II). Suppose further that it can be
   proved on other (geometric or topological) grounds
   that there are no non-trivial static vacuum
   solutions arising from blow-ups $g_{i}'  = \rho
   (x_{i})^{-2}\cdot  g_{i}.$ (As an example of this,
   we mention that it can be shown that there are no
   blow-up limits $(N, g_{o})$ which are non-trivial
   static vacuum solutions with either smooth or
   compact boundary $\partial N,$ (c.f. \S 5),
   provided the sequence $(M, g_{i})$ satisfies a
   uniform Sobolev inequality, as for instance in
   (3.97), and $M$ is irreducible, i.e. every
   2-sphere in $M$ bounds a 3-ball).
 
  Under these conditions, Theorem A immediately
  implies that one has the following strong estimate
	 \begin{equation} \label{e4.23} 
	 \rho (x) \geq c\Bigl(\frac{|u(x)|}{\sup|u|}\Bigr), 
	 \end{equation}
  whenever $u(x) \neq  0.$ If further (4.8) holds,
  then Proposition 4.2 implies that (4.23) can be
  strengthened to a bound $\rho (x) \geq  c >  0$
  everywhere on $(M, g_{i}),$ i.e. there is no
  degeneration.  \end{remark} 
 
 \section{\bf Completeness of the Blow-up Limits.}
 \setcounter{equation}{0}

  Theorem 3.10 proves under certain conditions the
  existence of non-flat blow-up limits for the
  sequence $\{g_{i}\}$ which model the degeneration
  of $\{g_{i}\}$ in a neighborhood of the 0-level of $u.$ 
  In particular, these blow-ups $(D, g' )$ are
  non-trivial solutions to the static vacuum Einstein
  equations. However, such solutions have only been
  defined locally. We now address the issue of the
  completeness of the solutions, together with some
  aspects of the behavior of the potential function
  $\Roof{u}{\bar}$ at the boundary and the
  renormalizations of $\{u_{i}\}$ used in the
  construction of $\bar{u}.$
 
  We first discuss the issue of completeness. Recall
  from Theorem 3.10 that the metric $g' $ is a limit
  of rescalings of the metrics $g_{i},$ i.e. 
 $g'  = \lim g_{i}' $ based at $y_{i},$ where
 \begin{equation} \label{e5.1} 
  g_{i}'  = \rho (y_{i})^{-2}\cdot  g_{i}, 
 \end{equation} Further,
 $u_{i}$ is renormalized to $\Roof{u}{\bar}_{i}$
 satisfying 
 \begin{equation} \label{e5.2}
 \bar{u}_{i}(y_{i}) = 1, 
 \end{equation} c.f.
 (3.69). For a fixed $\mu > 0$ small, let
 	\begin{equation} \label{e5.3} N_{i}(\mu ) =
 	\{x_{i}\in (M, g_{i}' ): \rho' (x_{i}) \geq
 	\mu\cdot \rho' (y_{i}) = \mu\}, 
 	\end{equation}
  where $\rho' $ is the $L^{2}$ curvature radius
  w.r.t. $g_{i}' .$ Let $N_{i}(\mu , y_{i})$ be the
  component of $N_{i}(\mu )$ containing the base
  point $y_{i}.$ By Theorem 1.5 and the non-collapse
  assumption (3.43), as previously discussed in \S 3.4,
 the pointed sequence of manifolds $\{(N_{i}(\mu , y_{i}),
  g_{i}' , y_{i})\}$ has a subsequence converging,
  weakly in $L^{2,2}$ on compact subsets, to a
  non-flat connected limit domain $(N(\mu , y' ), g'
  , y' ).$ The convergence is in the strong $L^{2,2}$
  topology away from the locus
  $\{|\bar{u}_{i}| \leq \varepsilon\},\ \varepsilon $
  small, and away from the metric boundary of $N_{i}(\mu, y_i)$, i.e. at any given small distance away from these sets. This limit
  clearly extends the limit $(D, g' , y' )$ given by
  Theorem 3.10. The $L^{2,2}$ limit metric $g' $ on
  $N(\mu , y' )$ induces a natural embedding of any
  domain $K$ with smooth and compact closure in $N(\mu , y' )$
  into $M,$ c.f. \S 1.4.
 
  Now choose a sequence $\mu  = \mu_{j} \rightarrow  0$ 
  and consider the double sequence $N_{i}(\mu_{j},
  y_{i}).$ There is a diagonal subsequence $\{(i, j_{i})\}$ 
  of the double sequence $\{(i, j)\},$ with
  $j_{i} \rightarrow  0$ sufficiently slowly as $i
  \rightarrow  \infty ,$ so that
  $\{N_{i}(\mu_{j_{i}},y_{i}), g_{i}' , y_{i}\}$
  converges, uniformly on compact subsets in the weak
  $L^{2,2}$ topology, to a limit $(N, g' , y'
  ).$ Although the limit may apriori have more than
  one component, (the connected domains
  $N_{i}(\mu_{j_{i}}, y_{i})$ may pinch off into
  several limit components), we assume that $(N, g' ,
  y' )$ is the component of the limit containing the
  base point $y' .$ The limit $N$ is contained in the
  union of the inclusions $N(\mu_{k}, y' ) \subset
  N(\mu_{k+1}, y' ),$ for a sequence $\mu_{k}
  \rightarrow  0,$ and is an open Riemannian manifold
  with $L^{2,2}$ metric $g' .$ Any smooth open
  domain $K$, with smooth boundary, and properly contained in the connected limit $N$, is embedded
  as an open domain in $M$ via the metric $g' .$ It
  is clear from thw work in \S 3 that the limit
  functions $\bar{u}_i$ converge, uniformly
  on compact subsets, to a non-constant limit potential
  function $\bar{u}$ on $N.$
 
  We call $(N, g', y',\Roof{u}{\bar})$ the {\it maximal solution}  
  of the static vacuum equations
  associated to the buffered sequence $\{y_{i}\}.$
  Observe that this maximal solution determined by the (convergent) sequence $\{y_{i}\}$ is unique, (up to
  isometry), since any maximal limit must contain the
  initial domain $(D, g' )$ and static vacuum
  solutions have unique extensions. Hence each domain
  $N_{i}(\mu , y_{i})$ converges in the pointed
  Hausdorff topology, (c.f. [G, Ch.5A]), to its limit
  domain in $N.$
 
  Let $\partial N$ denote the metric boundary of $(N, g' ),$ 
  i.e. the set of (ideal) limits of $g'$-Cauchy sequences, 
   non-convergent in $N.$ The
  boundary $\partial N$ might be empty, for instance
  perhaps if $N$ is the complete (isometrically
  doubled, c.f. Remark 5.3(ii)) Schwarzschild metric
  (0.17), but in `most' cases one expects $\partial N$ 
  to be non-empty. In fact one expects $\partial N$ 
  to often be singular in the sense that the
  Riemannian metric $g' $ does not extend to
  $\partial N.$ Of course the completion
  $\Roof{N}{\bar} = N\cup\partial N$ is a complete
  metric space w.r.t. the length metric induced by
  $g' .$
 
  In case $\partial N \neq  \emptyset ,$ we need to
  describe $\partial N$ also as a set of (ideal)
  limits of points of $\{(M, g_{i}' )\}.$ For $p\in
  \Roof{N}{\bar} = N\cup\partial N,$ and $p_{i}\in
  M,$ define 
  \begin{equation} \label{e5.4} 
p_{i} \rightarrow  p, \end{equation} 
if, for all $\delta > 0,$ there exists 
  $q = q_{\delta}\in N,\  \mu  = \mu_{\delta} > 0$ and 
  a sequence $q_{i}\in (N_{i}(\mu , y_{i}), g_{i}' )$ 
  with $\dist_{g_{i}'}(p_{i}, q_{i}) \leq  \delta ,$ 
 such that $q_{i} \rightarrow  q,$ in
 the $L^{2,2}$ convergence of $N_{i}(\mu , y_{i})
 \rightarrow  N(\mu , y' )$ discussed above. It follows from the convergence and the definition of metric completion that $\dist_{g'}(p, q_{\delta}) \leq \delta$. In particular, since $\delta$ is arbitrarily small, $p$ in (5.4) is uniquely determined, in that there do not exist two distinct points $p$ and $p'$ to which the sequence $p_i$ converges.

  A sequence of subsets $Z_{i} \subset  (M, g_{i}' )$ is said
  to converge in the {\it  Hausdorff topology
  based at} $y_{i} $ to a set $Z \subset
  \Roof{N}{\bar} = N\cup\partial N$ if, for each $z\in
  Z,$ there exist $z_{i}\in Z_{i}$ such that $z_{i}
  \rightarrow  z$ in the sense of (5.4), and
  conversely any bounded sequence of points $z_{i}\in
  Z_{i},$ (i.e. $\dist_{g_{i}'}(z_{i},y_{i}) \leq  C$
  for some $C$), has a subsequence converging to a
  limit $z\in Z$ in the sense of (5.4). This
  definition extends the notion of Hausdorff
  convergence of the domains $N(\mu , y_{i})$ above
  to the boundary $\partial N.$
 
  Observe that $\partial N$ is formed from base
  points of higher order curvature concentration than
  $\{y_{i}\}.$ Thus, $p\in\partial N$ if and only if
  there is a sequence $p_{i}\in (M, g_{i}),$ with
  $p_{i} \rightarrow  p$ in the sense of (5.4) with
  $\rho' (p_{i}) \rightarrow  0,$ or equivalently,
  $\rho (p_{i})  <  < \rho (y_{i}).$
 
  The next result shows that, in a certain sense, $N$
  is defined and complete at least up to the
  Hausdorff limits of the $\varepsilon$-levels
  $L^{\varepsilon}$ of $\Roof{u}{\bar}_{i}.$

 \begin{theorem}[\bf Completeness] \label{t 5.1. (Completeness).}
   Let $(N, g' , y' , \Roof{u}{\bar})$ be the maximal
   solution of the static vacuum Einstein equations
   associated to the $(\rho ,c)$ buffered sequence
   $\{y_{i}\}$ in (5.1); (hence the non-collapse
   assumption (3.43) is assumed). Let
   $\Roof{U}{\bar}_{i}^{\varepsilon}$ be the component
   of $\{x_{i}\in (M, g_{i}' ): \Roof{u}{\bar}(x_{i})
   >  \varepsilon\}$ containing the base point
   $y_{i},$ for $\Roof{u}{\bar}_{i}$ normalized as in
   (5.2).
 
 
  Let $N^o \subset N$ be the maximal domain on which
  the potential function $\bar{u} > 0.$ 
  Then there is a sequence $\varepsilon_i \rightarrow 0$
  such that $N^o$ is contained in the $y_i$-based Hausdorff
  limit of the domains $\bar{U}_i^{\varepsilon_i}\subset (M,g_i).$
  Further the metric boundary $\partial N^o \equiv \bar{L}^o$
  of $N^o$ w.r.t. $g'$ is contained in the $y_i$-based
  Hausdorff limit of the $\varepsilon_i$-levels 
  $\bar{L}_{i}^{\varepsilon_{i}} \equiv \partial\bar{U}_{i}^{\varepsilon_i}$
  of $\bar{u}_{i}.$ The union $\bar{N}^o \equiv N^o \cup \partial N^o 
  \subset \bar{N}$ is complete w.r.t. the metric $\bar{g}$
  but may be singular at the boundary. 

  Let $t(x) = \dist_{g'}(x, \partial N^{o}),$ for
    $x\in N^{o}.$ Then, as in (3.16), there is an absolute constant $K  <   \infty $ 
    such that
 	\begin{equation} \label{e5.5} 
 		|z|(x) 
 		\leq
 		\frac{K}{t^{2}(x)},
 		\ \ \ \ \ 
 		\Roof{u}{\bar}^{-1}|\nabla \bar{u}|(x) \leq
 		\frac{K}{t(x)} 
 	\end{equation}
 \end{theorem} 
\begin{pf}
  In the notation above, consider any sequence of
  points $\{z_{i}\}\in \Roof{U}{\bar}_{i}^{\varepsilon}$
  with $\dist_{g_{i}'}(z_{i}, y_{i}) \leq  R,$ for
  some arbitrary (large) constant $R  <   \infty .$
  Suppose further that 
	 \begin{equation} \label{e5.6}
	\dist_{g_{i}'}(z_{i},\Roof{L}{\bar}_{i}^{\varepsilon}) \geq  \delta , 
	 \end{equation}
  for an arbitrary (small) constant $\delta > 0,$ 
  where $\Roof{L}{\bar}_{i}^{\varepsilon} = \partial
  \Roof{U}{\bar}_{i}^{\varepsilon}.$
 
  It follows from Theorem 3.3 that there is a
  constant $a = a(\varepsilon ) > 0$ such that
	 \begin{equation} \label{e5.7} 
           	 \rho' (z_{i}) \geq a\cdot \delta , 
	 \end{equation} 
 Thus, the curvature
 of $\{g_{i}'\}$ does not blow-up (in $L^{2})$ within
 bounded distance to $\{y_{i}\},$ provided one stays
 a fixed distance away from the level set
 $\Roof{L}{\bar}_{i}^{\varepsilon}$ in
 $\Roof{U}{\bar}_{i}^{\varepsilon}.$ Note also that the
 assumption (3.43) prevents collapse in these balls
 $B_{y_{i}}' (\rho' (z_{i})).$
 
  Let $N_{i}'  = N_{i}' (\varepsilon ,\delta )$ be the
  connected component of the set $\{z_{i}\in (M,
  g_{i}' ): z_{i}\in \Roof{U}{\bar}_{i}^{\varepsilon},
  \dist(z_{i},
  \Roof{L}{\bar}_{i}^{\varepsilon})\geq\delta\}$
  containing the base point $y_{i}.$ Using Theorem
  1.5 and Theorem 3.4 in the usual way, the pointed
  Riemannian manifolds $\{(N_{i}' , g_{i}' ,
  y_{i})\}$ (sub)-converge in the strong $L^{2,2}$
  topology based at $\{y_{i}\}$ to a limit smooth
  metric $g' ,$ based at $y' ,$ defined on an open
  domain $N' (\varepsilon ,\delta );$ the convergence is
  uniform on compact subsets. The functions
  $\Roof{u}{\bar}_{i}$ converge to a limit harmonic
  function $\bar{u}$ on $N' (\varepsilon ,\ \delta
  ), \Roof{u}{\bar} \geq  \varepsilon ,$ so that the
  pair $(g' , \Roof{u}{\bar})$ is a non-flat solution
  of the static vacuum Einstein equations on 
  $N' (\varepsilon , \delta ).$ Of course 
  $y'\in N' (\varepsilon ,\delta )$ 
  for any $\varepsilon ,\delta $ small and $\bar{u}(y' ) = 1.$ 
  Given sequences $\varepsilon_{j} \rightarrow  0$ and 
  $\delta_{j}\rightarrow  0,$   there exist
  suitable diagonal subsequences $j=j_i$ of the double
  sequence $N_i' (\varepsilon_j,\delta_j),$
  (with $\varepsilon_{j_i}$ and $\delta_{j_i} \rightarrow 0$
  sufficiently slowly as $i\rightarrow \infty),$
  which converge, uniformly on compact subsets, to a limit.
  As above in the construction of $N$, we consider only the 
  component $(N^o,g')$ of the limit containing $y'.$
  The connected limit $(N^o, g', y')$ is an is an open Riemannian
  manifold, and is contained in the union of the domains 
  $N' (\varepsilon_k,\delta_k) \subset N' (\varepsilon_{k+1},\delta_{k+1})$
  for some sequence $k \rightarrow 0$.
  The limit potential function $\bar{u}$
  on $N^o$ satisfies $\bar{u} \geq 0$ on $N^o,$
  and $\bar{u}(y')=1$. Thus, $\bar{u}>0$
  by the minimum principle for harmonic functions.
  As discussed in \S1.3, $(N^o, g')$ is hence 
  a smooth solution of the static vacuum equations.

  The construction of $N$ and $N^{o}$ implies that
 \begin{equation} \label{e5.8} 
 (N^{o}, g' ) \subset (N, g' ),  
 \end{equation} 
 (for any choice of
 subsequence above) and $N^{o}$ is unique. It is
 possible, although not necessary, that $N^{o} = N.$
 It follows immediately from the definition preceding
 Theorem 5.1, (c.f. also the argument below), that
 $(N^{o}, g' )$ is contained in the $y_{i}$-based
 Hausdorff limit of the domains
 $(\Roof{U}{\bar}^{\varepsilon_{i}}, g_{i}' ),$ where
 $\varepsilon_{i} = \varepsilon_{j_{i}}$ is defined above.
 
  Let $\partial N^{o}$ be the metric boundary of
  $N^{o}$ w.r.t. $g' $ and let $\Roof{N}{\bar}^{o} =
  N^{o}\cup\partial N^{o}$ be the metric completion.
  Of course $\partial N^{o} \subset  \Roof{N}{\bar}.$
 
  We claim that $\partial N^{o}$ is contained in
  the $y_{i}$ based Hausdorff limit of the levels
  $\Roof{L}{\bar}^{\varepsilon_{i}}, \ \varepsilon_{i} =
  \varepsilon_{j_{i}}.$ To see this, suppose
  $p\in\partial N^{o},$ so that $p = \lim p_{j},$
  for $\{p_{j}\}$ a Cauchy sequence in $N^{o}.$ Thus,
  $p_{j}\in N' (\varepsilon_{j}, \delta_{j}),$ for some
  $\varepsilon_{j}, \delta_{j} >  0.$ But $N'
  (\varepsilon_{j}, \delta_{j}) =
  \lim_{i\to\infty} N_{i}' (\varepsilon_{j},
  \delta_{j})$ and hence for any $j,$ there exist
  sequences $q_{i,j} \in  N_{i}'
  (\varepsilon_{j},\delta_{j})$ with $q_{i,j}
  \rightarrow  p_{j}$ as $i \rightarrow  \infty ,$ as
  required.
 
  Suppose first $\Roof{N}{\bar}^{o} \subset  N,$ (so
  that in particular $\partial N^{o}\cap\partial N =
  \emptyset  ).$ Then the metric $g' $ extends as an
  $L^{2,2}$ metric past $\Roof{N}{\bar}^{o}$ and the
  convergence $g_{i}'  \rightarrow  g' $ is in the
  weak $L^{2,2}$ topology in an open region $R$
  satisfying $\Roof{N}{\bar}^{o} \subset  R \subset
  N(\mu , y' ) \subset  N$, for some $\mu  >  0.$
  (More precisely, $\mu $ may depend also on the
  distance to the base point $y' $ if $\partial
  N^{o}$ is non-compact). Since the functions
  $\Roof{u}{\bar}_{i}$ also converge in $L^{2,2}$ in
  $R$, (in fact strongly in $L^{2,2}$ by (3.35)), it
  follows that $\partial N^{o}$ is identified with
  the 0-level set of $\Roof{u}{\bar},$ i.e. with the
  event horizon $\Sigma $ of $(N, g' ).$
 
  In this case, since $g' $ is $L^{2,2}$ in a
  neighborhood of the event horizon $\Sigma ,$ as
  noted in \S 1.3, (c.f. the discussion following Theorem 1.1), elliptic regularity arguments show
  that in fact $g' $ is $C^{\infty}$ smooth across
  $\Sigma $ and $\Sigma $ is a collection of smooth totally geodesic
  surfaces forming the (smooth) topological boundary of
  $N^{o}.$ Regardless of this fact, the weak maximum principle for $L^{2,2}$ harmonic functions, c.f. [GT, Thm. 8.1], implies that $\bar{u} < 0$ in $N \setminus \bar{N}^{o}$.
 
  More generally, since the discussion above is
  local, it holds for the parts of $\partial N^{o}$
  which are contained in $N,$ i.e. those points of
  $\partial N^{o}$ which admit an open neighborhood
  contained in $N.$  Now both $\partial N$ 
  and $\partial N^{o}$ are closed, so the remaining
  points of $\partial N^{o}$ are those contained in
  $\partial N.$ Thus, the domain $N^{o}$ is the
  maximal domain in $N$ on which $\bar{u} > 0.$ 
  The set $Z = \partial N^{o}\cap\partial N$
  need not be empty, and one would frequently expect
  $\partial N^{o} = \partial N.$
 
  To conclude, the estimates (5.5) follow immediately
  from Theorem 3.2(II), (applied to any smooth
  subdomain in $N^{o}),$ and the fact that
  $\bar{u} >  0$ on $N^{o},$ c.f. also
  (A.22) and the discussion following it at the end
  of the Appendix. \\ 
\end{pf}
 
  Note that the limit function satisfies
  $\bar{u}(y' ) = 1,$ since
  $\bar{u}(y_{i}) = 1$ and, since $\rho'
  (y_{i}) = 1,$ the convergence to the limit is
  controlled near $y_{i}.$ In particular, $\partial
  N$ does not intersect $B_{y'}(1)$ and if $\partial
  N^{o}$ intersects $B_{y'}(1),$ it does so smoothly,
  as noted above.
 
  However, it is apriori possible that all levels
  $\Roof{L}{\bar}_{i}^{\varepsilon}$ approach for
  instance the level $\Roof{L}{\bar}_{i}^{1}$ as $i
  \rightarrow  \infty $ away from $B_{y_{i}}' (1).$
  Thus, although the domains $N_{i}' (\varepsilon ,
  \delta )$ contain $N_{i}' (1,\delta )$ for
  $\varepsilon < 1,$ they may apriori give rise to
  the same limit domains, outside some region
  containing $B_{y'}(1).$ In other words, the
  functions $\roof{u}{\bar}_{i}$ might descend from
  the value 1 (for example) to a value arbitrarily
  near 0 in arbitrarily short $g_{i}' $ distances as
  $i \rightarrow  \infty .$ Of course, this can only
  occur at points where $\rho'  \rightarrow  0,$ i.e.
  at $\partial N.$ For the same reasons, we do not
  assert that $\roof{u}{\bar}(x_{j})$ approaches 0
  whenever $x_{j}\in N^{o}$ approaches $\partial
  N^{o} \equiv  \roof{L}{\bar}^{o}.$
 
  Thus, $\roof{L}{\bar}^{o}$ cannot necessarily be
  identified as the 0-level set of the potential
  $\roof{u}{\bar},$ i.e. with the event horizon.
  Further remarks on the structure of the maximal
  limits $(N^{o}, g' )$ and $(N, g' )$ follow below.

 \begin{remark} \label{r 5.2.}
   We conjecture that for the blow-up limits $(N, g' , y' )$ 
   given by Theorem 5.1, one has
 	\begin{equation} \label{e5.9} \partial N^{o} =
 	\{\roof{u}{\bar} = 0\} = \Sigma , 
 	\end{equation}
  so that $\partial N^{o}$ is identified with the
  event horizon of the static vacuum solution. This
  amounts to proving that $\roof{u}{\bar}(x_{k})
  \rightarrow  0$ whenever $x_{k}$ converges to a
  point in $\partial N^{o}.$
  \end{remark}
 
 
 
 
 \begin{remark} \label{r 5.3.(i).}
   {\bf (i).} In general relativity, it is usually assumed that
   a static vacuum solution is complete up to the
   event horizon $\Sigma  = \{\roof{u}{\bar} = 0\},$
   and in particular that $\roof{u}{\bar}(x)
   \rightarrow  0$ as $x \rightarrow  \Sigma .$ Such
   static vacuum solutions are the most natural
   physically. As noted in \S 1.3, the static vacuum
   equations are formally degenerate at the locus
   $\Sigma ,$ but are formally non-degenerate away
   from $\Sigma .$ By Theorem 3.2(I), there are no
   non-trivial complete solutions with empty $\Sigma .$
 
  There are however examples of static vacuum
  solutions which are not complete (or defined) up to
  $\Sigma$, c.f. [An4, \S 2]. For instance, let $\nu $ be a bounded
  harmonic function on ${\Bbb R}^{3}\setminus B(1),$ with
  axially symmetric, but discontinuous or non-smooth
  boundary values on $\partial B(1).$ Then $\nu $ is
  axially symmetric and hence defines a Weyl solution
  as in (1.17), which does not extend everywhere past
  $B(1)$. For suitable boundary values, the solution
  will not extend anywhere into $B(1)$, and hence one
  has a solution complete away from $B(1)$, but not
  defined up to the event horizon, (where $\nu  = -\infty ).$
 
  Note that such examples may have $\max \nu ,$ and so
  $\max u,$ occurring on $\partial B(1),$ while $\min \nu $ 
  occurs at infinity. In this case, such solutions
  formally have `negative mass', c.f. \S 1.3.
 
  Of course, from Remark 5.2, we conjecture that such
  solutions cannot arise as blow-up limits of unit
  volume Yamabe metrics.
 
  (In the opposite direction, there are also
  solutions which are defined and smooth up to
  $\Sigma ,$ but which are not complete away from
  $\Sigma ,$ for example the A2 or B2 solutions, c.f.
  [EK, \S 2-3.6]).
 
 {\bf (ii).}
  While there are many static vacuum solutions which
  are singular at $\Sigma ,$ but smooth and complete
  away from $\Sigma ,$ the Schwarzschild metric
  (0.17) is of course smooth up to and at $\Sigma .$
  Here $\Sigma $ is given by a round, totally
  geodesic 2-sphere $S^{2}.$ Thus, the Schwarzschild
  metric $g_{s}$ may be isometrically doubled across
  $\Sigma $ giving a complete, smooth metric on
  $S^{2}\times {\Bbb R} ,$ asymptotically flat at both
  ends. This is clearly the maximal (abstract) smooth
  extension of $g_{s},$ call it $\roof{N}{\tilde}.$
  However, it is not necessarily the case that
  $\tilde{N}$ agrees with the maximal
  extension $N$ defined preceding Theorem 5.1.
  Apriori, it is possible that the curvature of
  $g_{i}' $ blows up near the event horizon $\Sigma
  ,$ (where one loses strong convergence), even
  though the limit in this particular case is smooth
  across $\Sigma .$
 
  Basically because of this, we will usually consider
  only the behavior of limits up to the event
  horizon, (in case this is defined), i.e. within the
  domain $N^{o} \subset  N;$ c.f. however \S 6.5.
 
 {\bf (iii).}
  Although stated only for buffered sequences
  $\{y_{i}\}$ and the associated limit, it is clear
  from the discussion preceding Theorem 5.1 and its
  proof that Theorem 5.1 remains valid for arbitrary
  maximal limits $(N, g' , x),$  
 $g'      = \lim g_{i}' ,$   
 $g_{i}'  = \rho (x_{i})^{-2}\cdot  g_{i},$ 
 for which
  $\rho (x_{i}) \rightarrow  0$ and say
  $(u_{i}/T_{i})(x_{i}) \geq  u_{o},$ for some
  arbitrary constant $u_{o} > 0.$ Of course in
  this case, the limits may be flat solutions to the
  static vacuum equations.  \end{remark} 
 
  The following result, valid for general static
  vacuum solutions, shows that the estimate (5.5) can
  be improved in certain regions, c.f. also (A.27).
 \begin{lemma} \label{l 5.4.}
    Suppose $(N, g, u)$ is an (arbitrary) static vacuum
    solution, $u > 0$ on $N,$ and $u$ is bounded
    above. Let $\{x_{j}\}$ be a maximizing sequence
    for $u$ in N, 
    i.e. $u(x_{j}) \rightarrow \  \sup \ u
     <   \infty ,$ 
     and 
     $t(x) = \dist(x, \partial N),$ 
     where $\partial N$ is the metric boundary of $N.$ 
   Then 
     \begin{equation} \label{e5.10} 
     |z|(x_{j}) 
     \leq
       \frac{\mu_{j}}{t^{2}(x_{j})} , 
       \ \ \ \ \
      |\nabla log u|(x_{j}) 
     \leq
      \frac{\mu_{j}}{t(x_{j})}, 
      \end{equation} 
      where
 $\mu_{j} = \mu_{j}(x_{j}) \rightarrow  0$ as $j
 \rightarrow  \infty .$
 \end{lemma} 
\begin{pf}
  Consider the natural rescalings of $ (N, g)$  based at
  $\{x_{j}\},$ i.e. the sequence of metrics $  (N, g_{j}, x_{j}),$ 
  where $g_{j} = t(x_{j})^{-2}\cdot
  g_{o}.$ Let $t_{j} = t/t(x_{j})$ be the distance to
  $\partial N$ w.r.t. the metric $g_{j},$ so that
  $t_{j}(x_{j}) = 1.$ Now of course the curvature of
  $\{g_{j}\}$ will blow up at any fixed base point
  $x_{o}\in N,$ (with $z(x_{o}) \neq  0),$ if
  $t(x_{j}) \rightarrow  \infty ,$ and hence there is
  no limit in this region. Similarly, for example the
  curvature of $g_{j}$ may blow up when $t(x_{j})
  \rightarrow  0.$
 
  However, given any $\delta > 0,$ the
  curvature, (w.r.t. $g_{j}),$ of the domain
  $N_{\delta}(j) = \{p_{j}\in (N, g_{j}):
  t_{j}(p_{j}) \geq  \delta\},$ containing $x_{j}$
  remains uniformly bounded by the (scale-invariant)
  estimate (3.16). If the sequence $\{(N_{\delta}(j),
  g_{j}, x_{j})\}$ is non-collapsing at $x_{j},$ it
  follows from the (local) Cheeger-Gromov theory, (or
  Theorem 1.5), that a subsequence of
  $\{B_{x_{j}}(\frac{1}{2}), g_{j}, x_{j})\}$
  converges to a limit $(B_{\infty}, g_{\infty},
  x_{\infty}).$ By Lemma A.2, the convergence is
  smooth and uniform on compact subsets, and the
  limit is again a smooth solution of the static
  vacuum equations. If the sequence instead collapses
  at $x_{j},$ then the collapse may be unwrapped by
  passing to the universal cover
  $(\widetilde{B}_{x_j}(\tfrac{1}{2}), g_{j});$ this
  sequence no longer collapses, and hence
  sub-converges to a limit as above. (We refer to the
  Appendix - Lemma A.2 and Corollary A.3 - for more
  details of this type of standard argument).
 
  On the other hand, it is clear that $u(x_{\infty})
  = \sup u,$  and so by the maximum principle, the
  limit harmonic function $u$ is constant on
  $B_{\infty}$ and $g_{\infty}$ is hence flat. Since
  the convergence of $N_{\delta}(j)$ to the limit is
  smooth, the curvature of $g_{j}$ is almost 0, and
  $u$ is almost constant, in $N_{\delta}(j),$ away
  from the boundary. This implies (5.10) by
  scale-invariance.  \\
\end{pf}

  Next, we prove that the domain $N^{o} \subset  N$
  is large in a natural sense; it particular it is
  unbounded.
 \begin{proposition} \label{p 5.5.}
   Let $(N, g' , y' )$ be a maximal non-flat limit
   solution as in Theorem 5.1, with domain 
   $N^{o} \subset  N$ on which $\bar{u} > 0.$
   There is a smooth curve $\gamma : {\Bbb R}^{+}
   \rightarrow  N^{o},$ parametrized by arclength,
   with $\gamma (0) = y' ,$ and positive constants
   $\delta_{1}, \delta_{2}$ such that 
   \begin{equation}
        \label{e5.11} \dist(\gamma (s), y' ) 
        \geq \delta_{1}\cdot  \dist(\gamma (s), \partial N^{o}) ,
   \end{equation} 
 and such that the cone 
 $V = V(d_{2}) = \{x: \dist\,(x, \gamma (s)) \leq  \delta_{2}\cdot  s\}$
 over $\gamma $ satisfies 
 	\begin{equation} \label{e5.12} 
 	V \subset  N^{o}, 
 	\ \ \ 
 	i.e. \ \ \ 
 	V\cap\partial
 	N^{o} = \emptyset  .  
 	\end{equation} 
 In particular, the function $t(x) = \dist_{g'}(x,\partial N^{o})$
 has linear, and hence unbounded, growth in $V.$
 \end{proposition} 
\begin{pf}
  This follows from the descent construction of the
  limit $N$ in Theorem 3.10. Thus, recall that the
  buffered sequence $y_{i} = x_{i}^{k_{o}}$ in $(M, g_{i}' )$, $k_{o} = k_{o}(i) \rightarrow \infty$, as $i \rightarrow \infty$, 
  has predecessors 
  $x_{i}^{j'}$, $j' = j'(i,j) = k_{o} - j$, for any fixed $j' > 0$. Thus, we have relabeled so that $y_{i} = x_{i}^{j'}$, for $j' = 0$.  By the construction, these have 
  $\bar{u}$-values at least
  $(\frac{1}{2})^{j'},$ c.f. (3.73). 
  In particular, these points are in $N_{i}' (\varepsilon ,1),$ 
  for $\varepsilon $ sufficiently small,
  (depending on $j'$), c.f. the proof of Theorem 5.1;
  compare also with Remark 3.15(iv). Further, in the
  scale $g_{i}' $ associated to $y_{i}$ these points
  have very large curvature radius; in fact $\rho'
  (x_{i}^{j'}) \geq  (d_{1})^{-j'},$ for all  
  $j > 0,$ while both
  $\dist_{g_{i}'}(x_{i}^{j'}, y_{i})$ and
  $\dist_{g_{i}'}(x_{i}^{j'},
  \roof{L}{\bar}^{\varepsilon})$ are also on the order
  of $(d_{1})^{j'},$ (for $\varepsilon $ sufficiently
  small, depending on $j'$, c.f. (3.66)-(3.67)). Note that from (1.27), for
  any $q\in B_{x_{i}^{j'}}(\rho' (x_{i}^{j'})),$ one
  has $\rho' (q) \geq \dist_{g_{i}'}
              (q , \partial B_{x_{i}^{j'}} (\rho'(x_{i}^{j'}))).$

  Thus, let $\gamma (s)$ be a path joining the limit
  points $x^{j'} = \lim x_{i}^{j'}$ in $N^{o},$
  approximating a minimizing geodesic joining each
  $x^{j'+1}$ to $x^{j'}$ within $B^{j'} =
  B_{x_{j'}}(\rho' (x^{j'})).$ One may then define 
  \begin{equation*}
    \displaystyle V = 
     \bigcup_{j'\geq 0}  B_{x^{j'}}(\rho' (x^{j'})) \subset N^{o},
  \end{equation*} 
   and the result follows.  \\
\end{pf}
 
  The curvature estimate (5.5) on $N^{o}$ can be
  improved in the region $V$ in (5.11), since the
  curvature radius has linear growth in $V,$ and the
  predecessors $x_{i}^{j'}$ are not $(\rho ,c)$
  buffered. Thus, in fact
 \begin{equation} \label{e5.13} 
 |z|(x) 
 \leq
 \frac{\kappa}{t^{2}(x)},
 \ \ \ \ \ 
 \bar{u}^{-1}|\nabla \bar{u}|(x) 
 \leq
 \frac{\kappa}{t(x)} 
 \end{equation}
  for $x\in V,$ where $\kappa  = \kappa (c)$ may be
  made small by choosing the buffer constant $c$
  sufficiently small. This follows since (5.13) holds
  in an $L^{2}$ sense, from the fact that $x$ is not
  $(\rho ,c)$ buffered, together with elliptic
  regularity for the static vacuum equations, which
  gives an $L^{\infty}$ bound in terms of an $L^{2}$
  bound.
 
  It seems possible that Proposition 5.5 may not hold
  for arbitrary `complete' static vacuum solutions.
  (For example, consider the positive measure $\mu $
  formed by placing positive multiples of the Dirac
  measure at all integer lattice points in ${\Bbb
  R}^{3},$ weighted so that the total mass is finite.
  Let $v$ be the Newtonian potential of $\mu .$ Then
  there may possibly be static vacuum solutions whose
  potential $u$ resembles the geometry of $v$).

 \begin{remark} \label{r 5.6.}
   We conclude with some remarks on the issue of the
   renormalization of $\{u_{i}\}$ and the relation of
   $\roof{u}{\bar}$ with the initial sequence
   $\{u_{i}\}.$
 
  Recall that the static vacuum solution constructed
  in Theorem 3.10 may well live in the region of
  $(M, g_{i}' )$ where $u_{i}$ is converging uniformly to
  0; the potential function $\roof{u}{\bar}$ of the
  static limit is obtained by renormalizing,
  (possibly infinitely many times), the original
  sequence $\{u_{i}\}.$ This will occur for instance
  if $u_{i}$ behaved as the function $v_{i} =
  t^{\delta_{i}},$ with $\delta_{i} \rightarrow  0$
  as $i \rightarrow  \infty ,$ and $t(x) = \dist(x, L^{o}).$ 
  Namely, if $x_{\varepsilon}$ satisfies
  $v_{i}(x_{\varepsilon}) = \varepsilon $ for any fixed
  $\varepsilon  >  0,$ and the distance between
  $L^{\varepsilon}$ and $L^{o}$ is scaled to size 1,
  then $v_{i}$ approaches the constant function
  $\varepsilon $ uniformly near $x_{\varepsilon}$ as
  $\delta_{i} \rightarrow  0.$
 
  In fact, we point out that there are static vacuum
  solutions $(N, g, u)$ and points $x_{j}\in N$ with
  $t(x_{j})~=~1,$ $(t(x)~=$ distance to event horizon), such that 
	  \begin{equation} \label{e5.14}
		   u(z_{j}) \geq  t^{\delta_{j}}(z_{j}), 
	  \end{equation}
 for a given sequence $\delta_{j} \rightarrow  0,$
 for all points $z_{j}\in B_{x_{j}}(2).$ Namely,
 consider for example the Weyl solution (1.17) with
 potential function $\nu $ determined by the Riesz
 measure
 	\begin{equation} \label{e5.15} d\mu_{\zeta} =
 	\frac{1}{1+|\zeta|}\, dL, 
 	\end{equation}
  where $dL$ is Lebesgue measure on the axis $A,$ and
  $\zeta $ is a parameter for $A,$ c.f. (1.18). This
  solution is complete up to the event horizon, so
  $\bar{L}^{o} = \Sigma  = \{u = 0\};$ the
  event horizon corresponds to the full axis $A.$ If
  $x_{j}$ diverges to infinity in $(N, g, u)$ at any
  fixed distance to $\Sigma ,$ then the potential
  function $u = e^{\nu}$ satisfies $u(x_{j})
  \rightarrow  1\  = \ \sup \, u$ and has the property
  (5.14).
 
  Thus, in addition to the discussion in Remark 5.9, this shows that the infinite descent down the $u$-levels in
  Theorem 3.10 is necessary.
 
  It is worth noting that in this example, $\Sigma $
  is highly singular, and the Riesz measure $d\mu $
  converges to 0 at infinity on the axis $A \sim
  \Sigma .$ In fact, any Weyl solution whose Riesz
  measure has regions of arbitrarily long length on
  $A,$ with arbitrarily small, but non-zero, density
  w.r.t. Lebesgue measure on $A$ will have behavior
  similar to (5.14).
 
  On the other hand, although one must carry out the
  inductive construction in Theorem 3.10 arbitrarily
  many times, i.e. $k_{o}(i) \rightarrow  \infty ,$
  as $i \rightarrow  \infty ,$ this does not {\it
  necessarily}  imply an infinite descent down the
  levels $L^{k}$ of $u.$ Recall that $u(x_{i}^{k})
  \geq  2^{-k},$ and it may well happen that that in
  fact $u(x_{i}^{k}) \geq  2^{-1},$ for all $k \leq k_{o},$ 
  (if $c$ is chosen sufficiently small). This
  is seen specifically for example in the
  Schwarzschild metric, as in Example 3.9. Thus,
  choose initial base points $\{x_{i}^{1}\}$ going to
  infinity in the Schwarzschild metric, so that
  $u(x_{i}^{1}) \rightarrow  1.$ If one carries out
  the inductive process of Theorem 3.10, then it is
  necessary to take $k_{o}(i) \rightarrow  \infty $
  in order to obtain a $(\rho ,c)$ buffered sequence
  $\{y_{i}\}$ from the sequence $\{x_{i}^{1}\}.$
  However, one easily sees that $u(x_{i}^{k}) \geq
  2^{-1}$ for all $k \leq  k_{o},$ so that one has
  $u(y_{i}) \geq  2^{-1}$ also, (for a suitable
  choice of $c$). In this example, the inductive
  process in the proof of Theorem 3.10 just
  recaptures the original Schwarzschild metric,
  (up to a bounded scale change).
 
  Similar behavior holds for Weyl solutions whose
  Riesz measure on $A$ has density w.r.t. the Lebesgue
  measure on $A$ either 0 or uniformly bounded below.
 \end{remark} 
 
 \section{\bf Construction of Yamabe Sequences with Singular Limits.}

  In this section, we discuss in detail some
  constructions of sequences of Yamabe metrics which
  illustrate the sharpness of the results in \S 3. In
  particular, these constructions exhibit the
  possibility that blow ups of Yamabe metrics may
  give rise only to trivial (i.e. flat), or
  super-trivial solutions of the static vacuum
  equations, see (3.15). (Discussions with R.
  Hamilton and R. Schoen were helpful to me in
  clarifying some aspects of the construction in
  Example 1). In \S 6.3, we analyse the behavior of
  the splittings (2.6) and (2.10) of $z$ and $g$
  respectively on the (singular) limits of these
  sequences.
 
  On the other hand, in \S 6.4 and \S 6.5, we
  construct examples of sequences of Yamabe metrics
  which do satisfy all the hypotheses of Theorem
  A(II)/Theorem 3.10, (so that this result is
  non-vacuous). The existence and basic properties of
  these four classes of examples are explained from a
  somewhat more general perspective in \S 7.
 
\subsection{\bf Example 1} 
 \setcounter{equation}{0}
   Let $(M, g_{o})$ be a hyperbolic 3-manifold.
   Clearly, $g_{o}$ is a critical point of
   $v^{-1/3}\cdot {\cal S} $ on ${\Bbb M} ,$ and thus
   of $v^{2/3}\cdot  s$ on ${\cal C} .$ In fact,
   $g_{o}$ is a local maximum of $v^{2/3}\cdot  s$ on
   ${\cal C} ,$ c.f. [Bes, 4.60].
 
  Let $N$ be any closed, oriented 3-manifold with
  $\sigma (N) > 0,$ for instance $N = S^{3},\ 
  S^{3}/\Gamma ,\ S^{2}\times S^{1},$ or a connected sum of
  such manifolds. In the following, we will construct
  Yamabe metrics on the manifold $N\#M,$ which are
  geometrically close to the original hyperbolic
  manifold $(M, g_{o}).$
 
  Let $\gamma $ be any smooth metric (not necessarily
  Yamabe) of positive scalar curvature on $N.$ Such a
  metric admits a positive Green's function
  $G_{y}(x),$ with pole at $y$, for the conformal
  Laplace operator $-
  8\Delta_{\gamma}+s_{\gamma},$ c.f.
  [LP, Thm2.8]. Consider the conformally related
  metric 
 \begin{equation} \label{e6.1} 
           g = G_{y}^{4}(x)\cdot \gamma , 
 \end{equation} 
 for any fixed $y\in N.$ This metric is a complete,
 scalar-flat metric on $N\setminus\{y\},$ which is
 asymptotically flat in the sense that, outside a
 large compact set, (thus in a small neighborhood of
 $ \{y\}),$
 \begin{equation} \label{e6.2} 
  g_{ij} = \left(1+\frac{2m}{r}\right)\delta_{ij} + O(r^{-2}),
 \end{equation}
  c.f. (1.15). Here $m$ is the mass of the metric
  satisfying $m \geq  0,$ with equality if and only
  if $g$ is flat, by the positive mass theorem [SY].
  From [Sc1], the metric $g$ is flat only when $(N,
  \gamma )$ is the canonical constant curvature
  metric on $S^{3}.$ Note that the curvature tensor
  $R_{g}$ of $g$ satisfies the decay condition
 \begin{equation} \label{e6.3} 
 |R_{g}|(x) 
 \leq
 c\cdot  r(x)^{-3}, 
 \end{equation}
  as $r(x) \rightarrow  \infty .$
 
  Choose any fixed value of $\varepsilon ,$ with 
  $0  <   4\varepsilon   <   i_{o} \equiv \inj_{g_{o}}(M),$ 
  where  $\inj$ denotes the injectivity
  radius. Given any $x\in M,$ we may then glue in the
  metric $g$ above to $B_{x}(2\varepsilon )
  \subset (M, g_{o})$ as follows. Given a fixed center point
  $p\in N\setminus \{y\},$ let $B(R)$ denote the geodesic $R$-ball
  in $(N\setminus \{y\}, g)$ centered at $p$. Scale $g|_{B(R)}$ to
  size $\varepsilon ,$ i.e. define 
	 \begin{equation}
	 \label{e6.4} g_{\varepsilon} =
	 \left(\tfrac{\varepsilon}{R}\right)^{2}\cdot  g , 
	 \end{equation} so
 that $g_{\varepsilon}|_{B_{x}(\varepsilon )}$ is
 homothetic to $g|_{B_{p}(R)}.$ From (6.2) and (6.3),
 it follows that in the geodesic annulus $A(\varepsilon
 ,2\varepsilon ),$ the sectional curvature of
 $g_{\varepsilon}$ is on the order of $|K_{\varepsilon}| =
 O(\varepsilon^{-2}R^{-1}).$ Hence the curvature is
 bounded in $A(\varepsilon ,2\varepsilon )$ if
 $\varepsilon^{2}R$ is bounded away from 0, while the
 metric $g_{\varepsilon}$ is almost flat in this band if
 $\varepsilon^{2}R >>  1.$ For simplicity, we
 assume from now on that $\varepsilon^{2}R >>  1.$
 Note that the unit ball $B_{p}(1)$ in $(N, g)$ is of
 radius $\varepsilon /R$ in $g_{\varepsilon};$ for
 $\varepsilon^{2}R >  1,$ $\varepsilon /R >> \varepsilon^{3}.$
 
  In order to bend the metric $g_{\varepsilon}$ on
  $B(2\varepsilon )$ so it has almost constant scalar
  curvature $-6=s_{g_{o}},$ define 
  \begin{equation} \label{e6.5} 
   \tilde{g}_{\varepsilon} = \psi^{4}\cdot  g_{\varepsilon}, 
   \end{equation} where
 $\psi  = 2^{1/2}(1-\tau^{2})^{-1/2}$ on $B(2\varepsilon)$ 
 and where $\tau $ will be determined below. Note
 that if $g_{\varepsilon}$ were flat on $B(2\varepsilon )$
 and $\tau  = t = \dist_{g_{\varepsilon}}(x, \cdot  ),$
 the metric $\roof{g}{\tilde}_{\varepsilon}$ is the
 hyperbolic metric of scalar curvature $-6$. Since
 $g_{\varepsilon}$ is scalar-flat a simple computation
 using (1.12) shows that \begin{equation}
 \label{e6.6} \Roof{s}{\tilde}_{\varepsilon} = -
 (1-\tau^{2})\Delta\tau^{2} -
 6\tau^{2}|d\tau|^{2}.  \end{equation} We choose
 $\tau $ so that $\Delta\tau^{2} = 6,$ and
 $\tau $ is close to the distance function $t(z) =
 \dist_{g_{\varepsilon}}(z, x).$ 
 Hence (6.6) becomes
 \begin{equation} \label{e6.7}
       \tilde{s}_{\varepsilon} = - 6[1+\tau^{2}(|d\tau|^{2}- 1)].  
 \end{equation}
 
  To construct such functions, return to the
  asymptotically flat manifold $(N\setminus \{y\}, g).$ We claim
  there is a function $\alpha $ on $N,$ with
  $\Delta\alpha  = 6 $ and $\alpha $
  asymptotic to $r^{2}.$ To define $\alpha $ consider
  the function $r^{2}+mr+2mlnr;$ this is defined and
  smooth outside a compact set $K$ of $(N\setminus \{y\}, g),$
  and off $K,$ one computes that 
  $\Delta (r^{2}+mr+2mlnr) = 6+4mr^{-3}+0(r^{-4}).$  
  Let $h$ be a smooth extension of $r^{2}+mr+2mlnr$ to all of
  $N\setminus \{y\};$ then $\Delta h = f,$ where the
  function $f$ satisfies $|6-f| = 0(r^{-3}).$ We
  define $\alpha $ as $\alpha  = h+\phi ,$ where
  $\phi (q) = -\int G(q,z)(6-f)dz$ and $G$ is the
  (positive) Green's function for $\Delta $
  on $N\setminus \{y\}.$ Since the product $\dist(q,z)\cdot
  G(q,z)$ is bounded below and above, $\phi $ is well
  defined and decays at infinity as $0(r^{-1}).$ It
  follows that the function $\alpha $ satisfies
  $\Delta\alpha  = \Delta (h+\phi )
  = 6$ on $(N\setminus \{y\}, g).$ By adding a suitable
  constant, we may assume $\alpha >0.$
 
  Rescaling to the metric $g_{\varepsilon},$ we then
  define $\tau  = (\varepsilon /R)\alpha^{1/2}$ on
  $B_{x}(2\varepsilon ).$ Since the metrics $g$ and
  $g_{\varepsilon}$ are homothetic, scaling properties
  imply $\Delta_{g_{\varepsilon}}\tau^{2} = 6$
  and $|d\tau|_{g_{\varepsilon}} = |d\alpha^{1/2}|_{g}
  \approx  1$ in $A(\varepsilon ,2\varepsilon ),$ while
  $|d\tau|_{g_{\varepsilon}}$ is bounded everywhere.
  Note also that $\tau  \rightarrow  0$ uniformly in
  $B_{x}(2\varepsilon )$ as $\varepsilon  \rightarrow  0.$
 
  It follows from (6.7) that the metric
  $\roof{g}{\tilde}_{\varepsilon}$ has scalar curvature
  converging to $-6$ in $B(2\varepsilon )$ as $\varepsilon
  \rightarrow  0$ and $\varepsilon^{2}R \rightarrow
  \infty .$ Further, from the remarks following (6.5)
  the metric $\roof{g}{\tilde}_{\varepsilon}$ approaches
  the hyperbolic metric in $A(\varepsilon ,2\varepsilon )$
  in the $C^{2}$ topology, as $\varepsilon  \rightarrow
   0$ with $\varepsilon^{2}R \rightarrow  \infty .$
  Thus, the metric $\roof{g}{\tilde}_{\varepsilon}$ may
  be perturbed a small amount in $C^{2}$ in
  $A(\varepsilon , 2\varepsilon )$ to match with the
  hyperbolic metric $g_{o}$ at $S(2\varepsilon )$ to
  give a smooth metric, again called
  $\roof{g}{\tilde}_{\varepsilon}$ on $N\#M,$ with scalar
  curvature close to $-6$ everywhere. Note that the
  conformal factor $\psi $ converges pointwise to the
  conformal factor $\psi_{o}$ bending the flat metric
  to the hyperbolic metric. In particular, on
  $(B_{x}(2\varepsilon ), g_{\varepsilon}),$ we have
 \begin{equation} \label{e6.8}
  \psi  = \psi_{\varepsilon} \rightarrow  2^{1/2}, 
  \ \ \ \ \ {\rm as}\ \ \
  \varepsilon  \rightarrow  0, 
   \ \ \ \varepsilon^{2}R \rightarrow  \infty .  
  \end{equation}
 
  Finally, let ${\bar{g}}_{\varepsilon}$ be the
  Yamabe metric in the conformal class of
  ${\tilde{g}}_{\varepsilon}$ with the same volume.
  Thus 
    \begin{equation} \label{e6.9}
          {\bar{g}}_{\varepsilon} = w^{4}\cdot {\tilde{g}}_{\varepsilon}, 
    \end{equation}
  where $w$ satisfies 
 	\begin{equation} \label{e6.10}
 	w^{5}  \bar{s}_{\varepsilon} = 
 	- 8\Delta w + {\tilde {s} }_{\varepsilon}w.  
 	\end{equation}
  Noting that $\min \ w <  1 <   \max \ w$ and
  evaluating (6.10) at points realizing the minimum
  and maximum of $w,$ gives the estimates
 \begin{equation} \label{e6.11} 
         \min \ w \geq
             \inf  \left|  \frac{{\tilde{s}}_\varepsilon}
                 {{\bar{s}}_\varepsilon}
            \right|,  \ \ \ 
    \max \ w \leq
     \sup\left|
         {   {\tilde{s}_\varepsilon}  \over  {\bar{s}_\varepsilon}    }
        \right|.  
 \end{equation}
  Since ${\tilde{s}}_{\varepsilon}$ converges to $- 6$ 
   in the $C^{o}$ norm as $\varepsilon  \rightarrow  0 $ 
  and $\varepsilon^{2}R \rightarrow  \infty ,$ each
  ratio in (6.11) must converge to $1,$ so that
 	\begin{equation} \label{e6.12} 
             w \rightarrow  1
 	\end{equation} 
 in the $C^{o}$ topology. Thus, the
 family of Yamabe metrics $\{{\bar{g}}_{\varepsilon}\}$ on 
 $N\#M$ have scalar curvatures converging to $-6 = s_{g_{o}},$ 
 and clearly converge smoothly to the hyperbolic metric
 on $M$ away from the point $x.$ The manifold $N$ is
 being crushed to the point $x$ under
 $\{{\bar{g}}_{\varepsilon}\},$ as $\varepsilon \rightarrow  0.$

  It is easy to see that
  $\{{\bar{g}}_{\varepsilon}\}$ is a degenerating
  family in the sense of (0.8), that is
 $$\int_{N \# M}
     |z_{{\bar{g}}_{\varepsilon}}|^{2}\ dV_{\roof{g}{\bar}_{\varepsilon}} 
  \rightarrow  \infty , \ \ \ {\rm as}\  \varepsilon  \rightarrow  0, 
 $$
  provided $(N, \gamma ) \neq  (S^{3}, g_{can}).$
  Namely, in this case, note that the metric $g$ on
  $N\setminus \{y\}$ satisfying (6.2) 
  has a definite amount of curvature, say 
 	 \begin{equation} \label{e6.13} 
 	    0 <  \int_{B_{p}(R)}  |z_{g}|^{2}\   dV_{g} 
 	       = \kappa  < \infty , 
 	 \end{equation}
  for $R$ large. The scaling properties of curvature then imply that
    \begin{equation} \label{e6.14}
          \int_{ N \# M }
              \left| z_{ 
                    \tilde{g}_{\varepsilon}
                  } 
               \right|^{2} 
            \   dV_{\tilde{g}_\varepsilon}
    \geq  
         \tfrac{1}{2} \displaystyle
         \int_{B( \varepsilon )}
        |z_{g_{\varepsilon}} |^{2}      \     dV_{g_{\varepsilon}} 
     =
          \frac{\kappa}{2}   \cdot \frac{R}{\varepsilon}
          \rightarrow  \infty ,\ \ \ \ {\rm  as }\ \varepsilon  \rightarrow  0,
    \end{equation}
  where the first inequality follows from (6.8) and
  elliptic regularity applied to the equation
  defining $\psi .$ 
  (The factor $\frac{1}{2}$ can be replaced by any 
  number $ <   1).$ 
  Similarly, the estimate (6.14) then also holds for the Yamabe
  metrics $\{ \bar{g}_\varepsilon \},$ using
  (6.12) and elliptic regularity for the equation
  (6.10) defining $w$. 
 
  Observe however that all blow-ups by the $L^{2}$
  curvature radii of $\{ \bar{g}_{\varepsilon}\},$
  as $\varepsilon  \rightarrow  0,$ \  
   $\varepsilon^{2}R \rightarrow  \infty ,$ 
   have either flat limits, or are exactly the 
   limit metric $g$, given by (6.1)
  above. It is clear that this metric can only be a
  super-trivial solution of the static vacuum
  equations, with potential function $u$ identically
  zero. For instance, since the metric $g$ is
  asymptotically flat and smooth everywhere, Theorem
  1.1 implies that any non-trivial static vacuum
  solution would have to be the Schwarzschild metric, which is
  obviously not isometric to $g.$ Thus all blow-up
  limits of $\{{\bar{g}}_{\varepsilon}\}$ are
  either trivial, i.e. flat, or are super-trivial. It
  is easily seen that hypothesis (i) and, using
  Theorem 2.10, hypothesis (ii) of Theorem A(II) are
  satisfied. In particular, this shows that the
  assumption (iii) in Theorem A(II) is necessary.
 
  We make several further remarks on this
  construction, c.f. also \S 6.3.
 
 \begin{remark} \label{r 6.1. (i).}
   {\bf (i).}
   It is obvious that this construction can be
   carried out on any collection of disjoint
   $4\varepsilon $ balls in $(M, g_{o}).$ Choosing a
   sequence $R_{i} \rightarrow  \infty $ and
   $\varepsilon_{i} \rightarrow  0,$ so that
   $\varepsilon_{i}^{2}\cdot  R_{i} \rightarrow  \infty
   ,$ one may construct sequences of Yamabe metrics
   $\{g_{i}\}$ on manifolds of the form 
$$ M_{k_{i}} =\left( \overset{k_i}{\underset{1}{\Huge\#}} N_{j}\right)\# M,$$
  where $N_{j}$ is any closed oriented manifold 
  with $\sigma (N_{j}) >  0$ and $\{k_{i}\}$ 
  either bounded or divergent to $\infty .$  
  Such sequences satisfy $v^{2/3}s(g_{i}) \rightarrow  v^{2/3}s(g_{o}),$ 
  and may be arranged that their curvature blows up 
  on a progressively denser set in any prescribed domain 
  $(\Omega ,g_{o}) \subset  (M,g_{o}).$
 
 {\bf (ii).}
  One may also glue in metrics of the form (6.1),
  where the Green's function based at $y$ is replaced
  by a finite sum of Green's functions based at a
  finite number of poles $\{y_{i}\}\in N.$ This gives
  rise to a complete scalar-flat metric with a finite
  number of asymptotically flat ends $E_{i},$ each of
  which may be glued as above into a small ball of a
  hyperbolic manifold $M_{i}.$ The limit as
  $\varepsilon  \rightarrow  0$ is then a finite number
  of hyperbolic manifolds, glued together at one
  point.
 
 {\bf (iii).}
  It is easy to see that this construction does not
  require that the gluing be into a hyperbolic
  manifold $(M, g_{o}).$ The same procedure is valid
  on any 3-manifold $(M, g)$  with Yamabe metric $g$ of
  scalar curvature $s_{g}  <   0.$
 
  Namely, it is a standard result that given any such
  $(M, g)$ and $x\in M,$ there is a $\delta  = \delta (M, g, x) >  0$ 
  such that the metric $g$ may
  be deformed to a metric $\hat{g}$ within
  $B_{x}(\delta )$ so that $\hat{g}$ is the
  constant curvature metric of scalar curvature
  $s_{g}$ in $B_{x}(\delta^{2})$ and such that the
  scalar curvature of $\hat g$ is arbitrarily
  close, (depending only on $\delta ),$ to the
  constant $s_{g}.$ Such deformations are local, in
  the sense that $\hat{g} = g$ outside
  $B_{x}(\delta ),$ and further small, in the sense
  that $|\hat{g}- g|_{C^{o}} \leq  \delta'  =
  \delta ' (\delta )$ in $B_{x}(\delta ).$ A proof of
  this fact is given for instance in [Kb1, Lemma
  3.2].
 
  As an example, for any $\varepsilon > 0,$ this
  glueing construction may be applied to the manifold $(N\# M,
  \bar{g}_{\varepsilon})$ constructed above;
  thus, one may choose $\varepsilon_{1} <  < \varepsilon $ 
  and glue another manifold $N_{1}$ into a small
  $\varepsilon_{1}$-ball in $(N, {\bar{g}}_{\varepsilon}) 
            \subset  (N \# M, \bar{g}_{\varepsilon})$ 
   for instance, giving Yamabe
  metrics on $N_{1}\# N \# M.$ Again, this process may be
  repeated inductively. Thus, one produces examples
  with curvature going to infinity at different
  scales or rates.  \end{remark} 

  These variations of the construction have all
  blow-up limits which are either super-trivial or
  trivial solutions of the static vacuum equations,
  c.f. Remark 6.8 for further discussion; (note the
  one possible exception there however). Hence the
  geometry of blow-up limits near the locus $\{u=0\}$
  can be quite arbitrary, in contrast to the
  geometric structure of blow-ups obtained when the
  hypotheses of Theorem A(II) are satisfied. The only
  common structure of the blow-up limits here is that
  they are scalar-flat and asymptotically flat.

\subsection{\bf Example 2} 
  Let $(M, g_{o})$ be as in Example 1. Here, we will
  construct degenerating sequences of Yamabe metrics
  on connected sums, whose blow-up limits are the
  (doubled) Schwarzschild metric.
 
  As indicated in the discussion in \S 3.1, static
  Einstein metrics are closely related to solutions,
  possibly only locally defined, of the equation
  $L^{*}u = 0.$ Consider then static Einstein
  manifolds of the form $X = \Omega \times_{h}S^{1},$ with
  scalar curvature $-12.$ Analogous to the
  (scalar-flat) Schwarzschild metric, consider the
  metric $g$ on $\Omega  = {\Bbb R}^{+}\times S^{2},$
  given by 
 \begin{equation} \label{e6.15} 
   g = dt^{2} + f^{2}(t)ds^{2}_{S^{2}}, 
 \end{equation} where $f$
 is the solution to the following initial value
 problem; $f(0) = a >  0,$ $f'  \geq  0,$
 	\begin{equation} \label{e6.16} 
         (f' )^{2} = 1+f^{2} - (a^{3}+a)f^{-1}, 
 	\end{equation} see [Bes, Ch. 9J].
 The scalar curvature of $g$ is $-6$. Setting $h = f' ,$ 
 one computes that 
 \begin{equation} \label{e6.17}
        L^{*}h = 0, 
 \end{equation}
  and $X$ is Einstein, with $Ric_{\X} = -3\cdot g_{\X}$ 
  by Proposition 3.0. The 2-sphere 
   $ S = \{t=0\}=\{f=a\}=\{h=0\}$ 
  is totally geodesic and of
  constant curvature $a^{-2},$ while the metric $g$
  is asymptotic to the hyperbolic metric $H^{3}(-1).$
  In fact, as $a \rightarrow  0,$ the metric $g =
  g_{a}$ converges to the hyperbolic metric
  $H^{3}(-1)$ off $S$, while $S$ is being crushed to a
  point. The curvature goes to infinity in small
  neighborhoods of $S,$ and if one rescales $g_{a}$ by
  the $L^{2}$ curvature radius, the blow-ups converge
  to the Schwarzschild metric.

\medskip

  Thus, as in Example 1, choose a small ball
  $B_{x}(\varepsilon )\subset (M, g_{o}).$ The smallest
  choice for $\varepsilon $ is $\varepsilon  \sim
  a^{1/3}.$ For any $\varepsilon  \geq  a^{1/3},$ the
  neighborhood $N_{\varepsilon /2,\varepsilon} = \{x:\dist(x, S) 
  \in  (\varepsilon /2,\varepsilon )\}$ of $S$
  in $(\Omega , g_{a})$ has uniformly bounded
  curvature as $a \rightarrow  0,$ and for $\varepsilon
  >>  a^{1/3}$ is almost isometric to the
  hyperbolic annulus $A(\varepsilon /2,\varepsilon )
  \subset  B_{x}(\varepsilon ).$ As in Example 1, we
  assume for simplicity that $\varepsilon  >> a^{1/3},$ 
  but $\varepsilon  \rightarrow  0$ as 
  $a \rightarrow  0.$ Thus, by a small smooth
  perturbation, one may glue on $N_{0,\varepsilon} =
  I\times S^{2}$ with metric $g_{a}$ onto $A(\varepsilon
  /2,\varepsilon )$ to obtain a metric $g_{a}$ on
  $ M\setminus (3$-ball), with scalar curvature almost $-6,$ and
  with totally geodesic, constant curvature boundary
  $S = S^{2}(a).$
 
  This process may be performed on any other
  hyperbolic manifold $(M' , g_{o}),$ and by matching
  at the isometric boundaries $S$, one obtains a smooth
  metric $g_{a}$ on $M' \# M,$ with scalar curvature
  almost $-6,$ and with $g_{a}$ isometric to $g_{o}$
  outside the $\varepsilon$ -neighborhood of $S.$ The
  metric $g_{a}$ may then be conformally deformed to
  a Yamabe metric ${\bar{g}}_{a}.$ One proves as
  in (6.10)-(6.12) that ${\bar{g}}_{a}$ is
  almost isometric to $g_{a};$ the conformal factor
  $w$ converges to 1 in the $C^{o}$ topology as 
  $ a \rightarrow  0.$
 
  Here, the limit of $(M' \#M, {\bar{g}}_{a})$ as
  $a \rightarrow  0$ is the 1-point union, at $x,$ of
  the hyperbolic manifolds $M$ and $M'.$ The limit
  of blow-ups by the $L^{2}$ curvature radius at $x,$
  or at points where $h_{a} = \frac{1}{2}$ for
  example, is the Schwarzschild metric (0.17) doubled
  isometrically across the event horizon.
 
  Consider briefly the behavior of the function $h =
  h_{a},$ as $a \rightarrow  0.$ Clearly, on any
  fixed interval $[t_{o},\mu ),$ for $t_{o} >  0,$
  the functions $h_{a}$ converge smoothly to the
  function $\cosh t >  1.$ However, from (6.16),
  $h_{a}(0) = 0,$ $h_{a}' (0) = a^{-1},$ so that
  $h_{a}$ increases very rapidly from 0 to 1; compare
  with Propositions 3.18 or 3.19. In particular, the
  analogue of the descent construction of Theorem
  3.10, with the function $h_{a}$ in place of
  $u_{i},$ gives rise to blow-ups limits given by the
  Schwarzschild metric with limit potential function
  $h = (1- 2mt^{-1})^{1/2},$ c.f. (0.17).

\medskip

  However, for the metrics ${\bar{g}}_{a}$ on
  $M' \#M,$ there is no non-trivial descent down the
  $u$-levels, that is the function $u$ does not satisfy
  the hypothesis (iii) of Theorem A(II), or
  Proposition $3.18/3.19.$ In fact, the limit
  function $u$ is identically $0$ on the (doubled)
  Schwarzschild blow-up limit, giving a super-trivial
  solution as in Example 1.
 
  To see this, if $u$ were non-zero on the
  Schwarzschild blow up, one must have
 \begin{equation} \label{e6.18} u = c\cdot  h,
 \end{equation} for some constant $c \neq  0.$ (Since
 the limit is the Schwarzschild metric, as noted in
 Example 3.9 and Remark 5.6, the descent construction
 in Theorem 3.10 terminates at levels where $u_{a}$
 is bounded away from 0, so that the renormalization
 to $ \bar{u}_{a}$ is not necessary).
 
  Now, by Proposition 2.9, $\Delta u$ is
  uniformly bounded in $L^{2}$ (as $a \rightarrow  0).$ 
 Thus, 
	  \begin{equation} \label{e6.19}
	 \Delta u \rightarrow  0 \ \ \ \ \ 
         {\rm in} \  L^{2},
	 \end{equation} 
 in the blow-ups, (as in the proof of Proposition 3.1). 
 Thus, in passing from $M' $ into
 $M$ through S, the function $u$ must change sign; in
 the blow-up Schwarzschild limit, $u$ is {\it odd}  w.r.t.
 reflection in the core~$S,$ and is asymptotic to $- c$
 or $+c$ at either end of the doubled Schwarzschild
 metric, corresponding to a region in $M' $ or $M$
 respectively. This implies that the functions $u =
 u_{a}$ for the metrics ${\bar{g}}_{a}$ converge
 to a limit function $u$ on the hyperbolic manifolds
 $M' \setminus \{x\}$ and $M\setminus \{x\}$ with 
 \begin{equation} \label{e6.20} 
  u(y) \rightarrow  - c \ \ \ \ \ {\rm as}\ 
  y \rightarrow  x \ {\rm in}\ M' ,{\rm  while }
 \end{equation}
 $$
  u(y) \rightarrow  +c \ \ \ \ \ {\rm as}\ y \rightarrow  x \ {\rm in \ }  M. 
 $$

  It turns out that the behavior (6.20) is never
  possible when $c \neq  0,$ c.f. Remark 6.8,
  although the proof is non-trivial. We prove the
  impossibility of (6.20) here in a simple special
  case, namely when $M = M' .$ In this situation,
  there is an isometry $\iota $ of $(M \# M,$
  ${\bar{g}}_{a})$ defined by reflection in the
  core 2-sphere $S = S^{2}(a).$ Let $u_{a}'  =
  u_{a}\circ\iota .$ From the invariance of the
  $L^{2}$ metric (1.3), the $L^{2}$ orthogonality of
  $L^{*}u$ and $\xi ,$ and the invariance of
  ${\bar{g}}_{a}$ under $\iota ,$ it follows
  easily that $L^{*}(u_{a}- u_{a}' ) = 0$ on 
  $(M\# M, {\bar{g}}_{a}).$ 
  But $\Ker L^{*}= 0,$ (c.f. the
  beginning of \S 2), so that $u_{a} = u_{a}' .$
  Thus, $u_{a}$ is an even function w.r.t. reflection
  in $S$. This property passes to the limit, (or the
  blow-up limit), and hence shows (6.20) is
  impossible.

 \begin{remark} \label{r 6.2.(i).}
   {\bf (i)}
   With only minor changes, this construction can be
   performed more generally on any pair of manifolds
   $M, N,$ with Yamabe metrics of negative and equal
   scalar curvatures, using the local deformation to
   hyperbolic metrics in Remark 6.1(iii). We point
   out that both the parameters $\varepsilon $ and $a$
   must be chosen sufficiently small, (in addition to
   satisfying $a^{1/3} <  < \varepsilon ),$ depending
   on the choice of Yamabe metrics on $M$ and $N.$
 
  To compare with Example 1, suppose $g$ is any
  Yamabe metric on a manifold $M$ with scalar
  curvature~$-6$  say. Note that any manifold $N$ with
  $\sigma (N) > 0$ carries unit volume Yamabe
  metrics with scalar curvature $-\delta ,$ for any
  $\delta > 0.$ Such a metric may be rescaled so
  that the scalar curvature is $-6$ , so that it then
  has volume on the order of $\delta^{3/2}$ and
  diameter on the order of $\delta^{1/2}.$ The
  construction in Example 2 may then be carried out
  on the pair $M$ and $N,$ (here $\varepsilon $ will be
  much smaller than $\delta ).$ This gives rise to a
  sequence of Yamabe metrics $\{g_{i}\}$ on $N\#M,$
  which converge smoothly on $M\setminus \{pt\}$ to $g_{o},$
  while crushing $N$ to a point. Thus, this sequence
  has the same basic features as the sequence in
  Example 1, except that blow-up limits are given by
  the isometrically doubled Schwarzschild metric. The
  triviality of the limit potential function $u$ is
  discussed in Remark 6.8.
 
 {\bf (ii).}
  In fact, the construction in Example 2 is more
  general in most respects than that of Example 1,
  since one may form in addition connected sums with
  all manifolds $N$ satisfying $\sigma (N) \geq  0,$
  (for example all graph manifolds). Namely, any such
  $N$ has unit volume Yamabe metrics with scalar
  curvature $-\delta .$ These may be scaled to make
  the scalar curvature $- 6$ and volume $\sim
  \delta^{3/2}$ and the construction on $N\#M$ proceeds
  as above. When $\sigma (N) = 0,$ the diameter of
  these rescaled metrics may remain large however, so
  that $N$ may no longer be crushed to a point, but
  be (volume) collapsed, away from $\{pt\},$ to a
  possibly arbitrarily long lower dimensional space.
  (The volume collapse may not necessarily be with
  uniformly bounded curvature).
 
 {\bf (iii).}
  Similarly, all the constructions mentioned in
  Remark 6.1 can also be recaptured by modifications
  of the construction in Example 2, as above, except
  the construction in Remark 6.1(ii). Regarding this
  construction, we note that there are no smooth and
  complete static vacuum solutions which have more
  than 2 ends; for instance there is no complete
  static vacuum solution on a $k$-punctured 3-sphere,
  with $k \geq  3.$ The Schwarzschild metric is a
  complete, conformally flat and scalar-flat metric
  on a 2-punctured 3-sphere.  \end{remark} 

\subsection{} 
  In this subsection, we analyse in somewhat greater
  detail the structure of the degeneration in
  Examples 1 and 2, in particular the structure
  induced on the limit (singular) manifold.
 
  Throughout this subsection, let 
 \begin{equation} \label{e6.21} 
   M_o = \bigcup_1^n M_k
 \end{equation} 
 be a finite collection of connected
 closed hyperbolic manifolds, identified at a point
 $x,$ with hyperbolic metric $g_{o}.$ Let $\{g_{i}\}$
 be a sequence of Yamabe metrics on a closed
 connected manifold $M,$ converging in the
 Gromov-Hausdorff topology to $(M_{o}, g_{o})$ and
 converging smoothly to the hyperbolic metric $g_{o}$
 on $M_{o}\setminus \{x\}.$ For example, $M$ may be the
 connected sum of the components $M_{k},$ or $M$ may
 be of the form $M = N\#M_{o},$ where $M_{o}\setminus \{x\}$ is
 connected, (i.e. $n = 1),$ and $\sigma (N) > 0,$
 so that the sequence $\{g_{i}\}$ crushes $N$ to the
 point $\{x\}.$
 
  We also assume throughout \S 6.3, (until Remark 6.8
  at the end), that $\{g_{i}\}$ has no blow-ups which
  are non-trivial solutions to the static vacuum
  equations with potential $u \not \equiv 0$. By the various constructions in \S 6.1
  and \S 6.2, such sequences exist at least for many
  configurations $M_{o}.$
 
  From Theorem 2.10, it is clear that the $L^{2}$
  norm of $z^{T}$ is uniformly bounded for the
  sequence $\{g_{i}\},$ and hence by (2.35), $\lambda
  $ is bounded away from 0 and so $\delta  =
  1-\lambda $ is bounded away from 1. Further, from
  Proposition 4.2 combined with Theorem A, the function $f = f_{i}$ does not
  converge to any constant in $L^{2}.$ In particular
  we see that, for some subsequence,
 \begin{equation} \label{e6.22} 
    \delta  = \delta_{i} \rightarrow  \delta_{o} >  0.  
  \end{equation} 
 (If $\delta_{i} \rightarrow  0,$ then from (2.35),
 $z^{T} \rightarrow  0$ in $L^{2},$ so that from the
 proof of Proposition 4.2, it follows that $u
 \rightarrow  1$ or $f \rightarrow  0$ in $L^{2},$
 which is impossible by the statement above).
 
  The convergence of the metrics $g_{i}$ to the limit
  $g_{o}$ is smooth away from $x.$ Hence for instance the trace-less Ricci curvature $z_{i}$ converges smoothly to the limit $z = 0$
  away from $x.$ Since $f_{i}$ is uniformly bounded in
  $T^{2,2}(M, g_{i})$ by Proposition 2.9, $f_{i}$
  converges at least weakly in $T^{2,2}$ away from
  $x$ to a limit function $f \in  T^{2,2}(M_{o},
  g_{o}),$ while $\xi_{i}$ converges weakly in
  $L^{2}$ to a limit $\xi  \in  L^{2}(M_{o}, g_{o}).$
  In fact, the convergence of $f_{i}$ and $\xi_{i}$
  to their limits is smooth away from $x,$ as one sees
  by applying elliptic regularity to the equation
  (2.29). Since $\Delta f \in
  L^{2}(M_{o}\setminus \{x\}, g_{o}),$ and each component
  $(M_{k}\setminus \{x\}, g_{o})$ of $(M_{o}\setminus \{x\}, g_{o})$
  extends smoothly to the hyperbolic metric on
  $M_{k},$ elliptic theory implies that $f$ is
  bounded in $L^{2,2}$ and hence is a $C^{1/2}$
  function on each $M_{k},$ c.f. [GT, Thm.8.12].
 
  The next result shows that $f$ extends continuously
  through $x$ on $M_{o}.$

 \begin{lemma} \label{l 6.3.}
    On the limit $(M_{o}\setminus \{x\}, g_{o}),$ we have
 \begin{equation} \label{e6.23} 
  f(y) \rightarrow  - 1 \ \ \ \ \ {\ as}\ y 
       \rightarrow  x, 
  \end{equation} 
 in any component
 $M_{k}$ of $M_{o}\setminus \{x\}.$
 \end{lemma} 
\begin{pf}
  To see this, first note that since the curvature of
  $\{g_{i}\}$ blows up near $x,$ Theorem 3.3 implies
  that $u_{i}$ goes to 0 somewhere near $x.$ Since
  there are no non-trivial blow-up limit solutions,
  Proposition 3.18 and Theorem 3.10 imply that for
  any $\varepsilon  >  0,$ with $B_{i} =
  (B_{x}(\varepsilon ), g_{i}),$
 	\begin{equation} \label{e6.24} 
        \osc_{B_{i}} u_{i} \leq  \delta  = \delta (\varepsilon ), 
        \end{equation}
  where $\delta (\varepsilon ) \rightarrow  0$ as
  $\varepsilon  \rightarrow  0.$ (Note also that $T_{i}$
  is uniformly bounded by the remark following the
  proof of Proposition 4.2). It follows that (6.24)
  holds on the limit $(M_{o}, g_{o}),$ which gives
  (6.23). \\
\end{pf}
 
  The basic identity (2.6), that is 
 $$z_{i} = L^{*}f_{i} + \xi_{i}, $$ 
 on $(M, g_{i}),$ passes to the limit and gives the equation 
             \begin{equation} \label{e6.25} 
              0 = L^{*}f + \xi , 
            \end{equation}
  on $(M_{o}\setminus \{x\}, g_{o}).$
 
  Of course there are no non-trivial smooth solutions
  to (6.25) on all of $(M_{o}, g_{o}),$ since the
  terms are $L^{2}$ orthogonal. However, this is not
  the case in the presence of the singular point
  $\{x\}.$

 \begin{lemma} \label{l 6.4.}
   The terms $L^{*}f$ and $\xi $ in (6.25) are not
   identically 0 on each component $M_{k}\setminus \{x\}$ of
   $M_{o}\setminus \{x\},$ so that the equation (6.25) is
   non-trivial on each $M_{k}.$
 \end{lemma} 
\begin{pf}
  Suppose $\xi = 0$ on $M_{k}\setminus \{x\}$ for some $k$. The
  trace equation (2.12), valid for $f = f_{i}$ on 
  $(M, g_{i}),$ passes to the limit away from $\{x\},$ so
  that we would then have 
          $$2\Delta f + sf = 0, $$ 
 on $M_{k}\setminus \{x\}.$ Since $f$ is bounded, a
 standard (Riemann) removable singularity result
 implies that $f$ extends smoothly over $\{x\}$ to give
 a smooth solution to this equation on $M_{k}.$
 However, since $s = - 6  <   0,$ the operator
 $-\Delta  -  \frac{s}{2}$ 
  is positive, and thus one must have $f \equiv  0$ on $M_{k}.$ This
 contradicts Lemma 6.3. \\
\end{pf}
 
  In particular, on each $(M_{k}\setminus \{x\}, g_{o}),$ we
  have a non-trivial solution of the trace equation
 	\begin{equation} \label{e6.26} 
 	2\Delta f + sf = \tr \xi .  
 	\end{equation}

  One sees easily that if $\xi $ were bounded in a
  neighborhood of $\{x\}$ in $(M_{k}\setminus \{x\}, g_{o}),$
  then (6.25) holds weakly across $\{x\},$ which
  implies that $\xi $ vanishes; more precisely, one
  uses a simple cutoff function argument to prove
  this. Thus $\xi $ is unbounded near $\{x\},$ as is
  $\tr\xi ,$ on each component $M_{k}\setminus \{x\}.$
 
  To analyse the behavior of $\tr \xi ,$ since $L(\xi) = 0,$ 
  $\delta\xi  = 0$ and $z = 0$ on
  $(M_{o}\setminus \{x\}, g_{o}),$ we have from the defining
  equation (2.2) \begin{equation} \label{e6.27}
 \Delta \tr\xi  + \frac{s}{3}\tr\xi  = 0.
 \end{equation}
  Exactly as in the proof of the non-triviality of
  (6.25), note that this equation has no bounded
  non-zero (weak) solutions on compact manifolds.
  Since $\tr\xi\in L^{2}(M_{o}\setminus \{x\}),$ and $\tr\xi $
  is smooth away from $x,$ it follows that, (up to a
  multiplicative constant), the unique non-zero
  solution to (6.27) is given by the Green's function
  $G_{x}(y)$ for the positive operator
  $-\Delta  -  \frac{s}{3} = -\Delta + 2.$
 
  Thus, $\tr\xi $ is asymptotically of the form
  $c_{k}/t,\  t = \dist_{g_{o}}(\cdot  , x),$ as 
  $t \rightarrow  0$ on each $M_{k},$ i.e.
 	\begin{equation} \label{e6.28} 
 		 t\cdot  \tr \xi \rightarrow  c_{k}, 
 	\end{equation}
  for some constant $c_{k} \neq  0.$ The next result
  evaluates the constant $c_{k},$ c.f. also Lemma 2.3.
 \begin{lemma} \label{l 6.5.}
   On each component $M_{k}\setminus \{x\}$ of $M_{o}\setminus \{x\},$
   we have
 \begin{equation} \label{e6.29} 
  \int_{M_{k}}|\xi|^{2}
 = -\frac{s}{3}\int_{M_{k}}\tr \xi  = 4\pi c_{k} > 0.  
 \end{equation}
 \end{lemma}
\begin{pf}
  Let $ \eta = \eta_k$ be a cutoff function with $ \eta \equiv  1$
  on $ M_{k}\setminus B_{x}(\varepsilon ), $ 
  $ \eta \equiv  0$ on $B_{x} (\varepsilon /2),$  
 for $\varepsilon$  small. Multiplying the
  equation (6.25) by  $\eta \xi$ gives
 \begin{equation} \label{e6.30)}
     \int \eta |\xi |^2 = \int \langle L^* f, \eta \xi\rangle =
     \int f L(\eta \xi) = f(x) \int L (\eta \xi) + O(\varepsilon),
 \end{equation}
  where $O(\varepsilon ) \rightarrow  0$ as $\varepsilon
  \rightarrow  0.$ The last equality follows from
  Lemma 6.3. As above in (6.27), from (2.2) and the
  facts that $L(\xi ) = 0$ and $\delta\xi  = 0,$ one
  computes 
\begin{equation} \label{e6.31}
   L(\eta \xi) =- \tr\xi\Delta \eta-2\langle d\tr \xi,d\eta\rangle
                          +\langle D^2\eta,\xi\rangle.
 \end{equation} 
 Thus, the divergence theorem gives 
 \begin{equation} \label{e6.32}
  \int L(\eta\xi) = \int \eta \Delta\tr\xi 
             + \int \langle d\eta, \delta\xi\rangle = \int \eta\Delta\tr\xi
             = -{s\over 3} \int \eta \tr \xi, 
 \end{equation}
  where the last inequality follows from (6.27).
  Hence, combining (6.30) and (6.32) and letting
  $\varepsilon  \rightarrow  0$ gives 
 $$
 \int_{M_{k}}|\xi|^{2} = - \frac{f(x)s}{3}  \int_{M_{k}} \tr \xi , 
 $$
  and the first equality in (6.29) follows from Lemma 6.3.
 
  For the second equality, apply the divergence
  theorem to (6.27) to obtain 
    $$-\frac{s}{3}
          \int_{M_{k}\setminus B_{x}(\varepsilon )}\tr \xi  =
           \int_{S_{x}(\varepsilon )}  \langle d\tr\xi , \nu\rangle , $$
 where $\nu $ is the outward unit normal. Since 
 $\nu = -\frac{d}{dt},$ use (6.28) to obtain
 	$$\lim_{\varepsilon\to 0}  
 	       \int_{S_{x}(\varepsilon)}
                \langle d\tr\xi , \nu\rangle  = 4\pi c_{k}. $$\\
\end{pf}
 
  The maximum principle applied to (6.27) implies
  that $\tr\xi $ does not change sign on any component
  $M_{k}\setminus \{x\}$ of $M_{o}\setminus \{x\}.$ Hence, from Lemma
  6.4, on each $M_{k},$ 
 \begin{equation} \label{e6.33} 
            \tr \xi  >  0.  
 \end{equation}
 
  Returning to (6.26), it follows from Lemma 6.3 and
  (6.28) that near $x$ on $M_{k}\setminus \{x\}, f$ has the
  expansion 
 \begin{equation} \label{e6.34} 
  f = - 1 + a_{k}t + o(t), \ \ {\rm as}\  t \rightarrow  0, 
 \end{equation}
 where $a_{k} = c_{k}/4 > 0.$ In particular, $f$
 extends to a Lipschitz function on each $(M_{k}, g_{o}).$
 
  Finally, from (6.25), we have 
 \begin{equation}
 \label{e6.35} LL^{*}f = 0, 
 \end{equation} on
 $(M_{o}\setminus \{x\}, g_{o}).$  Since $z = 0$ on
 $M_{o}\setminus \{x\},$ from (2.27), (6.35) is of the form
 	\begin{equation} \label{e6.36}
 	2\Delta\Delta f +
 	\frac{5}{3}s\Delta f + \frac{s^{2}}{3}f = 0.  
 	\end{equation}
  It is easy to see, again from standard removable
  singularity results, that there are no non-zero
  $C^{2}$ solutions to (6.36) on $M_{k},$ so that $f$
  does not extend to a $C^{2}$ function on $M_{k}.$
  Observe that since $\tr \xi > 0$ by (6.33), the
  maximum principle applied to the trace equation
  (6.26) implies 
 	\begin{equation} \label{e6.37} 
 	  f  <   0 \ \ \ \ \ {\rm on} \ M_{o}.  
 	\end{equation}

  Recall that $\tr\xi $ 
  is, (up to a multiplicative constant), 
  the positive Green's function or
  fundamental solution of the operator
  $\Delta +\frac{s}{3}$ on $M_{k},$ with pole
  at $x.$ Analogously for $f,$ we have:
 	\begin{proposition} \label{p 6.6.}
 	  The function $f$ satisfying (6.35) and (6.34) is,
 	  up to a multiplicative constant, the fundamental
 	  solution of the elliptic operator $LL^{*}$ on
 	  $M_{k},$ with singularity at $\{x\}.$
 
  Further, the function $f$ is the unique solution of
  (6.35) on $M_{o},$ smooth on $M_{o}\setminus \{x\}$ with $|f(y) + 1| = O(t)$ as $t \rightarrow 0$, where $t(y) = dist(y, x)$ in $M_{o}.$ 
\end{proposition} 
\begin{pf}
  The leading order term of $LL^{*}$ is the
  bi-Laplacian $\Delta^{2},$ whose
  fundamental solution $F_{x}(p)$ based at $x$ is
  asymptotic to $c\cdot  t(p) = c\cdot  \dist(x, p),$
  for some $c \neq  0.$ Since $f$ is smooth away from
  $\{x\},$ this implies the first statement.  
 
  The second statement follows from the fact that the
  only solution of the equation 
   $  LL^{*}h = 0 $, with  $|h| = O(t)$, as $t \rightarrow 0$,
  is $h \equiv  0.$ To see this, standard elliptic regularity applied to the equation (6.36) away from $x$ implies that
$$|\Delta h| \cdot t^{2} = o(t) \  {\rm and} \ |\nabla h| = O(t),$$
as $t \rightarrow 0$. Now pair (6.36) with $\eta \cdot h$, where $\eta$ is a cutoff function as in Lemma 6.5, with supp$|\nabla \eta| \subset A(\varepsilon, 2\varepsilon)$, $|\nabla \eta| \leq c/{\varepsilon}$, $|D^{2} \eta | \leq c/{\varepsilon^{2}}$. Integrating by parts gives
$$\int{\eta (\Delta h)^2 + \eta |s||\nabla h|^{2} + \eta s^{2}h^{2}} =
-\int{\Delta h [h \Delta \eta + 2<\nabla h, \nabla \eta>]} +
s \int{h<\nabla h, \nabla \eta>}.$$

\noindent
The estimates above on $h$ and its derivatives imply that the right side of this equation goes to $0$, as $\varepsilon \rightarrow 0$. Since all terms on the left are non-negative, it follows that $h = 0$. 

\end{pf}
 
  Since the function $f$ exists on $M_{o}$, Proposition 6.6 implies that (conversely),
  independent of how $f$ was constructed as a limit
  of $\{f_{i}\}$ from the geometry of $(M, g_{i}),$
  on any component $M_{k}$ of $M_{o}$ there exists a
  unique solution to (6.35) or (6.36), smooth away
  from $x$, which approaches the value $-1$ at $x$ linearly in $t$.
  Such a solution has the expansion (6.34), and
  hence the constant $a_{k} = c_{k}/4 > 0$ in
  (6.34) is an invariant of the punctured manifold
  $(M_{k},x).$ It might be considered as a kind of
  mass for the operator $LL^{*}$ at $\{x\}.$ It follows
  that the terms $L^{*}f$ and $\xi $ in (6.25) exist
  and are determined by the geometry of each
  $(M_{k}, x),$  independent of the approximating
  Yamabe sequence.

\medskip

  The fact that $\xi $ and $f$ are non-zero on the
  limit $(M_{o}\setminus \{x\}, g_{o})$ shows that they detect
  the `bad' convergence of $\{g_{i}\}$ near the
  singular point $\{x\},$ even on the limit. If
  $\{\gamma_{i}\}$ is a sequence of Yamabe metrics on
  a connected hyperbolic manifold $M_{o} = M,$
  converging smoothly to $(M_{o}, g_{o})$ everywhere,
  then of course $\xi  \rightarrow  0,$ $z^{T}
  \rightarrow  0$ and $f \rightarrow  0$ smoothly.
  When $M \neq  M_{o},$ it is obvious though that the
  metric $g_{o}$ itself extends smoothly over $\{x\}$
  to give the hyperbolic metric on the closed
  manifold $M_{o}.$ This is related to the following
  consequence of Lemma 6.5.

 \begin{corollary} \label{c 6.7.}
   For $(M, g_{i})$ as above, we have
 \begin{equation} \label{e6.38}
 \lim_{i\to\infty}\int_{M_{k}}|\xi_{i}|^{2}dV_{g_{i}}
 = \int_{M_{k}}|\xi|^{2}dV_{g_{o}} 
 \end{equation} and
 \begin{equation} \label{e6.39}
 \lim_{i\to\infty}\int_{M_{k}}|L^{*}u_{i}|^{2}dV_{g_{i}}
 = \int_{M_{k}}|L^{*}u|^{2}dV_{g_{o}}.
 \end{equation}
 \end{corollary} 
\begin{pf}
  Theorem 2.2, (2.10), and the continuity of the
  scalar curvature and volume in the convergence to
  the limit $(M_{o}, g_{o})$ show that (6.38) and
  (6.39) are equivalent. Fatou's theorem, (lower
  semi-continuity), implies that, for each $k,$
 	\begin{equation} \label{e6.40}
 	\int_{M_{k}}|\xi|^{2}dV_{g_{o}} =
 	\int_{M_{k}\setminus \{x\}}|\xi|^{2}dV_{g_{o}} 
         \leq
 	\liminf_{i\to\infty }\int_{M_{k}\setminus \{x\}}|\xi_{i}|^{2}dV_{g_{i}},
 	\end{equation}
  where for the integral on the right we consider
  $M_{k}\setminus \{x\} \subset  M$ via a Gromov-Hausdorff
  approximation. \\
  Summing over $k$ then gives
 		  $$\int_{M_{o}}|\xi|^{2}dV_{g_{o}} \leq
 		\liminf_{i\to\infty}\int_{M}|\xi_{i}|^{2}dV_{g_{i}}
		= \liminf_{i\to\infty}-\frac{s}{3}\int_{M}\tr\xi_{i}\ dV_{g_{i}}
 		= -\frac{s}{3}\int_{M_{o}}\tr\xi\  dV_{g_{o}}
 		=\int_{M_{o}}|\xi|^{2}dV_{g_{o}}. 
 	$$
  Here, the first equality uses (2.14), the second
  uses the fact that $\xi_{i} \rightarrow  \xi $
  weakly in $L^{2}$ and the third follows from Lemma
  6.5. It follows that equality holds in (6.40) and
  the $\liminf$ may be replaced by limit.\\
\end{pf}
 
  Of course the $L^{2}$ norm of $L^* f_{i}$ and of
  $z_{i}$ diverges to $\infty $ as $i \rightarrow
  \infty .$ In other words, Corollary 6.6 shows that
  the splitting of the metric $g_{i}$ in (2.10) is
  continuous in the limit $i \rightarrow  \infty ,$
  but the splitting of the curvature $z_{i}$ is far
  from continuous.
 
  As a curiosity at this stage, we note that an easy
  computation, using the relation (6.29), implies
 \begin{equation} \label{e6.41}
 \frac{d}{dt}(v^{2/3}\cdot  s)(g_{o}+t\xi )|_{t=0} =
 -  v^{-1/3}\int_{M_{o}\setminus \{x\}}|\xi|^{2}dV.
 \end{equation}
  Thus deforming the metric $g_{o}$ on $M_{o}\setminus \{x\}$
  in the direction $-\xi $ increases the Yamabe
  constant $v^{2/3}\cdot  s$ of~$g_{o}.$ Of course
  since $\xi $ is unbounded, this is not a smooth
  perturbation of the hyperbolic metric on~$M_{o}.$

 \begin{remark} \label{r 6.8.}
   In this remark, we consider briefly the converse
   of the previous discussion in \S 6.1-\S 6.3,
   namely to what extent blow-up limits of Yamabe
   metrics converging to a singular hyperbolic limit
   as above {\it  must}  be super-trivial solutions
   of the static vacuum equations.
 
  Thus, let $(M, g_{i})$ be any sequence of Yamabe
  metrics converging to $M_{o}\cup X,$ where $M_{o}$
  is as in (6.21) and $X,$ (possibly empty) is a lower
  dimensional space of zero volume, compare with Remark
  6.2(ii). Suppose further that the blow-up limit at
  the singular point $x\in M_{o}$ is a complete,
  smooth, non-flat and asymptotically flat manifold
  $(N, g_{b}),$ i.e. all ends are asymptotically
  flat. All of the constructions and their variations
  in \S 6.1 and \S 6.2 satisfy these assumptions,
  (assuming the construction is not iterated).
 
  By Theorem A(I), the blow-up limit $(N, g_{b})$ is
  a solution of the static vacuum equations, which is
  non-flat; $(x$ is essentially the point of maximal
  curvature concentration). Hence, for instance by
  Theorem 1.1, either $(N, g_{b})$ is a super-trivial
  solution or it is the (isometrically doubled)
  Schwarzschild metric, with non-vanishing potential
  function $u$. Thus, we are interested in considering
  only if the latter possibility occurs.
 
  If the blow-up limit $(N, g_{b})$ is the
  Schwarzschild metric, then of course one must have
  $M = M_{1}\#M_{2},$ with the Schwarzschild neck
  $S^{2}\times{\Bbb R} $ joining the two factors. This
  situation occurs in all versions of the
  construction in Example 2 in \S 6.2. However, for
  the construction in \S 6.1, it occurs only in the
  situation of Remark 6.1(ii), where a sum of two
  Green's functions is used w.r.t. the conformal
  Laplacian on the standard $(S^{3}, g_{can}).$
 
  We may assume that $M_{2}$ is hyperbolic, and
  $M_{1}$ is either hyperbolic or is a manifold
  satisfying $\sigma (M_{1}) \geq  0.$ In the former
  case, $M_{o}$ is the one point union at $x$ of
  $M_{1}$ and $M_{2},$ while in the latter case,
  $M_{o} = M_{2}$ and $M_{1}$ is collapsed to $X.$
  (Note that one may have $X = \{x\};$ for instance,
  this is necessarily the case for the construction
  in Remark 6.1(ii), which requires $N = M_{1}$ with
  $\sigma (M_{1}) > 0).$
 
  It turns out that if $M_{1}$ is hyperbolic, then
  the Schwarzschild neck $(N, g_{b})$ must be a
  super-trivial solution, for any sequence
  $\{g_{i}\}$ satisfying the hypotheses above. On the
  other hand, if $M_{1}$ is any graph manifold, (in
  particular $\sigma (M_{1}) \geq  0),$ and
  $\{g_{i}\}$ collapses $M_{1}$ with bounded
  curvature away from the neck to $X,$ then the limit
  solution is non-trivial. The limit potential
  function $u$ is given by the Schwarzschild
  potential, odd under reflection through the event
  horizon, as in (6.20). The construction in Theorem
  3.10 produces this non-trivial blow-up limit.
 
  For the remaining possible situations, where
  $\sigma (M_{1})$ is not a graph manifold, or
  $\{g_{i}\}$ volume collapses $M_{1}$ to a lower
  dimensional space, but the collapse is with
  unbounded curvature away from the neck, the
  triviality or non-triviality of the solutions is
  unknown. (Note this situation includes the
  construction in \S 6.1, Remark 6.1(ii)).
 
  We do not give a complete proof of these statements
  here, since essentially no new ideas are needed,
  and since these statements will not be used. Instead,
  we just sketch the ideas involved.
 
  When $M_{1}$ is a graph manifold and $\{g_{i}\}$
  collapses $M_{1}$ with bounded curvature, this
  statement can be deduced quite easily from the work
  in \S 6.5, in particular from Proposition 6.9. In
  this case, in fact $\xi  \rightarrow  0,$ $z^{T}
  \rightarrow  0$ and $u \rightarrow  1$ in $L^{2}(M, g_{i});$ 
  c.f. Remark 6.10(ii). Note that $M_{1}$
  becomes invisible, in terms of volume and total
  scalar curvature, in $\{(M, g_{i})\}.$
 
  On the other hand, if $M_{1}$ is assumed
  hyperbolic, then $u$ cannot converge to 1 in
  $L^{2},$ essentially for the reasons discussed
  following (6.19). One may then prove that $u_{i}$
  is approximately 0 for $i$ large on the
  Schwarzschild neck by using the minimizing property
  Proposition 2.7 for $L^{*}u,$ together with the
  existence of non-trivial solutions $\xi,\ L^{*}f$
  to (6.25), (6.27) and (6.34)-(6.35) on
  $M_{o}\setminus \{x\},$ as discussed following Proposition
  6.6.  \end{remark}

\subsection{\bf Example 3} 
   In this next example, we let $M = S^{2}\times S^{1},$
   and consider a maximizing family $\{g_{i}\}$ of
   Yamabe metrics on $M$ with $\vol_{g_{i}}M = \vol S^{3}(1).$ 
  It is known that 
 \begin{equation} \label{e6.42} 
 \sigma (S^{2}\times S^{1}) = \sigma (S^{3}).  
 \end{equation} In fact, see [Kb1,2],
 [Sc2], there are conformally flat, spherically (i.e.
 $S^{2})$ symmetric Yamabe metrics $g_{i},$ thus of
 the form 
 \begin{equation} \label{e6.43} 
           g_{i} = dt^{2} + f_{i}^{2}(t)ds^{2}_{S^{2}}, 
 \end{equation}
 with $s_{g_{i}} \rightarrow  \sigma (S^{3}).$ As
 noted in [La], and as in Example 2, setting $h_{i} =
 f_{i}' ,$ one computes 
 \begin{equation}
        \label{e6.44} L^{*}(h_{i}) = 0, 
 \end{equation}
 w.r.t. the metrics $g_{i}$ globally on $M.$ In fact,
 c.f. [La], $$\Ker L^{*} = <h_{i}>, $$ on 
 $(M, g_{i}).$ Again from Proposition 3.0, the
 $4$-manifolds $X = M\times_{h_{i}}S^{1}$ are Einstein, with
 scalar curvature $s_{\X} \rightarrow  12,$ the scalar
 curvature of $S^{4}(1).$ (These metrics have cone
 singularities along two totally geodesic 2-spheres
 in $X$, corresponding to the locus where $h_{i} = 0).$
 
  As $i \rightarrow  \infty ,$ the manifolds
  $(S^{2}\times S^{1}, g_{i})$ converge to the space 
  $(Z, g_{o}),$ where $Z = S^{3}/\{p\}\sim\{-p\}$ is the
  3-sphere with two antipodal points identified, and
  $g_{o}$ is the canonical metric of constant
  curvature and volume 1, see [An1], [Sc2]. The
  convergence is smooth on the complement of the
  point $\{p\}\sim\{-p\},$ and the functions $h_{i},$
  for an appropriate normalization, converge to an
  eigenfunction $h$ of the Laplacian on $S^{3},$ with
  eigenvalue 3, with $h_{i}(p) \rightarrow  1,$
  $h_{i}(-p)\rightarrow  - 1.$
 
  Now consider the equation (2.10) on $(M, g_{i}),$
  i.e.  \begin{equation} \label{e6.45} L^{*}u + \xi
  = -\frac{s}{3}g.  \end{equation} The function $u =
  u_{i}$ is only determined up to functions in $\Ker L^{*}.$ 
  Using suitable multiples of the functions
  $h_{i}$ in (6.44), we may arrange that
  \begin{equation} \label{e6.46} 
    ||u_{i}||_{L^{2}} \rightarrow  \infty , 
     \ \ \ \ \ {\rm as }\ i \rightarrow  \infty .
   \end{equation} 
 (We note that it may well be
 necessary to use the functions $h_{i}$ in this way;
 there may be representatives $u'\in\{u+\Ker L^{*}\}$
 which are uniformly bounded in $L^{2}).$ In case
 (6.46) holds, as discussed in \S 3.2-\S 3.4, the
 construction of the buffered blow-up in Theorem 3.10
 requires that $u = u_{i}$ be renormalized by its
 maximum. Thus, setting ${\bar{u}} =
 \frac{u}{\sup u},$ we obtain from (6.45) that
 \begin{equation} \label{e6.47}
 L^{*}{\bar{u}}_{i} \rightarrow  0.
 \end{equation} Then in fact ${\bar{u}}_{i}$
 approaches the function $h_{i}$ in (6.44) and
 $\{{\bar{u}}_{i}\}$ limits on the eigenfunction
 $h$ of the Laplacian as above. If one performs the
 descent construction of Theorem 3.10 on
 $\{{\bar{u}}_{i}\},$ one obtains as blow-up
 limit the Schwarzschild metric (0.17) doubled across
 the totally geodesic boundary $S^{2}.$ In the
 blow-up, the functions ${\bar{u}}_{i}$ limit on
 the potential function $u$ of the Schwarzschild
 metric, and the harmonic function $u$ is odd w.r.t
 reflection in $S^{2}.$ At one end, $u \rightarrow
 +1$ while $u \rightarrow  - 1$ at the other end,
 corresponding to the two points \{p\} and $\{- p\}$
 respectively in the limit $(Z, g_{o}).$
 
  Thus in this simple example, the blow-up limit
  obtained from Theorem A(II) mirrors exactly the
  degeneration of the sequence $(M, g_{i}).$

\noindent
\subsection{\bf Example 4} 
  The discussion in \S 6.1 and \S 6.2 raises the
  question of what are the simplest kinds of
  degenerations of Yamabe metrics which {\it  are}
  modeled on non-trivial solutions of the static
  vacuum equations, (leaving Example 3 aside). Thus,
  let $M_{1}$ and $M_{2}$ be closed hyperbolic
  3-manifolds and $M = M_{1} \#M_{2}.$ Based on a
  modification of Example 2, we construct a sequence
  $\{g_{i}\}$ of Yamabe metrics on $M$ which crush
  the essential 2-sphere in $M,$ converge smoothly
  almost everywhere w.r.t. volume to the hyperbolic
  metrics on $M_{1}$ and $M_{2},$ and for which
  blow-up limits obtained via Theorem 3.10 are
  non-trivial solutions to the static vacuum
  equations. It is interesting to compare the
  construction below with the work of O. Kobayashi in
  [Kb1, Thm.2].
 
  First, return to the construction in Example 2 on
  each manifold $M_{k},$ $ k = 1,2.$ This gives a metric
  $g_{a}$ on $M_{k}\setminus B_{x_{k}}(a),$ with totally
  geodesic, constant curvature boundary $S =
  S^{2}(a),$ which is hyperbolic outside the ball
  $B_{x_{k}}(\varepsilon )$ and is ``Schwarzschild-like''
  in the neck $A(a,\varepsilon /2).$ As in Examples 1,2
  note that the curvature of $g_{a}$ is uniformly
  bounded in $A(\varepsilon /2, \varepsilon )$ provided
  $\varepsilon  \geq  a^{1/3}.$ As previously, we assume
  $\varepsilon > > a^{1/3},$ so that the
  sectional curvature of $g_{a}$ is very close,
  (depending only on $\varepsilon /a^{1/3})$ to $- 1$ on
  $A(\varepsilon /2,\varepsilon ),$ and the scalar
  curvature of $g_{a}$ is almost $- 6$ everywhere.
 
  Now instead of identifying 
     $M_{k}\setminus B_{x_{k}}(a)$
  along their isometric boundaries, we extend $g_{a}$
  past $S^{2}(a)$ with the metric $g$ on $S^{2}\times{\Bbb
  R}^{+}$ in (6.15). The resulting metric, still
  called $g_{a},$ is a complete metric on
  $\tilde{M}_{k} \equiv  M_{k}\setminus B_{x_{k}}(a)\cup
  S^{2}\times{\Bbb R}^{+}.$ Note that this metric is
  conformally equivalent, in a natural way, to a
  metric on $M_{k}\setminus \{x_{k}\}$ which extends smoothly
  over $\{x_{k}\},$ so that
  $(\tilde{M}_{k}\cup\{x_k\}, g_{a})$ is
  conformally compact. 
  (Of course $(M_{k}\setminus B_{x_k}(a), g_{a})$ 
  cannot be conformally compactified by adding a point). 
  The blow-up limit of $({\tilde{M}}_{k}, g_{a})$ 
   at or near $x_{k}$ as  $a \rightarrow  0$ is the complete,
  isometrically doubled Schwarzschild metric. The end
  of $({\tilde{M}}_{k}, g_{a})$ is asymptotically
  hyperbolic. 
   In fact, just as before on the other
  `outer' side of $S^{2}(a)$ in ${\tilde{M}}_{k},$
  the `inner' annulus 
  $A_{i} = A(\varepsilon /2, \varepsilon ) \subset  (S^{2}\times {\Bbb R}^{+}, g_a)$
  is almost isometric to an annulus
  $A_{-1}(\varepsilon /2, \varepsilon)$ in $H^3 (-1)$ 
  for $\varepsilon >> a^{1 /3}.$  In particular, the curvature
  of $g_a$ is bounded in both the inner annulus
  $A_i\subset S^2 \times \Bbb R^+$  and the (isometric) `outer'
  annulus $A_o\subset {M}_k.$
 
  To join the two inner annuli $A_{i},$ we choose a
  metric $\gamma_{a}$ on $G_{a} \equiv 
  S^{3}\setminus \ (B_{1}\cup B_{2}),\  B_{j}$ a 3-ball, so that
  $\gamma_{a}$ is isometric to $(A_{i}, g_{a})$ near
  its boundary. Thus, the metrics
  $({\tilde{M}}_{1}, g_{a})$ and
  $({\tilde{M}}_{2}, g_{a})$ can be glued along
  $(S^{3}\setminus B_{1}\cup B_{2}), \gamma_{a})$ to give a
  smooth family of metrics, still denoted $g_{a},$ on
  $M = {\tilde{M}}_{1}\cup G_{a}\cup
  {\tilde{M}}_{2} = M_{1}\#M_{2}.$ The gluing
  metric $\gamma_{a}$ is chosen to have curvature
  $z_{\gamma_{a}}$ uniformly bounded as  
  $a \rightarrow  0,$ scalar curvature 
  $s_{\gamma_{a}} \rightarrow  - 6$ 
  everywhere as  $a \rightarrow  0,$
  and to be {\it  volume collapsing,}  so that
 \begin{equation} \label{e6.48} 
       \vol_{\gamma_{a}}G_{a} \rightarrow  0, 
       \ \ \  {\rm as\ } a \rightarrow  0.  
 \end{equation}
  Of course since the metrics $g_{a}$ are not
  collapsing near the boundaries $\partial B_{i}$ of
  $G_{a},$ the collapse of $\gamma_{a}$ takes place
  on a scale much larger than $\varepsilon .$ We sketch
  the construction of such metrics.

\bigskip

 \noindent
 {\bf Construction of glueing metrics:} View $S^{3}$ as
  the union of two solid tori $D^{2}\times S^{1},$ glued
  along the torus boundary $T^{2}$ by interchanging
  the two circles in $T^{2}.$ Consider a complete
  warped product metric on $D^2\times S^1$ of the form
 \begin{equation} \label{e6.49} 
     h = g_{\D^{2}} +f^{2}d\theta^{2}, 
 \end{equation} 
 where $f$ is a
 positive function on $D^{2}.$ We assume that outside
 a compact set, $h$ is isometric to a rank 2
 hyperbolic cusp, so that in particular $f = f(r) =
 e^{-r}$ for $r$ large.
 
  It is a standard fact that $h$ may be chosen to
  satisfy 
       \begin{equation} \label{e6.50} 
    \int s_{h}\ dV_{h} <  0, 
       \end{equation} 
  c.f.  [Bes, Thm.4.32] 
  for example. Now let $\bar{h}$
 be a Yamabe metric conformal to $h$, so that
  $\bar{h} = \psi^{4}\cdot  h.$ It is also quite
 standard, (c.f. [AM, Thm.C] for a proof), that
 (6.50) and the assumptions on the asymptotic form of
 $h$ imply that ${\bar{h}}$ exists and is
 complete, with $\psi $ uniformly bounded away from $0$
 and $\infty .$ Possibly after a rescaling, the
 scalar curvature ${\bar{s}}$ of $\bar{h}$
 satisfies $\bar{s} = - 6.$
 
  Now we claim that the metric $\bar h$ is also
  a warped product, i.e. is invariant under the
  $S^{1}$ action on $D^2\times S^{1}.$ To see this, let
  $F = F_{\theta}$ be an isometry of $h$, (from the
  $S^{1}$ family). We have 
   $$F^*{\bar{h}} 
   = (F^{*}\psi )^{4}\cdot  (F^{*}h) 
   = (F^{*}\psi)^{4}\cdot  h 
   = v^{4}\cdot  {\bar{h}},  
  $$ 
  where $v = (F^{*}\psi )/\psi $ is a bounded function. The
 metrics $F^*{\bar{h}}$ and ${\bar{h}}$ are
 thus conformal Yamabe metrics, of the same scalar
 curvature $- 6$ and volume. Hence $v$ satisfies the
 equation (1.12), i.e.  $$- 6v^{5} = - 8\Delta v -  6v, $$ 
 w.r.t. the
 ${\bar{h}}$ metric. The maximum principle then
 implies that $v \equiv  1.$ Thus, $F$ is also an
 isometry of~$\bar{h}.$
 
  It follows that $\psi $ is a function on $D^{2}$
  and the defining equation (1.12) becomes, on
  $D^2,$ 
 \begin{equation} \label{e6.51}
          8\Delta_{\D^2}\psi  + 8 < \nabla \psi\, ,\, \nabla log f>  
          = 6\psi^{5} + s\psi .  
 \end{equation}
 Since $\psi $ is bounded and $s = - 6$ outside a
 compact set, the maximum principle applied to (6.51)
 shows that $\psi (x) \rightarrow  1 $ as $x
 \rightarrow  \infty $ in $D^{2}\times S^{1}.$
 
  It follows that ${\bar{h}}$ is also of the form
  (6.49), with $f$ replaced by ${\bar{f}},$ and is
  asymptotic to a rank 2 hyperbolic cusp at
  infinity.
 
  Now observe that the full curvature of
  ${\bar{h}}$ is unchanged when the length of the
  $S^{1}$ factor is changed. In fact, for any
  $\alpha  >  0 $ (small), the metrics
	 \begin{equation} \label{e6.52} 
		   {\bar{h}}_{\alpha} = {\bar{g}}_{\D^{2}} 
		  +(\alpha {\bar{f}})^{2}\ d\theta^{2}, 
	 \end{equation} 
  are locally isometric to ${\bar{h}}.$ Note that the
  length of the $S^{1}$ fiber is asymptotic to 
  $\alpha e^{-r}.$
 
  Next, for $R$ sufficiently large depending on the
  choice of $\alpha ,$ on the annulus $A(R,2R),$ we may
  bend the metric 
  ${\bar{g}}_{\D^{2}} \sim dr^{2}+e^{-2r}d\phi^{2}$ 
  slowly on $A(R,2R)$ so that
  on $A(2R- 1,2R),$ it has the form 
  $dr^{2}+(\alpha {\bar{f}})^{2}\ d\phi^{2}.$ 
  Here we assume that both $\theta $ and $\phi $ are parameters in
  $[0,2\pi ].$ The resulting metric
  ${\tilde{h}}_{\alpha}$ has uniformly bounded
  curvature, independent of $\alpha $ and $R,$ scalar
  curvature arbitrarily close, depending only on $R$,
  to $- 6,$ and is almost isometric to a rank 2
  hyperbolic cusp on $A(2R-1, 2R).$
 
  Finally, on the annulus $A(2R,4R),$  the function
  $\alpha \bar{f}$ is again bent slowly so that
  at the boundary $S(4R),$ $\alpha {\bar{f}}(4R)
  \sim  \alpha e^{-4R},\  \bar{f}' (4R) = 0$ and
  $\bar{f}$ extends smoothly as an even function
  under reflection through $4R.$ Again this may be done
  so that the curvature is uniformly bounded and the
  scalar curvature is arbitrarily close to $- 6.$
  Call the resulting metric on $D^{2}\times S^{1}$ again
  $\tilde{h}_{\alpha}.$
 
  Now the boundary $S(4R)$ is a totally geodesic, flat,
  and square torus $T^{2},$ in the sense that the
  generators $d\phi $ and $d\theta $ are orthogonal
  and of the same length. Hence the metric
  $\tilde{h}_{\alpha}$ may be doubled across the
  boundary $T^{2},$ with an interchange of the
  $S^{1}$ factors, to give a smooth metric on
  $S^{3}.$ This is the metric $\gamma_{a},$ where
  $\alpha$ in (6.52), essentially the maximal length of the
  $S^{1}$ fiber in (6.52), is relabeled to $a~=~a(\alpha ),$ 
  the size of the core $S^{2}$ in
  $(M_{k}, g_{a})$ above. We require  $a <    < \alpha ,$ 
  in fact $\varepsilon   <  <   a,$ but $a \rightarrow  0$ 
  implies $\alpha  \rightarrow  0.$
 
  By choosing $R = R_{\alpha}$ suitably, the
  resulting family of smooth metrics $\gamma_{a}$ on
  $S^{3}$ has scalar curvature converging smoothly to
  $- 6$ as  $a\rightarrow  0,$ has uniformly bounded
  curvature, and is volume collapsing in the sense of
  (6.48). Observe that this volume collapse requires
 \begin{equation} \label{e6.53}
           \diam_{\gamma_{a}}S^{3} \sim  2R \rightarrow  \infty , 
           \ \ \ \ \ {\rm  as\ } a \rightarrow  0.  
  \end{equation}
  Now the inner annuli $(A_{i}(\varepsilon /2, \varepsilon
  ), g_{a})$ in each ${\tilde{M}}_{k}$ are almost
  isometric to hyperbolic $\varepsilon $-annuli. For
  $\varepsilon   <  <   \alpha ,$ the metrics
  $\gamma_{a}$ on the geodesic ball $B_{p}(\varepsilon
  )$ contained in each $D^{2}\times S^{1} \subset  S^{3}$
  are fixed, independent of $a$ or $\alpha ;$ here $p =
  (p_{o}, \theta ),$ where $p_{o}$ is the `center' of
  $D^{2}.$ (The metric on $D^{2}$ is fixed, only the
  length of the fiber $S^{1}$ changes with $a).$ Hence,
  just as with the first or outer glueing as in
  Example 2, the metric $\gamma_{a}$ may be smoothly
  matched to $(A_{i}, g_{a}),$ keeping the curvature
  uniformly bounded and the scalar curvature
  arbitrarily close to $- 6,$ (c.f. the local
  deformation in Remark 6.1(iii)). This completes the
  construction of the glueing metrics.
 
\bigskip

  Finally, as in Examples 1, 2, let
  ${\bar{g}}_{a}$ be the Yamabe metric conformal
  to $g_{a}$ with the same volume on $M,$ so that
  $\roof{g}{\bar}_{a} = w^{4}\cdot  g_{a}.$ Since the
  scalar curvature of $g_{a}$ approaches $- 6$
  everywhere as $a \rightarrow  0,$ as in (6.12), $w
  \rightarrow  1$ pointwise as  $a\rightarrow  0.$
  Thus, the metric $\roof{g}{\bar}_{a}$ is $C^{o}$
  close to $g_{a}.$ Further, from the defining
  equation (1.12), $|\Delta w| \rightarrow  0$ pointwise as  
  $a\rightarrow 0.$ Using $L^{2}$ estimates 
  for the elliptic equation (1.12) on any
  geodesic ball $B$ about $q,$ of radius $\rho (q),$
  one obtains w.r.t. the $g_{a}$ metric,
 \begin{equation} \label{e6.54} 
   \int_{B}|D^{2}w|^{2}  <  <   \max\, \Bigl( 1, 
                   \int_{B}|z|^{2}
                          \Bigr)  
  \ \ \  {\rm as \ }  a \rightarrow  0, 
  \end{equation}
  provided $\rho (q) \leq  1.$ Of course
  $\roof{g}{\bar}_{a}$ is smoothly close to $g_{a}$
  on $M\setminus B_{x_{1}}(\varepsilon )\cup B_{x_{2}}(\varepsilon
  )).$ This completes the construction of the Yamabe
  metrics $\{\roof{g}{\bar}_{a}\}$ on $M.$
 
  The justification for this construction comes from
  the following result.
 \begin{proposition} \label{p 6.9.}
    On the family $(M, \bar{g}_{a}),$
 \begin{equation} \label{e6.55} 
   ||z^{T}||_{L^{2}} \rightarrow  0 \ \ \ \ \  
   {\rm as\ } a \rightarrow  0.  
 \end{equation}
 In particular, $\xi  \rightarrow  0$ in $L^{2},$ and
 $\delta  \rightarrow  0.$  \end{proposition} 
 
  From Proposition 4.1 for instance, it follows that
 \begin{equation} \label{e6.56} u = u_{a}
 \rightarrow  1  
 \ \ \ \ \  {\rm in}\  L^{2}, \end{equation} 
 and hence,
 (c.f. Proposition 4.2), one may carry out the
 descent construction in Theorem 3.10 to obtain
 non-trivial blow-up limit solutions of the static
 vacuum equations. In fact for an initial sequence of
 points $\{x_{i}\}$ in either $M_{k}\setminus B_{x_{k}}(a),$
 with $x_{i} \rightarrow  x_{k},\  u_{i}(x_{i})
 \rightarrow  1,$ and $a = a(i) \rightarrow  0,$
 the maximal blow-up limit $(N, g_{o}),$ c.f. \S 5.1,
 based at a buffered sequence $\{y_{i}\}$ obtained by
 $u$-descent from $\{x_{i}\}$ is the complete doubled
 Schwarzschild metric. The limit potential function
 is $u = \pm (1 - \frac{2m}{t})^{1/2},$ c.f. (0.17).
 
\medskip

\begin{pf*}{\bf Proof of Proposition 6.9.}
  To prove (6.55), we use the characterization (2.41)
  of $||z^{T}||_{L^{2}},$ i.e.
 	\begin{equation} \label{e6.57} \int_{M}|z^{T}|^{2}
 	=  \inf_{\int\phi =0 }\int_{M}|z -  L^{*}\phi|^{2}.
 	\end{equation}
  Of course (6.57) is taken w.r.t. the
  $\roof{g}{\bar}_{a}$ metric. However, using the
  fact that $\roof{g}{\bar}_{a}$ is $C^{o}$ close to
  $g_{a},$ and smoothly close away from
  $M\setminus ( B_{x_{1}}(\varepsilon )\cup B_{x_{2}}(\varepsilon)),$ 
  together with (6.54), it suffices to estimate
  the right hand side of (6.57) w.r.t. the $g_{a}$
  metric. To see this, refer to the expansion (2.42).
  Using standard formulas, c.f. [Bes, Ch.1J], under
  the conformal change from $g_{a}$ to
  $\roof{g}{\bar}_{a},$ all terms in (2.42) differ by
  arbitrarily small quantities as  $a \rightarrow  0,$
  except the term containing $|z|^{2},$ which differs
  by a term on the order of $|D^{2}w|^{2},$ where $w$
  is the conformal factor. By (6.54), this is small
  in $B_{x_{1}}(\varepsilon )\cup B_{x_{2}}(\varepsilon )$
  compared with the dominant $|z|^{2}$ term in
  (2.42).
 
  Of course, since the metric $g_{a}$ is quite
  explicit, it is easier to estimate (6.57) w.r.t.
  $g_{a}$ than w.r.t. $\roof{g}{\bar}_{a}.$ All the
  estimates to follow are then on $(M, g_{a}).$
 
  The function $\phi  = \phi_{a}$ is defined on the
  various pieces of $(M, g_{a})$ as follows. First,
  on each $M_{k}\setminus B_{x_{k}}(\varepsilon ),$ since $z = 0, $
  we set $\phi  \equiv  0,$ so that
 \begin{equation} \label{e6.58}
 \int_{M_{k}\setminus 
      B_{x_{k}}(\varepsilon )}|z - L^{*}\phi|^{2} = 0 .  
 \end{equation} 
 Next, on the (doubled) Schwarzschild-like neck $N_{a}$ joining
 the inner and outer annuli $A(\varepsilon /2,\varepsilon)$ 
 on each $\tilde{M}_{k},$ let
 \begin{equation} \label{e6.59} 
                  \phi  = h- 1,
 \end{equation} 
 where $h = h_{a}$ is the function from (6.16)-(6.17). 
 Note that $h$ is odd under reflection through the core 
 $S = S^{2}(a)$ in each $M_{k}.$ Since $L^{*}h = 0$ 
 w.r.t. the metric $g_{a},\  L^{*}\phi  = r$ on $N_{a}.$ Hence, since 
 $\vol N_{a} \rightarrow  0$ as $a \rightarrow  0,$ and
 the scalar curvature is uniformly bounded,
	 \begin{equation} \label{e6.60} 
	   \int_{N_{a}}|z - L^{*}\phi|^{2} \rightarrow  0, 
           \ \ \ {\rm as\ } a \rightarrow  0.
	 \end{equation} 
 To examine the behavior on the outer
 or first glueing annulus $A_{o} = A(\varepsilon /2,
 \varepsilon ) \subset  M_{k},$ the metric $g_{a}$ has
 uniformly bounded curvature on $A_{o}$ as 
 $ a \rightarrow  0.$ Now from the discussion in Example
 2, in $A_{o},$ 
	 \begin{equation} \label{e6.61} 
	 h \sim  \cosh t = 1 + \frac{1}{2}t^{2} + o(t^{2}),
	 \end{equation} 
 so that $\phi  \sim  0.$  Further,
 since $h \sim  h(t)$, 
 	 \begin{equation} \label{e6.62}
 	|D^{2}h| \sim  |h' D^{2}t + h'' (dt)^{2}| \leq  C,
 	\end{equation} 
 uniformly in $A_{o}$ as 
 $a \rightarrow  0,$ by (6.61). Hence the function
 $\phi $ may be smoothed in $A_{o}$ so that $\phi
 \equiv  0$ near the outer boundary $S(\varepsilon )$ of
 $A_{o},$ keeping $D^{2}\phi $ uniformly bounded.
 Since $\vol A_{o} \rightarrow  0$ as $a \rightarrow  0,$ 
 it follows that (6.60) holds also over $A_{o}.$
 
  Finally, in the glueing region $G_{a}$ joining
  $\tilde{M}_{1}$ and $\tilde{M}_{2},$ and
  containing the inner annuli $A_{i},$ the curvature
  of $g_{a}$ is uniformly bounded while 
  $\vol G_{a} \rightarrow  0,$ so that 
       \begin{equation}
       \label{e6.63} 
        \int_{G_{a}}|z -  L^{*}\phi|^{2}
       \leq  \delta_{a} + 2\int_{G_{a}}|L^{*}\phi|^{2}
       \leq  \delta_{a}+ c\int_{G_{a}}(|D^{2}\phi|^{2}+
       \phi^{2}), 
  \end{equation} where $\delta_{a}
 \rightarrow  0$ as $a \rightarrow  0.$ Now in the
 inner annulus $A_{i} = A(\varepsilon /2, \varepsilon )
 \subset  S^{2}\times{\Bbb R}^{+},$ by considerations
 similar to (6.61), $\phi  \sim  - 2,$ while as in
 (6.62), $|D^{2}\phi|$ is uniformly bounded. Hence,
 using (6.53), we may extend $\phi $ from a
 neighborhood of both components of $\partial G_{a}$
 into $G_{a}$ - for instance as a function of a
 suitable smoothing of the distance function to
 $x_{1}$ (or $x_{2})$ - so that $\phi $ and
 $D^{2}\phi $ remain uniformly bounded and so that,
 (most importantly), 
  \begin{equation} \label{e6.64}
       \int_{M}\phi dV_{g_{a}} \sim\int_{M}\phi dV_{\roof{g}{\bar}_{a}} = 0.  
 \end{equation} 
 Again, since $\vol G_{a} \rightarrow  0,$ it follows that
 the right side of (6.63) goes to $0$ as $a \rightarrow 0.$ 
 These estimates together then prove (6.55). \\
\end{pf*} 

  The limit of $(M, \roof{g}{\bar}_{a})$ as 
  $a \rightarrow  0$ is the union of the two hyperbolic
  manifolds $M_{1}$ and $M_{2},$ glued along a lower
  dimensional length space of infinite diameter and 
  zero volume. In particular, $M_{1}$ and $M_{2}$ are
  infinitely far apart in the limit. The pointed
  Hausdorff limits of $(M, \roof{g}{\bar}_{a},
  q_{a})$ as $a \rightarrow  0$ are complete subsets
  of this structure, based at some limit point 
  $q = \lim q_{a}.$
 
  Of course, the limit here is quite different from
  the limits of the sequences in \S 6.1-\S 6.3.

 \begin{remark} \label{r 6.10.(i). }
  {\bf (i)}
  Observe that in obtaining (6.55), rather strong use
  has been made of the fact that the limit metrics on
  $M_{1}$ and $M_{2}$ are hyperbolic. Although the
  construction in Example 2 is valid quite generally,
  (c.f. Remark 6.2(i)), and does not require the
  limits to be hyperbolic, (6.55) will no longer hold
  for constructions as above with non-hyperbolic
  limits.

 {\bf (ii).} In addition to taking connected sums of
  hyperbolic manifolds, one may carry out this
  construction on manifolds of the form $N \# M,$ where
  $N$ is any closed oriented graph manifold, (and
  hence $\sigma (N) \geq  0),$ compare with Remark
  6.2(ii).
 
  Namely, carry out the construction above on 
  $M_{2} = M$ and replace the glueing region 
  $G_{a}=S^{3}\setminus(B_{1}\cup B_{2})$ 
  by $N\setminus B.$ Just as $S^{3}$
  admits metrics $\gamma_{a}$ satisfying (6.48) and
  the curvature conditions preceding (6.48), so does
  the manifold $N.$ The remainder of the argument then
  proceeds as before. Proposition 6.9 holds as before
  also.
 
  As in Remarks 6.1 and 6.2, since the construction
  is essentially local, it can be iterated
  (arbitrarily) many times.  \end{remark} 
 \begin{remark} \label{r 6.11.}
   It is worthwhile to mention explicitly that all
   known constructions or examples of degenerating
   sequences of Yamabe metrics on a given 3-manifold
   $M$ have blow-up limits which do have one common
   feature; namely, they are all scalar-flat and
   asymptotically flat.  
 \end{remark} 
 
 \section{\bf Palais-Smale Sequences for Scalar Curvature Functionals.}
 \setcounter{equation}{0}
  The discussion through \S 5 applies to quite
  general sequences of Yamabe metrics, and has made
  no assumptions that the sequence $\{g_{i}\}$ is a
  maximizing sequence for ${\cal S}|_{{\cal C}}$ or
  even a sequence for which ${\cal S} (g_{i})$
  approaches a critical value of ${\cal S}|_{{\cal
  C}}.$ Of course, we are most interested in
  understanding the degeneration of a maximizing
  sequence of Yamabe metrics, and one would expect
  that such a sequence may have some restrictions on
  its behavior not valid for general sequences.

  In this section, we examine more closely the
  structure of the space of metrics ${\Bbb M} $ and
  its completeness properties w.r.t. natural metrics.
  This is then used to understand the existence of
  Palais-Smale sequences for natural functionals on
  ${\Bbb M} $ or its subvarieties. In particular, the differences in the behaviors of the examples in \S 6 can be understood from this viewpoint.
 
\subsection{} 
 
  The space ${\Bbb M} $ is, of course, not complete
  with respect to the $L^2$ or $L^{2,2}$ metrics. We study its completion w.r.t. the $L^{2,2}$ norm, (and later the $T^{2,2}$ norm). The $L^{2,2}$ norm (1.4) on the collection $\{T_{g}{\Bbb M}\}$ of tangent spaces generates a length metric $L^{2,2}$ on ${\Bbb M} ,$ by defining
  $L^{2,2}(g_{1}, g_{2})$ to be the infimum of the lengths
  of curves joining $g_{1}$ to $g_{2}.$
 
  Let ${\Bbb M}_{L^{2,2}}$ be the Cauchy-completion of ${\Bbb M} $
  w.r.t. the $L^{2,2}$ metric, so an (ideal) point in ${\Bbb M}_{L^{2,2}}$ is a limit point of a Cauchy sequence w.r.t. the $L^{2,2}$ metric. Let 
 \begin{equation} \label{e7.1}
  \bar{\Bbb M} = {\Bbb M}_{L^{2,2}}\cap C^{0}({\Bbb M} ), 
 \end{equation} 
where $C^{0}({\Bbb M} )$ is the
 space of continuous Riemannian metrics on $M.$
 
  By work in [E], the space $\bar{\Bbb M}$ may be identified
  with the space of Riemannian metrics on $M$ which,
  with respect to a fixed smooth coordinate atlas, have
  local expressions which are $L^{2,2}$ functions of the
  local coordinates. The tangent space $T_{g}
  \bar{\Bbb M}$ is naturally identified with the
  Hilbert space of symmetric bilinear forms on $M,$
  locally in $L^{2,2}$, with inner product given as in
  (1.4).
 
  Let ${\Bbb M}_{o} \subset  {\Bbb M} $ be the space
  of metrics having a fixed volume form $dV,$  of
  total volume 1, and ${\Bbb M}_{1} \subset  {\Bbb M} $ 
  the space of metrics of total volume 1.
  Similarly, let $\bar{\Bbb M}_{o},$ and
  $\bar{\Bbb M}_{1}$ be the intersection of
  the $L^{2,2}$-completion of ${\Bbb M}_{o}$ or ${\Bbb
  M}_{1}$ with $C^{0}({\Bbb M}_{o})$ or $C^{0}({\Bbb
  M}_{1})$ respectively. The spaces 
  $\bar{\Bbb M},\ \bar{\Bbb M}_{o}$ and 
  $\bar{\Bbb M}_{1}$ are infinite dimensional Hilbert
  manifolds in the topology generated by the $L^{2,2}$ norm, c.f. [E].
  
  It is important to note that the spaces 
  $\bar{\Bbb M},\ \bar{\Bbb M}_{o}, \ \bar{\Bbb M}_{1}$ 
  need not apriori be (Cauchy) complete in the $L^{2,2}$ metric. 
  By definition, of course ${\Bbb M}_{L^{2,2}}$ is complete. 
  However, the symmetric bilinear forms in such a completion 
  may not be positive definite, and thus not metrics. Thus, the
  restriction that $\bar{\Bbb M}$ be contained
  in $C^{0}({\Bbb M} )$ implies that 
  $\bar{\Bbb M}$ might not apriori be complete.

  To understand the completeness of these spaces, following the discussion in [E], we first decompose the space of smooth metrics ${\Bbb M}$ as 
\begin{equation} \label{e7.2} 
{\Bbb M}  = \Vol(M)\times {\Bbb M}_{o},\ \ \ \ 
 g = \phi\cdot  h
 \end{equation}
where $\Vol(M)$ is the space of volume
 forms on $M,$ $h$ is an arbitrary metric in ${\Bbb
 M}_{o}$ and $\phi $ is the ratio of the volume forms
 of $g$ and $h$ to the power $2/n = 2/3.$ Note that
 the metric $g$ in (7.2) is conformal to $h$, so that
 ${\Bbb M}_{o}$ gives a representation of the space
 of conformal structures on $M.$ This should be
 compared with the splitting of ${\Bbb M} $ into the
 space of conformal classes and the space ${\cal C} $
 of Yamabe metrics. Here we recall a well-known result of Moser [Mo] that any metric in ${\Bbb M}_1$ may be pulled back by a diffeomorphism of $M$ to a metric in ${\Bbb M}_o$. In particular, for any unit
 volume Yamabe metric $g\in{\cal C}_{1},$ there is a
 diffeomorphism $\psi $ of $M$ such that $\psi^{*}g
 \in  {\Bbb M}_{o}.$

  Given a fixed background metric $g_{o}\in{\Bbb
  M}_{o},$ any metric $h\in{\Bbb M}_{o}$ may be
  written as
 \begin{equation} \label{e7.3} 
     h = g_{o}\cdot e^X, 
 \end{equation}
  where $X$ is trace-free w.r.t. $g_{o};$ here we are
  using the metric $g_{o}$ to identify bilinear forms
  with linear maps.
 \begin{lemma} \label{l 7.1.}
      The space $\bar{\Bbb M}_{o}$ is (Cauchy) complete.

 \end{lemma} 
\begin{pf}
  Let $\{g_{i}\}$ be a Cauchy sequence in 
  $\bar{\Bbb M}_{o}$ w.r.t. the $L^{2,2}$ metric. Then by
  definition, c.f. (7.1), the sequence $\{g_{i}\}$
  converges to an element $g\in{\Bbb M}_{L^{2,2}},$ i.e. a
  symmetric bilinear form on $M$ with local component
  functions in the $L^{2,2}$-completion of $C^{\infty}(M).$
  By Sobolev embedding $g_{i} \rightarrow  g$ in the
  $C^{0}$ topology. Since the volume forms
  $dV_{g_{i}}$ are fixed, it follows that $dV_{g} =
  dV_{g_{i}}.$ In particular, $dV_{g}$ is a
  (continuous) positive $3$-form on $M.$  However, the
  symmetric bilinear forms in ${\Bbb M}_{L^{2,2}}$ are
  necessarily positive semi-definite. It follows that
  $g$ is a positive definite form and hence in 
  $\bar{\Bbb M}_{o}.$\\
\end{pf}

   It follows from the Hopf-Rinow theorem in infinite dimensions that $\bar{\Bbb M}_{o}$ is also geodesically complete w.r.t. the $L^{2,2}$ metric, i.e. $L^{2,2}$ geodesics in $\bar{\Bbb M}_{o}$ exist for all time. We note that ${\Bbb M}_{o}$ itself is geodesically complete, (but of course not Cauchy-complete), w.r.t the $L^2$ metric, c.f. [E, Thm.8.9]; in fact the $L^2$ geodesics in ${\Bbb M}_{o}$ are given by the simple expression, c.f. (7.3),
\begin{equation} \label{e7.4} 
g(t) = g\cdot  e^{tA} \subset  {\Bbb M}_{o}, 
\end{equation}
where $A\in T_{g}{\Bbb M}_{o},$ so that $A$ is a smooth symmetric bilinear
  form, trace-free w.r.t. $g$. In particular, these geodesics exist for all time. Further, the $L^2$ exponential map is a diffeomorphism onto ${\Bbb M}_{o}.$
 
  On the other hand, the full space $\bar{\Bbb M}$
  is not complete with respect to the $L^{2,2}$ metric, nor geodesically complete w.r.t. the $L^2$ metric. With regard to the latter, in
  contrast to (7.4), the geodesics of ${\Bbb M} $
  with respect to the $L^{2}$ metric do not exist for
  all time, but may leave the space ${\Bbb M} $ in
  finite time. In [FG], the $L^{2}$ geodesic of
  ${\Bbb M} ,$ with initial position $(\mu ,g)$ and
  initial velocity $(w, A),$ c.f. (7.3) and (7.4), is
  calculated to be given by
 \begin{equation} \label{e7.5} g(t) =
               (v(t)^{2}+n^{2}t^{2})^{2/3}g\cdot
               \exp\left(\frac{\tan^{-1}(nt/v)}{n}A\right), 
 \end{equation}
  where $v(t) = 1+\frac{1}{2}(\frac{w}{\mu})t$ and
  $n=\frac{1}{4}(3\tr A^{2})^{1/2}.$
 
  A brief inspection, as noted in [FG], shows that
  these curves escape from ${\Bbb M} $ in finite
  time, (i.e. the form $g(t)$ is no longer positive
  definite), if the initial velocity matrix $A$
  vanishes somewhere on $M.$ For example, in the volume
  or conformal directions, when $A \equiv  0,$ it is
  apparent from (7.5) that $g(t)$ becomes
  degenerate, i.e. non-positive definite, in finite
  time. 

  Although we will not carry it out here,  it
  is not difficult to verify that there are curves
  $g(t),\ t \in  [0,1)$, of the form (7.5), with 
  $A \equiv  0,$ of finite length in the $L^{2,2}$, which degenerate as
  $t \rightarrow  1.$ In particular, the $L^{2}$ norm
  of the curvature blows up as $t \rightarrow  1.$
 
  On the other hand, if the matrix $A$ never vanishes
  on $M,$  i.e., if $A(p) \neq  0,$ for all $p\in M,$
  then the geodesic (7.5) with initial velocity $(w, A)$
  continues in ${\Bbb M} $ for all time, for any
  initial metric $(\mu , g).$ Thus, in most all directions in $\bar{\Bbb M}$, the $L^2$ geodesic (7.5) exists for all time.
 
   Given this background, we now discuss the main result of this section. Let ${\cal F} : \bar{\Bbb M}_{o} \rightarrow
  {\Bbb R}^{+}$ be a $C^{1}$ functional on 
  $\bar{\Bbb M}_{o}$, which is thus bounded below, We only consider 'natural' functionals, in the sense that ${\cal F}(\phi^{*}g) = {\cal F}(g)$, for any $\phi \in$ Diff($M$). For example, it is easy to see that the total scalar curvature functional ${\cal S}$, or the $L^{p}$ norm of the scalar curvature ${\cal S}^p$, for $1 < p \leq 2$, extends to a $C^1$ functional on $\bar{\Bbb M}_{o}$ or $\bar{\Bbb M}_{1}$.

 \begin{theorem} \label{t 7.2.}
   Suppose ${\cal F}: \bar{\Bbb M}_{o} \rightarrow {\Bbb R}_{+}$ is a $C^1$ natural functional, as above, and let $\{\gamma_{i}\}\in{\Bbb M}_{o}$ be a
    minimizing sequence for ${\cal F} .$ 
   Given any $\varepsilon > 0,$ there is another minimizing
    sequence $\{g_{i}\}\in{\Bbb M}_{o}$ for ${\cal F} ,$ 
    with $L^{2,2}(g_{i}, \gamma_{i}) \leq  \varepsilon ,$
    such that $\nabla{\cal F} (g_{i}) \rightarrow  0,$ 
    as $i \rightarrow  \infty $ in the dual $(L^{2,2})^{*}$ norm, i.e.
 	\begin{equation} \label{e7.6} 
 	     |\nabla{\cal F}|^{*}(g_{i}) 
 	     \equiv
 	     \sup_{|\alpha|_{L^{2,2}=1}}
 	    \left|
 	    \int_{M} 
              \langle \nabla_{g_{i}}{\cal F} , \alpha\rangle_{g_{i}}  
                 dV_{g_{i}}
 	    \right| \rightarrow  0,
 	    \ \ \ \ \ {\rm as}\  i \rightarrow  \infty , 
 	\end{equation} 
 for $\alpha = \alpha_{i}\in T_{g_{i}}{\bar{\Bbb M}}_{o}.$
 \end{theorem} 
\begin{pf}
  For otherwise, given any $\varepsilon  >  0,$ there
  exists $c = c(\varepsilon ) >  0$ such that, (for
  some subsequence), 
 	\begin{equation} \label{e7.7}
 	|\nabla{\cal F}|^{*}(g_{i}) \geq  c, \ \ \ \ \ {\rm as} \  
          i \rightarrow  \infty .  
 	\end{equation} 
 for all $g_{i} \in  \bar{\Bbb M}_{o}$ such that 
 	\begin{equation}
 	\label{e7.8} L^{2,2}(g_{i}, \gamma_{i}) \leq  \varepsilon .
 	\end{equation}

  Let $\gamma_{i}(t)$ be a smooth curve in 
  $\bar{\Bbb M}_{o}$, with $\gamma_{i}(0) = \gamma_{i},\ 
  |\frac{d}{dt}\gamma_{i}(t)| = 1,$ 
   (i.e.  $\gamma_{i}(t)$ is parameterized by arclength in
  the $L^{2,2}$ metric), and 
 	 \begin{equation} \label{e7.9} 
                 \int_{M}\langle \nabla_{\gamma_{i}(t)}{\cal F} , 
                 \frac{d}{dt}\gamma_{i}(t)
                         \rangle  dV_{\gamma_{i}(t)}
 		\leq  -  \frac{c}{2}.  
 	 \end{equation} 
 Since $\bar{\Bbb M}_{o}$ is an infinite dimensional manifold
 and $\nabla{\cal F} $ is continuous, it is obvious
 from (7.7) that such curves exist for
 $t$ sufficiently small, say for $0  <   t \leq t_{o} = t_{o}(i);$ 
 one could take for instance
 $\gamma_{i}(t)$ to be the $L^2$ geodesic (7.4) in
 the direction $-\nabla_{\gamma_{i}}{\cal F}$, but reparametrized w.r.t. $L^{2,2}$ arclength.
 
  If $t_{o}  <   \varepsilon ,$ we just repeat the process
 above starting at $\gamma_{i}(t_{o}).$ Since the $L^{2,2}$ metric is complete on  $\bar{\Bbb M}_{o}$, it is clear that one obtains in
  this way a continuous piecewise smooth curve
  $\gamma_{i}(t)$ in $\bar{\Bbb M}_{o},$
  satisfying (7.9) for  $0 \leq  t \leq  \varepsilon .$
  By definition, (7.9) then gives
	 \begin{equation} \label{e7.10}
	     \frac{d}{dt}{\cal F} (\gamma_{i}(t)) =
	     \int_{M}\langle \nabla_{\gamma_{i}(t)}{\cal F} ,
	     \frac{d}{dt}\gamma_{i}(t)\rangle  dV_{\gamma_{i}(t)}
	     \leq  -  \frac{c}{2}.  
	 \end{equation}
  for all $t \in $ [0, $\varepsilon ].$ By integration of (7.10), it
  follows that
 \begin{equation} \label{e7.11} 
    {\cal F} (\gamma_{i}(t)) \leq  {\cal F} (\gamma_{i}) -
    \frac{c}{2}\cdot  t.  
  \end{equation}
  The fact that (7.11) is valid for all $t \leq \varepsilon ,$ 
  contradicts the fact that
  $\{\gamma_{i}\}$ is a minimizing sequence for
  ${\cal F} .$ Hence (7.7) cannot hold in
  conjunction with (7.8) and the result follows.\\
\end{pf}

  A sequence of metrics $\{g_{i}\}$ satisfying (7.6)
  will be called a {\it  Palais-Smale}  sequence for
  ${\cal F} ,$ w.r.t. the $L^{2,2}$ norm. Note that this
  definition corresponds to one of the hypotheses in
  the well-known Condition $C$ of Palais-Smale [PS].
  Of course, the condition (7.6) should not be
  confused with Condition $C$---a compactness
  condition which is much too strong to be of any use
  in this context.
 
  Thus Theorem 7.2 implies that within any $\varepsilon$-neighborhood, w.r.t. the $L^{2,2}$ metric, of a minimizing sequence for ${\cal F}$, there exists a Palais-Smale minimizing sequence.

  Since, by Moser's theorem mentioned above, ${\Bbb M}_{1}$ is obtained from ${\Bbb M}_{o}$ by the action of the diffeomorphism group, and ${\cal F}$ together with (7.6) is diffeomorphism invariant, it follows that Theorem 7.2 holds for ${\cal F}$ considered as a functional ${\cal F}: \bar{\Bbb M}_{1} \rightarrow {\Bbb R}^{+}$, with $\alpha \in T{\Bbb M}_{1}$.

\begin{remark} \label{r 7.3}
   It is easy to see that Theorem 7.2 also holds for
  functionals on other natural submanifolds of 
  ${\Bbb M}_{1}$, obtained by diffeomorphims from ${\Bbb M}_{o}$. For instance, at least in case
  $\sigma (M) \leq  0,$ so that ${\cal C}_{1} $ is an
  infinite dimensional submanifold of ${\Bbb M}_{1} ,$
  one may verify without difficulty that Theorem 7.2
  holds for ${\cal S}|_{{\cal C}_{1}},$ (or $-{\cal
  S}|_{{\cal C}_{1}}).$ 
 \end{remark}

\subsection{}
  Since the $L^{2,2}$ metric is rather strong, the corresponding dual $(L^{2,2})^{*} = L^{-2,2}$ metric is rather weak. It is thus of interest to understand if the results above can be strengthened by use of a weaker norm than the $L^{2,2}$ norm. Of course, the main issue in \S 7.1 is the Cauchy completeness of $\bar {\Bbb M}_{o}$.

  In this respect, the following might be useful.

\begin{lemma} \label{l 7.4}
  The completion of ${\Bbb M}$ w.r.t. the $T^{2,2}$ norm (1.9), as in (7.1), is the same as its completion $\bar{\Bbb M}$ w.r.t. the $L^{2,2}$ norm.
\end{lemma}
\begin{pf}
  Let ${g_i}$ be a $T^{2,2}$ Cauchy sequence in ${\Bbb M}$, converging to an element $g$ in ${\Bbb M}_{T^{2,2}} \cap C^{0}(M)$. Thus $g$ is a continuous metric on $M$, and the local components $g_{kl}$ of $g$ in a smooth atlas on $M$ are $T_{g}^{2,2}$ functions of the local coordinates. It follows that $\Delta_{g}(g_{kl})$ is locally in $L^2$ on $M$. Hence, using the continuity of $g$, elliptic regularity theory, c.f. [GT, Thms.8.8,9.15], implies that $g_{kl}$ is locally in $L_{g}^{2,2}$. Thus $g \in \bar {\Bbb M}$, which gives the result.
\end{pf}

  It is now easy to verify that all of the results above, in particular Theorem 7.2, hold also w.r.t. the weaker $T^{2,2}$ norm; hence there exist (many) Palais-Smale sequences for ${\cal F}$ as above w.r.t. the $T^{2,2}$ norm.

  Let us apply this to the functional ${\cal S}|_{{\cal C}_{1}}$. If ${g_i}$ is a Palais-Smale sequence of Yamabe metrics in ${\cal C}_{1}$, then
\begin{equation} \label{e7.12} 
 \| \nabla_{g_i}  ( {\cal S}|_{{\cal C}_1 })
 \|_{T^{-2,2}(T{\cal C}_1 )} 
        = \sup_{|\alpha|_{T^{2,2}} = 1} 
          \left|   \int_{M} \langle z^{T}, \alpha\rangle dV_{g_{i}}
          \right|  
       \rightarrow  0,
 \end{equation}
\noindent
  where $z^{T} = z^{T}_{g_i}$ and $\alpha \in T_{g_{i}}{\cal C}_{1}$.
  Now (unfortunately), the condition that 
  $\alpha  = \alpha_{i}\in T_{g_{i}}{\cal C}_{1} $ is a global
  condition on $\alpha ,$ (as a tensor on $(M, g_{i}))$, 
   and hence it is difficult to obtain local information 
  on $z^{T}$ from (7.12). If (7.12) could be strengthened to hold
  for all $\alpha\in T_{g_{i}}{\Bbb M}_{1}$, i.e.
  if $\{g_{i}\}$ is  strongly  Palais-Smale in the sense of (4.6), 
  then the  estimate (7.12) does provide non-trivial local control. 
  In particular, Propositions 4.1 and 4.2 hold in this situation. 
  We do not know however if there exist such strongly Palais-Smale 
  sequences on general 3-manifolds $M$.

  With regard to the examples of \S 6, it is easily
  verified that all sequences in Example 4 are
  Palais-Smale for ${\cal S}|_{{\cal C}_1}$ w.r.t. the
  $T^{2,2}$ norm. In fact, such
  sequences are strongly Palais-Smale, in the sense
  of (4.6), since the full $L^{2}$ norm of $z^{T}$
  goes to 0 by Proposition 6.9. Note here that of
  course $||z^{T}||_{T^{-2,2}}
  \leq  ||z^{T}||_{L^{2}}$ on $\{g_{i}\}.$ Similarly,
  the sequences in Example 3 are strongly
  Palais-Smale, since the same, but in this case much
  simpler, reasoning proving Proposition 6.9 holds
  here also. On the other hand, Examples 1 and 2 are
  clearly not strongly Palais-Smale for ${\cal
  S}|_{{\cal C}_1}$ w.r.t. $T^{2,2},$
  since for instance
  $\lim_{i\to\infty}(z^{T})_{g_{i}}$ is
  non-zero on the limit $(M_{o}, g_{o}),$ c.f. \S
  6.3. It seems non-trivial to prove, although we
  certainly conjecture, that the sequences in Example
  1 and 2 are not Palais-Smale for ${\cal S}|_{{\cal C}_1}.$

 \section{\bf Appendix} \setcounter{equation}{0}

  In this Appendix, we prove Theorem 3.2. We break
  the proof into several parts. The first result
  states that the event horizon $\Sigma $ must be
  non-empty.

\appendix 
\setcounter{section}{1}
 \begin{theorem} \label{t A.1. }
  The only complete solution $(N, g, u)$ to the static
  vacuum Einstein equations with $u >  0$
  everywhere is the flat solution on ${\Bbb R}^{3}$
  or a quotient ${\Bbb R}^{3}/\Gamma ,$ with $u = const.$
 \end{theorem}
\begin{pf}
Since $u$ is positive on the complete
  manifold $N,$  the associated 4-manifold $X^{4} =
  N^{3}\times_{u}S^{1}$ is Ricci-flat and complete, c.f.
  (1.14). Define the function $h$ by $$h = \log u. $$
 A straightforward calculation gives \begin{equation}
 \label{eA.1} \Delta_{X}h = 0.
 \end{equation}

  Suppose first that either $u$ is bounded above, or
  bounded below away from 0. Thus, $h$ has either a
  uniform upper or a uniform lower bound. The
  Liouville theorem of Yau [Yu] then implies that $h = const.,$ 
  and thus $u = const. >  0.$ The
  equations then immediately imply that the metric
  $g$ is flat.
 
  In general, since $X^{4}$ is Ricci flat, the
  infinitesimal Harnack inequality of Cheng-Yau
  [CY, Theorem~6] gives the bound 
 $$
  \sup_{\scriptscriptstyle B^{4}(r)} |\nabla(\log v)| \leq  C\cdot \frac{1}{r},
  $$
  for any positive harmonic function $v$ defined on a
  geodesic ball $B^{4}(r)\subset X.$ In particular,
  by integration, this leads to the usual Harnack
  inequality 
 \begin{equation} \label{eA.2}
 \sup_{\scriptscriptstyle B^{4}(r)}v \leq  C\cdot   \inf_{\scriptscriptstyle B^{4}(r)}v.
 \end{equation}

  Now suppose $u$ is not constant. Then applying the
  estimate (A.2) to the functions $h - a$ and $b - h,$
  where $a$ and $b,$ depending on $r,$ are chosen to
  make these functions positive on $B^{4}(r)$
  implies, by a well known technique due to Moser,
  c.f. [GT, Theorem 8.22], that $u$ has oscillation
  growing at a definite power of $r,$ as 
  $r \rightarrow  \infty .$ 
  In fact, an observation of
  Cheng [Cg] is that if $h$ is non-constant, then $h$
  must have at least linear growth, i.e. there is a
  constant $C  <   \infty $ such that, for $r \geq  1,$
 \begin{equation} \label{eA.3} 
  \osc_{\scriptscriptstyle B^{4}(r)}  h = 
      \sup_{\scriptscriptstyle B^{4}(r)}h 
      -\inf_{\scriptscriptstyle B^{4}(r)}h  \geq  C\cdot r, 
 \end{equation} 
 where as above $B^{4}(r)$ is the
 geodesic $r$-ball about some fixed point $x_{o}\in X.$
 Suppose first that 
        \begin{equation} \label{eA.4}
        \inf_{\scriptscriptstyle B^{4}(r)}h  <   - c_{1}r, 
 \end{equation}
  for some $c_{1} >  0.$ It follows again from the
  Harnack inequality that for any $r,$ there are points
  $x_{r}\in S^{3}_{\pi (x_{o})}(r)\subset N$ such
  that $u(x) \leq  c_{2}e^{- c_{1}\cdot  r},$ for all
  $x\in B_{x_{r}}^{3}(1), \ r$ large, for some constant
  $c_{2} >  0;$ here $\pi : X^{4} \rightarrow
  N^{3}$ is the projection on the first factor and
  $S^{3} (B^{3})$ denote geodesic spheres (resp.
  balls) in~$N.$ We use this to estimate the volume of
  regions in $X^{4}.$ Setting $U_{x_{r}}(1) =
  \pi^{-1}(B_{x_{r}}^{3}(1)),$ we have
 	\begin{equation} \label{eA.5}
 	\vol_{X^{4}}U_{x_{r}}(1) = \int_{B_{x_{r}}(1)}u\ dV
 	\leq  c_{3}\cdot  \vol B_{x_{r}}^{3}(1)
 	e^{-c_{1}\cdot  r}. 
 	\end{equation}

  To estimate $\vol B_{x_{r}}^{3}(1),$ it follows from
  the Harnack inequality as above that there is a
  constant $C_{o}$, independent of $r,$ s.t.
 \begin{equation} \label{eA.6}
   C_{o}^{-1} \leq     \frac{\sup_{B} \ u(x)}{\inf_{B}\ u(x)}   \leq  C_{o}.  
 \end{equation} 
 where we set $B = B_{x_{r}}^{3}(1) \subset  U_{x_{r}}(1);$ here
 we are using the fact that $u$ is invariant in the
 fiber or $S^{1}$ factor. Consider the conformally
 equivalent metric $\tilde{g} = u^{2}g$ from
 (1.19). Thus, if $u_{o} =  \inf_{B} u(x),$ then
     $\vol_g B \leq  u_{o}^{-3}\ \vol_{\tilde{g}} B.$
 On the other hand, since distances in
 $\tilde{g}$ are at most $C_{o}u_{o}$ times
 distances in $g,$ it follows that $B \subset
 B_{\tilde{g}} (C_{o}\cdot  u_{o}).$ Since from
 (1.20), $Ric_{\tilde{g}} \geq  0,$ the volume
 comparison theorem for Ricci curvature implies
 $\vol_{\tilde{g}} B(C_{o}\cdot  u_{o}) 
    \leq \frac{4}{3}\pi C_{o}^{3}u_{o}^{3}.$ It follows there
 is a constant $D,$ independent of $r,$ s.t.
 	\begin{equation} \label{eA.7}
 	 \vol B_{x_{r}}^{3}(1) \leq  D, 
 	\end{equation} so that
 from (A.5), we obtain 
 	\begin{equation} \label{eA.8}
 	\vol_{X^{4}}U_{x_{r}}(1) \leq  c_{4} e^{- c_{1}\cdot r} , 
 	\end{equation}
  for some constant $c_{4}.$ However, (A.8)
  contradicts the (relative) volume comparison
  principle. Namely, for Ricci flat 4-manifolds, and
  $r \geq  1,$ one has the bound
 $$ 
 \frac{\vol B_{x_{r}}^{4}(r)}{r^{4}} \leq  \vol
 B_{x_{r}}^{4}(1) = \vol U_{x_{r}}(1), 
 $$
  where the equality follows from the fact that the
  fiber $S^{1}$ has exponentially small length near
  $x_{r}\in N\subset X.$ Together with (A.8), this
  implies
 $$\vol B_{x_{r}}^{4}(r) \leq  c_{4}r^{4} e^{-
 c_{1}\cdot  r} \rightarrow  0 
 \ \ \ \ \ {\rm as}\  r \rightarrow \infty , 
 $$
  which is of course impossible.
 
  On the other hand, if (A.4) does not hold, then by
  (A.3), one has the estimate
       $$ \sup_{\scriptscriptstyle B^{4}(r)} h > c_{5}\cdot   r, $$ 
 for some $c_{5} >  0.$ Hence,
 there are points $x_{r}\in N\subset X$ near which
 the fiber $S^{1}$ has exponentially large growth.
 Arguing in exactly the same way as above, from the
 volume comparison theorem for non-negative Ricci
 curvature, $\vol B_{x_{r}}(1) \geq  c_{6}r^{-3}$ and
 hence,
 	$$\vol B_{x_{r}}^{4}(r) \geq  c_{7}r^{- 3}
 	e^{c_{5}\cdot  r} , 
         $$
  which is again impossible by the relative volume
  comparison theorem on $X.$
 
  It follows that $h$ must be constant and thus 
  $(N, g)$ is flat.\\
\end{pf}
 
  We note the following compactness principle for
  solutions of the static equations.
 
 \begin{lemma} \label{l A.2}
   Let $(g_{i}, u_{i})$ be a sequence of solutions
  to the static vacuum equations (0.16), defined on a
  geodesic ball $B_{i} = B_{x_{i}}(1),$ (w.r.t. the
  metric $g_{i}).$ Suppose
 \begin{equation} \label{eA.9} 
                 \vol B_{i} \geq  v_{o} >0,
                 \ \ \ \ 
                  |r_{i}| \leq  \Lambda <  \infty ,
 \end{equation} 
 and 
 \begin{equation} \label{eA.10}
 u_{i}(x_{i}) \geq  u_{o} >  0, 
 \ \ \ \ \   u_{i}(y_{i}) > 0, 
 \ \  \ \forall y_{i}\in B_{x_{i}}(1),
 \end{equation}
  for some constants $v_{o},\  \Lambda ,\  h_{o}.$ Then,
  for any $d >  0,$ a subsequence of the
  Riemannian manifolds $(B_{x_{i}}(1- d), g_{i})$
  converges, in the $C^{\infty}$ topology, modulo
  diffeomorphisms, to a limit manifold $(B_{x}(1-
  d),$ g), with limit function $u.$ The triple
  $(B_{x}(1- d), g, u)$ is a smooth solution of the
  static vacuum equations.
 \end{lemma}
\begin{pf}
  The Cheeger-Gromov compactness theorem, c.f.
  Theorem 1.3 and references there, implies that a
  subsequence of $(B_{i}, g_{i})$ converges, in the
  $C^{1,\alpha}$ topology, to a limit manifold $(B,
  g).$  The bounds (A.9) and (A.10) imply further that
  the positive harmonic functions $u_{i}$ satisfy a
  uniform Harnack inequality in $B_{i}(1- d).$ If
  there is a uniform bound $u_{i}(x_{i}) \leq  u^{o}
   <   \infty ,$ it follows that a subsequence of
  $\{u_{i}\}$ converges in the $C^{o}$ topology to a
  limit positive function $u.$ If not, so $u_{i}(x_{i})
  \rightarrow  \infty ,$ we renormalize $u_{i}$ by
  setting $\bar{u}_{i} = u_{i}/u_{i}(x_{i}).$ 
  Then again the fact that
  $\bar{u}_{i}$ is harmonic and the Harnack
  inequality imply that $\{\bar{u}_{i}\}$ has
  a subsequence converging to a positive limit
  harmonic function $u$ in $B(1- d).$ As noted in \S
  1.3, since the metrics $g_{i}$ satisfy an elliptic
  equation (0.16), the regularity theory for elliptic
  equations implies that one has $C^{k,\alpha}$ bounds
  for solutions, in terms of $C^{1,\alpha}$ bounds,
  i.e. all covariant derivatives of the curvature
  are bounded in terms of the bounds (A.9)-(A.10).
  Equivalently, one can pass to the Ricci flat
  4-manifold $N$ and use the Einstein equation
  to obtain $C^{k,\alpha}$ regularity, for any $k.$
  These estimates imply convergence in the $C^k$ topology,
  and thus $C^\infty$ convergence.

  Using elliptic regularity on the equations (0.16),
  the $L^{\infty}$ bound on $|r|$ in (A.9) may be
  replaced by weaker bounds; for instance, it
  suffices to assume an $L^{2}$ bound on $r.$\\
\end{pf}
 
  The following Corollary gives an apriori estimate
  for the curvature of a solution away from the event
  horizon $\Sigma .$ The proof is a standard
  consequence of Theorem A.1 and Lemma A.2.
 
 \begin{corollary} \label{c A.3.}
   Let $(N, g, u)$ be a solution to static vacuum
   equations (0.16), and $U \subset  N$ a domain with
   smooth boundary on which $u >  0.$ Then there
   is an (absolute) constant $K  <   \infty ,$
   independent of $(N, g)$ and $U,$ such that for all
   $x\in U,$ 
 	\begin{equation} \label{eA.11} 
 	|z|(x) \leq  \frac{K}{t(x)^{2}}, 
 	\end{equation} 
 where $t(x) = \dist(x, \partial U).$
 \end{corollary} 
\begin{pf}
  The proof is by contradiction. Thus, assume that
  (A.11) does not hold. Then there are static vacuum
  solutions $(N_{i}, g_{i}, u_{i}),$ smooth domains
  $U_{i}\subset N_{i}$ on which $u_{i} >  0$ and
  points $x_{i}\in U_{i}$ such that
 \begin{equation} \label{eA.12}
	 t^{2}(x_{i})\cdot |z_{i}|(x_{i}) \rightarrow  \infty ,
	 \ \ \ 
	 {\rm as}\ i \rightarrow  \infty .  
 \end{equation}
  Let $t_{i} = t(x_{i}).$ Since it may not be
  possible to choose the points $x_{i}$ so that they
  maximize $|z_{i}|$ (over large domains), we shift
  the base points $x_{i}$ as follows: choose
  $s_{i}\in [0,t_{i})$ such that
 	\begin{equation} \label{eA.13} 
 	     s_{i}^{2} \sup_{B_{x_{i}}(t_{i}-s_{i})} |z_{i}|   
 	      =   \displaystyle \sup_{s\in [0,t_{i})} 
                 \bigl(s^{2} 
 	     \cdot    
 	     \sup_{B_{x_{i}}(t_{i}- s)}|z_{i}|\bigr) \rightarrow\infty , 
 	       \ \ \ {\rm as}\  i \rightarrow \infty , 
 	\end{equation}
  where the last estimate follows from (A.12), (set
  $s_{i} = t_{i}).$ Let $y_{i}\in
  B_{x_{i}}(t_{i}-s_{i})$ be points such that
 	\begin{equation} \label{eA.14} 
 	|z_{i}|(y_{i}) = \sup_{B_{x_{i}}(t_{i}-s_{i})}|z_{i}|.  
 	\end{equation}
 Further, setting $s = s_{i}(1-\frac{1}{k}),\  k >  1,$ 
 in (A.13), one obtains the estimate
 	\begin{equation} \label{eA.15}
 	s_{i}^{2}|z_{i}|(y_{i}) 
        \geq
 	s_{i}^{2}  \left(1-\tfrac{1}{k}\right)^{2}   \cdot
 	\sup_{B_{x_{i}}(t_{i}- s_{i}(1-\frac{1}{k}))}|z_{i}|
        \geq    
         s_{i}^{2} \left(1-\tfrac{1}{k}\right)^{2}  \cdot
 	\sup_{B_{y_{i}}(s_{i}/k)}|z_{i}| , 
 \end{equation} 
 so that 
        \begin{equation} \label{eA.16}
 	   \sup_{B_{y_{i}}(s_{i}/k)}|z_{i}|  
           \leq
 	   \left(1-\tfrac{1}{k}\right)^{-2} |z_{i}|(y_{i}),
 	\end{equation}

  Now rescale or blow-up the metric so that
  $|\tilde{z}_{i}|(y_{i}) = 1$ by setting
  $\tilde{g}_{i} = |z_{i}|(y_{i})\cdot  g,$
  and consider the pointed sequence $(U_{i},
  \tilde{g}_{i}, y_{i}).$ We have
 	\begin{equation} \label{eA.17}
 	|\tilde{z}_{i}|(y_{i}) = 1, 
 \end{equation} 
 and by (A.13) and scale invariance, 
 	\begin{equation}
 		\label{eA.18} \dist_{\tilde{g}_{i}}(y_{i},
 	\partial U_{i}) \rightarrow  \infty , 
 	\ \ \ \ \  {\rm as}\  i
 	\rightarrow\infty .  
 	\end{equation} 
 Also, it follows
 from (A.16) that 
 		\begin{equation} \label{eA.19}
 		| \tilde{z_i}|(x) \leq
 		C(\dist_{\tilde{g}_{i}}(x, y_{i})).
 	\end{equation}
  We also normalize $u$ by setting \begin{equation}
 \label{eA.20} \tilde{u}_{i}(x) =
 \frac{u(x)}{u(y_{i})}, \end{equation} and note by
 construction that $\tilde{u}_{i} >  0 $ on
 $U_{i}.$
 
  Suppose first that there exists $v >  0$ such
  that 
 \begin{equation} \label{eA.21}
 \vol_{\tilde{g}_{i}}B_{y_{i}}(1) \geq  v.
 \end{equation} 
 Then the volume comparison
 principle for bounded curvature and (A.21) 
 imply that $\vol_{{\tilde{g}}_i}B_{x_i}(1)\geq v(x_{i}),$ 
  where $v(x_{i})$ depends only on
 $\dist_{\tilde{g}_{i}}(x_{i}, y_{i}).$ By the
 compactness of solutions, Lemma A.2, it follows that
 a subsequence converges, in the $C^{\infty}$
 topology, to a limit solution $(U_{\infty}, g_{\infty}, u_{\infty}),$ 
 which is complete and satisfies $u_{\infty} >  0$ everywhere. 
 (The minimum principle for harmonic functions implies
 that $u_{\infty}$ cannot vanish anywhere). By
 Theorem A.1, $g_{\infty}$ must be flat and
 $u_{\infty}$ constant. However, the smooth
 convergence guarantees that the equality (A.17)
 passes to the limit, contradicting the fact that
 $g_{\infty}$ is flat.
 
  If (A.21) is not satisfied, so that
  $\vol_{\tilde{g}_{i}}B_{y_{i}}(1) \rightarrow
   0,$ as $i \rightarrow  \infty ,$ it follows that
  the sequence $(U_{i}, \tilde{g}_{i}, y_{i})$
  is collapsing in the sense of Cheeger-Gromov on
  balls $(B_{y_{i}}(R_{i}), \tilde{g}_{i}),$
  where $R_{i} \rightarrow  \infty $ as $i
  \rightarrow  \infty .$ In dimension 3, the
  structure theory of collapse implies that the
  collapse is along an injective F-structure. More
  precisely, the balls $(B_{y_{i}}(R_{i}),
  \tilde{g}_{i})$ have the structure of a
  Seifert fibration, with fibers that are injective
  in the fundamental group, c.f. [An2, \S 2,3] or
  [An II, \S 2]. Thus, one may pass to the universal
  cover of $B_{y_{i}}(R_{i}).$ Since the universal
  covers no longer collapse, i.e. (A.21) is
  satisfied, one may apply the discussion above to
  the universal covers, and obtain a contradiction in
  the same manner.\\
\end{pf}
 
  We note that exactly the same proof can be used to
  prove also that
 	\begin{equation} \label{eA.22} 
 		|\nabla \log u|(x) \leq
 		\frac{K}{t(x)}.  
 	\end{equation}
  In particular, these results together prove Theorem 3.2.
 
  Note that since $K$ is independent of the domain $U$,
  (A.11) holds for $t$ the distance to the event
  horizon $\Sigma ,$ (provided this is defined), even
  if $\Sigma $ is singular. More generally but for
  the same reasons, (A.11) holds for $t$ the distance
  to $\partial N^{o},$ where $N^{o}$ is the maximal
  domain on which the potential $\bar{u}$ is
  positive, c.f. Theorem 5.1. This is because $N^{o}$
  is the Hausdorff limit of an exhaustion of $N^{o}$
  by smooth subdomains.

   \begin{remark} 
   {\bf (i).} Similarly, using elliptic regularity associated to
   the static vacuum equations, one may show in the
   same way that for any $k \geq  1,$
 \begin{equation} \label{eA.23} 
 |\nabla^{k}z|(x)
 \leq  \frac{C(k)}{t^{2+k}(x)},  
  \ \ \ \ \   |D^{k}\log u|(x)
 \leq  \frac{C(k)}{t^{k}(x)}.  
 \end{equation} 
 {\bf (ii).}
  Under the same hypotheses as Corollary A.3, suppose
  also that there is an end $E \subset  N$ such that
 \begin{equation} \label{eA.24} 
  u(x) \rightarrow const. > 0, 
 \end{equation} as $x \rightarrow
 \infty $ in $E$. Then essentially the same proof shows
 that there is a function $\mu  = \mu (t)$ such that
 \begin{equation} \label{eA.25} 
 |z|(x) 
 \leq
 \frac{\mu (t(x))}{t^{2}(x)}, 
 \ \ \ \ \  
  |\nabla \log u|(x) \leq
 \frac{\mu (t(x))}{t(x)}, 
 \end{equation}
  where $t(x)$ is the distance to a fixed base point in
  $E$. Namely, if (A.25) does not hold, then one
  derives the same estimates as (A.12), (A.13), (A.18),
  with $\infty $ replaced by some constant $c > 0 $
  and concludes the proof as before using the fact
  that static vacuum solutions with $u = const.$ are flat.
   \end{remark}
  \bibliographystyle{plain}

  \bigskip

\begin{center}
May, 1998; revision, November, 1998
\end{center}
\bigskip
\address{Department of Mathematics\\
S.U.N.Y. at Stony Brook\\
Stony Brook, N.Y. 11794-3651}\\
E-mail: anderson@math.sunysb.edu

 \end{document}